\documentclass[11pt,twoside]{article}
\usepackage{fancyhdr}
\usepackage{amsfonts,epsfig,graphicx}
\usepackage{afterpage}
\usepackage{amsmath,amssymb,amsthm,dsfont} 
\usepackage{fullpage}
\usepackage[T1]{fontenc} 
\usepackage{epsf} 
\usepackage{graphics} 
\usepackage{amsfonts,amsmath}
\usepackage[ruled,boxed,commentsnumbered]{algorithm2e}
\usepackage[]{natbib} 
\usepackage{psfrag,xspace}
\usepackage{color,etoolbox}
\usepackage{booktabs}

\RequirePackage[OT1]{fontenc}
\RequirePackage{amsthm,amsmath}
\RequirePackage[]{natbib}
\usepackage{amsfonts}
\usepackage{amssymb}
\usepackage{multirow}
\usepackage{graphicx}
\usepackage{mathtools}
\usepackage[colorlinks,citecolor=blue,urlcolor=blue,linkcolor=blue,bookmarks=false]{hyperref}

\newtheorem{theorem}{Theorem}[section]
\newtheorem{definition}{Definition}[section]
\newtheorem{lemma}{Lemma}[section]
\newtheorem{proposition}{Proposition}[section]
\newtheorem{corollary}{Corollary}[section]

\theoremstyle{remark}
\newtheorem{remark}{Remark}[section]
\newtheorem{example}{Example}[section]

\newcommand{\mP}{\mathbb{P}}
\newcommand{\mE}{\mathbb{E}}
\newcommand{\mV}{\mathbb{V}}

\newcommand{\tr}{\text{tr}}
\newcommand{\ind}{\mathds{1}}  % Indicator

\newcommand{\convP}{\overset{p}{\longrightarrow}}
\newcommand{\convD}{\overset{d}{\longrightarrow}}

\DeclarePairedDelimiter\floor{\lfloor}{\rfloor}
\DeclarePairedDelimiter\ceil{\lceil}{\rceil}

\begin{document}

\title{}

\begin{center} {\Large{\bf{Robust Multivariate Nonparametric Tests via Projection-Averaging}}}

\vspace*{.3in}

{\large{
\begin{tabular}{ccccc}
Ilmun Kim$^\dagger$&&Sivaraman Balakrishnan$^\dagger$ && Larry Wasserman$^\dagger$ \\
\end{tabular}

\vspace*{.1in}

\begin{tabular}{ccc}
Department of Statistics and Data Science$^{\dagger}$ \\
\end{tabular}

\begin{tabular}{c}
Carnegie Mellon University \\
Pittsburgh, PA 15213
\end{tabular}
}}

\vspace*{.2in}

\begin{abstract}
In this work, we generalize the Cram{\'e}r-von Mises statistic via projection-averaging to obtain a robust test for the multivariate two-sample problem. The proposed test is consistent against all fixed alternatives, robust to heavy-tailed data and minimax rate optimal against a certain class of alternatives. Our test statistic is completely free of tuning parameters and is computationally efficient even in high dimensions. When the dimension tends to infinity, the proposed test is shown to have comparable power to the existing high-dimensional mean tests under certain location models. As a by-product of our approach, we introduce a new metric called {\em the angular distance} which can be thought of as a robust alternative to the Euclidean distance. Using the angular distance, we connect the proposed method to the reproducing kernel Hilbert space approach. In addition to the Cram{\'e}r-von Mises statistic, we demonstrate that the projection-averaging technique can be used to define robust, multivariate tests in many other problems.
\end{abstract}
\end{center}

\vskip 1em

\section{Introduction}

Let $X$ and $Y$ be random vectors defined on a common probability space $(\Omega, \mathcal{A}, \mP)$ with distributions $P_X$ and $P_Y$, respectively. Given two mutually independent samples $\mathcal{X}_m = \{X_1,\ldots,X_m\}$ and $\mathcal{Y}_n =
\{Y_1,\ldots, Y_n\}$ from $P_X$ and $P_Y$, we want to test 
\begin{align} \label{Eq: Two-Sample Hypotheses}
H_0 : P_X = P_Y\quad \text{versus} \quad H_1 : P_X \neq P_Y.
\end{align}
This fundamental problem has received considerable attention in statistics with a wide range of applications \citep[see e.g.][for
a review]{thas2010comparing}. A common statistic for the univariate two-sample testing is the Cram{\'e}r-von Mises (CvM) statistic \citep{anderson1962distribution}:
\begin{align*} %\label{Eq: Univaraite CvM-statistic}
\frac{mn}{m+n} \int_{\infty}^{\infty}  \big( \widehat{F}_X(t) - \widehat{F}_Y(t) \big)^2 d \widehat{H}(t),
\end{align*} 
where $\widehat{F}_X(t)$ and $\widehat{F}_Y(t)$ are the empirical distribution functions of $\mathcal{X}_m$ and $\mathcal{Y}_n$, respectively, and $(m+n)\widehat{H}(t) = m \widehat{F}_X(t) + n \widehat{F}_Y(t)$. Another approach is based on the energy statistic, which is an estimate of the
squared energy distance \citep{szekely2013energy}:
\begin{align*}
E^2 = 2 \mE\normalsize[ | X_1 - Y_1 | \normalsize] - 
\mE \normalsize[ | X_1 - X_2 | \normalsize] - \mE \normalsize[ | Y_1 - Y_2 | \normalsize].
\end{align*}
The energy distance is well-defined assuming a finite first moment and it can be written in a form that is similar to Cram{\'e}r's distance~\citep{cramer1928composition}, namely,
\begin{align*}
E^2 = 2 \int_{-\infty}^{\infty}  \big( F_X(t) - F_Y(t) \big)^2 dt,
\end{align*}
where $F_X(t)$ and $F_Y(t)$ are the distribution functions of $X$ and $Y$, respectively.

The CvM-statistic has several advantages over the energy statistic for univariate two-sample testing. For instance, the CvM-statistic is distribution-free under $H_0$ \citep{anderson1962distribution} and its population counterpart is well-defined without any moment assumptions. It also has an intuitive probabilistic interpretation in terms of probabilities of concordant and discordance of four independent random variables \citep[][]{baringhaus2017cramer}. Nevertheless, the CvM-statistic has rarely been studied for multivariate testing. A primary reason is that the CvM-statistic is essentially rank-based, which leads to a challenge to generalize it in a multivariate space. In contrast, the energy statistic can be easily applied in arbitrary dimensions as in \cite{baringhaus2004new} and \cite{szekely2004testing}. Specifically, they defined the squared multivariate energy distance by
\begin{align} \label{Eq: Definition of Energy Distance}
E_d^2(P_X,P_Y) = 2\mE \normalsize[ \| X_1 - Y_1 \| \normalsize] - 
\mE \normalsize[ \| X_1 - X_2 \| \normalsize] - \mE \normalsize[ \| Y_1 - Y_2 \| \normalsize],
\end{align}
where $\| \cdot \|$ is the Euclidean norm in $\mathbb{R}^d$. The multivariate energy distance maintains the characteristic property that it is always non-negative and equal to zero if and only if $P_X = P_Y$. It can also be viewed as the average of univariate Cram{\'e}r's distances of projected random variables \citep{baringhaus2004new}:
\begin{align} \label{Eq: Projection-type expression for Energy distance}
E^2_d(P_X,P_Y)  = \frac{\sqrt{\pi}(d-1)\Gamma(\frac{d-1}{2})}{\Gamma(\frac{d}{2})}
\int_{\mathbb{S}^{d-1}} \int_{\mathbb{R}} \left(F_{\beta^\top X}(t) - F_{\beta^\top Y}(t)\right)^2 dt d\lambda(\beta),
\end{align}
where $\lambda$ represents the uniform probability measure on the $d$-dimensional unit sphere $\mathbb{S}^{d-1} = \{x \in \mathbb{R}^d :\|x\| = 1\}$ and $\Gamma(\cdot)$ is the gamma function.

Although the multivariate energy distance can be easily estimated in any dimension, it still requires the finite moment assumption as in
the univariate case. When the underlying distributions violate this moment condition with potential outliers, the resulting energy test might suffer from low power.  Given that outlying observations arise frequently in practice with high-dimensional data, there is a need to develop a robust counterpart of the energy distance. The primary goal of this work is to introduce a robust, tuning parameter free, two-sample testing procedure that is easily applicable in arbitrary dimensions and consistent against all fixed alternatives. Specifically, we modify the univariate CvM-statistic to generalize it to an arbitrary dimension by averaging over all one-dimensional projections. In detail, the proposed test statistic is an unbiased estimate of the squared multivariate CvM-distance defined as follows:
\begin{align} \label{Eq: Definition of the multivariate CvM-distance}
W^2_d(P_X,P_Y) = \int_{\mathbb{S}^{d-1}} \int_{\mathbb{R}} \left(F_{\beta^\top X}(t) - F_{\beta^\top Y}(t)\right)^2 dH_\beta(t)d\lambda(\beta),
\end{align}
where $H_\beta(t) = \vartheta_X F_{\beta^\top X}(t) + \vartheta_Y F_{\beta^\top Y}(t)$ and $\vartheta_X$ is a fixed value in $(0,1)$ and $\vartheta_Y = 1 - \vartheta_X$.  For simplicity and when there is no ambiguity, we may omit the dependency on $P_X,P_Y$ and write $W_d(P_X,P_Y)$ as $W_d$. % \textcolor{red}{TO DO} $P_X,P_Y$ distributions / $F_X,P_Y$ distribution functions

%Depending on the underlying distributions, we may write $W_d(F_X,P_Y)$ or $W_d(P,Q)$

% For simplicity and when there is no ambiguity, we may omit the dependency on $F_X,P_Y$ and write $W_d(P_X,P_Y)$ or $W_d(P,Q)$ as $W_d$. %It is also worth mentioning that $W_d$ is invariant to the choice of $\vartheta_X$ whenever $\beta^\top X$ and $\beta^\top Y$ have continuous distribution functions for $\lambda$-almost all $\beta \in \mathbb{S}^{d-1}$ (Theorem~\ref{Theorem: Multivaraite CvM Expression}). 
%As shall we see in Theorem~\ref{Theorem: Multivaraite CvM Expression}, $W^2_d$ is invariant to the choice of $\vartheta_X$ and $\vartheta_Y$ whenever $\beta^\top X$ and $\beta^\top Y$ have continuous distribution functions for all $\beta \in \mathbb{S}^{d-1}$, $\lambda$-almost surely. 

Throughout this paper, we refer to the process of averaging over all projections as {\em projection-averaging}.

\subsection{Summary of our results}

The proposed multivariate CvM-distance shares some appealing properties of the energy distance while being robust to heavy-tailed
distributions or outliers. For example, $W_d$ satisfies the characteristic property (Lemma~\ref{Lemma: Multivaraite CvM Characteristic Property}) and is invariant to orthogonal transformations. More importantly, it is straightforward to estimate $W_d$ without using any tuning parameters (Theorem~\ref{Theorem: Multivaraite CvM Expression}). Based on an unbiased estimate of $W_d^2$, we apply the permutation test procedure to determine a critical value of the test statistic. Although the permutation approach has been standard in practical implementations of two-sample testing, its theoretical properties have been less explored beyond simple cases \citep[e.g.][]{pesarin2001multivariate}. Indeed, previous studies usually consider asymptotic tests in their theory section whereas their actual tests are calibrated via permutations. We bridge the gap between theory and practice by presenting both theoretical and empirical results on the permutation test under various scenarios. Our main results regarding the CvM-distance are summarized as follows:

\begin{itemize}
	\itemsep 1pt
	\parskip 0pt
	\item \textbf{Closed form expression} (Section~\ref{Section: Projection Averaging-Type Cramer-von Mises Statistics}):~Building on \cite{escanciano2006consistent} and \cite{zhu2017projection}, we show that the test statistic has a simple closed-form expression.
	\item \textbf{Asymptotic power} (Section~\ref{Section: Projection Averaging-Type Cramer-von Mises Statistics}):~We prove that the permutation test based on the proposed statistic has the same asymptotic power as the oracle test against fixed and contiguous alternatives.
	\item \textbf{Robustness} (Section~\ref{Section: Robustness}):~We show that the permutation test based on the proposed statistic maintains good power in the contamination model, while the energy test becomes completely powerless in this setting.
	\item \textbf{Minimax optimality} (Section~\ref{Section: Minimax Optimality}):~We analyze the finite-sample power of the proposed permutation test and prove its minimax rate optimality against a class of alternatives that differ from the null in terms of the CvM-distance. We also show that the energy test is not optimal in our context.
	\item \textbf{HDLSS behavior} (Section~\ref{Section: High Dimension, Low Sample Size Analysis}):~We consider a \emph{high-dimension, low-sample size} (HDLSS) regime where the dimension tends to infinity while the sample size is fixed. Under this regime, we establish sufficient conditions under which the power of the proposed test converges to one. In addition, we show that the proposed test has comparable power to the high-dimensional mean tests introduced by \cite{chen2010two} and \cite{chakraborty2017tests} under certain location models.
	\item \textbf{Angular distance} (Section~\ref{Section: Connection to Generalized Energy Distance and MMD}):~We introduce the angular distance between two vectors and use this to show that the multivariate CvM-distance is a special case of the generalized energy
	distance~\citep{sejdinovic2013equivalence}. Furthermore, the CvM-distance is the maximum mean discrepancy~\citep{gretton2012kernel} associated with the angular distance.
\end{itemize}

Beyond the CvM-statistic, the projection-averaging technique can be widely applicable to other nonparametric statistics. In the second part of this study, we revisit some famous univariate sign- or rank-based statistics and propose their multivariate counterparts via projection-averaging. Although there has been much effort to extend univariate sign- or rank-based statistics in a multivariate space \citep[see e.g.][]{hettmansperger1998affine,oja2004multivariate,liu2006data,oja2010multivariate}, they are either computationally expensive to implement or less intuitive to understand. Our projection-averaging approach addresses these issues by providing a tractable calculation form of statistics and by having a direct interpretation in terms of projections. In Section~\ref{Section: Multivariate Extensions via Projection Averaging}, we demonstrate the generality of the projection-averaging approach by presenting multivariate extensions of several existing univariate statistics.

\subsection{Literature review}
There are a number of multivariate two-sample testing procedures available in the literature. We list some fundamental methods and recent developments. \cite{anderson1994two} proposed the two-sample statistic based on the integrated square distance between two kernel density estimates. The energy statistic was introduced by \cite{baringhaus2004new} and \cite{szekely2004testing} independently. \cite{biswas2014nonparametric} modified the energy statistic to improve the performance of the previous test for the high-dimensional location-scale and scale problems. \cite{gretton2012kernel} introduced a class of distances between two probability distributions, called the maximum mean discrepancy (MMD), based on a reproducing kernel Hilbert approach. \cite{sejdinovic2013equivalence} showed that the energy distance is a special case of the MMD associated with the kernel induced by the Euclidean distance. Recently, \cite{pan2018ball} proposed a new metric, named the ball divergence, between two probability distributions and connected it to the MMD. A further review of kernel-based two-sample tests can be found in \cite{harchaoui2013kernel}.

Another line of work is based on graph constructions. \cite{schilling1986multivariate} and \cite{henze1988multivariate} introduced a multivariate two-sample test based on the $k$ nearest neighbor (NN) graph. \cite{mondal2015high} pointed out that the previous NN test may suffer from low power for the high-dimensional location-scale problem and provided an alternative that addresses this limitation. Another variant of the NN test, which is tailored to imbalanced samples, can be found in \cite{chen2013ensemble}. \cite{friedman1979multivariate} considered minimum spanning tree (MST) to present a generalization of the univariate run test in \cite{wald1940test}. The MST test proposed by \cite{friedman1979multivariate} has recently been modified by \cite{chen2017new} and \cite{chen2018weighted} to improve power under scale alternatives and imbalanced samples, respectively. \cite{rosenbaum2005exact} proposed a distribution-free test in finite samples based on cross-matches. More recently, \cite{biswas2014distribution} introduced another distribution-free test based on the shortest Hamiltonian path. A general theoretical framework for graph-based tests has been established by \cite{bhattacharya2015distribution,bhattacharya2015general}. Other recent developments include \cite{liu2011triangle}, \cite{kanamori2012f}, \cite{bera2013smooth}, \cite{lopez2016revisiting}, \cite{zhou2017two}, \cite{mukhopadhyay2018nonparametric}, among others.

%The aforementioned studies have their own strengths, but also weaknesses. First, most of them employ asymptotic tests to theoretically justify their approach. It is not trivial and largely unknown whether the permutation tests possess the same properties as the asymptotic tests for large samples. Second, some of them require the choice of tuning parameters (e.g. Gaussian kernel bandwidth, the number of nearest neighbors) for their implementation. This makes the testing procedure more complicated. Third, the previous work mainly focus on light-tailed distributions such as normal distributions. The validity of their approach to heavy-tailed data has been neglected. Lastly, certain types of graph-based tests require an algorithm with high computational complexity. We tackle these issues with the proposed test statistic and demonstrate its efficacy with both theoretical and empirical results.

The projection-averaging approach to CvM-type statistics can be found in other statistical problems. For example, \cite{zhu1997estimated} and \cite{cui2002average} considered the CvM-statistic using projection-averaging to investigate one-sample goodness-of-fit tests for multivariate distributions. \cite{escanciano2006consistent} proposed the CvM-based goodness-of-fit test for parametric regression models. To the best of our knowledge, however, this is the first study that investigates the CvM-statistic for the multivariate two-sample problem via projection-averaging. 

Our technique to obtain a closed-form expression for projection-averaging statistics is based on  \cite{escanciano2006consistent}. The same principle has been exploited by \cite{zhu2017projection} in the context of testing for multivariate independence. We further extend the result of \cite{escanciano2006consistent} to more general cases and provide an alternative proof using orthant probabilities for normal distributions.

\vskip 1em

\paragraph{Outline.} The rest of this paper is organized as follows. In Section~\ref{Section: Projection Averaging-Type Cramer-von Mises Statistics}, we introduce our test statistic and the permutation test procedure. We then study their limiting behaviors under the conventional fixed dimension asymptotic framework. In Section~\ref{Section: Robustness}, we compare the power of the CvM test with that of the energy test and highlight the robustness of the CvM test. Section~\ref{Section: Minimax Optimality} establishes minimax rate optimality of the proposed test against a certain class of alternatives associated with the CvM-distance. In Section~\ref{Section: High Dimension, Low Sample Size Analysis}, we study the asymptotic power of the CvM test in the HDLSS setting. We introduce the angular distance between two vectors in Section~\ref{Section: Connection to Generalized Energy Distance and MMD} to show that the CvM-distance is the generalized energy distance built on the introduced metric. In Section~\ref{Section: Multivariate Extensions via Projection Averaging}, the projection-averaging technique is applied to other sign- or rank-based statistics and this allows us to provide new multivariate extensions. Simulation results are reported in Section~\ref{Section: Simulations} to demonstrate the competitive power performance of the proposed approach with finite sample size. All proofs not contained in the main text are in the supplementary material.

\vskip 1em

\paragraph{Notation.} For $U_1,U_2 \in \mathbb{R}^d$, we denote the angle between $U_1$ and $U_2$ by 
$\mathsf{Ang}(U_1,U_2) = \text{arccos} \big\{ U_1^\top U_2 / (\|U_1\| \| U_2\|) \big\}$ where the symbol $\top$ stands for the transpose operation. For $1 \leq q \leq p$, we let $(p)_q = p(p-1) \cdots (p-q+1)$. Let $\mP_0$ and $\mP_1$ be the probability measures under $H_0$ and $H_1$, respectively. Similarly $\mE_0$ and $\mE_1$ stand for the expectations with respect to $\mP_0$ and $\mP_1$. For any two real sequences $\{a_n\}$ and $\{b_n\}$, we use $a_n \asymp b_n$ if there exist constants $C,C^\prime > 0$ such that $C < |a_n / b_n | < C^\prime$ for large $n$. We write $a_n = O(b_n)$ if there exists $C > 0$ such that $|a_n| \leq C |b_n|$ for large $n$. For any given $c>0$, if $|a_n| \leq c |b_n|$ holds for large $n$, we write $a_n = o(b_n)$. For a sequence of random variables $X_n$, we write $X_n = O_{\mP}(a_n)$ if, for any $\epsilon >0$, there exists $M>0$ such that  $\mP(|X_n/a_n| > M) < \epsilon$ for large $n$. The acronym \emph{i.i.d.} stands for independent and identically distributed and we use the symbol $X_1,\ldots,X_n \overset{i.i.d.}{\sim} P$ to represent that $X_1,\ldots,X_n$ are $i.i.d.$ samples from distribution $P$. We denote the $d \times d$ identity matrix by $I_d$. The symbol $\ind(\cdot)$ is used for indicator functions. We write summation over the set of all $k$-tuples drawn without replacement from $\{1,\ldots,n\}$ by $\sum_{i_1,\ldots,i_k = 1}^{n,\neq}.$ Throughout this paper, we assume that all vectors are column vectors and $m,n \geq 2$.

\vskip 1em

\section{Projection Averaging-Type Cram{\'e}r-von Mises Statistics}  \label{Section: Projection Averaging-Type Cramer-von Mises Statistics}

In this section, we start with the basic properties of the CvM-distance. We then introduce our test statistic and study its limiting behavior. We end this section with a description of the permutation test and its large sample properties. Throughout this section, we consider the conventional asymptotic regime where the dimension is fixed and
%In this section, we start with the basic properties of the CvM-distance and introduce our test statistic. We then investigate the limiting behaviors of the proposed statistic as well as permutation test. Throughout this section (and only this section), we consider the conventional asymptotic regime where the dimension is fixed and
\begin{align} \label{Eq: Limit of Sample Proportion}
\frac{m}{m+n} \rightarrow \vartheta_X \in(0,1) \quad \text{and} \quad \frac{n}{m+n} \rightarrow \vartheta_Y \in (0,1) \quad \text{as} \quad  N=m+n \rightarrow \infty.
\end{align}

Let us first establish the characteristic property of the CvM-distance, meaning that $W_d$ is nonnegative and equal to zero if and only if $P_X = P_Y$. %As a consequence, the resulting test is consistent against all fixed alternatives.

\begin{lemma}\label{Lemma: Multivaraite CvM Characteristic Property}
	$W_d$ is nonnegative and has the characteristic property:
	\begin{align*}
	W_d(P_X, P_Y) = 0 \quad \text{if and only if} \quad P_X = P_Y.
	\end{align*}
\end{lemma}

Note that $W_d$ involves integration over the unit sphere. One way to approximate this integral is to consider a subset of $\mathbb{S}^{d-1}$, namely  $\{\beta_1,\ldots,\beta_k\}$, and then to take the sample mean over $k$ different univariate CvM-statistics \citep[see e.g.][]{zhu1997estimated}. However, this approach has a clear trade-off between accuracy and computational time depending on the choice of $k$. Our approach does not suffer from this issue by explicitly calculating the integral over $\mathbb{S}^{d-1}$. The explicit form of the integration is mainly due to \cite{escanciano2006consistent} who provided the following lemma:
\begin{lemma}\citep{escanciano2006consistent} \label{Lemma: Integration over Unit Sphere (2 terms)}
	For any two non-zero vectors $U_1,U_2 \in \mathbb{R}^{d}$,
	\begin{align*}
	\int_{\mathbb{S}^{d-1}} 
	\ind(\beta^\top U_1 \leq 0) \ind(\beta^\top U_2 \leq 0) d \lambda(\beta) = \frac{1}{2} - \frac{1}{2\pi}\mathsf{Ang}\left( U_1, U_2 \right).
	\end{align*}
%	\begin{proof}
%		The original paper proved this result using the volume of a spherical wedge. We provide an alternative proof by connecting this problem to the calculation of orthant probabilities for normal distributions. Suppose that $Z \in \mathbb{R}^d$ is normally distributed with mean zero and identity covariance, i.e. $Z \sim N(0,I_d)$. Then for given $U_1,U_2$, we observe that $(Z^\top U_1, Z^\top U_2)$ follows a bivariate normal distribution with mean zero and correlation matrix $R_{ij}$ where $R_{11}=R_{22}=1$ and $R_{12} = R_{21} = U_1^\top U_2 / \{ \|U_1\| \|U_2\|\}$. Also observe that $Z/\|Z\|$ is uniformly distributed on $\mathbb{S}^{d-1}$. Thus $\ind(\beta^\top U_1 \leq 0) \ind(\beta^\top U_2 \leq 0)$ is equally distributed as $\ind(Z^\top U_1 \leq 0) \ind(Z^\top U_2 \leq 0)$. The explicit formula for the expectation of the latter quantity has been well-established \citep[e.g.][]{childs1967reduction}. Using this result, we can complete the proof. The same trick can be applied to a more general integration involving more than two identity functions. See Lemma~\ref{Lemma: Extension of Escanciano (2006) with three arguments} and the supplementary material. 
%	\end{proof}
\end{lemma}

\vskip .5em 

\begin{remark} \label{Remark: Extension}
	\cite{escanciano2006consistent} proved Lemma~\ref{Lemma: Integration over Unit Sphere (2 terms)} using the volume of a spherical wedge. In the supplementary material, we provide an alternative proof of this result based on orthant probabilities for normal distributions. We also extend this result to integration involving three or more than three indicator functions in Lemma~\ref{Lemma: Extension of Escanciano (2006) with three arguments} and the supplementary material, respectively.
\end{remark}

\vskip .5em

%\begin{remark} \label{Remark: Continuity Assumption}
%	We would like to note that the arccosine is indeterminate when either $U_1$ or $U_2$ is a zero vector. To avoid this issue, we simply assume that each of $X$ and $Y$ contains at least one continuous random variable so that $\beta^\top X$ and $\beta^\top Y$ have continuous distribution for all $\beta \in \mathbb{S}^{d-1}$ almost surely.
%\end{remark}

%\begin{remark} \label{Remark: Continuity Assumption}
%	We would like to note that the arccosine is indeterminate when either $U_1$ or $U_2$ is a zero vector. To avoid this issue, we simply assume that $\beta^\top X$ and $\beta^\top Y$ have continuous distribution for all $\beta \in \mathbb{R}^{d-1}$ $\lambda$-almost surely for the rest of this paper.	
%\end{remark}

Based on Lemma~\ref{Lemma: Integration over Unit Sphere (2 terms)}, we give another representation of $W^2_d$ in terms of the expected angle involving three independent random vectors. Here and hereafter, we assume that 
\begin{align} \label{Eq: continuity assumption}
\text{\emph{$\beta^\top X$ and $\beta^\top Y$ have continuous distribution functions for $\lambda$-almost all $\beta \in \mathbb{S}^{d-1}$.}}
\end{align}
This continuity assumption greatly simplifies the alternative expression for $W_d^2$ and avoids the possibility that $\textsf{Ang}(\cdot,\cdot)$ is not well-defined when one of the inputs is a zero vector. This issue may be handled by defining $\textsf{Ang}(\cdot,\cdot)$ differently for those exceptional cases, but we do not pursue this direction here. 
%The corresponding result without the continuity assumption can be found in the supplementary material. 

\begin{theorem}[Closed form expression] \label{Theorem: Multivaraite CvM Expression}
	Suppose that $X_1,X_2 \overset{i.i.d.}{\sim} P_X$ and, independently, $Y_1,Y_2 \overset{i.i.d.}{\sim} P_Y$. Then the squared multivariate CvM-distance can be written as
	\begin{align*}
	W_d^2(P_X,P_Y) ~=~ \frac{1}{3}  - \frac{1}{2\pi}\mE \left[ \mathsf{Ang} \left( X_1 - Y_1, X_2 -Y_1 \right)   \right] - 
	\frac{1}{2\pi}\mE \left[ \mathsf{Ang} \left( Y_1- X_1, Y_2 - X_1 \right) \right].
	\end{align*}
	\begin{proof}
		After expanding the square term in $W_d^2$, we may get several pieces including 
		\begin{align*}
		\vartheta_Y \int_{ \mathbb{S}^{d-1}} \int_{\mathbb{R}} \big( F_{\beta^\top X}(t) \big)^2 dF_{\beta^\top Y}(t) d\lambda(\beta).
		\end{align*}
		By Fubini's theorem, the above term can be written as
		\begin{align*}
		\vartheta_Y \mE\bigg[\int_{ \mathbb{S}^{d-1}} \ind \big\{ \beta^\top(X_1-Y_1) \leq 0 \big\} \ind \big\{\beta^\top(X_2-Y_1) \leq 0 \big\} d\lambda(\beta) \bigg].
		\end{align*}
		We then apply Lemma~\ref{Lemma: Integration over Unit Sphere (2 terms)} to have an expression that involves the angle between $X_1-Y_1$ and $X_2-Y_1$. Applying the same principle to the other terms and simplifying them by using the continuity assumption, we may obtain the desired expression. The details can be found in the supplementary material. 
	\end{proof}
\end{theorem}

\vskip .8em

\begin{remark}
	Theorem~\ref{Theorem: Multivaraite CvM Expression} highlights that $W_d(P_X,P_Y)$ is invariant to the choice of $\vartheta_X$ and $\vartheta_Y$ under the continuity assumption (\ref{Eq: continuity assumption}).
\end{remark}

\vskip 1em

\subsection{Test Statistic and Limiting Distributions}
Theorem~\ref{Theorem: Multivaraite CvM Expression} leads to a natural empirical estimate of $W^2_d$ based on a $U$-statistic. Consider the kernel of order two:
\begin{equation}
\begin{aligned} \label{Eq: U-stat kernel}
h_{\text{CvM}}(x_1,x_2; y_1,y_2) ~= ~ \frac{1}{3} - \frac{1}{2\pi} \mathsf{Ang} \left( x_1-y_1, x_2 - y_1 \right) - 
\frac{1}{2\pi} \mathsf{Ang} \left( y_1-x_1, y_2 - x_1 \right).
\end{aligned}
\end{equation}
Then we define our test statistic as follows:
\begin{align}  \label{Eq: Definition of U-statistic}
U_{\text{CvM}} & = \frac{1}{(m)_2 (n)_2} \sum_{i_1,i_2 =1}^{m,\neq} \sum_{j_1,j_2 =1}^{n,\neq} h_{\text{CvM}}(X_{i_1},X_{i_2}; Y_{j_1},Y_{j_2}).
\end{align}

Leveraging the basic theory of $U$-statistics \citep[e.g.][]{lee1990u}, it is clear that $U_{\text{CvM}}$ is an unbiased estimator of $W_d^2$. Additionally, $U_{\text{CvM}}$ is a degenerate $U$-statistic under the null hypothesis as we proved in the supplementary material. Hence we can apply the asymptotic theory for a degenerate two-sample $U$-statistic \citep[Chapter 3 of][]{bhat1995theory} to obtain the following result.

%In the next theorem, we explore the asymptotic distribution of the proposed test statistic under the null hypothesis. %Under the null hypothesis, the proposed test statistic is a degenerate $U$-statistic, which converges weakly to a mixture of chi-square random variables.

\begin{theorem}[Asymptotic null distribution of $U_{\text{CvM}}$] \label{Theorem: Asymptotic Null Distribution}
	Let $\lambda_{k}$ be the eigenvalue with the corresponding eigenfunction $\phi_k$ satisfying the integral equation
	\begin{align} \label{Eq: Eigenvalue and Eigenfucntion}
	\mE \Big\{ \mE \Big[ \widetilde{h}_{\text{\emph{CvM}}}(x_1,X_2;Y_1,Y_2) \big| X_2 \Big] \phi_k(X_2) \Big\} = \lambda_{k} \phi_k(x_1) \quad \text{for} \ k = 1,2,\ldots,
	\end{align}
	where $\widetilde{h}_{\text{\emph{CvM}}}(x_1,x_2;y_1,y_2) = h_{\text{\emph{CvM}}}(x_1,x_2;y_1,y_2)/2 + h_{\text{\emph{CvM}}}(x_2,x_1;y_2,y_1)/2$. Then $U_{\text{\emph{CvM}}}$ has the limiting null distribution under the limiting regime (\ref{Eq: Limit of Sample Proportion}) given by
	\begin{align*}
	& N  U_{\text{\emph{CvM}}} \convD  \vartheta_X^{-1} \vartheta_Y^{-1} \sum_{k=1}^\infty \lambda_k (\xi_k^2 - 1), 
	\end{align*}
	where $\xi_k \overset{i.i.d.}{\sim} N(0,1)$ and $\convD$ stands for convergence in distribution.
\end{theorem}

\vskip .8em

\begin{remark}
	The eigenvalues $\{\lambda_i\}_{i=1}^\infty$ may depend on the underlying distribution, which implies that the test statistic is not distribution-free even asymptotically. Nevertheless, for the univariate continuous case, explicit expressions for the eigenvalues and the eigenfunctions are available as $\lambda_i = 2/(i\pi)^2$ and $\phi_i(x) = \sqrt{2} \text{cos}(i \pi x)$ for $i=1,2,\ldots$ \citep[e.g.][]{chikkagoudar2014limiting}.
\end{remark}

\vskip .8em

Under a fixed alternative hypothesis where $P_X$ and $P_Y$ do not change with $m$ and $n$, the proposed test statistic converges weakly to a normal distribution. We build on Hoeffding's decomposition of a two-sample $U$-statistic \citep[e.g. page 40 of][]{lee1990u} to prove the following result.

\begin{theorem}[Asymptotic distribution of $U_{\text{CvM}}$ under fixed alternatives] \label{Theorem: Asymptotic distribution under fixed alternatives}
	Let us define
	\begin{align*}
	& \sigma_{h_X}^2 = \mV \Big\{ \mE \Big[ \widetilde{h}_{\text{\emph{CvM}}}(X_1,X_2;Y_1,Y_2) \big| X_1 \Big] \Big\}, \\[.5em]
	& \sigma_{h_Y}^2 = \mV \Big\{ \mE \Big[ \widetilde{h}_{\text{\emph{CvM}}}(X_1,X_2;Y_1,Y_2) \big| Y_1 \Big] \Big\}.
	\end{align*}
	Then under the limiting regime (\ref{Eq: Limit of Sample Proportion}) and fixed alternative $P_X \neq P_Y$, we have
	\begin{align*}
	\sqrt{N} (U_{\text{\emph{CvM}}} - W_d^2) \convD N\left(0, 4 \vartheta_X^{-1} \sigma_{h_X}^2 + 4\vartheta_Y^{-1} \sigma_{h_Y}^2\right).
	\end{align*}
	%\begin{proof}
	%	 Under the fixed alternative, $U_{\text{\emph{CvM}}}$ is a non-degenerate $U$-statistic. Hence the result follows by Theorem 2 in page 165 of \cite{kowalski2008modern}.
	%\end{proof}
\end{theorem}

\vskip .8em

%From the previous results, it is seen that $N U_{\text{CvM}}$ is stochastically bounded under the null hypothesis whereas it diverges to infinity under the alternative hypothesis. Thus any reasonable testing procedure based the proposed test statistic is consistent against all fixed alternatives. 

%The large-sample power of many reasonable tests tend to one against any fixed alternative, which may be of less interest. 

The problem of distinguishing two fixed distributions becomes too easy in large sample situations and may be of less interest. We therefore turn now to a more challenging scenario where a distance between $P_X$ and $P_Y$ diminishes as the sample size increases. To this end, we make a standard assumption that the underlying distributions belong to quadratic mean differentiable (QMD) families \citep[e.g.][]{bhattacharya2015general}.

%Many nonparametric tests have been shown to be pointwise consistent in power against any fixed alternative, which may be of less interest. We therefore turn now to a more challenging scenario where a distance between $P_X$ and $P_Y$ diminishes as the sample size increases. To this end, we make a standard assumption that the underlying distributions belong to the quadratic mean differentiable (QMD) families \citep[e.g.][]{bhattacharya2015general}. \textcolor{red}{(revised)}

%In particular, we focus on a (local) contiguous alternative

%We therefore turn now to a (local) contiguous alternative where the distance between $P_X$ and $P_Y$ diminishes as the sample size increases. In order to invoke

%In order to obtain a nontrivial asymptotic power, we 

%so that it becomes much more difficult to distinguish the alternative from the null. In particular, we focus on a class of contiguous alternatives 

%In the following, we study the large sample behavior of the test statistic under a (local) contiguous alternative where the distance between $P_X$ and $P_Y$ diminishes as the sample size increases. To this end, we make a standard assumption that the underlying distributions belong to the quadratic mean differentiable (QMD) families \citep[e.g.][]{bhattacharya2015general}. 

\begin{definition}\citep[Quadratic Mean Differentiable Families, page 484 of][]{lehmann2006testing}
	Let $\{P_{\theta}, \theta \in \Omega \}$ be a family of probability distributions on $(\mathbb{R}^d, \mathcal{B})$ where $\mathcal{B}$ is the Borel $\sigma$-field associated with $\mathbb{R}^d$. Assume each $P_{\theta}$ is absolutely continuous with respect to Lebesgue measure and set $p_\theta(t) = dP_{\theta}(t) /dt$. The family $\{ P_{\theta},\theta \in \Omega\}$ is quadratic mean differentiable at $\theta_0$ if there exists a vector of real-valued functions $\eta(\cdot,\theta_0) = (\eta_1(\cdot,\theta_0),\ldots,\eta_k(\cdot,\theta_0))^\top$ such that 
	\begin{align*}
	\int_{\mathbb{R}^d} \Big[ \sqrt{p_{\theta_0+b}(t)} - \sqrt{p_{\theta_0}(t)} - \langle \eta(t,\theta_0),b \rangle  \Big]^2 d t = o(\|b\|^2),
	\end{align*}
	as $\|b\| \rightarrow 0$.
\end{definition}

\vskip .8em

The QMD families include a broad class of parametric distributions such as exponential families in natural form. By focusing on the QMD families, we are particularly interested in asymptotically non-degenerate situations where the limiting sum of the type I and type II errors of the optimal test is non-trivial, i.e.~bounded by zero and one. It has been shown that when $P_{\theta_0}$ and $P_{\theta_N}$ belong to the QMD families, this non-degenerate situation occurs when $\|\theta_0 - \theta_N \| \asymp N^{-1/2}$ \citep[Chapter 13.1 of][]{lehmann2006testing}. Hence, we consider a sequence of contiguous alternatives where $\theta_N = \theta_0 + b N^{-1/2}$ for some $b \in \mathbb{R}^k$ and establish the asymptotic behavior of $U_{\text{CvM}}$ under the given scenario. Our result builds on the prior work by \cite{chikkagoudar2014limiting} and extends it to multivariate cases.

%Under this sequence of contiguous alternatives, we establish the asymptotic behavior of $U_{\text{CvM}}$. 
%Our result builds on the prior work by \cite{chikkagoudar2014limiting} and extends it to multivariate cases. \textcolor{red}{(revised)}

\begin{theorem}[Asymptotic distribution of $U_{\text{CvM}}$ under contiguous alternatives] \label{Theorem: Asymptotic distribution under contiguous alternatives}
	Assume $\{P_{\theta}, \theta \in \Omega \}$ is quadratic mean differentiable at $\theta_0$ with derivative $\eta(\cdot,\theta_0)$ and $\Omega$ is an open subset of $\mathbb{R}^k$. Define the Fisher Information matrix to be the matrix $I(\theta)$ with $(i,j)$ entry 
	\begin{align*}
	I_{i,j}(\theta) = 4 \int_{\mathbb{R}^d} \eta_i(t,\theta) \eta_j(t,\theta) dt,
	\end{align*}
	and assume that $I(\theta)$ is nonsingular. Suppose we observe $\mathcal{X}_m \overset{i.i.d.}{\sim} P_{\theta_0}$ and $\mathcal{Y}_n \overset{i.i.d.}{\sim} P_{\theta_0+bN^{-1/2}}$ for $b \in \mathbb{R}^k$. Then under the limiting regime (\ref{Eq: Limit of Sample Proportion}), 
	\begin{align*}
	& N  U_{\text{\emph{CvM}}} \convD  \vartheta_X^{-1} \vartheta_Y^{-1} \sum_{k=1}^\infty \lambda_k \{(\xi_k+ \vartheta_X^{1/2}a_k)^2 - 1 \},
	\end{align*}
	where
	\begin{align*}
	a_k = \int_{\mathbb{R}^d}  \big\langle b, 2\eta(x,\theta_0)p_{\theta_0}^{-1/2}(x) \big\rangle \phi_k(x) dP_{\theta_0}(x).
	\end{align*}
	\begin{proof}
		We provided a more general result in Lemma~\ref{Lemma: Contiguous alternative} and this is a direct consequence of Lemma~\ref{Lemma: Contiguous alternative} with $r=2$. 
	\end{proof}
\end{theorem}

\begin{remark}
	As can be seen by putting $b=0$, Theorem~\ref{Theorem: Asymptotic Null Distribution} is a special case of Theorem~\ref{Theorem: Asymptotic distribution under contiguous alternatives} for the QMD families. Theorem~\ref{Theorem: Asymptotic distribution under contiguous alternatives} also shows that if there exists $k \geq 1$ such that $a_k \neq 0$  and $\lambda_k >0$, the oracle test and the permutation test considered later in Theorem~\ref{Theorem: Asymptotic Power of Permutation Tests} have asymptotic power greater than $\alpha$ \citep[see, page 615 of][]{lehmann2006testing}.
\end{remark}
%\begin{remark}
%	Another statistic can be based on a $V$-statistic (see Lemma~\ref{Lemma: V-statistic as a plug-in estimator} in the appendix). Although the $V$-statistic can be a more natural extension of the univariate CvM-statistic, we prefer to use the $U$-statistic since it is an unbiased estimator of $W_d^2$ and it avoids the issue that the arccosine is indeterminate.
%\end{remark}
\vskip 1em

\subsection{Critical Value and Permutation Test} \label{Section: Critical Value and Permutation Test}
%Unfortunately, when $d \geq 2$, the null distribution of $U_{\text{CvM}}$ depends on the underlying distributions $P_X$ and $P_Y$; therefore, the resulting test is not distribution-free. 
%A common approach to determining the critical value of a test is to compute the quantile of the limiting distribution under the null. However, there are several drawbacks associated with this asymptotic approach in our setting. First, the limiting null distribution of $U_{\text{CvM}}$ depends on an infinite number of the nuisance parameters that need to be estimated in multivariate cases. Second, $U_{\text{CvM}}$, which is a degenerate $U$-statistic, may have different limiting behaviors in different asymptotic regimes \citep[e.g.][]{hall1984central}. This means that the validity of the test is questionable even when the sample size goes to infinity. A permutation approach, on the other hand, does not suffer from these issues and presents a valid $p$-value for any finite sample size  \citep[e.g.~Chapter 15 of][]{lehmann2006testing}. 

We next describe the permutation test based on $U_{\text{CvM}}$ and examine its large sample properties under the conventional asymptotic regime. Let us start by introducing the oracle test and then compare it to the permutation test. Suppose that the mixture distribution $\vartheta_X P_X + \vartheta_Y P_Y$ is known. Then the critical value of the oracle test can be defined as follows:

%In this subsection, we first describe the permutation test based on $U_{\text{CvM}}$ and then investigate its limiting properties under the conventional asymptotic regime. To define a test, one could also use the limiting null distribution studied in Theorem~\ref{Theorem: Asymptotic Null Distribution}. However, this asymptotic approach is undesirable in our setting for several reasons. First, the limiting null distribution of $U_{\text{CvM}}$ depends on an infinite number of the nuisance parameters that need to be estimated in multivariate cases. Second, $U_{\text{CvM}}$, which is a degenerate $U$-statistic, may have different limiting behaviors in different asymptotic regimes \citep[e.g.][]{hall1984central}. This means that the validity of the test is questionable even when the sample size goes to infinity. A permutation approach, which is the main focus of this paper on the other hand, does not suffer from these issues and presents a valid $p$-value for any finite sample size  \citep[e.g.~Chapter 15 of][]{lehmann2006testing}. 

%We first introduce the oracle test and compare it to the permutation test. Suppose that the mixture distribution $\vartheta_X P_X + \vartheta_Y P_Y$ is known where $m/N \rightarrow \vartheta_X$ and $n/N \rightarrow \vartheta_Y$ as $N \rightarrow \infty$. Then the critical value of the oracle test can be decided as follows:
\paragraph{$\bullet$ Oracle Test}
\begin{enumerate}
	\item Consider new $i.i.d.$ samples $\{\widetilde{Z}_1,\ldots,\widetilde{Z}_{N}\}$ from the mixture $\vartheta_X P_X + \vartheta_Y P_Y$.
	\item Let $T_{m,n}(\widetilde{Z})$ be the test statistic of interest calculated based on $\widetilde{\mathcal{X}}_m = \{\widetilde{Z}_1,\ldots,\widetilde{Z}_{m}\}$ and $\widetilde{\mathcal{Y}}_n = \{\widetilde{Z}_{m+1},\ldots, \widetilde{Z}_{N}\}$. 
	\item Given a significance level $0< \alpha <1$, return the critical value $c^\ast_{\alpha, m,n}$ defined by
	\begin{align} \label{Eq: Critical Value of Oracle Test}
	& c^\ast_{\alpha,m,n} := \inf \Big\{ t \in \mathbb{R} :  1-\alpha \leq \mP \Big( T_{m,n}(\widetilde{Z}) \leq t \Big) \Big\}.
	\end{align}
\end{enumerate}

\begin{remark}
	It is worth pointing out that $T_{m,n}(\widetilde{Z})$ has the same distribution as the test statistic based on the original samples under $H_0$, but not necessarily under $H_1$. Hence the oracle test based on $c^\ast_{\alpha,m,n}$ is exact under $H_0$ and can be powerful under $H_1$. 
\end{remark}

The critical value of the permutation test can be obtained without knowledge of the mixture distribution $\vartheta_X P_X + \vartheta_Y P_Y$ as follows:

\paragraph{$\bullet$ Permutation Test}
\begin{enumerate}
	\item Let $\{Z_1,\ldots,Z_{N}\} = \{ X_1,\ldots,X_m,Y_1,\ldots,Y_n\}$ be the pooled samples and $Z_{\varpi} = \{Z_{\varpi(1)},\ldots,\\ Z_{\varpi(N)}\}$ where $\varpi = \{ \varpi(1),\ldots,\varpi(N)\}$ is a permutation of $\{1,\ldots,N\}$.
	\item Let $T_{m,n}(Z_\varpi)$ be the test statistic of interest calculated based on $\mathcal{X}_m^\varpi = \{Z_{\varpi(1)},\ldots,Z_{\varpi(m)}\}$ and $\mathcal{Y}_n^\varpi = \{Z_{\varpi(m+1)},\ldots, Z_{\varpi(N)}\}$. 
	\item Given a significance level $0< \alpha <1$, return the critical value $c_{\alpha, m,n}$ defined by
	\begin{align}  \label{Eq: Critical Value of Permutation Test}
	& c_{\alpha,m,n} := \inf \Big\{ t \in \mathbb{R} :  1-\alpha \leq \frac{1}{N!} \sum_{\varpi \in \mathcal{S}_N} \ind\Big( T_{m,n}(Z_\varpi) \leq t \Big) \Big\},
	\end{align}
	where $\mathcal{S}_N$ is the set of all permutations of $\{1,\ldots, N\}$. 
\end{enumerate}

%In the next theorem, we show that the difference between $c^\ast_{\alpha,m,n}$ and $c_{\alpha,m,n}$ of the proposed statistic is asymptotically negligible under both the null and alternative hypotheses. From this result (combined with Theorem~\ref{Theorem: Asymptotic Null Distribution} and Theorem~\ref{Theorem: Asymptotic distribution under fixed alternatives}), it can be easily seen that the asymptotic power of the oracle test and the permutation test are the same and converge to one against any fixed alternative. \textcolor{red}{(TO DO)}

In the next theorem, we show that the difference between $c^\ast_{\alpha,m,n}$ and $c_{\alpha,m,n}$ for the proposed statistic is asymptotically negligible under both the null and alternative hypotheses. In doing so, we develop a general asymptotic theory for the permutation distribution of a two-sample degenerate $U$-statistic under $H_0$. This general result is established based on Hoeffding's conditions~\citep{hoeffding1952large} and extended to $H_1$ via the coupling argument~\citep{chung2013exact}. The details can be found in Appendix~\ref{Section: Permutation Tests}.

\begin{theorem}[Asymptotic behavior of the critical values] \label{Theorem: Critical value of Permutation test}
	Consider the conventional limiting regime in (\ref{Eq: Limit of Sample Proportion}). Let $c^\ast_{\alpha,\text{\emph{CvM}}}$ and $c_{\alpha,\text{\emph{CvM}}}$ be the critical values of the oracle test and the permutation test based on the scaled CvM-statistic, that is $N U_{\text{\emph{CvM}}}$, as described in (\ref{Eq: Critical Value of Oracle Test}) and (\ref{Eq: Critical Value of Permutation Test}), respectively. Then under both the null and (fixed or contiguous) alternative hypotheses,
	\begin{align*}
	c^\ast_{\alpha,\text{\emph{CvM}}} - c_{\alpha,\text{\emph{CvM}}}  \convP 0.
	\end{align*}
	Here $\convP$ stands for convergence in probability.
\end{theorem}
%Beyond the consistency, we prove that the asymptotic power of the oracle test and the permutation test are the same under the contiguous alternative hypothesis.

Leveraging the previous result combined with Slutsky's theorem, we prove that the asymptotic power of the oracle test and the permutation test are identical against any fixed and contiguous alternatives. This clearly highlights an advantage of the permutation test as it is exact under $H_0$ and asymptotically as powerful as the oracle test under $H_1$. More importantly, the permutation test does not require any prior information on the underlying distributions. 

%Thus the permutation test is asymptotically as powerful as the oracle test and controls the type I error in any finite sample size, but it does not require  

\begin{theorem}[Asymptotic equivalence of power]\label{Theorem: Asymptotic Power of Permutation Tests}
	The oracle test and the permutation test control the type I error under the null hypothesis as
	\begin{align*}
	\mP_0\big( N U_{\text{\emph{CvM}}} > c^\ast_{\alpha,\text{\emph{CvM}}} \big)  \leq \alpha \quad \text{and} \quad \mP_0 \big( N U_{\text{\emph{CvM}}} > c_{\alpha,\text{\emph{CvM}}} \big) \leq \alpha.
	\end{align*}
	On the other hand, under the fixed or contiguous alternative hypotheses considered in Theorem~\ref{Theorem: Asymptotic distribution under fixed alternatives} and Theorem~\ref{Theorem: Asymptotic distribution under contiguous alternatives}, we have that
	\begin{align*}
	\mP_1 \big( N U_{\text{\emph{CvM}}} > c^\ast_{\alpha,\text{\emph{CvM}}} \big) - \mP_1 \big( N U_{\text{\emph{CvM}}} > c_{\alpha,\text{\emph{CvM}}} \big)  \rightarrow 0 ~ \text{as} ~ N \rightarrow \infty.
	\end{align*}
\end{theorem}

\vskip 0.8em 

%We end this section by discussing a computational aspect of the permutation approach.

\begin{remark} \label{Remark: Monte Carlo Permutations}
	Except for small sample sizes, it may not be feasible to implement the permutation procedure as in (\ref{Eq: Critical Value of Permutation Test}) due to computational cost. A common approach to alleviate this computational issue is to use Monte Carlo sampling of random permutations and approximate the exact permutation $p$-value. In more detail, note first that the permutation test function can be written as $\ind(\widehat{p}_{\text{{CvM}}} \leq \alpha)$ where $\widehat{p}_{\text{{CvM}}}$ is the permutation $p$-value given by
	\begin{align*}
	\widehat{p}_{\text{{CvM}}} = \frac{1}{N!} \sum_{\varpi \in \mathcal{S}_N} \ind\{ U_{\text{{CvM}}}(Z_\varpi) \geq U_{\text{{CvM}}} \}.
	\end{align*}
	Let $\varpi^{(1)}, \ldots, \varpi^{(B)}$ be independent and uniformly distributed on $\mathcal{S}_N$. Then the Monte Carlo version of the permutation $p$-value is computed by
	\begin{align*}
	\widehat{p}_{\text{{CvM}}}^{(B)} = \frac{1}{B+1} \left[ \sum_{i = 1}^B \ind\{ U_{\text{{CvM}}}(Z_{\varpi^{(i)}}) \geq U_{\text{{CvM}}} \} + 1 \right].
	\end{align*}
	It is well-known that $\ind (\widehat{p}_{\text{{CvM}}}^{(B)} \leq \alpha)$ is also a valid level $\alpha$ test for any finite sample size and $\widehat{p}_{\text{{CvM}}} - \widehat{p}_{\text{{CvM}}}^{(B)} \convP 0$ as $B \rightarrow \infty$ \citep[e.g. page 636 of][]{lehmann2006testing}. Throughout this paper, we also adapt this approach for our simulation studies.
\end{remark}
% Lehmann Example 11.2.13

%\begin{remark}
%	In Appendix~\ref{Section: Permutation Tests}, we establish general conditions under which the permutation distribution based on a degenerate $U$-statistic is asymptotically equivalent to the corresponding unconditional distribution. This general result includes the energy statistic, Gaussian MMD and Proposition~\ref{Theorem: Critical value of Permutation test} as a corollary.
%\end{remark}

\section{Robustness} \label{Section: Robustness}

Recall that the energy distance and the CvM-distance can be represented by integrals of the $L_2^2$-type difference between two distribution functions. In view of this, the main difference between the energy distance and the CvM-distance is in their weight function. More precisely, the energy distance is defined with $dt$, which gives a uniform weight to the whole real line. On the other hand, the CvM-distance is defined with $dH_\beta(t)$, which gives the most weight on high-density regions. As a result, the test based on the CvM-distance is more robust to extreme observations than the one based on the energy distance. It is also important to note that the CvM-distance is well-defined without any moment conditions, whereas the energy distance is only well-defined assuming a finite first moment. When the moment condition is violated or there exist extreme observations, the test based on the energy distance may suffer from low power. The purpose of this section is to demonstrate this point both theoretically and empirically by using contaminated distribution models.

\subsection{Theoretical Analysis}

Suppose we observe samples from an $\epsilon$-contamination model:
\begin{align} \label{Eq: Contaminated Model}
X \sim P_{X,N}  := (1-\epsilon) Q_{X} + \epsilon G_{N} \quad \text{and} 
\quad Y \sim  P_{Y,N} := (1-\epsilon) Q_{Y} + \epsilon G_{N},
\end{align}
where $G_{N}$ can change arbitrarily with $N$ and $\epsilon \in (0,1)$. Suppose that $Q_X$ and $Q_Y$ are significantly different so that a given test has high power to distinguish between $Q_X$ and $Q_Y$ without contaminations. Then it is natural to expect that the power of the same test would not decrease much for the contamination model when $\epsilon$ is close to zero. In other words, an ideal test would maintain robust power against any choice of $G_{N}$ as long as $Q_X$ and $Q_Y$ are different and $\epsilon$ is small. Unfortunately, this is not the case for the energy test. As we shall see, for any arbitrary small (but fixed) $\epsilon$, there exists a heavy-tail contamination $G_{N}$ such that the energy test becomes asymptotically powerless under mild moment conditions for $Q_X$ and $Q_Y$. On the other hand, the CvM test is uniformly powerful over any choice of $G_{N}$ as sample size tends to infinity.

%With this model assumption, we would like to assess whether a test can maintain adequate power even when $G_{N}$ behave in favor of the null hypothesis \textcolor{red}{(revised)}. We mainly focus on statistical power to study robustness because one can always employ the permutation procedure to control the type I error under $H_0 : P_{X,N} = P_{Y,N}$.

Let us consider the energy statistic based on a $U$-statistic:
\begin{equation}
\begin{aligned} \label{Eq: Energy U-statistic}
U_{\text{Energy}} & =  \frac{2}{mn} \sum_{i=1}^{m} 
\sum_{j=1}^{n} \| X_{i} - Y_{j}\| - \frac{1}{(m)_2} \sum_{i_1,i_2=1}^{m,\neq} \| X_{i_1} - X_{i_2}\| \\
&  - \frac{1}{(n)_2} \sum_{j_1,j_2=1}^{n,\neq} \| Y_{j_1} - Y_{j_2}\|. 
\end{aligned}
\end{equation}
Then the main result of this subsection is stated as follows.

\begin{theorem}[Robustness under contaminations] \label{Theorem: Robustness}
	Suppose we observe samples $\mathcal{X}_m$ and $\mathcal{Y}_n$ from the contaminated model in (\ref{Eq: Contaminated Model}) with an arbitrary small but fixed contamination ratio $\epsilon$. Assume that $Q_X$ and $Q_Y$ are fixed but $Q_X \neq Q_Y$ while $N$ changes. In addition, assume that $Q_X$ and $Q_Y$ have their finite second moments. Consider the tests based on $U_{\text{\emph{CvM}}}$ and $U_{\text{\emph{Energy}}}$ given by
	\begin{align*} 
	\phi_{\emph{CvM}} := \ind(U_{\text{\emph{CvM}}} > c_{\alpha,\text{\emph{CvM}}}) \quad \text{and} \quad  \phi_{\emph{Energy}} := \ind(U_{\text{\emph{Energy}}} > c_{\alpha,\text{\emph{Eng}}}),
	\end{align*}
	where $c_{\alpha,\text{\emph{CvM}}}$ and $c_{\alpha,\text{\emph{Eng}}}$ are $\alpha$ level permutation critical values of $U_{\text{\emph{CvM}}}$ and $U_{\text{\emph{Energy}}}$ respectively. Then for any $(Q_X,Q_Y)$, there exists a certain $G_{N}$ such that the energy test becomes asymptotically powerless under the asymptotic regime in (\ref{Eq: Limit of Sample Proportion}). On the other hand, the CvM test is asymptotically powerful uniformly over all possible $G_{N}$, that is
	\begin{align}
	&  \lim_{m,n \rightarrow \infty} \inf_{G_{N}} \mE_1 \left[ \phi_{\emph{Energy}} \right] \leq \alpha ~~ \text{and} ~~  \lim_{m,n \rightarrow \infty} \inf_{G_{N}} \mE_1 \left[ \phi_{\emph{CvM}} \right] = 1 .
	\end{align}
	\begin{proof}
		We sketch the proof of the negative result for the energy test. The details can be found in the supplementary document. Assume that $G_{N}$ is a multivariate normal distribution with zero mean vector and 
		covariance matrix $\sigma_{N}^2 I_d$ where $\sigma_N^2 \in \mathbb{R}$ is a positive sequence that tends to infinity as $N \rightarrow \infty$. Let us define the truncated random vectors $\widetilde{X}$ and $\widetilde{Y}$ coupled with $X$ and $Y$ as
		\begin{align*}
		\widetilde{X} = 
		\begin{cases}
		(0,\ldots,0)^\top, \quad & \text{if} ~ X \sim Q_X, \\
		X / \sigma_N, \quad & \text{if} ~ X \sim G_N,
		\end{cases} 
		\quad \text{and} \quad
		\widetilde{Y} = 
		\begin{cases}
		(0,\ldots,0)^\top, \quad & \text{if} ~ Y \sim Q_Y, \\
		Y / \sigma_N, \quad & \text{if} ~ Y \sim G_N.
		\end{cases} 
		\end{align*}
		By the construction, it is clear that $\widetilde{X} $ and $\widetilde{Y}$ have the same mixture distribution as
		\begin{align*}
		\widetilde{X} ,\widetilde{Y}  \sim  \widetilde{P}:=(1-\epsilon) Q_{\delta_0} + \epsilon \widetilde{G},
		\end{align*}
		where $Q_{\delta_0}$ is the degenerate distribution at $(0,\ldots,0)^\top$ and $\widetilde{G}$ is the standard multivariate normal distribution, i.e. $N((0,\ldots,0)^\top, I_d)$. Now we consider the two energy statistics: one based on the original samples and the other based on the corresponding truncated samples. Denote these two statistics by $U_{\text{Energy}}$ and $\widetilde{U}_{\text{Energy}}$, respectively. In the supplementary material, we show that $N\sigma_N^{-1} U_{\text{Energy}}$ and $N \widetilde{U}_{\text{Energy}}$ are asymptotically the same under a certain choice of $\sigma_N^2$. We also show that these two statistics have the same permutation distribution in large sample scenarios. Since the power of the permutation test based on $N \widetilde{U}_{\text{Energy}}$ cannot exceed $\alpha$, this implies that the permutation test based on $N\sigma_N^{-1} U_{\text{Energy}}$ becomes asymptotically powerless. This completes the proof.
	\end{proof}
\end{theorem}

\begin{remark}
	In Theorem~\ref{Theorem: Robustness}, we made the assumption that $Q_X$ and $Q_Y$ are fixed and have finite second moments. We also assumed the asymptotic regime in (\ref{Eq: Limit of Sample Proportion}). These assumptions are mainly for the energy test and are not necessary for the CvM test. In fact, the same result can be derived for the CvM test given that there is a positive sequence $b_{m,n} \rightarrow \infty$ increasing arbitrary slowly with $m,n$ such that $W_d(Q_X,Q_Y) \geq b_{m,n} (1/\sqrt{m} + 1/\sqrt{n})$ (see Theorem~\ref{Theorem: Upper Bound}). 
\end{remark}

\begin{remark} \label{Remark: Integral representations}
	From the integral representations in (\ref{Eq: Projection-type expression for Energy distance}) and (\ref{Eq: Definition of the multivariate CvM-distance}), it is seen that $E_d(P_{X,N},P_{Y,N}) = (1-\epsilon) E_d(Q_X,Q_Y)$ and $W_d(P_{X,N},P_{Y,N}) \geq (1-\epsilon) W_d(Q_X,Q_Y)$, which are positive provided that $Q_X \neq Q_Y$. This explains that the poor performance of the energy test is not because of lack of signal in the contamination model but because of non-robustness of the energy test statistic. 
\end{remark}

\begin{remark}
	We mainly focus on statistical power to study robustness because one can always employ the permutation procedure to control the  type I error under $H_0 : P_{X,N} = P_{Y,N}$.
\end{remark}

\vskip 1em

\subsection{Empirical Analysis}

\begin{figure}[t!]
	\begin{center}		
		\begin{minipage}[b]{0.45\textwidth}
			\includegraphics[width=\textwidth]{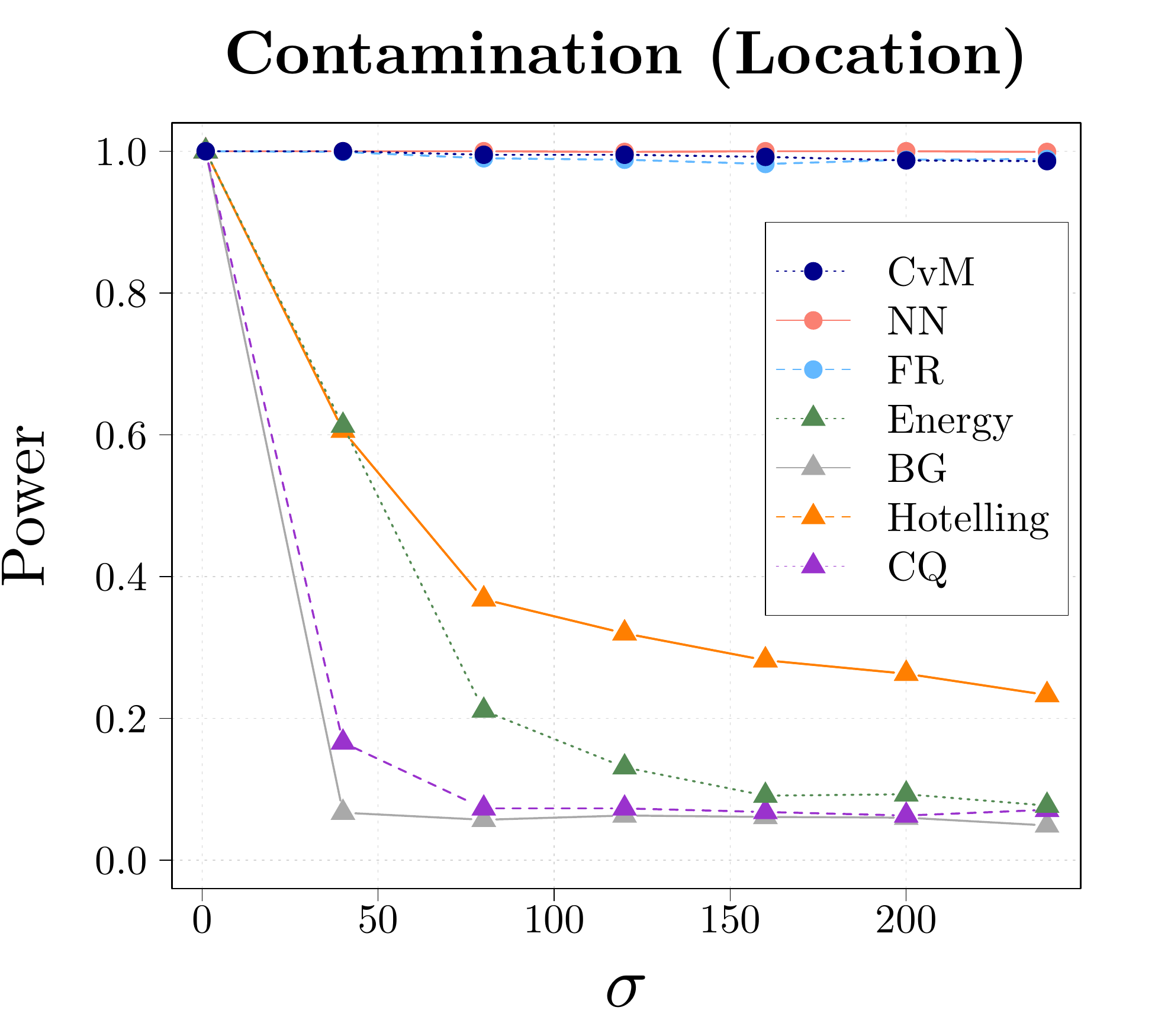}
		\end{minipage} 
		\hskip 1em
		\begin{minipage}[b]{0.45\textwidth}
			\includegraphics[width=\textwidth]{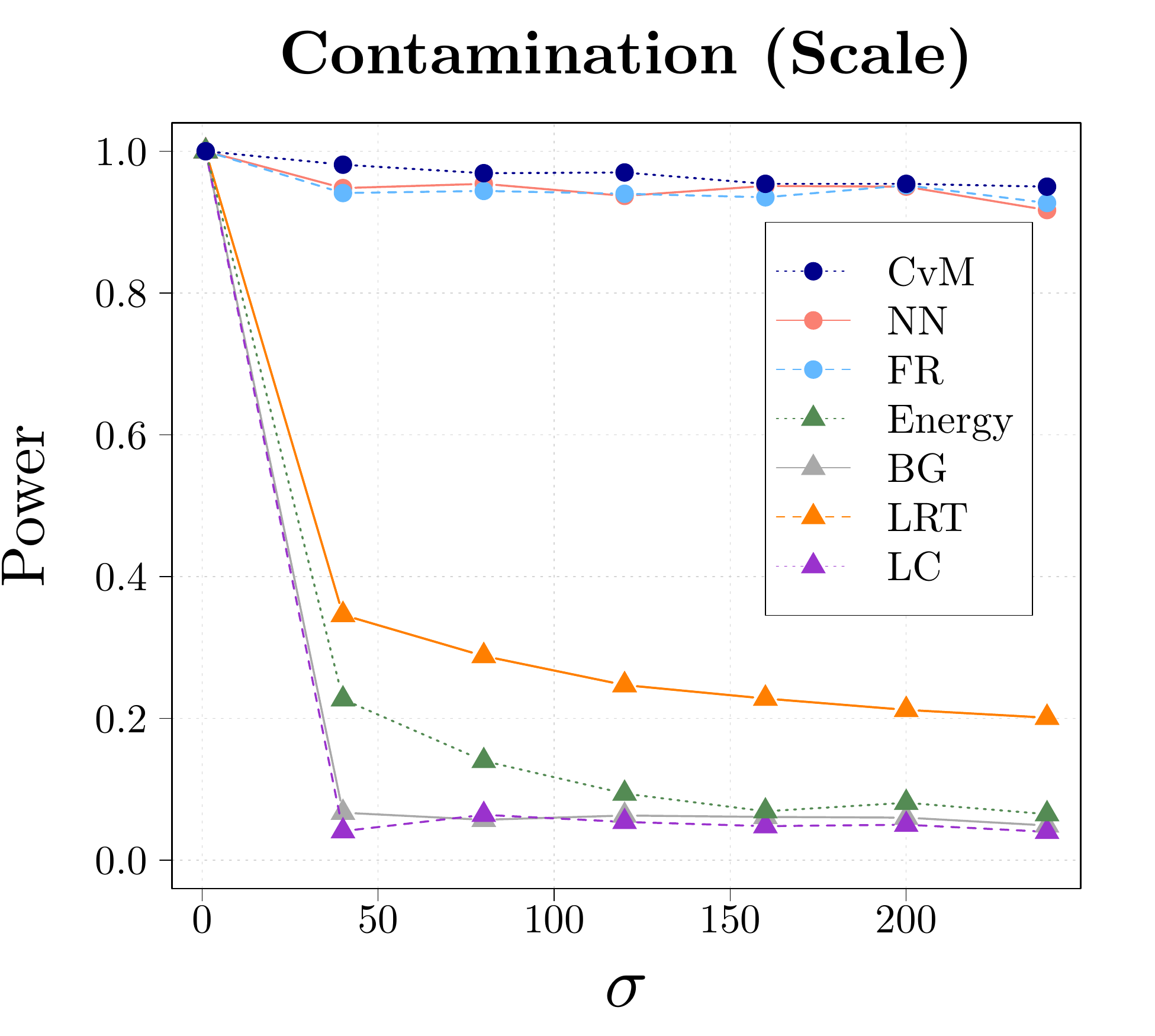}
		\end{minipage}
		\caption{\small Empirical power of NN, FR, Energy, BG, Hotelling, CQ, LRT, LC and CvM tests under the contamination models with $\epsilon = 0.05$. See Example~\ref{Proposition: Location model for robustness} and \ref{Proposition: Scale model for robustness}
			for details.} \label{Figure: Robust Study}
	\end{center}
\end{figure}

To illustrate Theorem~\ref{Theorem: Robustness} with finite sample size, we carried out simulation studies using the contamination model in (\ref{Eq: Contaminated Model}). In our simulation, we take $Q_X$ and $Q_Y$ to have multivariate normal distributions with different location parameters or different scale parameters. In both examples, we take $G_N$ to have a multivariate normal distribution given by
\begin{align*}
G_N := N((0,\ldots,0)^\top, \sigma^2 I_d),
\end{align*}
where $\sigma$ controls the degree of heavy-tailedness.

\begin{example}[Location difference] \label{Proposition: Location model for robustness} 
	For the location alternative, we compare two multivariate normal distributions, where the means are different but the covariance matrices are identical. Specifically, we set
	\begin{align*}
	& Q_X = N((-0.5,\ldots,-0.5)^\top,I_d), \quad \text{and} \quad Q_Y = N((0.5,\ldots,0.5)^\top, I_d),
	\end{align*}
	with $\epsilon = 0.05$. We then change $\sigma = 1,40,80,120,160,200$ and $240$ to investigate the robustness of the tests against heavy-tail contaminations.
\end{example}

\begin{example}[Scale difference] \label{Proposition: Scale model for robustness}
	Similar to the location alternative, we again choose multivariate normal distributions which differ in their scale but not in their location parameters. In detail, we have
	\begin{align*}
	& Q_X = N((0,\ldots,0)^\top ,0.1^2 \times I_d), \quad \text{and} \quad Q_Y = N((0,\ldots,0)^\top, I_d),
	\end{align*}
	with $\epsilon = 0.05$. Again, we change $\sigma = 1,40,80,120,160,200$ and $240$ to assess the effect of heavy-tail contaminations.
\end{example}

In addition to the energy test, we further considered three nonparametric tests in our simulation studies, namely, the $k$-nearest neighbor test by \cite{schilling1986multivariate} with $k=3$, the MST test proposed by \cite{friedman1979multivariate} and the inter-point distance test by \cite{biswas2014nonparametric}. For future reference, we refer to them as the NN test, the FR test and the BG test, respectively. We also added the high-dimensional mean test by \cite{chen2010two} and Hotelling's $T^2$ test \citep[e.g.~page 188 of][]{anderson2003introduction} for the location alternative and the high-dimensional covariance test by \cite{li2012two} and the conventional likelihood ratio test \citep[e.g.~page 412 of][]{anderson2003introduction} for the scale alternative. We refer to them as the CQ test, Hotelling's test, the LC test and the LRT test, respectively.

Experiments were run $1,000$ times to estimate the power of different tests with $m=n=40$ and $d=10$ at significance level $\alpha=0.05$. The $p$-value of each test was computed using $500$ permutations as in Remark~\ref{Remark: Monte Carlo Permutations}. As can be seen from Figure~\ref{Figure: Robust Study}, the power of the CvM test is consistently robust to the value of $\sigma,$ which supports our theoretical result. The power of the energy test, on the other hand, drops down significantly as $\sigma$ increases for both location and scale differences. As explained in the proof of Theorem~\ref{Theorem: Robustness}, this poor performance was attributed to the fact that the energy statistic is very much dominated by extreme observations from $G_N$ when $\sigma$ is large. The graph-based tests, i.e.~the NN and FR tests, also show a robust power performance against the contamination models. Intuitively speaking, they perform robust under the given scenarios as their test statistics, which count the number of edges in a graph, do not vary a lot even in the presence of outliers; but as far as we know, there is no theoretical support for this result in the current literature. The other four tests (Hotelling's test, the LRT test, the LC test and the CQ test) perform poorly for large $\sigma$, which may be explained similarly as to why the energy test has low power in these examples.

\vskip 1em

\section{Minimax Optimality} \label{Section: Minimax Optimality}

Although our choice of the $U$-statistic was a natural one to estimate $W_d^2$, it remains unclear whether one can come up with a better test statistic for testing whether $H_0: W_d=0$ or $H_0: W_d > 0$. One might also wonder whether there exists a testing procedure that leads to significantly higher power than the permutation test while controlling the type I error. In this section, we shall show that the answer is negative from a minimax point of view. In particular, we prove that the permutation test based on $U_{\text{CvM}}$ is minimax rate optimal against a class of alternatives associated with the CvM-distance.

To formulate the minimax problem, let us define the set of two multivariate distributions which are at least $\epsilon$ far apart in terms of the CvM-distance, i.e.
\begin{align*}
\mathcal{F}(\epsilon) := \big\{ (P_X,P_Y)  :  W_d(P_X,P_Y) \geq \epsilon \big\}.
\end{align*}
For a given significance level $\alpha \in (0,1)$, let $\mathds{T}_{m,n}(\alpha)$ be the set of measurable functions $\phi : \{\mathcal{X}_m, \mathcal{Y}_n\} \mapsto\{0,1\}$ such that
\begin{align*}
\mathds{T}_{m,n}(\alpha) = \{\phi : \mP_0(\phi = 1) \leq \alpha \}.
\end{align*}
We then define the minimax type II error as follows:
\begin{align}
1 - \beta_{m,n} (\epsilon) = \inf_{\phi \in \mathds{T}_{m,n}(\alpha)} \sup_{P_X,P_Y \in \mathcal{F}(\epsilon)} \mP_1 (\phi = 0).
\end{align}
Our primary interest is in finding the minimum separation rate $\epsilon_{m,n}$ satisfying
\begin{align*}
\epsilon_{m,n} = \inf \big\{ \epsilon : 1 -\beta_{m,n}(\epsilon) \leq \zeta \big\},
\end{align*}
for some $0 < \zeta < 1 - \alpha$.

\subsection{Lower Bound}
We begin by presenting a lower bound of the multivariate CvM-distance.
\begin{lemma} \label{Lemma: Lower Bound of CvM-distance}
	The multivariate CvM-distance is lower bounded by
	\begin{align} \label{Eq: Lower Bound of W_d}
	W_d(P_X,P_Y)  \geq   \int_{\mathbb{S}^{d-1}} \Big| \frac{1}{2} - \mP \left( \beta^\top X \leq \beta^\top Y \right)  \Big| d\lambda(\beta).
	\end{align}
\end{lemma}

Consider two independent random vectors $X^\ast$ and $Y^\ast$ such that their first coordinates have normal distributions as $\xi_1 \sim N(\mu_{X^\ast},1)$ and $\xi_2 \sim N(\mu_{Y^\ast},1)$ and the other coordinates have the degenerate distribution at zero, i.e.
\begin{align*}
X^\ast := (\xi_1,0,\ldots,0)^\top \quad \text{and} \quad Y^\ast := (\xi_2,0,\ldots,0)^\top.
\end{align*}
Given $\beta = (\beta_1,\ldots,\beta_d)^\top \in \mathbb{S}^{d-1}$, we have $\beta^\top X^\ast \sim N(\beta_1\mu_{X^\ast}, \beta_1^2 )$ and $\beta^\top Y^\ast \sim N(\beta_1\mu_{Y^\ast}, \beta_1^2 )$; therefore $\beta^\top X^\ast$ and $\beta^\top Y^\ast$ have continuous distributions for $\lambda$-almost all $\beta \in \mathbb{S}^{d-1}$. Under this setting, the multivariate CvM-distance is lower bounded as follows:
\begin{lemma} \label{Lemma: Least Favorable Distribution}
	Consider independent random vectors $X^\ast$ and $Y^\ast$ described above with $\mu_{X^\ast} = c m^{-1/2}$ and $\mu_{Y^\ast} = -c n^{-1/2}$ for some constant $c>0$. Let us denote the corresponding distributions by $P_{X^\ast}$ and $P_{Y^\ast}$. Then there exists another constant $C>0$ independent of the dimension satisfying
	\begin{align*}
	W_d(P_{X^\ast},P_{Y^\ast}) \geq C \left( \frac{1}{\sqrt{m}} + \frac{1}{\sqrt{n}} \right).
	\end{align*} 
	Furthermore, the lower bound is tight up to constant factors.
	\begin{proof}
		From Lemma~\ref{Lemma: Lower Bound of CvM-distance}, it is enough to show
		\begin{align*}
		\int_{\mathbb{S}^{d-1}} \Big| \frac{1}{2} - \mP \left( \beta^\top X^\ast \leq \beta^\top Y^\ast \right)  \Big| d\lambda(\beta) \geq C \left( \frac{1}{\sqrt{m}} + \frac{1}{\sqrt{n}} \right).
		\end{align*}
		For any fixed $\beta \in \mathbb{S}^{d-1}$, we have $\beta^\top (X^\ast -Y^\ast) \sim N(\beta_1( \mu_{X^\ast} - \mu_{Y^\ast}), 2\beta_1^2 )$. Let $\Phi(\cdot)$ and $\varphi(\cdot)$ denote the cumulative distribution function and the density function of the standard normal distribution respectively. Then
		\begin{align*}
		\Big| \frac{1}{2} - \mP \left( \beta^\top X^\ast \leq \beta^\top Y^\ast \right)  \Big| & = \bigg| \frac{1}{2} - \Phi \left( - \text{sign}(\beta_1) \cdot \frac{c}{\sqrt{2}} \left( \frac{1}{\sqrt{m}} + \frac{1}{\sqrt{n}} \right) \right) \bigg| \\[.5em]
		& \geq \frac{c}{\sqrt{2}} \left( \frac{1}{\sqrt{m}} + \frac{1}{\sqrt{n}} \right) \cdot \varphi\left(\frac{c}{\sqrt{2}} \left( \frac{1}{\sqrt{m}} + \frac{1}{\sqrt{n}} \right) \right) \\[.5em]
		& \geq \frac{c}{\sqrt{2}} \left( \frac{1}{\sqrt{m}} + \frac{1}{\sqrt{n}} \right) \cdot \varphi\left(\frac{c}{2\sqrt{2}} \right), 
		\end{align*}
		This lower bound holds for $\lambda$-almost all $\beta \in \mathbb{S}^{d-1}$ and thus the result follows. To have an upper bound, notice that 
		\begin{align*}
		W_d^2 (P_{X^\ast}, P_{Y^\ast}) \leq ~ & \int_{ \mathbb{S}^{d-1}} \sup_{t \in \mathbb{R}}\left(  F_{\beta^\top X}(t) - F_{\beta^\top Y}(t) \right)^2   d\lambda(\beta) \\[.5em]
		\overset{(i)}{\leq} ~ & \frac{1}{2} \int_{ \mathbb{S}^{d-1}} \mathsf{KL}\big(N(\beta_1\mu_{X^\ast},\beta_1^2), N(\beta_1 \mu_{Y^\ast}, \beta_1^2) \big) d\lambda(\beta) \\[.5em]
		= ~ & \frac{c^2}{2} \left( \frac{1}{\sqrt{m}} + \frac{1}{\sqrt{n}} \right)^2,
 		\end{align*}
 		where $\mathsf{KL}(\cdot,\cdot)$ is the Kullback-Leibler divergence between two distributions and we used the Pinsker's inequality for $(i)$ \citep[e.g.~Lemma 2.5 of][]{tsybakov2009introduction}. This shows the tightness of the lower bound. 
	\end{proof}
\end{lemma}

The previous result combined with Neyman-Pearson lemma establishes a lower bound for the minimum separation rate in the next theorem.

\begin{theorem}[Lower Bound]  \label{Theorem: Lower Bound}
	For $0 < \zeta < 1 - \alpha$, 
	there exists some constant $b = b(\alpha, \zeta)$ independent of the dimension such 
	that $\epsilon_{m,n} = b (m^{-1/2} + n^{-1/2})$ and the minimax type II error is lower bounded by $\zeta$, i.e.
	\begin{align*}
	1 - \beta_{m,n} \left(\epsilon_{m,n} \right) \geq \zeta.
	\end{align*}
\end{theorem}

\subsection{Upper Bound}
According to Theorem \ref{Theorem: Lower Bound}, \emph{no} test can have considerable power against all alternatives when $\epsilon_{m,n}$ is of order $m^{-1/2} + n^{-1/2}$. Therefore it presents a lower bound for the minimum separation rate. We now prove that this lower bound is tight by establishing a matching upper bound. In particular, the upper bound is obtained by the permutation test based on $U_{\text{CvM}}$, highlighting that the proposed approach is minimax rate optimal.

\begin{theorem}[Upper Bound] \label{Theorem: Upper Bound}
	Recall the CvM test $\phi_{\text{\emph{CvM}}} $ given in Theorem~\ref{Theorem: Robustness}. For a sufficiently large $c > 0$, let $\epsilon_{m,n}^\star$ be the radius of interest defined by 
	\begin{align} \label{Eq: Radius of interest}
	\epsilon_{m,n}^\star := c \left( \frac{1}{\sqrt{m}} + \frac{1}{\sqrt{n}} \right).
	\end{align}
	Then there exists $\zeta \in (0, 1- \alpha)$ such that the type II error of $\phi_{\text{\emph{CvM}}}$ is uniformly bounded by $\zeta$, i.e.
	\begin{align*}
	\sup_{P_X,P_Y \in \mathcal{F}(\epsilon_{m,n}^\star)} \mP_1 \left( \phi_{\text{\emph{CvM}}} = 0  \right) < \zeta.
	\end{align*}
	\begin{proof}
		Note that the permutation critical value $c_{\alpha,\text{CvM}}$ is a random quantity depending on $\mathcal{X}_m$ and $\mathcal{Y}_n$. To control the randomness from $c_{\alpha,\text{CvM}}$, we use a similar idea in \cite{fromont2013two} \citep[see also][]{albert2015tests} where they considered the quantile of a permutation critical value. Specifically, let $c^\ast_{\zeta/2}$ be the upper $\zeta/2$ quantile of the distribution of $c_{\alpha,\text{CvM}}$, and let $\mV_1$ be the variance under $H_1$. Then it suffices to show that 
		\begin{align} \label{Eq: Minimax Power Sufficient Condition}
		\mE_1\left[ U_{\text{CvM}} \right] \geq  c_{\zeta/2}^\ast  +  \sqrt{\frac{2}{\zeta} \mV_1(U_{\text{CvM}})}
		\end{align}
		uniformly over $P_X,P_Y \in \mathcal{F}(\epsilon_{m,n}^\star)$ by choosing a sufficiently large $c$. In detail, we have
		\begin{align*}
		& \mP_1 \left( U_{\text{CvM}} < c_{\alpha,\text{CvM}} \right)  \\[.5em]
		 = ~ & \mP_1 \left( U_{\text{CvM}} < c_{\alpha,\text{CvM}}, ~ c_{\alpha,\text{CvM}} > c^\ast_{\zeta/2} \right) + \mP_1 \left( U_{\text{CvM}} < c_{\alpha,\text{CvM}}, ~ c_{\alpha,\text{CvM}} \leq c^\ast_{\zeta/2} \right) \\[.5em]
		  \leq ~& \mP_1 \left( c_{\alpha,\text{CvM}} > c^\ast_{\zeta/2}  \right)  + \mP_1 \left( U_{\text{CvM}} \leq c^\ast_{\zeta/2} \right) \\[.5em]
		 \leq ~&  \frac{\zeta}{2} + \mP_1 \left( U_{\text{CvM}} \leq c^\ast_{\zeta/2} \right),
		\end{align*}
		where the second inequality is by the definition of $c^\ast_{\zeta/2}$. To control the second term, we apply Chebyshev's inequality
		\begin{align*}
		\mP_1 \left( U_{\text{CvM}} \leq c^\ast_{\zeta/2} \right) & = ~ \mP_1 \left( \frac{U_{\text{CvM}} - \mE_1 \left[ U_{\text{CvM}} \right] }{\sqrt{\mV_1 \left( U_{\text{CvM}} \right)}}  \leq \frac{c^\ast_{\zeta/2}  - \mE_1 \left[ U_{\text{CvM}} \right] }{\sqrt{\mV_1(U_{\text{CvM}})}} \right) \\[.5em]
		& = ~ \mP_1 \left( \frac{-U_{\text{CvM}} + \mE_1 \left[ U_{\text{CvM}} \right] }{\sqrt{\mV_1 \left( U_{\text{CvM}} \right)}}  \geq \frac{\mE_1 \left[ U_{\text{CvM}} \right] - c^\ast_{\zeta/2}}{\sqrt{\mV_1(U_{\text{CvM}})}} \right) \\[.5em]
		& \leq ~ \frac{\mV_1 \left( U_{\text{CvM}} \right) }{\left( \mE_1 \left[U_{\text{CvM}}\right]  - c^\ast_{\zeta/2} \right)^2} \\[.5em]
		& \leq ~ \frac{\zeta}{2},
		\end{align*}
		where the last inequality uses (\ref{Eq: Minimax Power Sufficient Condition}). Indeed, (\ref{Eq: Minimax Power Sufficient Condition}) holds and the details can be found in the supplementary document. Hence, the result follows.
	\end{proof}
\end{theorem}

\begin{remark}
	We would like to emphasize that no assumption has been made in Theorem~\ref{Theorem: Upper Bound} regarding the ratio of the sample sizes. This implies that the proposed test can be consistent against general alternatives even when the two sample sizes are highly unbalanced as $m/n \rightarrow 0$ or $m/n \rightarrow \infty$. In addition, our minimax result is based on the permutation test, which tightly controls the type I error. This is in contrast to the previous studies \citep[see e.g.][]{arias2018remember} that employed a loose cut-off value to prove minimax rate optimality.
\end{remark}

There are computationally more efficient ways of estimating $W_d^2$. For example, one can use the linear-type statistic defined as
\begin{align} \label{Eq: Linear Statistic}
L_{\text{CvM}} = \frac{1}{M}\sum_{i=1}^{M} \frac{1}{2}\left[ h_{\text{CvM}}(X_{2i-1},X_{2i}; Y_{2i-1},Y_{2i}) + h_{\text{CvM}}(X_{2i},X_{2i-1}; Y_{2i},Y_{2i-i})\right],
\end{align} 
where $M = \floor{n/2}$ and $m=n$ for simplicity. While $L_{\text{CvM}}$ is also an unbiased estimator of $W_d^2$ and can be computed in linear time, the test based on $L_{\text{CvM}}$ is notably sub-optimal in terms of minimax power. In detail, we show that the oracle test based on $L_{\text{CvM}}$ can have full power only against alternatives shrinking slower than $N^{-1/4}$ rate, whereas the minimax optimal rate is $N^{-1/2}$ when $m=n$. We build on the observation that $L_{\text{CvM}}$ converges to a normal distribution under both $H_0$ and $H_1$ to prove the following result.

\begin{proposition}[Non-optimality of the linear time test] \label{Proposition: Power of Linear Statistic}
	Let $c_{\alpha,\text{\emph{linear}}}$ be the $\alpha$ level critical value of the oracle test (see Section~\ref{Section: Critical Value and Permutation Test}) based on $L_{\text{\emph{CvM}}}$ in (\ref{Eq: Linear Statistic}) and define the corresponding test function by
	\begin{align*}
	\phi_{L_{\text{\emph{CvM}}}} := \ind(L_{\text{\emph{CvM}}} > c_{\alpha,\text{\emph{linear}}}).
	\end{align*} 
	Consider a sequence of alternatives such that
	\begin{align*}
	W_d(P_{X},P_{Y}) ~  \asymp ~ N^{-\varepsilon} \quad \text{where} \quad \varepsilon > 1/4.
	\end{align*}
	Then for $0 < \alpha < 1/2$,
	\begin{align*}
	\lim_{m,n \rightarrow \infty} \mP_1 (\phi_{L_{\text{\emph{CvM}}}} = 1) \leq 1/2.
	\end{align*}
\end{proposition}

As a straightforward consequence of Theorem~\ref{Theorem: Robustness}, we also show that the energy test, which is our main competitor, is not minimax rate optimal in our context.

\begin{proposition}[Non-optimality of the energy test] \label{Proposition: Energy Test Not Optimal}
	Recall the energy test $\phi_{\emph{Energy}}$ given in Theorem~\ref{Theorem: Robustness}. Then there exists a pair of distributions that belongs to $\mathcal{F}(\epsilon_{m,n}^\star)$ such that the energy test becomes asymptotically powerless, i.e.
	\begin{align*}
	\lim_{m,n \rightarrow \infty} \inf_{P_X, P_Y \in  \mathcal{F}(\epsilon_{m,n}^\star)} \mP_1 (\phi_{\emph{Energy}} = 1) \leq \alpha.
	\end{align*}
	\begin{proof}
		Consider $P_{X,N}= (1-\epsilon)Q_X + \epsilon G_N, P_{Y,N} = (1-\epsilon)Q_Y + \epsilon G_N$ in (\ref{Eq: Contaminated Model}) where $Q_X$ and $Q_Y$ are fixed but $Q_X \neq Q_Y$ and they have their finite second moments. Then as noted in Remark~\ref{Remark: Integral representations}, there exists a constant $\delta > 0$ such that $W_d(P_{X,N},P_{Y,N}) >  \delta$. In other words, $P_{X,N},P_{Y,N} \in \mathcal{F}(\epsilon^\star_{m,n})$. Then the result follows by Theorem~\ref{Theorem: Robustness}. 
	\end{proof}
\end{proposition}

\vskip 1em

\section{High Dimension, Low Sample Size Analysis} \label{Section: High Dimension, Low Sample Size Analysis}
We now turn our attention to the asymptotic regime where the sample size is fixed and the dimension tends to infinity. This HDLSS regime has received increasing attention in recent years and has been frequently employed to give statistical insights into high-dimensional two-sample testing \citep[e.g.][]{biswas2014nonparametric,biswas2014distribution,mondal2015high,chakraborty2017tests}. 

The goal of this section is twofold: Firstly, we provide sufficient conditions under which the proposed test is consistent in HDLSS situations. %Secondly, we show that $U_{\text{CvM}}$ has an adaptivity property that it becomes asymptotically equivalent to the existing high-dimensional mean test statistics against certain location models. From this result, we further prove that the asymptotic test based on $U_{\text{CvM}}$ has the same power as the high-dimensional mean tests under certain conditions. adapts to the high-dimensional location models having the same asymptotic behavior as the high- dimensional mean test statistics 
Secondly, we show that $U_{\text{CvM}}$ has the same asymptotic behavior as the high-dimensional mean test statistics proposed by \cite{chen2010two} and \cite{chakraborty2017tests} under certain location models. Along with these mean test statistics, we further establish the equivalence among $U_{\text{CvM}}$, the energy statistic and the MMD statistic with the Gaussian kernel. The latter connection was motivated by \cite{ramdas2015adaptivity} who showed that the energy statistic, the MMD statistic and the mean test statistic by \cite{chen2010two} are asymptotically equivalent under different scenarios.

%We also observe that the energy statistic and the MMD statistic with the Gaussian kernel are asymptotically equivalent to the aforementioned high-dimensional mean test statistics. as already discussed by \cite{ramdas2015adaptivity} under different settings. 

Let us denote $\mE(X) = \mu_X$, $\mE(Y) = \mu_Y$, $\mV(X) = \Sigma_X$ and $\mV(Y) = \Sigma_Y$ where $\Sigma_X$ and $\Sigma_Y$ are positive definite matrices. To begin we state the two assumptions.
\begin{flalign*}
\textbf{(A1).} \quad & \mV(\|Z_1^\ast - Z_2^\ast \|^2) = O(d), ~ \text{and} ~  \mV \{(Z_1^\ast-Z_3^\ast)^\top (Z_2^\ast - Z_3^\ast) \} = O(d), \\[.5em] 
& \text{where $Z_1^\ast,Z_2^\ast,Z_3^\ast$ are independent and each $Z_i^\ast$ follows either $P_X$ or $P_Y$.}\\[.5em]
\textbf{(A2).} \quad & d^{-1}\tr(\Sigma_X) \rightarrow \overline{\sigma}_X^2,  \ d^{-1}\tr(\Sigma_Y) \rightarrow \overline{\sigma}_Y^2, \ d^{-1}\|\mu_X - \mu_Y\|_2^2 \rightarrow \overline{\delta}_{XY}^2 \\[.5em] 
& \text{where $0 < \overline{\sigma}_X^2, \overline{\sigma}_Y^2 < \infty$ and $0 \leq \overline{\delta}_{XY}^2 <\infty$.} &&
\end{flalign*}

Assumption \textbf{(A1)} implies that component variables are weakly dependent. Under the distributional assumptions (including multivariate normal distributions) made in \cite{bai1996effect} and \cite{chen2010two}, \textbf{(A1)} is satisfied when
\begin{align} \label{Eq: CQ assumption}
(\mu_X - \mu_Y)^\top (\Sigma_X + \Sigma_Y) (\mu_X - \mu_Y)  = O(d) \quad \text{and} \quad  \tr \{ (\Sigma_X + \Sigma_Y )^2\} = O(d).
\end{align}
The details of this derivation can be found in the supplementary material. Assumption \textbf{(A2)} is common in the HDLSS literature \citep[e.g.][]{hall2005geometric} and facilitates the analysis. Under these conditions, the following theorem establishes the HDLSS consistency of the proposed test. 

\begin{theorem}[HDLSS consistency] \label{Theorem: HDLSS Consistency}
	Suppose \textbf{\emph{(A1)}} and \textbf{\emph{(A2)}} hold. Assume that $\overline{\sigma}_{X}^2 \neq \overline{\sigma}_{Y}^2$ or $\overline{\delta}_{XY}^2 >0$. Then for $\alpha > 1/ \{(m+n)! /(m!n!)\}$ when $m \neq n$ and for $\alpha > 2 / \{(m+n)! /(m!n!)\}$ when $m = n$, the permutation test based on $U_{\text{\emph{CvM}}}$ is consistent under the HDLSS regime, that is $\lim_{d \rightarrow \infty} \mE_1[ \phi_{\emph{CvM}} ] = 1$. 
	\begin{proof}
		Let $U_{\text{CvM}}^\varpi$ be the CvM-statistic calculated based on $\mathcal{X}_m^\varpi = \{Z_{\varpi(1)}, \ldots, Z_{\varpi(m)} \}$ and $\mathcal{Y}_m^\varpi = \{Z_{\varpi(m+1)}, \ldots, Z_{\varpi(N)} \}$ and let $\varpi_0 = \{1,\ldots,N\}$. A a high-level, the proof follows by showing that $U_{\text{CvM}}^{\varpi_0}$ achieves the maximum among other permuted test statistics under $H_1$ as $d \rightarrow \infty$. If we choose a permutation critical value such that it becomes less than $U_{\text{CvM}}^{\varpi_0}$ in the limit, then the power will converges to one as $d \rightarrow \infty$. This proof requires a careful analysis of the order among the limit values of $U_{\text{CvM}}^\varpi$ and we defer the details in the supplementary document. 
	\end{proof}
\end{theorem}

Next we focus on mean difference alternatives with equal covariance matrices. There are many types of high-dimensional mean inference procedures in the literature \citep[][for a recent review]{hu2016review}. For example, \cite{chen2010two} suggested the test statistic based on an unbiased estimator of $\|\mu_X - \mu_Y\|^2$. Specifically, their test statistic is given by
\begin{align*}
U_{\text{CQ}}=  \frac{1}{(m)_2 (n)_2} \sum_{i_1,i_2 =1}^{m,\neq} \sum_{j_1,j_2 =1}^{n,\neq} (X_{i_1} - Y_{j_1})^\top (X_{i_2} - Y_{j_2}).
\end{align*}
More recently, \cite{chakraborty2017tests} defined the test statistic based on spatial ranks as
\begin{align*}
U_{\text{WMW}} = \frac{1}{(m)_2 (n)_2} \sum_{i_1,i_2 =1}^{m,\neq} \sum_{j_1,j_2 =1}^{n,\neq} \frac{(X_{i_1} - Y_{j_1})}{\|X_{i_1} - Y_{j_1}\|}^\top \frac{(X_{i_2} - Y_{j_2})}{\|X_{i_2} - Y_{j_2}\|}.
\end{align*}
They proved that $U_{\text{CQ}}$ and $U_{\text{WMW}}$ are asymptotically equivalent under a certain HDLSS setting. Independently, the equivalence between $U_{\text{CQ}}$, $U_{\text{Energy}}$ and the MMD statistic with the Gaussian kernel was established by \cite{ramdas2015adaptivity} under different settings. Let us denote the MMD statistic with the Gaussian kernel by
\begin{align*} \label{Eq: MMD U-statistic}
U_{\text{MMD}} & =   \frac{1}{(m)_2} \sum_{i_1,i_2=1}^{m,\neq} \exp \Big( -\frac{1}{2\varsigma_d^2}\| X_{i_1} - X_{i_2}\|^2 \Big)  + \frac{1}{(n)_2} \sum_{j_1,j_2=1}^{n,\neq} \exp \Big( -\frac{1}{2\varsigma_d^2} \| Y_{j_1} - Y_{j_2}\|^2 \Big) \\[.5em] 
& - \frac{2}{mn} \sum_{i=1}^{m} \sum_{j=1}^{n} \exp \Big(-\frac{1}{2\varsigma_d^2} \| X_{i} - Y_{j}\|^2 \Big),
\end{align*}
where $\varsigma_d^2$ is the bandwidth parameter. Here we combine and further extend these results by presenting sufficient conditions under which $U_{\text{CvM}}$, $U_{\text{Energy}}$, $U_{\text{MMD}}$, $U_{\text{CQ}}$ and $U_{\text{WMW}}$ are asymptotically equivalent. To establish the result, we need two more assumptions. 
\begin{flalign*}
\textbf{(A3).} \quad & \mV\{(Z_1^\ast-Z_2^\ast)^\top (Z_3^\ast - Z_4^\ast) \} = O(d), \text{where $Z_1^\ast,Z_2^\ast,Z_3^\ast,Z_4^\ast$ are independent and} \\[.5em]
& \text{each $Z_i^\ast$ follows either $P_X$ or $P_Y$.} \\[.5em]
\textbf{(A4).} \quad & \Sigma_X = \Sigma_Y ~ \text{and} ~ \| \mu_X - \mu_Y \|^2 = O(\sqrt{d}). &&
\end{flalign*}

Assumption \textbf{(A3)} is required for studying $U_{\text{CQ}}$ and $U_{\text{WMW}}$. As Assumption \textbf{(A1)}, \textbf{(A3)} is satisfied under (\ref{Eq: CQ assumption}). Notice that $U_{\text{CQ}}$ and $U_{\text{WMW}}$ are only sensitive to location parameters whereas $U_{\text{CvM}}$, $U_{\text{Energy}}$ and $U_{\text{MMD}}$ are sensitive to both location and scale parameters. This suggests that the equal covariance assumption in \textbf{(A4)} is crucial for our result and cannot be dropped. The condition $\| \mu_X - \mu_Y \|^2 = O(\sqrt{d})$ is also important for our analysis and it was also considered in \cite{chakraborty2017tests}. Under the given assumptions, we make repeated use of Taylor expansions to establish the equivalence among the test statistics stated as follows.

\begin{theorem}[HDLSS equivalence] \label{Theorem: HDLSS Equivalence}
	Suppose \textbf{\emph{(A1)}}, \textbf{\emph{(A2)}}, \textbf{\emph{(A3)}} and \textbf{\emph{(A4)}} hold. Let $\varpi$ be an arbitrary permutation of $\{1,\ldots,N\}$ and $\overline{\sigma}_d^2 = d^{-1}\emph{\tr}(\Sigma_X)$. We denote by $U_{\text{\emph{CvM}}}^{\varpi}$, $U_{\text{\emph{Enregy}}}^{\varpi}$, $U_{\text{\emph{MMD}}}^{\varpi}$, $U_{\text{\emph{CQ}}}^{\varpi}$ and $U_{\text{\emph{WMW}}}^{\varpi}$, the CvM, Energy, MMD, CQ, and WMW test statistics, respectively, calculated based on $\mathcal{X}_m^\varpi = \{Z_{\varpi(1)},\ldots,Z_{\varpi(m)}\}$ and $\mathcal{Y}_n^\varpi = \{Z_{\varpi(m+1)},\ldots,Z_{\varpi(N)}\}$. Assume that the bandwidth parameter of the Gaussian kernel satisfies $\varsigma_d^2 \asymp d$. Then under the HDLSS asymptotics, we have that 
	\begin{equation}
	\begin{aligned}  \label{Eq: HDLSS Equivalence}
	 & \sqrt{d}U_{\text{\emph{CvM}}}^{\varpi} = \frac{1}{2 \pi \sqrt{3d} \overline{\sigma}_d^2}U_{\text{\emph{CQ}}}^{\varpi} + O_{\mathbb{P}}(d^{-1/2}), \quad	 U_{\text{\emph{Energy}}}^{\varpi} = \frac{1}{\sqrt{2d} \overline{\sigma}_d} U_{\text{\emph{CQ}}}^{\varpi} + O_{\mathbb{P}}(d^{-1/2}), \\[.5em]
	 & \sqrt{d}U_{\text{\emph{WMW}}}^{\varpi}  = \frac{1}{\sqrt{d}\overline{\sigma}_d^2 } U_{\text{\emph{CQ}}}^{\varpi} + O_{\mathbb{P}}(d^{-1/2}), \quad \sqrt{d} U_{\text{\emph{MMD}}}^{\varpi} =  \frac{\sqrt{d}}{\varsigma_d^2} e^{-d \overline{\sigma}_d^2/\varsigma_d^2} U_{\text{\emph{CQ}}}^{\varpi} + O_{\mathbb{P}}(d^{-1/2}).
	\end{aligned}
	\end{equation}
\end{theorem}

Note that the asymptotic equivalence established in (\ref{Eq: HDLSS Equivalence}) holds for any permutations. Leveraging this result, we show that the permutation critical values of the test statistics are asymptotically the same as well. 
\begin{corollary}[Permutation critical values] \label{Corollary: HDLSS Equivalence Critical Values}
	Consider the same assumptions made in Theorem~\ref{Theorem: HDLSS Equivalence}. Let $c_{\alpha,\text{\emph{CvM}}}$, $c_{\alpha,\text{\emph{Eng}}}$, $c_{\alpha,\text{\emph{MMD}}}$, $c_{\alpha,\text{\emph{CQ}}}$ and $c_{\alpha,\text{\emph{WMW}}}$ be the $1-\alpha$ quantile of the permutation distribution of $2\pi \sqrt{3d}\overline{\sigma}_d^2 U_{\text{\emph{CvM}}}$, $\sqrt{2}\overline{\sigma}_d U_{\text{\emph{Energy}}}$, $\varsigma_d^2 e^{-d\overline{\sigma}_d^2/ \varsigma_d^2} U_{\text{\emph{MMD}}} / \sqrt{d}$, $U_{\text{\emph{CQ}}}/\sqrt{d}$ and $\sqrt{d} \overline{\sigma}_d^2 U_{\text{\emph{WMW}}}$, respectively. Then
	\begin{align*}
	c_{\alpha,\text{\emph{CvM}}} & = c_{\alpha,\text{\emph{Eng}}} + O_{\mathbb{P}}(d^{-1/2}) = c_{\alpha,\text{\emph{MMD}}} + O_{\mathbb{P}}(d^{-1/2})  \\[.5em]
	& = c_{\alpha,\text{\emph{CQ}}}  + O_{\mathbb{P}}(d^{-1/2}) = c_{\alpha,\text{\emph{WMW}}} + O_{\mathbb{P}}(d^{-1/2}). 
	\end{align*}
	\begin{proof}
		We will only show that $c_{\alpha,\text{CvM}} = c_{\alpha, \text{CQ}} + O_{\mathbb{P}}(d^{-1/2})$. The remaining results follow similarly. From Theorem~\ref{Theorem: HDLSS Equivalence}, we know that 
		\begin{align*}
		2\pi\sqrt{3d} \overline{\sigma}_d^2(U_{\text{CvM}}^{\varpi_1},\ldots,U_{\text{CvM}}^{\varpi_{N!}}) = d^{-1/2}(U_{\text{CQ}}^{\varpi_1},\ldots,U_{\text{CQ}}^{\varpi_{N!}}) + O_{\mathbb{P}}(d^{-1/2})
		\end{align*}
		where $\varpi_i$ is an element of $\mathcal{S}_N$ for $i=1,\ldots,N!$. For simplicity, let us write $2\pi\sqrt{3d} \overline{\sigma}_d^2 U_{\text{CvM}}^{\varpi_i} = U_{\text{CvM},s}^{\varpi_i}$ and $d^{-1/2}U_{\text{CQ}}^{\varpi_i} = U_{\text{CQ},s}^{\varpi_i}$. Then $c_{\alpha,\text{CvM}}$ and $c_{\alpha,\text{CQ}}$ are the $\ceil{N!(1-\alpha)}$th order statistic of $\{ U_{\text{CvM},s}^{\varpi_1},\ldots, U_{\text{CvM},s}^{\varpi_{N!}}\}$ and $\{ U_{\text{CQ},s}^{\varpi_1},\ldots, U_{\text{CQ},s}^{\varpi_{N!}}\}$, respectively. It is well-known that the order statistic is a Lipschitz function \citep[e.g.~page 43 of][]{wainwright2019high}. More specifically, using Pigeonhole principle, it can be seen that 
		\begin{align*}
		|c_{\alpha,\text{CvM}} - c_{\alpha,\text{CQ}}| \leq \Bigg\{ \sum_{i=1}^{N!} (U_{\text{CvM},s}^{\varpi_i} - U_{\text{CQ},s}^{\varpi_i})^2 \Bigg\}^{1/2} = O_{\mathbb{P}}(d^{-1/2}).
		\end{align*}
		Hence the result follows. 
	\end{proof}
\end{corollary}

From the previous results, we may conclude that the considered permutation tests have comparable power in the limit as further illustrated by our simulation results in Section~\ref{Section: Simulations}. We would like to emphasize, however, that when the moment assumption is violated, the power of these tests can be entirely different. For instance, our simulation results in Section~\ref{Section: Simulations} demonstrate that the CQ, energy and MMD tests perform poorly when $X$ and $Y$ have Cauchy distributions with different location parameters. In contrast, the CvM and WMW tests maintain robust power against the same Cauchy alternative.

%As an implication of Theorem~\ref{Theorem: HDLSS Equivalence} and Corollary~\ref{Corollary: HDLSS Equivalence Critical Values}, the permutation tests based on the considered test statistics may have the same power in the limit as further illustrated by our simulation results in Section~\ref{Section: Simulations}. 

%The explicit expression for the power of a permutation test is difficult to determine in the HDLSS setting due to the randomness of its critical value. Instead, we consider the asymptotic tests based on the considered statistics and provide 

We end this section with an explicit expression for the limiting power function of the asymptotic tests based on the considered statistics. To this end, we need more restrictions on $X$ and $Y$ such as stationary $\rho$-mixing condition. Then we build on the asymptotic results established in \cite{chakraborty2017tests} combined with Theorem~\ref{Theorem: HDLSS Equivalence} to have the following corollary.

% We prove the following result by building on the central limit theorem result in HDLSS settings established by \cite{chakraborty2017tests}.

%To this end, we need more restrictions on $X$ and $Y$ such as stationary $\rho$-mixing condition as in \cite{chakraborty2017tests}.

%As another corollary of Theorem~\ref{Theorem: HDLSS Equivalence}, we show that the limiting power of the asymptotic tests based on these statistics are the same under additional restrictions. To this end, we employ asymptotic tests and consider stationary $\rho$-mixing sequences \citep[e.g.][]{zhengyan1997limit} as in \cite{chakraborty2017tests}.  We would like to emphasize, however, that when the moment assumption is violated, the power of these tests can be entirely different. For instance, our simulation results in Section~\ref{Section: Simulations} demonstrate that the CQ, energy and MMD tests perform poorly when $X$ and $Y$ have Cauchy distributions with different location parameters. In contrast, the CvM and WMW tests maintain robust power against the same Cauchy alternative.

\begin{corollary}[Power of asymptotic tests] \label{Corollary: HDLSS Power equivalence}
	Consider the same assumptions made in Theorem~\ref{Theorem: HDLSS Equivalence}. Assume that $X = \mu_X + V_X$ and $Y = \mu_Y + V_Y$ where $\mE(V_X)= \mE (V_Y) = 0$ and $V_X$ and $V_Y$ are mutually independent random vectors in $\mathbb{R}^d$. In addition, assume that the components of $V_X = (V_{X,1},V_{X,2},\ldots,)$ are strictly stationary and satisfy $\sum_{k=1}^\infty \rho_X(2^k) <\infty$ where $\rho_X(\cdot)$ is the $\rho$-mixing coefficient. The components of $V_Y = (V_{Y,1},V_{Y,2},\ldots,)$ are similarly defined with another mixing coefficient $\rho_Y(\cdot)$. Let $\{{X}_i\}_{i=1}^m$ be i.i.d. copies of $X$ and $\{{Y}_i\}_{i=1}^n$ be i.i.d. copies of $Y$. Denote
	\begin{align*}
	\psi_{m,n} = \emph{\tr}(\Sigma^2)\{ 2/m_{(2)} + 2/n_{(2)} + 4/(mn) \},
	\end{align*}
	and $\phi^\prime_{\emph{CvM}} = \ind( 2\pi\sqrt{3}d \overline{\sigma}^2 U_{\text{\emph{CvM}}} > z_{\alpha}\psi_{m,n}^{1/2} )$, $\phi^\prime_{\emph{Energy}} = \ind( \sqrt{2d} \overline{\sigma} U_{\text{\emph{Energy}}} > z_{\alpha} \psi_{m,n}^{1/2} )$, $\phi^\prime_{\emph{MMD}} = \ind(\varsigma_d^2 e^{-d\overline{\sigma}_d^2/ \varsigma_d^2} U_{\text{\emph{MMD}}} > z_{\alpha} \psi_{m,n}^{1/2})$, $\phi^\prime_{\emph{CQ}} = \ind(  U_{\text{\emph{CQ}}} > z_{\alpha} \psi_{m,n}^{1/2} )$ and $\phi^\prime_{\emph{WMW}} = \ind(  d \overline{\sigma}^2 U_{\text{\emph{WMW}}} > z_{\alpha} \psi_{m,n}^{1/2})$.
	Then under the HDLSS setting, 
	\begin{align*}
	& \lim_{d \rightarrow \infty} \mE[\phi^\prime_{\text{\emph{CvM}}} ] = \lim_{d \rightarrow \infty} \mE[\phi^\prime_{\text{\emph{Energy}}}] =  \lim_{d \rightarrow \infty} \mE[\phi^\prime_{\text{\emph{MMD}}}] = \lim_{d \rightarrow \infty} \mE[\phi^\prime_{\text{\emph{CQ}}}] = \lim_{d \rightarrow \infty} \mE[\phi^\prime_{\text{\emph{WMW}}}],
	\end{align*}
	which converges to
	\begin{align*}
	\Phi \Big(- z_\alpha + \psi_{m,n}^{-1/2}\|\mu_X - \mu_Y\|^2 \Big),
	\end{align*}
	where $z_\alpha$ is the upper $\alpha$ quantile of the standard normal distribution.
\end{corollary}

\vskip 2em

\section{Connection to the Generalized Energy Distance and MMD}  \label{Section: Connection to Generalized Energy Distance and MMD}

Recall that the energy distance is defined with the Euclidean distance under the finite first moment condition. By considering a semimetric space $(\mathds{Z}, \rho)$ of negative type, \cite{sejdinovic2013equivalence} generalized the energy distance by 
\begin{align*}
E_{\rho}^2 = 2 \mE [ \rho(X_1,Y_1) ] - \mE[ \rho(X_1,X_2)]- \mE[ \rho(Y_1,Y_2)].
\end{align*} 
They further established the equivalence between the generalized energy distance and the MMD with a kernel induced by $\rho(\cdot,
\cdot)$. Given a distance-induced kernel $k(\cdot,\cdot)$, the squared MMD is given by
\begin{align*}
\text{MMD}^2_k = \mE[ k(X_1,X_2)] +  \mE[ k(Y_1,Y_2) ]- 2 \mE[ k(X_1,Y_1)].
\end{align*}

In this section, we will show that the multivariate CvM-distance is a member of the generalized energy distance by the use of the angular distance and thus also a member of the MMD. Let $\mathcal{M}_X$ and $\mathcal{M}_Y$ be the support of $X$ and $Y$ respectively and let $\mathcal{M} = \mathcal{M}_X \cup \mathcal{M}_Y \subseteq \mathbb{R}^d$. Then we define the {\em angular distance} as follows:

\begin{definition}[Angular distance] \label{Definition: (Main) Angular Distance}
	Let $Z^\ast$ be a random vector having mixture distribution $(1/2) P_X + (1/2) P_Y$. For $z,z^\prime \in \mathcal{M}$, 
	denote the scaled angle between $z - Z^\ast$ and $z^\prime - Z^\ast$ by
	\begin{align*}
	\rho_{Angle}(z,z^\prime; Z^\ast) = \frac{1}{\pi} \mathsf{Ang} \left( z-Z^\ast, z^\prime - Z^\ast \right).
	\end{align*}
	The angular distance is defined as the expected value of the scaled angle:
	\begin{align} \label{Definition: Angle distance}
	\rho_{Angle}(z,z^\prime) = \mE \left[  \rho_{Angle}(z,z^\prime; Z^\ast)  \right].
	\end{align}
\end{definition}

The next lemma shows that $\rho_{Angle}$ is a metric of negative type defined on $\mathcal{M}$.

\begin{lemma} \label{Lemma: Angle distance is of negative type}
	For $\forall z,z^\prime, z^{\prime \prime} \in \mathcal{M}$ and $\rho_{\text{Angle}}: \mathcal{M} \times \mathcal{M} \mapsto [0,\infty)$, the following conditions are satisfied
	\begin{enumerate}
		\item $\rho_{\text{Angle}}(z,z^\prime) \geq 0$ and $\rho_{\text{Angle}}(z,z^\prime) = 0$ if and only if $z = z^\prime$. 
		\item $\rho_{\text{Angle}}(z,z^\prime) = \rho_{\text{Angle}}(z^\prime,z)$. 
		\item $\rho_{\text{Angle}}(z,z^\prime) \leq \rho_{\text{Angle}}(z,z^{\prime\prime}) + \rho_{\text{Angle}}(z^\prime,z^{\prime\prime})$.
	\end{enumerate}
	In addition, for $\forall n \geq 2$, $z_1,\ldots, z_n \in \mathcal{M}$, and $\alpha_1,\ldots,\alpha_n \in \mathbb{R}$, with $\sum_{i=1}^n \alpha_i = 0$,
	\begin{align*}
	\sum_{i=1}^n \sum_{j=1}^n \alpha_i \alpha_j \rho_{Angle}(z_i,z_j) \leq 0.
	\end{align*} 
\end{lemma}

By the use of the angular distance, we establish the identity between the generalized energy distance and the CvM-distance in the next proposition. As a result, we conclude that the multivariate CvM-distance is a special case of the generalized energy distance based on the angular distance. 

\begin{proposition}[Another view of the CvM-distance]  \label{Proposition: Identity between the generalized energy distance and CvM-distance}
	Let us consider the angular distance defined in (\ref{Definition: Angle distance}). Then
	\begin{align*}
	2W_d^2 = 2 \mE \left[ \rho_{Angle} (X_1,Y_1) \right] - \mE \left[ \rho_{Angle} (X_1,X_2) \right] - \mE \left[\rho_{Angle}(Y_1,Y_2)\right].
	\end{align*}
\end{proposition}

\vskip 1em

\begin{remark}\label{Remark: Generalization of Angular distance}
	The angular distance can be generalized by taking the expectation with respect to a different measure. For instance, when the expectation is taken with respect to Lebesgue measure, the generalized angular distance is proportional to the Euclidean distance, i.e.
	\begin{align*}
	\int_{\mathbb{R}^d} \rho_{Angle}(z,z^\prime;t) dt = \gamma_d \| z -  z^\prime \|,
	\end{align*}
	where $\gamma_d$ depends solely on the dimension (see the proof of Lemma~\ref{Lemma: Angle distance is of negative type} for more
	details). The main difference between the Euclidean distance and the proposed angular distance is that the latter takes into account information	from the underlying distribution and is less sensitive to outliers. In this aspect, the introduced angular distance can be viewed as a robust alternative for the Euclidean distance.
\end{remark}

\section{Other Multivariate Extensions via Projection-Averaging} 
\label{Section: Multivariate Extensions via Projection Averaging}

The projection-averaging approach used for the multivariate CvM-statistic can be applied to many other univariate robust statistics. In this section, we illustrate the utility of the projection-averaging approach by considering several examples including Kendall's tau, the coefficient by \cite{blum1961distribution} and the sign covariance~\citep{bergsma2014consistent}.  We begin by considering one-sample and two-sample robust statistics. Given a pair of random variables $(X,Y)$, define $Z = X-Y$. The univariate sign test statistic is an estimate of $\mathsf{T}_{\text{sign}} := \mP(Z>0) - 1/2$ and it is used to test whether
\begin{align*}
H_0: \mP(Z>0)  = 1/2 \quad \text{versus} \quad H_1: \mP(Z>0)  \neq 1/2.
\end{align*}
The projection-averaging technique extends $\mathsf{T}_{\text{sign}}$ to a multivariate case as follows: 
\begin{proposition}[One-sample sign test statistic] \label{Proposition: One-Sample Sign Test}
	For $i.i.d.$ random vectors $Z_1,Z_2$ from a multivariate distribution $P_Z$ where $Z \in \mathbb{R}^d$, the projection-averaging approach generalizes $\mathsf{T}_{\text{\emph{sign}}}$  as 
	\begin{align} \label{Eq: Projection-based one-sample sign statistic}
	\int_{\mathbb{S}^{d-1}} \bigg( \mP(\beta^\top Z_1 > 0) - \frac{1}{2} \bigg)^2 d\lambda(\beta) = \frac{1}{4} - \frac{1}{2\pi}\mE \left[ \mathsf{Ang} \left( Z_1,  Z_2\right) \right].
	\end{align} 
	\begin{proof}
		Given $\beta \in \mathbb{S}^{d-1}$, note that
		\begin{align*}
		\left( \mP(\beta^\top Z_1 > 0) - \frac{1}{2} \right)^2 = \frac{1}{4} - \mE \left[ \ind(\beta^\top Z_1 > 0) \right] + \mE \left[ \ind(\beta^\top Z_1 > 0) \ind( \beta^\top Z_2 > 0) \right].
		\end{align*}
		Applying Lemma~\ref{Lemma: Integration over Unit Sphere (2 terms)} with Fubini's theorem yields
		\begin{align*}
		& \mE \left[ \int_{ \mathbb{S}^{d-1}}  \ind(\beta^\top Z_1 > 0) d \lambda(\beta) \right] = \frac{1}{2}, \\[.5em]
		& \mE \left[ \int_{ \mathbb{S}^{d-1}}   \ind(\beta^\top Z_1 > 0) \ind(\beta^\top Z_2 > 0)  d \lambda(\beta) \right] = \frac{1}{2} - \frac{1}{2\pi}  \mE \left[ \mathsf{Ang} \left(Z_1, Z_2 \right) \right].
		\end{align*}
		This completes the proof.
	\end{proof}
\end{proposition}

\vskip 1em

Given univariate two samples $\mathcal{X}_m =  \{X_1,\ldots,X_m\}$ and $\mathcal{Y}_n = \{Y_1,\ldots,Y_n\}$, the Wilcoxon-Mann-Whitney test is designed for testing whether
\begin{align*}
H_0: \mP(X > Y) = 1/2 \quad \text{versus} \quad H_1: \mP(X > Y) \neq 1/2.
\end{align*}
Its test statistic is based on an estimate of $\mathsf{T}_{\text{WMW}} := \mP(X > Y) - 1 / 2$. The next proposition extends $\mathsf{T}_{\text{WMW}}$ to a multivariate case via projection-averaging. 

\begin{proposition}[Two-sample Wilcoxon-Mann-Whitney test statistic] \label{Proposition: Two-Sample Wilcoxon-Mann-Whitney Test}
	Let $X_1,X_2 \overset{i.i.d.}{\sim} P_X$ and, independently, $Y_1,Y_2 \overset{i.i.d.}{\sim} P_Y$ where $X_1,Y_1 \in \mathbb{R}^d$. The projection-averaging approach generalizes $\mathsf{T}_{\text{\emph{WMW}}}$ as
	\begin{align} \label{Eq: Projection-based two-sample WMW statistic}
	\int_{\mathbb{S}^{d-1}} \bigg( \mP(\beta^\top X_1 > \beta^\top Y_1) - \frac{1}{2} \bigg)^2 d\lambda(\beta) = \frac{1}{4} - \frac{1}{2\pi}\mE \left[ \mathsf{Ang} \left( X_1 - Y_1, X_2 - Y_2 \right) \right].
	\end{align} 
	\begin{proof}
		The result follows by replacing $Z_1$, $Z_2$ with $X_1 - Y_1$, $X_2 - Y_2$ in Proposition~\ref{Proposition: One-Sample Sign Test}.
	\end{proof}
\end{proposition}

\begin{remark} \label{Remark: Equivalence of PA to SP}
	The first order Taylor approximation of the inverse cosine function shows that the representations given in the right-side of (\ref{Eq: Projection-based one-sample sign statistic}) and (\ref{Eq: Projection-based two-sample WMW statistic}) are related to the spatial sign-statistics introduced by \cite{wang2015high} and \cite{chakraborty2017tests}, respectively. In fact, when $U$-statistics are used to estimate (\ref{Eq: Projection-based one-sample sign statistic}) and (\ref{Eq: Projection-based two-sample WMW statistic}), the projection-averaging statistics and the spatial sign-statistics are asymptotically equivalent under some regularity conditions (see Section~\ref{Section: Asymptotic Equivalences between Projection-Averaging and Spatial-Sign Statistics} in the supplementary document). We believe, however, that our projection-averaging-type statistics --- which can be viewed as the average of univariate statistics based on projected random variables --- is more intuitive to understand.	
\end{remark}

The same technique can be further applied to some robust statistics for independence testing. To test for independence between two random variables, Kendall's tau statistic is defined as an estimate of $\tau := 4\mP \left( X_1 < X_2, Y_1 < Y_2 \right) - 1$. We present a
multivariate extension of $\tau$ as follows:

\begin{theorem}[Kendall's tau] \label{Theorem: Kendall's tau}
	For i.i.d. pairs of random vectors $(X_1,Y_1), \ldots ,(X_4,Y_4)$ from a joint distribution $P_{XY}$ where $X \in \mathbb{R}^{p}$ and $Y \in \mathbb{R}^{q}$, the multivariate extension of $\tau$ via projection-averaging is given by
	\begin{align*}
	& \int_{\mathbb{S}^{p-1}} \int_{\mathbb{S}^{q-1}} \Big[ 4 \mP\left(\alpha^\top (X_1 - X_2) < 0 , \beta^\top (Y_1 - Y_2) < 0
	\right) - 1 \Big]^2 d \lambda(\alpha) d \lambda(\beta) \\[.5em] 
	= ~ &  \mE \left[ \left(2 - \frac{2}{\pi}\mathsf{Ang}\left( X_1 - X_2, X_3 - X_4\right) \right) \cdot \left( 2 - \frac{2}{\pi}\mathsf{Ang} \left( Y_1 - Y_2,
	Y_3 - Y_4\right) \right) \right] - 1.
	\end{align*}
\end{theorem}

Kendall's tau has been frequently used in practice due to its robustness, simplicity and interpretability. Nonetheless, the main
limitation of Kendall's tau is that it can be zero even when there exists a certain association between random variables. There have been
alternative approaches to resolve this issue in the literature. For a multivariate case, \cite{zhu2017projection} extended Hoeffding's
coefficient \citep{hoeffding1948non} via projection-averaging. Specifically, they defined the projection correlation between $X \in \mathbb{R}^{p}$ and $Y \in \mathbb{R}^q$ as
\begin{align} \label{Eq: Projection correlation}
\int_{\mathbb{S}^{p-1}} \int_{\mathbb{S}^{q-1}} \int_{\mathbb{R}^2} \left[ F_{\alpha^\top X, \beta^\top Y}(u,v) - F_{\alpha^\top X}(u) F_{\beta^\top Y}(v) \right]^2 d\omega_1 (u,v,\alpha,\beta),
\end{align}
where $d\omega_1 (u,v,\alpha,\beta) =d F_{\alpha^\top X, \beta^\top Y} (u,v) d\lambda(\alpha) d\lambda(\beta)$. Although the projection
correlation is more broadly sensitive than Kendall's tau is in detecting dependence among random variables, it can still be zero even
when $X$ and $Y$ are dependent.  A counterexample for the univariate case can be found in \cite{hoeffding1948non}.

On the other hand, the coefficient introduced by \cite{blum1961distribution} overcomes this issue by replacing $dF_{X,Y}$ with $d F_X dF_Y$. The univariate Blum-Kiefer-Rosenblatt (BKR) coefficient \citep{blum1961distribution} is defined by
\begin{align*}
\int_{\mathbb{R}^2} \left[ F_{XY}(u,v) - F_X(u) F_Y(v) \right]^2 dF_X(u) dF_Y(v).
\end{align*}

Next, we generalize the univariate BKR coefficient to a multivariate space via projection-averaging. 

\begin{theorem}[Blum-Kiefer-Rosenblatt (BKR) coefficient] \label{Theorem: BKR coefficient}
	Let us consider weight function $d \omega_2(u,v,\alpha,\beta) =
	dF_{\alpha^\top X}(u) dF_{\beta^\top Y}(v) d\lambda(\alpha)
	d\lambda(\beta)$. For i.i.d. random vectors $(X_1,Y_1),\ldots,(X_6,Y_6)$ from a joint distribution $P_{XY}$ where $X \in \mathbb{R}^p$ and $Y \in \mathbb{R}^q$, the univariate BKR coefficient can be extended to a multivariate case by
	\begin{align*}
	& \int_{\mathbb{S}^{p-1}} \int_{ \mathbb{S}^{q-1}}  \int_{\mathbb{R}^2} \left[ F_{\alpha^\top X, \beta^\top Y}(u,v) - F_{\alpha^\top X}(u) F_{\beta^\top Y}(v) \right]^2  d\omega_2(u,v,\alpha,\beta) \\[.5em]
	= ~ & \mE \left[ \left( \frac{1}{2} - \frac{1}{2 \pi} \mathsf{Ang}\left( X_1- X_3, X_2 - X_3\right) \right) \cdot \left( \frac{1}{2} - \frac{1}{2 \pi} \mathsf{Ang}\left( Y_1- Y_4, Y_2 - Y_4\right) \right)  \right] \\[.5em]
	+~  & \mE \left[ \left( \frac{1}{2} - \frac{1}{2 \pi} \mathsf{Ang} \left( X_1 - X_5, X_2 - X_5\right) \right)\cdot \left( \frac{1}{2} - \frac{1}{2 \pi} \mathsf{Ang} \left( Y_3- Y_6, Y_4 - Y_6\right) \right)  \right] \\[.5em] 
	-2   & \mE \left[ \left( \frac{1}{2} - \frac{1}{2 \pi} \mathsf{Ang}\left( X_1 - X_4, X_2 - X_4\right) \right) \cdot\left( \frac{1}{2} - \frac{1}{2 \pi} \mathsf{Ang}\left( Y_1- Y_5, Y_3 - Y_5\right) \right)  \right].
	\end{align*}
\end{theorem}

\vskip 1em

Recently, \cite{bergsma2014consistent} introduced a modification of Kendall's tau, which is zero if and only if random variables are independent under some mild conditions. Let us denote the univariate Bergsma-Dassios sign covariance by
\begin{align} \label{Eq: Univariate Tau Star}
\tau^\ast = \mE \left[ a_{\text{sign}}(X_1,X_2,X_3,X_4) \cdot a_{\text{sign}}(Y_1,Y_2,Y_3,Y_4) \right],
\end{align}
with $a_{\text{sign}}(z_1,z_2,z_3,z_4) = \text{sign} \left( |z_1 - z_2| + |z_3 - z_4| - |z_1 - z_3| - |z_2 - z_4| \right)$. Motivated by the projection-averaging approach, we propose the multivariate $\tau^\ast$ as follows:

\begin{definition}[Multivariate $\tau^\ast$] \label{Definition: Multivariate tau-star}
	Suppose $(X_1,Y_1), \ldots ,(X_4,Y_4)$ are i.i.d.~random vectors from a joint distribution $P_{XY}$ where $X \in \mathbb{R}^p$ and $Y \in \mathbb{R}^q$.  We define the multivariate $\tau^\ast$ by
	\begin{align*}
	\tau^\ast_{p,q} = \int_{ \mathbb{S}^{p-1}} \int_{ \mathbb{S}^{q-1}} \mE \big[ & a_{\text{\emph{sign}}}(\alpha^\top X_1,\alpha^\top X_2,\alpha^\top X_3,\alpha^\top X_4) \\[.5em]
	\times & a_{\text{\emph{sign}}}(\beta^\top Y_1,\beta^\top Y_2,\beta^\top Y_3,\beta^\top Y_4) \big] d \lambda(\alpha) d \lambda(\beta).
	\end{align*}
\end{definition}

Since the kernel of $\tau^\ast$ is sign-invariant, i.e. $a_{\text{sign}}(z_1,z_2,z_3,z_4) = a_{\text{sign}}(-z_1,-z_2,-z_3,-z_4)$, it is easy to see that $\tau_{p,q}^\ast$ becomes the univariate $\tau^\ast$ when $p=q=1$. Also, note that since $X$ and $Y$ are
independent if and only if $\alpha^\top X$ and $\beta^\top Y$ are independent for all $\alpha \in \mathbb{S}^{p-1}$ and $\beta \in
\mathbb{S}^{q-1}$, the characteristic property of $\tau_{p,q}^\ast$ follows by that of the univariate $\tau^\ast$.

To have an expression for $\tau_{p,q}^\ast$ without involving integrations over the unit sphere, we first generalize
Lemma~\ref{Lemma: Integration over Unit Sphere (2 terms)} with three indicator functions presented in Lemma~\ref{Lemma: Extension of Escanciano (2006) with three arguments}. Then based on this result, we provide an alternative expression for $\tau_{p,q}^\ast$ in
Theorem~\ref{Theorem: Multivariate tau star}.

\begin{lemma} \label{Lemma: Extension of Escanciano (2006) with three arguments}
	For arbitrary vectors $U_1, U_2, U_3 \in \mathbb{R}^d$, we have
	\begin{align*}
	\int_{\mathbb{S}^{d-1}} \prod_{i=1}^3 \ind(\beta^\top U_i \leq 0)  d\lambda (\beta)  =  \frac{1}{2} -\frac{1}{4\pi} \left[ \mathsf{Ang} \left( U_1, U_2 \right)  + \mathsf{Ang} \left( U_1, U_3\right) + \mathsf{Ang} \left( U_2, U_3\right) \right].
	\end{align*}
\end{lemma}

\vskip .5em 

%\begin{remark}
%	Although the expression of (\ref{Eq: integration over four indicators}) looks complicated, it simplifies the integration over the $d$-dimensional unit sphere to a single integration over the unit interval.
%\end{remark}

\vskip 1em

For $U_1,U_2,U_3 \in \mathbb{R}^d$, define $g_d(U_1, U_2, U_3)$ and $h_d(Z_1,Z_2,Z_3,Z_4)$ by
\begin{align*}
g_d(U_1, U_2, U_3) = \frac{1}{2} - \frac{1}{4\pi} \left[ \mathsf{Ang} \left( U_1, U_2 \right)  + \mathsf{Ang} \left( U_1, U_3\right) + \mathsf{Ang} \left( U_2, U_3\right) \right]
\end{align*}
and 
\begin{align*}
& h_d(Z_1,Z_2,Z_3,Z_4)   \\[.5em] 
= ~ & g_d(Z_1-Z_2, Z_2-Z_3, Z_3-Z_4) + g_d(Z_2 - Z_1, Z_1 - Z_3, Z_3 - Z_4)  \\[.5em]
+~  & g_d(Z_1 - Z_2, Z_2 - Z_4, Z_4 - Z_3) + g_d (Z_2 - Z_1, Z_1 - Z_4, Z_4 - Z_3).	
\end{align*}
Based on the kernel $h_d$, we present an alternative expression for $\tau_{p,q}^\ast$ as follows:

\begin{theorem}[Closed form expression for $\tau^\ast_{p,q}$]
	\label{Theorem: Multivariate tau star}
	For i.i.d.~random vectors $(X_1,Y_1), \ldots ,(X_4,Y_4)$ 
	from a joint distribution $P_{XY}$ where $X \in \mathbb{R}^p$ and $Y \in \mathbb{R}^q$, $\tau_{p,q}^\ast$ can be written as
	\begin{align*}
	\tau^\ast_{p,q}  ~ = ~ & \mE \left[ h_p(X_1, X_2, X_3, X_4) \cdot h_q(Y_1,Y_2,Y_3,Y_4)\right]  \\[.5em]
	+ &  \mE \left[ h_p(X_1, X_2, X_3, X_4) \cdot h_q(Y_3,Y_4,Y_1,Y_2)\right] \\[.5em]
	-2 & \mE \left[ h_p(X_1, X_2, X_3, X_4) \cdot  h_q(Y_1,Y_3,Y_2,Y_4)\right].
	\end{align*}
\end{theorem}

\vskip .8em

Theorem~\ref{Theorem: Multivariate tau star} leads to a straightforward empirical estimate of $\tau_{p,q}^\ast$ based on a
$U$-statistic. This is also true for the other multivariate generalizations introduced in this section. Using these estimates,
some theoretical and empirical properties of the proposed measures can be further investigated. These topics are reserved for future work.

\vskip .8em

%\begin{remark}
%	In order to estimate the projection correlation and the multivariate BKR coefficient, one needs at least five and six independent pairs of samples, respectively. On the other hand, an estimate of $\tau_{p,q}^\ast$ requires four independent pairs of samples; thereby $\tau_{p,q}^\ast$ has a computational advantage over the other robust measure of independence. 
%\end{remark}

%\begin{remark}
%	\cite{dhar2016study} studied the robustness of the univariate $\tau^\ast$ and contrasted it to the non-robustness of the distance covariance~\citep{szekely2007measuring}. One can similarly show that the multivariate $\tau^\ast_{p,q}$ also retains the robustness in terms of the maximum bias functional considered in \cite{dhar2016study}.
%\end{remark}

\vskip 1em

\section{Simulations} \label{Section: Simulations}
In this section, we report numerical results to support the argument in Section~\ref{Section: High Dimension, Low Sample Size Analysis} as well as to compare the performance of the CvM test with other competing nonparametric tests against heavy-tailed alternatives. Along with the energy, MMD, NN, FR and BG tests described before, we consider the cross-match test \citep{rosenbaum2005exact}, the multivariate run test \citep{biswas2014distribution}, the modified $k$-NN test \citep{mondal2015high} and the ball divergence test \citep{pan2018ball} for comparison. We refer to them as the CM test, run test, MBG test and ball test, respectively. In our simulations, we used the Gaussian kernel with the median heuristic \citep{gretton2012kernel} for the MMD test and we set the number of nearest neighbors as $k=3$ for both NN test and MBG test. Since finding the shortest Hamiltonian path for the run test is NP-complete, we employed Kruskal's algorithm \citep{kruskal1956shortest} as suggested by \cite{biswas2014distribution}.

Throughout our experiments, the significance level was set at 0.05 and the permutation procedure was used to determine the $p$-value of each test with $200$ permutations as in Remark~\ref{Remark: Monte Carlo Permutations}. The simulations were repeated 500 times to approximate the power of different tests. We set the sample size and the dimension by $m,n=20$ and $d=200$ for the balanced cases and by $m=35,n=5$ and $d=200$ for the imbalanced cases.
%Since the computational costs of $U$-statistics with permutation tests are expensive for large $n$ and $d$, we set the sample size and the dimension relatively small by $m=n=20$ ($m=35,n=5$ for the unbalanced case) and $d=200$, but which are sufficient to support our claim.

First, we consider several examples where the powers of the five tests (CvM, energy, MMD, CQ and WMW tests) in Section~\ref{Section: High Dimension, Low Sample Size Analysis} are approximately equivalent to each other. Specifically we use multivariate normal distributions with different means 
\begin{align*}
& \mu^{(0)}=(0,\ldots,0)^\top, \quad \mu^{(1)} = (0.15,\ldots,0.15)^\top \quad \text{and} \\
& \mu^{(2)} = \sqrt{0.045}(\underbrace{1,\ldots,1}_{d/2 \text{ elements}},\underbrace{0,\ldots,0}_{d/2 \text{ elements}})^\top
\end{align*}
and covariance matrices:
\begin{enumerate}
	\itemsep 1pt
	\parskip 0pt
	\item Identity matrix (denoted by $I$) where $\sigma_{i,i} = 1$ and $\sigma_{i,j} = 0$ for $i \neq j$. 
	\item Banded matrix (denoted by $\Sigma_{Band}$) where $\sigma_{i,i} = 1$, $\sigma_{i,j} = 0.6$ for $|i-j| = 1$, $\sigma_{i,j} = 0.3$ for $|i-j| = 2$ and $\sigma_{i,j} = 0$ otherwise.
	\item Autocorrelation matrix (denoted by $\Sigma_{Auto}$) where $\sigma_{i,i} = 1$ and $\sigma_{i,j} = 0.2^{|i-j|}$ when $i\neq j$.
	\item Block diagonal matrix (denoted by $\Sigma_{Block}$) where the $5\times 5$ main diagonal blocks $\mathbf{A}$ are defined by $a_{i,i} = 1$ and $a_{i,j} = 0.2$ when $i\neq j$, and the off-diagonal blocks are zeros.
\end{enumerate}
Then we generate random samples from $X \sim N(\mu^{(0)},\Sigma)$ and either $Y \sim N(\mu^{(1)},\Sigma)$ or $Y \sim N(\mu^{(2)},\Sigma)$.
The results are summarized in Table~\ref{Table: Equivalence}. As can be seen from the table, the empirical powers of the considered tests are very close under the given setting, which supports our theoretical results in Section~\ref{Section: High Dimension, Low Sample Size Analysis}. We also observe that the other nonparametric tests, not considered in Section~\ref{Section: High Dimension, Low Sample Size Analysis}, are significantly less powerful than the proposed test in all normal location alternatives.

\setlength{\abovetopsep}{0.5em}

\begin{table}[t!]
\begin{center}
	\renewcommand{\tabcolsep}{7pt}
	\renewcommand\arraystretch{1}
	\small
	\caption{\small Empirical power of the considered tests against the normal location models at $\alpha=0.05$.}
	\label{Table: Equivalence}
	\begin{tabular}{ccccccccc}
		\toprule
		\multicolumn{1}{l}{} & \multicolumn{2}{c}{$I_d$} & \multicolumn{2}{c}{$\Sigma_{{Band}}$} & \multicolumn{2}{c}{$\Sigma_{{Block}}$} & \multicolumn{2}{c}{$\Sigma_{{Auto}}$} \\  \cmidrule(lr){2-3} \cmidrule(lr){4-5} \cmidrule(lr){6-7} \cmidrule(lr){8-9} 
		 $m=20,n=20$             & $\mu^{(1)}$ & $\mu^{(2)}$ & $\mu^{(1)}$           & $\mu^{(2)}$          & $\mu^{(1)}$          & $\mu^{(2)}$          & $\mu^{(1)}$          & $\mu^{(2)}$         \\ \midrule
		CvM                  & \textbf{0.662}       & \textbf{0.646}       & \textbf{0.418}                 & \textbf{0.406}                & \textbf{0.572}                & \textbf{0.584}                & \textbf{0.452}                & \textbf{0.442}               \\
		Energy               & \textbf{0.656}       & \textbf{0.650}       & \textbf{0.420}                 & \textbf{0.408}                & \textbf{0.576}                & \textbf{0.584}                & \textbf{0.452}                & \textbf{0.444}               \\
		MMD                  & \textbf{0.658}       & \textbf{0.638}       & \textbf{0.412}                 & \textbf{0.398}                & \textbf{0.568}                & \textbf{0.570}                & \textbf{0.458}                & \textbf{0.444}               \\
		CQ                   & \textbf{0.656}       & \textbf{0.650}       & \textbf{0.416}                 & \textbf{0.412}                & \textbf{0.578}                & \textbf{0.580}                & \textbf{0.454}                & \textbf{0.448}               \\
		WMW                  & \textbf{0.668}       & \textbf{0.646}       & \textbf{0.420}                 & \textbf{0.402}                & \textbf{0.568}                & \textbf{0.580}                & \textbf{0.458}                & \textbf{0.444}               \\ \midrule
		NN                  & 0.288       & 0.288       & 0.164                 & 0.154                & 0.242                & 0.238                & 0.176                &  0.174             \\
		FR                   & 0.168       & 0.170       & 0.090                 & 0.084                & 0.158                & 0.116                & 0.112                & 0.088               \\
		MBG               & 0.050       & 0.050       & 0.050                 & 0.052                & 0.048                & 0.044                & 0.060                & 0.046               \\
		Ball                 & 0.240       & 0.254       & 0.186                 & 0.198                & 0.262                & 0.250                & 0.216                & 0.226				 \\
		CM           & 0.042       & 0.054       & 0.028                 & 0.040                & 0.052                & 0.050                & 0.038                & 0.034               \\
		BG                  & 0.070       & 0.060       & 0.074                 & 0.074                & 0.074                & 0.078                & 0.084                & 0.078               \\
		Run                 & 0.160       & 0.153       & 0.101                 & 0.105                & 0.146               & 0.128                & 0.110                & 0.102				 \\		
		\bottomrule              
	\end{tabular}
\end{center}
\end{table}

%\vspace{-2em}

\begin{table}[t!]
	\begin{center}
		\renewcommand{\tabcolsep}{7pt}
		\renewcommand\arraystretch{1}
		\small
		\caption{\small Empirical power of the considered tests against multivariate Cauchy distributions with $m=n=20$ at $\alpha=0.05$ where $\gamma, s$  represent the location and scale parameter, respectively. The three highest power estimates in each column are highlighted in boldface.}
		\label{Table: Cauchy alterantives Balanced Case}
		\begin{tabular}{ccccccccc}
			\toprule
			\multicolumn{1}{l}{} & \multicolumn{4}{c}{\emph{Location}}                      & \multicolumn{4}{c}{\emph{Scale}}     \\ \cmidrule(lr){2-5}\cmidrule(lr){6-9}
			\emph{$m=20,n=20$}              & $\gamma=2$ & $\gamma=3$ & $\gamma=4$ & $\gamma=5$ & $s=2$ & $s=3$ & $s=4$ & $s=5$ \\ \midrule
			CvM                  & 0.124      & 0.252      & 0.596      & 0.842      & \textbf{0.560} & \textbf{0.926} & \textbf{0.988} & \textbf{1.000} \\
			Energy              & 0.060      & 0.066      & 0.102      & 0.134     & 0.316 & 0.602 & 0.766 & 0.866 \\
			MMD                & 0.056     & 0.064      & 0.110     & 0.162      & 0.448 & 0.772 & 0.890 & 0.970 \\
			CQ                    & 0.138      & 0.268      & 0.360      & 0.456      & 0.046 & 0.070 & 0.042 & 0.068 \\
			WMW               & \textbf{0.324}      & \textbf{0.698}      & \textbf{0.912}      & \textbf{0.988}     & 0.052 & 0.064 & 0.062 & 0.056 \\ \midrule
			NN                   & \textbf{0.288}      & \textbf{0.662}      & \textbf{0.884}      & \textbf{0.976}      & 0.214 & 0.194 & 0.256 & 0.224 \\
			FR                    & \textbf{0.178}      & \textbf{0.462}      & \textbf{0.706}      & \textbf{0.888}      & 0.028 & 0.034 & 0.048 & 0.036 \\
			MBG                & 0.060     & 0.044      & 0.050      & 0.074      & \textbf{0.564} & \textbf{0.904} & \textbf{0.964} & \textbf{0.992} \\
			Ball                  & 0.064      & 0.064      & 0.076     & 0.098      & \textbf{0.606} & \textbf{0.936} & \textbf{0.994} &   \textbf{1.000}   \\ 		
			CM                  & 0.030      & 0.078      & 0.128      & 0.226      & 0.056 & 0.170 & 0.334 & 0.490 \\
			BG                   & 0.048      & 0.038      & 0.048      & 0.040      & 0.238 & 0.394 & 0.560 & 0.632 \\
			Run                 & 0.059       & 0.129       & 0.274                 & 0.422                & 0.220                & 0.506                & 0.767                & 0.864	 \\		
			\bottomrule     
		\end{tabular}
		
		\vskip 1em
		
		\renewcommand{\tabcolsep}{7pt}
		\renewcommand\arraystretch{1}
		\small
		\caption{\small Empirical power of the considered tests against multivariate Cauchy distributions with $m=35$ and $n=5$ at $\alpha=0.05$ where $\gamma, s$ represent the location and scale parameter, respectively. The three highest power estimates in each column are highlighted in boldface.}
		\label{Table: Cauchy alterantives Imbalanced Case}
		\begin{tabular}{ccccccccc}
			\toprule
			\multicolumn{1}{l}{} & \multicolumn{4}{c}{\emph{Location}}                      & \multicolumn{4}{c}{\emph{Scale}}     \\ \cmidrule(lr){2-5}\cmidrule(lr){6-9}
			\emph{$m=35,n=5$}              & $\gamma=5$ & $\gamma=6$ & $\gamma=7$ & $\gamma=8$ & $s=3$ & $s=4$ & $s=5$ & $s=6$ \\ \midrule
			CvM                 & \textbf{0.340}      & \textbf{0.498}      & \textbf{0.652}      & \textbf{0.758}      & \textbf{0.570} & \textbf{0.806} & \textbf{0.928} & \textbf{0.952} \\
			Energy             & 0.110      & 0.146     & 0.212      & 0.262     & \textbf{0.436} & \textbf{0.632} & \textbf{0.794} & \textbf{0.858} \\
			MMD                & 0.108     & 0.148      & 0.192     & 0.240      & \textbf{0.552} & \textbf{0.808} & \textbf{0.926} & \textbf{0.968} \\
			CQ                   & \textbf{0.284}      & \textbf{0.380}      & 0.454      & 0.544      & 0.178 & 0.210 & 0.262 & 0.290 \\
			WMW               & \textbf{0.796}      & \textbf{0.890}      & \textbf{0.942}      & \textbf{0.960}     & 0.110 & 0.126 & 0.134 & 0.148 \\ \midrule
			NN                   & 0.144      & 0.294      & 0.376      & 0.558     & 0.118 & 0.150 & 0.154 & 0.182 \\
			FR                    & 0.226      & 0.360      & \textbf{0.464}      & \textbf{0.588}      & 0.078 & 0.092 & 0.104 & 0.112 \\
			MBG                & 0.010     & 0.000      & 0.008      & 0.000      & 0.092 & 0.130 & 0.176 & 0.214 \\
			Ball                  & 0.072      & 0.088      & 0.098     & 0.122      & 0.238 & 0.406 & 0.594 &   0.762   \\ 		
			CM                  & 0.082      & 0.176      & 0.190      & 0.262      & 0.030 & 0.080 & 0.092 & 0.126 \\
			BG                   & 0.058      & 0.052      & 0.058      & 0.052      & 0.320 & 0.386 & 0.506 & 0.514 \\
			Run                 & 0.088       & 0.150       & 0.198                 & 0.228                & 0.106                & 0.174                & 0.248                & 0.326	 \\		
			\bottomrule     
		\end{tabular}
	\end{center}
\end{table}

In our second experiment, we consider several examples where the moment conditions are not satisfied. We focus on random samples generated from multivariate Cauchy distributions. Let $\text{Cauchy}(\gamma, s)$ refer to the univariate Cauchy distribution where $\gamma, s$ are the location parameter and the scale parameter, respectively. Let $X = (X^{(1)},\ldots, X^{(d)})$ and $Y=(Y^{(1)},\ldots,Y^{(d)})$ be random vectors where $X^{(i)} \overset{i.i.d.}{\sim} \text{Cauchy}(0,1)$ and $Y^{(i)} \overset{i.i.d.}{\sim} \text{Cauchy}(\gamma, s)$ for $i= 1,\ldots,d$. We first consider location differences where $\gamma$ is not zero but the scale parameters are identical, i.e. $s=1$. Similarly, we consider scale differences where the scale parameter $s$ changes, but the location parameters are identical, i.e. $\gamma=0$. 

From the results presented in Table~\ref{Table: Cauchy alterantives Balanced Case} and Table~\ref{Table: Cauchy alterantives Imbalanced Case}, it is seen that, unlike the multivariate normal cases, there are significant differences between power performance among CvM, energy, MMD, CQ and WMW tests. In particular, the tests based on the energy, MMD and CQ statistics have relatively low power against the heavy-tail location alternatives, whereas the tests based on the CvM and WMW statistics show better performance than the others. Turning to the scale problems, it can be seen that the CQ and WMW tests are not sensitive to detect scale differences, which makes sense because they are specifically designed for location problems. On the other hand, the CvM, energy and MMD tests perform reasonably well in these alternatives. Among the omnibus nonparametric tests, the MMD, energy and ball tests have competitive power against the scale differences, but not against the location differences in general. The MBG test is only powerful against the scale differences where the sample sizes are balanced. The CM and run tests are uniformly outperformed by the CvM test under all scenarios. The NN and FR tests perform strongly against the location alternatives especially for the balanced case, but not against the scale alternatives. When the sample sizes are unbalanced, the performance of the NN and FR tests are degraded a little bit, which can be explained by \cite{chen2013ensemble} and \cite{chen2018weighted}. The CvM test, on the other hand, performs consistently well against the heavy-tail location and scale alternatives and its performance appears immune to the sample proportion.

In summary, the proposed test has almost identical power as the high-dimensional mean tests against the light-tail location alternatives, whereas it outperforms many popular nonparametric competitors under the heavy-tail location and scale alternatives.

%As a side note, the problem of the NN and FR tests associated with the unbalanced sample sizes can be found in \cite{chen2013ensemble} and \cite{chen2018weighted}. 
%The performance of the NN and FR tests appear sensitive to the sample proportion, which coincides with the previous observations by \cite{chen2013ensemble} and \cite{chen2018weighted}. 

%\clearpage 

%\clearpage 

\vskip 1em

\section{Concluding Remarks}

In this work, we extended the univariate Cram{\'e}r-von Mises statistic for two-sample testing to the multivariate case using projection-averaging. The proposed statistic has a straightforward calculation formula in arbitrary dimensions and the resulting test has good statistical properties. Throughout this paper, we demonstrated its robustness, minimax rate optimality and high-dimensional power properties. In addition, we applied the same projection technique to other robust statistics and presented their multivariate extensions.

Beyond nonparametric testing problems, we believe that our approach can be used for other problems. For example, our work can be viewed as an application of the angular distance to the two-sample problem. The angular distance is closely connected to the Euclidean distance (Remark~\ref{Remark: Generalization of Angular distance}) but is more robust to outliers by incorporating information from the underlying distribution. Given that the use of distances is of fundamental importance in many statistical applications (including clustering, classification and regression), we expect that the angular distance can be applied to other statistical problems as a robust alternative for the Euclidean distance.

%\bibliographystyle{plainnat}
%\bibliography{bibtex}

\bibliographystyle{apalike}
\bibliography{reference}

%\clearpage

\appendix

\section{Permutation Tests} \label{Section: Permutation Tests}

%Permutation tests are among the most widely used conditional procedures of nonparametric inference. An essential feature of the
%permutation tests is that they are guaranteed to obtain an \emph{exact} rejection probability whenever the exchangeability condition is satisfied under $H_0$. However, except for the finite sample exactness, there has been little theory for permutation tests. The goal here is to enhance our understanding of permutation tests for the two-sample problem. Specifically, we would like to establish fairly general conditions under which the permutation distribution is asymptotically equivalent to the corresponding unconditional null distribution based on a two-sample $U$-statistic. We demonstrate our results using the proposed CvM-statistic.

In this section, we study the limiting behavior of the permutation distribution of a two-sample $U$-statistic under the conventional asymptotic framework (\ref{Eq: Limit of Sample Proportion}). Specifically, we establish fairly general conditions under which the permutation distribution of a two-sample $U$-statistic is asymptotically equivalent to the corresponding unconditional null distribution. We first focus on the large sample behavior of the permutation distribution under the null hypothesis in Section~\ref{Section: Asymptotic null behavior of permutation $U$-statistics} and then discuss how to generalize this result to the alternative hypothesis via coupling argument in Section~\ref{Section: The coupling argument}.

\subsection{Asymptotic null behavior of permutation $U$-statistics} \label{Section: Asymptotic null behavior of permutation $U$-statistics}

Let us start with some notation. For $r \geq 2$, consider a kernel $g(x_1,\ldots,x_r; y_1,\ldots, y_r)$ of degree $(r,r)$ such that
\begin{equation}
\begin{aligned} \label{Eq: Assumption - first/second moments}
& \mE\left[ g(X_1,\ldots,X_r; Y_1,\ldots, Y_r) \right] = \theta, \\[.5em]
& \mE \left[ \{ g(X_1,\ldots,X_r; Y_1,\ldots, Y_r) \}^2 \right] < \infty.
\end{aligned}
\end{equation}
Without loss of generality, we assume that $g(x_1,\ldots,x_r;y_1,\ldots,y_r)$ is symmetric in each set of arguments, which means that the value of the kernel is invariant to the order of the first $r$ arguments as well as the last $r$ arguments. The reason for this is that we can always redefine the kernel as
\begin{align} \label{Eq: Symmetrization}
\widetilde{g}(x_1,\ldots,x_r;y_1,\ldots,y_r) = \frac{1}{r!r!} \sum_{\varpi \in \mathcal{S}_r} \sum_{\varpi^\prime \in \mathcal{S}_r} g(x_{\varpi(1)},\ldots,x_{\varpi(r)}; y_{\varpi^\prime(1)},\ldots,y_{\varpi^\prime(r)}),
\end{align}
where $\mathcal{S}_r$ is the set of all permutations of $\{1,\ldots,r\}$.

Let us write the $U$-statistic based on the kernel $g$ by
\begin{align} \label{Eq: General Two-Sample U-statistic}
U_{m,n} = \frac{1}{\binom{m}{r}\binom{n}{r}} \sum_{\alpha_1,\ldots,\alpha_r}  \sum_{\beta_1,\ldots,\beta_r}g(X_{\alpha_1},\ldots,X_{\alpha_r}; Y_{\beta_1},\ldots, Y_{\beta_r}), 
\end{align}
where the sums are taken over all subsets $\{\alpha_1,\ldots,\alpha_r \} $ of $\{1,\ldots,m\}$ and $\{\beta_1,\ldots, \beta_r \}$ of $\{1,\ldots,n\}$ and $\binom{m}{r}$ and $\binom{n}{r}$ are the binomial coefficient defined by $m!/\{r!(m-r)!\}$ and $n!/\{r!(n-r)!\}$, respectively. For $0 \leq c,d \leq r$, let $g_{c,d}(x_1,\ldots,x_c;y_1,\ldots,y_d)$ be the conditional expectation given by
\begin{align} \label{Eq: definition of g_{c,d} }
g_{c,d}(x_1,\ldots,x_c;y_1,\ldots,y_d) :=  \mE \big[ g(x_1,\ldots,x_c,X_{c+1},\ldots, X_r;  y_1,\ldots,y_d,Y_{d+1}, \ldots, Y_r) \big] .
\end{align}
Further write the centered conditional expectation and its variance as
\begin{align}  \label{Eq: Definition of g ast}
& g^\ast_{c,d}(x_1,\ldots,x_c;y_1,\ldots,y_d) := g_{c,d}(x_1,\ldots,x_c;y_1,\ldots,y_d)  - \theta,  \\[0.5em]
& \sigma_{c,d}^2 := \mV \left[g_{c,d}(X_1,\ldots,X_c; Y_1,\ldots, Y_d)\right] = \mE\big[ \big\{ g^\ast_{c,d}(X_1,\ldots,X_c; Y_1,\ldots, Y_d) \big\}^2\big]. \label{Eq: Definition of xi_{cd}}
\end{align}

The kernel $g$ is \emph{non-degenerate} if both $\sigma_{0,1}$ and $\sigma_{1,0}$ are strictly positive, and \emph{degenerate} if $\sigma_{0,1} = \sigma_{1,0} = 0$. For the case where the kernel is non-degenerate, \cite{chung2016asymptotically} provided a sufficient condition under which the permutation distribution approximates the unconditional distribution of $U_{m,n}$. Their result, however, does not cover some important degenerate $U$-statistics including $U_{\text{CvM}}$, $U_{\text{Energy}}$ and $U_{\text{MMD}}$ in the main text. To fill this gap, we develop a similar result for the degenerate cases.

Consider the centered $U$-statistic scaled by $N=m+n$:
\begin{align*}
U^\ast_{m,n}(X_1,\ldots,X_m, Y_1,\ldots,Y_n) := N(U_{m,n} -\theta ),
\end{align*}
and let $\{Z_1,\ldots, Z_{m+n}\}= \{X_1,\ldots,X_m, Y_1,\ldots,Y_n \}$ be the pooled samples. Then the permutation distribution function of $U^\ast_{m,n} $ can be written as
\begin{align*}
\widehat{R}_{m,n}(t) = \frac{1}{N!} \sum_{\varpi \in \mathcal{S}_N} I\big\{  U^\ast_{m,n}(Z_{\varpi(1)}, \ldots, Z_{\varpi(N)}) \leq t \big\}.
\end{align*}  
Also, let $R(t)$ be the unconditional limiting null distribution of $U^\ast_{m,n}$. Then we present the following theorem.  

\begin{theorem} \label{Theorem: Two-Sample Degenerate Kernel}
	Suppose $g(x_1,\ldots,x_r; y_1,\ldots, y_r)$ is symmetric in each set of arguments and degenerate under $H_0$. Further assume that $\mE[g^2] < \infty$ and it satisfies
	\begin{enumerate}\setlength{\itemindent}{.7in}
		\item[Condition 1.] $g^\ast_{0,2}(z_1,z_2) = g^\ast_{2,0}(z_1,z_2)$ and $g^\ast_{1,1}(z_1,z_2) = \frac{1-r}{r}g^\ast_{0,2}(z_1,z_2)$,
		\item[Condition 2.] $\sigma_{0,1}^2 = \sigma_{1,0}^2 = 0$ and $\sigma_{0,2}^2, \ \sigma_{2,0}^2, \ \sigma_{1,1}^2 > 0$,
	\end{enumerate}
	Then under the conventional limiting regime (\ref{Eq: Limit of Sample Proportion}) and $H_0$,
	\begin{align} \label{Eq: Approximation of the permutation null distribution}
	\sup_{t \in \mathbb{R}} \Big| \widehat{R}_{m,n}(t) - R(t)  \Big| \convP 0.
	\end{align}
	\begin{proof}
		The proof can be found in Section~\ref{Section: Proof of Permutation Consistency}.
	\end{proof}
\end{theorem}

\subsection{The coupling argument} \label{Section: The coupling argument}

The proof of Theorem~\ref{Theorem: Two-Sample Degenerate Kernel} relies on the fact that $Z_{\varpi(1)}, \ldots, Z_{\varpi(N)}$ are $i.i.d.$ samples under the null hypothesis for any permutations. The main difficulty of generalizing this result to the alternative hypothesis is that the given samples are not identically distributed under $H_1$. We instead have $m$ samples $\{X_1,\ldots,X_m\}$ from $P_X$ and $n$ samples $\{Y_1,\ldots,Y_n\}$ from $P_Y$. In order to overcome such difficulty, we employ the coupling argument considered in \cite{chung2013exact}, which is summarized in Algorithm~\ref{Alg: Coupling}.

\begin{algorithm}
	\KwData{$\{ Z_1,\ldots,Z_{N} \}:= \{ X_1,\ldots,X_m, Y_1,\ldots,Y_n \}$ where $\{X_1,\ldots,X_m\} \overset{i.i.d.}{\sim}P_X$ and $\{Y_1,\ldots,Y_n\} \overset{i.i.d.}{\sim} P_Y$, a random permutation $\varpi_0$ of $\{1,\ldots,N\}$.}
	\KwResult{$\{ \overline{Z}_{\varpi_0(1)}, \ldots, \overline{Z}_{\varpi_0(N)} \}$.}
	\Begin{
		$B \sim \text{Binomial}(N,m/N)$\;
		\If{$B \geq m$}{
			\text{Generate $\{X_{m+1},\ldots, X_B \}$ $i.i.d.$ samples from $P_X$}\;
			\Return{$\{ \overline{Z}_{\varpi_0(1)}, \ldots, \overline{Z}_{\varpi_0(N)} \}:=\{ X_1,\ldots, X_m, Y_1,\ldots,Y_{N-B}, X_{m+1},\ldots,X_B \}$}\;
		}
		\Else{
			\text{Generate $\{Y_{n+1},\ldots, Y_{N-B} \}$ $i.i.d.$ samples from $P_Y$}\;
			\Return{$\{ \overline{Z}_{\varpi_0(1)}, \ldots, \overline{Z}_{\varpi_0(N)} \}:=\{ X_1,\ldots, X_{B}, Y_{n+1},\ldots, Y_{N-B} ,Y_1,\ldots,Y_{n}\}$}\;
		}
	}
	\caption{Coupling}
	\label{Alg: Coupling}
\end{algorithm}

Note that the output of Algorithm~\ref{Alg: Coupling} consists of $i.i.d.$ samples from $\frac{m}{N} P_X + \frac{n}{N} P_Y$. Also note that there are $D = |m-B|$ different observations between the original samples $\{Z_1,\ldots, Z_N\}$ and the coupled samples $\{ \overline{Z}_{\varpi_0(1)}, \ldots, \overline{Z}_{\varpi_0(N)} \}$. The main strategy of studying the permutation distribution under the alternative hypothesis is to establish that 
\begin{align} \label{Eq: Coupling Goal}
U^\ast_{m,n}(Z_{\varpi(1)},\ldots,Z_{\varpi(N)}) - U^\ast_{m,n}(\overline{Z}_{\varpi(\varpi_0(1))}, \ldots, \overline{Z}_{\varpi(\varpi_0(N))}) \convP 0.
\end{align}
If this is the case, then both statistics have the same limiting behavior, which means that we can still apply Theorem~\ref{Theorem: Two-Sample Degenerate Kernel}. We demonstrate this procedure by using the proposed CvM-statistic and prove Theorem~\ref{Theorem: Critical value of Permutation test} in the main text. The details can be found in the proof of Theorem~\ref{Theorem: Critical value of Permutation test}.

\begin{remark}
	The coupling argument in \cite{chung2013exact} requires the condition
	\begin{align} \label{Eq: Assumption of Prior}
	\frac{m}{N} - \vartheta_X = O\left( \frac{1}{\sqrt{N}} \right),
	\end{align}
	which turns out to be unnecessary in our application; we only need the assumption that $m/N \rightarrow \vartheta_X \in (0,1)$ and $n/N \rightarrow \vartheta_Y \in (0,1)$  as $N\rightarrow \infty$ without any further restriction. To remove the condition in (\ref{Eq: Assumption of Prior}), we first show that the test statistic based on permuted samples is close to that based on $i.i.d.$ samples from $\frac{m}{N} P_X + \frac{n}{N}P_Y$. Then we will show that the two test statistics --- one is based on $i.i.d.$ samples from $\frac{m}{N} P_X + \frac{n}{N}P_Y$ and the other one is based on $i.i.d.$ samples from $\vartheta_X P_X + \vartheta_Y P_Y$ --- have the same asymptotic behavior. 
\end{remark}

\allowdisplaybreaks

\vskip 1em

\section{Auxiliary Lemmas}
In this section, we collect some auxiliary lemmas used in our main proofs. We start with another expression for the CvM-distance. 

%To begin, we present another expression for the CvM-distance (Lemma~\ref{Lemma: Multivaraite CvM Expression}) as well as the proposed statistic (Lemma~\ref{Lemma: another expression for the U-statistic}). 
\begin{lemma}[Another expression for the CvM-distance] \label{Lemma: Multivaraite CvM Expression}
	Let $X_1,X_2,X_3 \overset{i.i.d.}{\sim} P_X$ and, independently, $Y_1,Y_2,Y_3 \overset{i.i.d.}{\sim} P_Y$. Furthermore, assume that $\beta^\top X_1$ and $\beta^\top Y_1$ have continuous distribution functions for $\lambda$-almost all $\beta \in \mathbb{S}^{d-1}$. Then the squared multivariate CvM-distance can be written as 
	\begin{align*}
	W^2_d(P_X,P_Y) ~ & =   \frac{1}{2\pi}\mE \left[ \mathsf{Ang}\left( X_1 - X_2, Y_1-X_2 \right) \right] + 
	\frac{1}{2\pi}\mE \left[ \mathsf{Ang} \left( X_1-Y_2, Y_1 - Y_2 \right) \right] \\[.5em]
	&~ - \frac{1}{4\pi} \mE \left[ \mathsf{Ang}\left( X_1-X_3, X_2 - X_3 \right) \right] - \frac{1}{4\pi}\mE \left[ \mathsf{Ang}\left( X_1-Y_1, X_2 - Y_1 \right)   \right] \\[.5em] 
	&~ - \frac{1}{4\pi} \mE \left[  \mathsf{Ang}\left( Y_1-Y_3, Y_2 - Y_3 \right)  \right] - \frac{1}{4\pi}\mE \left[  \mathsf{Ang}\left( Y_1-X_1, Y_2 - X_1 \right) \right].
	\end{align*}
	\begin{proof}
		Since the CvM-distance is invariant to the choice of $\vartheta_X$ and $\vartheta_Y$ (Theorem~\ref{Theorem: Multivaraite CvM Expression}), we may assume that $\vartheta_X=\vartheta_Y=1/2$ for simplicity. Then 
		\begin{align*}
		W_d^2  ~= ~ & \int_{\mathbb{S}^{d-1}} \int_{\mathbb{R}} \left(F_{\beta^\top X}(t) - F_{\beta^\top Y}(t)\right)^2 d\{ F_{\beta^\top X}(t)/2 + F_{\beta^\top Y}(t)/2\} d\lambda(\beta) \\[.5em]
		~= ~ & \mE \left[ \left( F_{\beta^\top X} (\beta^\top Z^\ast) \right)^2 \right] +  \mE_{\beta, Z^\ast} \left[ \left( F_{\beta^\top Y} (\beta^\top Z^\ast) \right)^2 \right] \\[.5em]
		& - 2 \mE \left[  F_{\beta^\top X} (\beta^\top Z^\ast) F_{\beta^\top Y} (\beta^\top Z^\ast)\right], \\[.5em]
		~ = ~ & (I) + (II) - 2 (III) \quad \text{(say),}
		\end{align*}
		where $Z^\ast \sim (1/2) P_X + (1/2) P_Y$. By the Fubini's theorem and the definition of $Z^\ast$, the first term $(I)$ has the identity
		\begin{align*}
		(I) & = ~\mE \left[ \ind(\beta^\top X_1 \leq \beta^\top Z^\ast, ~ \beta^\top X_2 \leq \beta^\top Z^\ast) \right] \\[.5em]
		& =~ \frac{1}{2} \mE \left[ \ind(\beta^\top X_1 \leq \beta^\top X_3, \beta^\top X_2 \leq \beta^\top X_3) \right]  + \frac{1}{2} \mE\left[ \ind(\beta^\top X_1 \leq \beta^\top Y_1, \beta^\top X_2 \leq \beta^\top Y_1) \right].
		\end{align*}
		Similarly, 
		\begin{align*}
		(II) & = ~ \mE \left[ \ind(\beta^\top Y_1 \leq \beta^\top Z^\ast, \beta^\top Y_2 \leq \beta^\top Z^\ast) \right] \\[.5em]
		& = ~ \frac{1}{2} \mE \left[ \ind(\beta^\top Y_1 \leq \beta^\top Y_3, \beta^\top Y_2 \leq \beta^\top Y_3) \right] + \frac{1}{2} \mE \left[ \ind(\beta^\top Y_1 \leq \beta^\top X_1, \beta^\top Y_2 \leq \beta^\top X_1) \right]
		\end{align*}
		and
		\begin{align*}
		(III) & = ~  \mE \left[ \ind(\beta^\top X_1 \leq \beta^\top Z^\ast, \beta^\top Y_1 \leq \beta^\top Z^\ast ) \right] \\[.5em]
		& = ~ \frac{1}{2} \mE \left[ \ind(\beta^\top X_1 \leq \beta^\top X_2, \beta^\top Y_1 \leq \beta^\top X_2) \right] + \frac{1}{2} \mE \left[ \ind(\beta^\top X_1 \leq \beta^\top Y_2, \beta^\top Y_1 \leq \beta^\top Y_2) \right].
		\end{align*}
		We then apply Lemma~\ref{Lemma: Integration over Unit Sphere (2 terms)} to obtain the desired result. 
	\end{proof}
\end{lemma}

\vskip 1em

Next we provide another expression for the CvM-statistic with a third-order kernel. 
\begin{lemma}[Another expression for the CvM-statistic] \label{Lemma: another expression for the U-statistic}
	Consider the kernel of order three
	\begin{align} \label{Eq: third-order kernel}
	& h^\star_{\text{\emph{CvM}}}(x_1,x_2,x_3; y_1,y_2,y_3) \\[.5em] \nonumber
	= ~& \frac{1}{2}\mE \big[ \{ \ind(\beta^\top x_1 \leq \beta^\top x_3 ) - \ind(\beta^\top y_1 \leq \beta^\top x_3 ) \} \cdot \{ \ind(\beta^\top x_2 \leq \beta^\top x_3 ) - \ind(\beta^\top y_2 \leq \beta^\top x_3 ) \} \big] \\[.5em] \nonumber
	+~ & \frac{1}{2}\mE \big[ \{ \ind(\beta^\top x_1 \leq \beta^\top y_3 ) - \ind(\beta^\top y_1 \leq \beta^\top y_3 ) \} \cdot \{ \ind(\beta^\top x_2 \leq \beta^\top y_3 ) - \ind(\beta^\top y_2 \leq \beta^\top y_3 ) \} \big]. \nonumber
	\end{align}
	Let us define the corresponding $U$-statistic by
	\begin{align*} 
	U^\star_{\text{\emph{CvM}}} & := \frac{1}{(m)_3 (n)_3} \sum_{i_1,i_2,i_3 = 1}^{m,\neq} \sum_{j_1,j_2,j_3 = 1}^{n,\neq} h^\star_{\text{\emph{CvM}}}(X_{i_1},X_{i_2},X_{i_3}; Y_{j_1},Y_{j_2},Y_{j_3}).
	\end{align*}
	Then $U^\star_{\text{\emph{CvM}}}$ is an unbiased estimator of $W_d^2$. Furthermore when $\beta^\top X$ and $\beta^\top Y$ are continuous for $\lambda$-almost all $\beta \in \mathbb{S}^{d-1}$, it is simplified as
	\begin{align} \label{Eq: another expression for U-statistic}
	U^\star_{\text{\emph{CvM}}} =  \frac{1}{(m)_2 (n)_2} \sum_{i_1,i_2 =1}^{m,\neq} \sum_{j_1,j_2 =1}^{n,\neq} h_{\text{\emph{CvM}}}(X_{i_1},X_{i_2}; Y_{j_1},Y_{j_2}).
	\end{align}
	\begin{proof}
		The unbiasedness property is trivial. We will show that (\ref{Eq: another expression for U-statistic}) holds under the given conditions. Since there is no tie with probability one, we have 
		\begin{align*}
		& \frac{1}{(m)_3} \sum_{i_1,i_2,i_3=1}^{m,\neq} \mE_{\beta} [ \ind(\beta^\top X_{i_1}  \leq \beta^\top X_{i_3})  \ind(\beta^\top X_{i_2}  \leq \beta^\top X_{i_3}) ] = \frac{1}{3}, \\[.5em]
		& \frac{1}{(n)_3} \sum_{j_1,j_2,j_3=1}^{n,\neq} \mE_{\beta} [ \ind(\beta^\top Y_{j_1}  \leq \beta^\top Y_{j_3})  \ind(\beta^\top Y_{j_2}  \leq \beta^\top Y_{j_3}) ] = \frac{1}{3}.
		\end{align*}
		Also the following identities are true
		\begin{align*}
		& \frac{2}{(m)_2 \cdot n} \sum_{i_1,i_2 = 1}^{m,\neq} \sum_{j=1}^n \mE_{\beta} [ \ind(\beta^\top X_{i_1}  \leq \beta^\top X_{i_2})  \ind(\beta^\top Y_{j}  \leq \beta^\top X_{i_2})]  \\[.5em]
		= ~ & 1 - \frac{1}{(m)_2 \cdot n }   \sum_{i_1,i_2 = 1}^{m,\neq} \sum_{j=1}^n \mE_{\beta} [ \ind(\beta^\top X_{i_1}  \leq \beta^\top Y_{j})  \ind(\beta^\top X_{i_2}  \leq \beta^\top Y_{j})] 
		\end{align*}
		and
		\begin{align*}
		& \frac{2}{m \cdot (n)_2} \sum_{i=1}^m \sum_{j_1,j_2 = 1}^{n,\neq} \mE_{\beta} [ \ind(\beta^\top Y_{j_1}  \leq \beta^\top Y_{j_2})  \ind(\beta^\top X_{i}  \leq \beta^\top Y_{j_2})]  \\[.5em]
		= ~ & 1 - \frac{1}{m \cdot (n)_2 }  \sum_{i=1}^m  \sum_{j_1,j_2 = 1}^{n,\neq}  \mE_{\beta} [ \ind(\beta^\top Y_{j_1}  \leq \beta^\top X_{i})  \ind(\beta^\top Y_{j_2}  \leq \beta^\top X_{i})]. 
		\end{align*}
		After expanding the terms in $h^\star_{\text{CvM}}$ and replacing the above identities, we can obtain
		\begin{align*}
		U^\star_{\text{CvM}}  & = ~  \frac{1}{(m)_2 \cdot n} \sum_{i_1,i_2 = 1}^{m,\neq} \sum_{j=1}^n \mE_{\beta} [ \ind(\beta^\top X_{i_1}  \leq \beta^\top Y_{j})  \ind(\beta^\top X_{i_2}  \leq \beta^\top Y_{j})]  \\[.5em]
		& +  \frac{1}{m \cdot (n)_2} \sum_{i=1}^m \sum_{j_1,j_2 =1}^{n,\neq} \mE_{\beta} [ \ind(\beta^\top Y_{j_1}  \leq \beta^\top X_{i})  \ind(\beta^\top Y_{j_2}  \leq \beta^\top X_{i})] - \frac{2}{3}, \\[.5em]
		& = ~	\frac{1}{(m)_2 (n)_2} \sum_{i_1,i_2 =1}^{m,\neq} \sum_{j_1,j_2 =1}^{n,\neq} h_{\text{CvM}}(X_{i_1},X_{i_2}; Y_{j_1},Y_{j_2}).
		\end{align*}
		Hence the result follows.
	\end{proof}
\end{lemma}

\vskip 1em

In the next lemma, we present an explicit expression for the variance of $U_{m,n}$, which will be used to bound the variance of the proposed statistic. 
\begin{lemma}[Theorem 2 of \cite{lee1990u} in Chapter 2] \label{Lemma: Variance of two-sample U-statistic}
	Let $U_{m,n}$ be a two-sample $U$-statistic based on a kernel having degrees $k_1$ and $k_2$. Then
	\begin{align*}
	\mV\left( U_{m,n} \right) = \sum_{c=0}^{k_1} \sum_{d=0}^{k_2} \frac{\binom{k_1}{c} \binom{k_2}{d} \binom{m-k_1}{k_1 -c} \binom{n_2 - k_2}{k_2 - d}}{\binom{n_1}{k_1} \binom{n_2}{k_2}}\sigma_{c,d}^2,
	\end{align*}
	where $\sigma_{c,d}^2$ is defined similarly as (\ref{Eq: Definition of xi_{cd}}).
\end{lemma}

\cite{hoeffding1952large} established a sufficient condition \citep[indeed the necessary condition proved by][]{chung2013exact} under which the permutation distribution approximates the corresponding unconditional distribution. The condition is stated as follows:

\begin{lemma}[Theorem 5.1 of \cite{chung2013exact}] \label{Lemma: Hoeffding's condition}
	Consider a sequence of random quantity $X^n$ taking values in a sample space $\mathcal{M}^n$ and suppose that $X^n$ has distribution $P^n$ in $\mathcal{M}^n$. Let $\mathcal{S}_N$ be a finite group of transformation from $\mathcal{M}^n$ onto itself. Let $T_n=T_n(X^n)$ be any real valued statistic and $\varpi_n$ be a random variable that is uniform on $\mathcal{S}_n$. Also, let $\varpi_n^\prime$ have the same distribution as $\varpi_n$, with $X^n$, $\varpi_n$ and $\varpi^\prime_n$ mutually independent. Suppose, under $P^n$,
	\begin{align} \label{Eq: Hoeffding's condition 1}
	(T_n(\varpi_nX^n), T_n(\varpi^\prime_n X^n)) \convD  (T, T^\prime),
	\end{align}
	where $T$ and $T^\prime$ are independent, each with common cumulative distribution function $R(\cdot)$. Here, $\varpi_n X^n$ denotes the composition of $X^n$ with $\varpi_n$ and $\varpi_n^\prime X^n$ is similarly defined. Let $\widehat{R}_n$ be the randomization distribution function of $T_n$ defined by 
	\begin{align*}
	\widehat{R}_n(t) = \frac{1}{\#|\mathcal{S}_n|} \sum_{\varpi_n \in \mathcal{S}_n} \ind\{ T_n(\varpi_n X^n) \leq t \},
	\end{align*}	
	where $\#|\mathcal{S}_n|$ denotes the cardinality of $\mathcal{S}_n$. Then, under $P^n$, 
	\begin{align} \label{Eq: Hoeffding's condition 2}
	\widehat{R}_n(t) \convP R(t),
	\end{align}
	for every $t$ which is a continuity point of $R(\cdot)$. Conversely, if (\ref{Eq: Hoeffding's condition 2}) holds for some limiting cumulative distribution function $R(\cdot)$ whenever $t$ is a continuity point, then (\ref{Eq: Hoeffding's condition 1}) holds.
\end{lemma}

\vskip 1em

\cite{chikkagoudar2014limiting} studied the limiting distribution of a two-sample $U$-statistic under contiguous alternatives for the univariate case \citep[see Theorem 3.1 therein and also][]{gregory1977large}. Here we extend their result to the multivariate case.

\vskip .5em

First we prepare for some notation. Let $P_{\theta_0}^N$ and $P_{\theta_0+bN^{-1/2}}^N$ denote  the joint distribution of the pooled samples $\{X_1,\ldots,X_m,Y_1,\ldots,Y_n\}$ under the null and contiguous alternative, respectively. Let $\lambda_{k,g}$ and $\phi_{k,g}(\cdot)$ be the eigenvalue and the corresponding eigenfunction satisfying the following integral equation
\begin{align*}
\mE[ g^\ast_{2,0}(x_1,X_2) \phi_{k,g}(X_2) ]= \lambda_{k,g} \phi_{k,g}(x_1) \quad \text{for} \ k = 1,2,\ldots,
\end{align*}
where $g^\ast_{2,0}(\cdot,\cdot)$ is defined in (\ref{Eq: Definition of g ast}) under the null hypothesis. For a sequence of random variables $Z_N$, we write $Z_N = o_{P_{\theta_0}^N}(1)$, if
\begin{align*}
\lim_{N \rightarrow \infty} P_{\theta_0}^N\big( |Z_N| \geq  \epsilon \big) = 0,
\end{align*}
for any $\epsilon > 0$. Then we have the following result.

\begin{lemma} \label{Lemma: Contiguous alternative}
	Recall the two-sample $U$-statistic, $U_{m,n}$, given in (\ref{Eq: General Two-Sample U-statistic}). Consider the same assumptions used in Theorem~\ref{Theorem: Asymptotic distribution under contiguous alternatives} and Theorem~\ref{Theorem: Two-Sample Degenerate Kernel}. Then under $P_{\theta_0+bN^{-1/2}}^N$ ,
	\begin{align*}
	N(U_{m,n} -\mE_{\theta_0}[U_{m,n}]) \convD ~ \frac{r(r-1)}{2\vartheta_X\vartheta_Y} \sum_{k=1}^\infty \lambda_{k,g} \{(\xi_k+ \vartheta_X^{1/2}a_{k,g})^2 - 1 \},
	\end{align*}
	where
	\begin{align*}
	a_{k,g} = \int_{\mathbb{R}^d} \big\langle b, 2\eta(x,\theta_0)p_{\theta_0}^{-1/2}(x) \big\rangle \phi_{k,g}(x) d P_{\theta_0}(x) .
	\end{align*}
	\begin{proof}
		Let us denote the likelihood ratio as
		\begin{align*}
		L_{N,h} = \frac{\prod_{i=1}^m p_{\theta_0}(X_i) \prod_{j=1}^n p_{\theta_0+ bN^{-1/2}}(Y_j)}{\prod_{i=1}^m p_{\theta_0}(X_i) \prod_{j=1}^n p_{\theta_0}(Y_j)} =  \frac{\prod_{j=1}^n p_{\theta_0+ bN^{-1/2}}(Y_j)}{\prod_{j=1}^n p_{\theta_0}(Y_j)}. 
		\end{align*}
		Then under the given conditions, one can establish
		\begin{align} \label{Eq: Likelihood ratio}
		\log L_{N,h} = \frac{1}{\sqrt{n}} \sum_{i=1}^n \langle h, \widetilde{\eta}(Y_i,\theta_0) \rangle  - \frac{1}{2} \langle h, I(\theta_0) h \rangle + o_{P_{\theta_0}^N}(1),
		\end{align}
		where $\widetilde{\eta}(x,\theta) = 2 \eta(x,\theta) / p_{\theta}^{1/2}(x)$ \citep[see Example 12.3.7 of][for details]{lehmann2006testing}. Then by Corollary 12.3.1 of \cite{lehmann2006testing}, $P_{\theta_0}^N$ and $P_{\theta_0+bN^{-1/2}}^N$ are mutually contiguous. 
		
		Without loss of generality, we assume that $\mE_{\theta_0} [U_{m,n}] = 0$ and denote the projection of $U_{m,n}$ under \emph{condition 2} in Theorem~\ref{Theorem: Two-Sample Degenerate Kernel} by 
		\begin{align*}
		\widehat{U}_{m,n} = & \frac{r(r-1)}{m(m-1)} \sum_{1 \leq i_1  < i_2 \leq m} g^\ast_{2,0}(X_{i_1},X_{i_2}) +  \frac{r(r-1)}{n(n-1)} \sum_{1 \leq j_1  < j_2 \leq m} g^\ast_{0,2}(Y_{j_1},Y_{j_2}) \\[.5em] 
		& + \frac{r^2}{mn} \sum_{i=1}^m \sum_{j=1}^n g^\ast_{1,1}(X_{i},Y_{j}). 
		\end{align*}
		Then as in Lemma 2.2 of \cite{chikkagoudar2014limiting}, it can be seen that 
		\begin{align*}
		N U_{m,n} = ~  N \widehat{U}_{m,n} + o_{P_{\theta_0}^N}(1),
		\end{align*}
		and the same approximation holds under $P_{\theta_0+ bN^{-1/2}}^N$ by contiguity. As a result, it is enough to study the limiting distribution of $N \widehat{U}_{m,n}$. 
		
		Now following the same steps in the proof of Theorem 3.1 in \cite{chikkagoudar2014limiting} and using (\ref{Eq: Likelihood ratio}), we can arrive at
		\begin{align*}
		 N \widehat{U}_{m,n}  \convD \frac{r(r-1)}{2\vartheta_X\vartheta_Y} \sum_{k=1}^\infty \lambda_{k,g} \{(\xi_k+ \vartheta_X^{1/2}a_{k,g})^2 - 1\},
		\end{align*}
		under $P_{\theta_0+ bN^{-1/2}}^N$. Hence the result follows. 
	\end{proof}
\end{lemma}

\vskip 2em

\section{Proofs} \label{Section: Proofs}

In addition to the notation given in the main text, we introduce further notation that will be used throughout this section.\\[1em]
\noindent \textbf{Notation.} 
We denote the probability measure under permutations by $\mP_{\varpi}$. The expectation and variance with respect to $\mP_{\varpi}$ are denoted by $\mE_{\varpi}$ and $\mV_{\varpi}$, respectively. We write the expectation with respect to the uniform probability measure $\lambda$ on $\mathbb{S}^{d-1}$ by $\mE_{\beta}$. The symbol $\#|A|$ stands for the cardinality of $A$. We denote the Kullback-Leibler divergence between two probability distributions $P$ and $Q$ by $\mathsf{KL}(P,Q)$. For $x,y \in \mathbb{R}$, we use $x\vee y$ and $x\wedge y$ to denote $\max\{x,y\}$ and $\min\{x,y\}$, respectively. Given a permutation $\varpi$ of $\{1,\ldots,N\}$ and the pooled samples $\{Z_1,\ldots,Z_{m+n}\} =\{X_1,\ldots,X_m,Y_1,\ldots,Y_n\}$, we may write $U_{\text{CvM}}(Z_{\varpi(1)},\ldots,Z_{\varpi(N)})$ or $U_{\text{CvM}}^{\varpi}$ to denote the CvM-statistic computed based on $\mathcal{X}_m = \{Z_{\varpi(1)},\ldots,Z_{\varpi(m)}\}$ and $\mathcal{Y}_n = \{Z_{\varpi(m+1)},\ldots,Z_{\varpi(m+n)} \}$. For the original permutation, which is $\varpi = \{1,\ldots,N\}$, we write $U_{\text{CvM}}$ or $U_{\text{CvM}}(Z_{1},\ldots,Z_{1})$ to denote the CvM-statistic computed based on $\mathcal{X}_m = \{Z_{1},\ldots,Z_{m}\}$ and $\mathcal{Y}_n = \{Z_{1},\ldots,Z_{m+n}\}$. The similar notation will be used for other test statistics. In general, we will write $\widetilde{h}$ to denote the symmetrized version of a kernel $h$ in the sense of (\ref{Eq: Symmetrization}). For any two real sequences $\{a_n\}$ and $\{b_n\}$, we write $b_n \gtrsim a_n$ or equivalently $a_n \lesssim b_n$ if there exists $C > 0$ such that $a_n \leq C b_n$ for each $n$. $c,C, C_0,C_1,C_2,C_3,C_4,C_5$ are some universal constants whose values may differ in different places of this section.

\vskip 2em

\subsection{Proof of Lemma~\ref{Lemma: Multivaraite CvM Characteristic Property}}
From the definition of $W_d^2$, it is clear to see that $W_d^2 \geq 0$ and it becomes zero if $P_X = P_Y$. For the other direction, we will show that if $W_d^2 = 0$, then $X$ and $Y$ have the same characteristic function:
\begin{align*}
\mE_X \left[ e^{it \beta^\top X}\right] = \mE_Y \left[ e^{it \beta^\top Y}\right] \quad \text{for all} \ (\beta, t) \in\mathbb{S}^{d-1} \times \mathbb{R},
\end{align*}
which implies $P_X = P_Y$.

\vskip 1em

\noindent \textbf{1. Univariate case} 

\noindent In the univariate case, $W^2= 0$ implies that $F_X(t) = F_Y(t)$ for $d\{\vartheta_X F_{X}(t) + \vartheta_Y F_{Y}(t) \}$-almost all $t$, hence we conclude $P_X = P_Y$ \citep[see also Lemma 4.1 of][]{lehmann1951consistency}.

%\noindent We begin with the univariate case. Suppose that $X$ and $Y$ are univariate random variables. We claim that the univariate CvM-distance has the characteristic property:
%\begin{align*}
%W^2 = \int_{\mathbb{R}} \left( F_{X}(t) - F_{Y}(t) \right)^2 d\{\vartheta_X F_{X}(t) + \vartheta_Y F_{Y}(t) \} = 0,
%\end{align*}
%if and only if  $F_{X}(t) = F_{Y}(t)$ for all $t \in \mathbb{R}$. Define $E$ to be the smallest set $E \subseteq \mathbb{R}$ such that
%\begin{align*}
%\int_E d F_X = \int_E dF_Y = 1.
%\end{align*}
%Since $\mP(X \in E^c) = \mP(Y \in E^c) = 0$, we only need to show that $W^2=0$ implies $F_X(t) = F_Y(t)$ for all $t \in E$. Suppose there exists $t_0 \in E$ such that $F_X(t_0) \neq F_Y(t_0)$. Suppose further that either $F_X$ or $F_Y$ is discontinuous at $t_0$. Then, there is a nonzero probability mass at $t_0$ and thus $W^2 > 0$. This contradicts the assumption. Next, suppose that both $F_X$ and $F_Y$ are continuous at $t_0$. Then $W^2 > 0 $ follows similarly from Lemma 4.1 of \cite{lehmann1951consistency}, which again contracts the assumption. Therefore, we conclude that $W^2  = 0$ implies that $X$ and $Y$ have the same distribution. 

\vskip 1em

\noindent \textbf{2. Multivariate case}

\noindent Recall that $\lambda(\cdot)$ is the uniform probability measure on $\mathbb{S}^{d-1}$. From the characteristic property of the univariate CvM-distance, $W_d^2 = 0$ implies that $\beta^\top X$ and $\beta^\top Y$ are identically distributed for $\lambda$-almost all $\beta \in \mathbb{S}^{d-1}$. Now, by continuity of the characteristic function, we conclude that 
\begin{align*}
\mE_X \left[ e^{it \beta^\top X}\right] = \mE_Y \left[ e^{it \beta^\top Y}\right] \quad \text{for all} \ (\beta, t) \in\mathbb{S}^{d-1} \times \mathbb{R}.
\end{align*}

\vskip 1em
\subsection{Proof of Lemma~\ref{Lemma: Integration over Unit Sphere (2 terms)}} \label{Section: Proof of Lemma: Integration over unit sphere}
%In the original paper by \cite{escanciano2006consistent}, the integration was with respect to $d \beta$ rather than $d \lambda(\beta)$:
%\begin{align*}
%\int_{\mathbb{S}^{d-1}} \ind(\beta^\top U_1 \leq 0) \ind(\beta^\top U_2 \leq 0) d\beta = c(d) \Bigg\{ \pi -  \text{{arccos}}\left( \frac{U_1^\top U_2}{\|U_1\| \|U_2\|} \right)\Bigg\},
%\end{align*}
%where $c(d) = \pi^{d/2-1} / \Gamma(d/2)$ as corrected in \cite{zhu2017projection}. Our expression in Lemma~\ref{Lemma: Integration over Unit Sphere (2 terms)} follows by the fact that the area of the surface of the $d$-dimensional unit ball is $2 \pi^{d/2} / \Gamma(d/2)$.

%\vskip 1em

%\noindent \textbf{Alternative proof}

%\vskip .5em 

Here we provide an alternative proof of Lemma~\ref{Lemma: Integration over Unit Sphere (2 terms)} based on the orthant probability for normal distribution. First we state a recent result on the bivariate normal distribution function presented by \cite{monhor2013inequalities}.
\begin{lemma}\citep[Theorem 4 of][]{monhor2013inequalities} \label{Lemma: Monhor}
	 Let $(\xi_1,\xi_2)^\top$ has a bivariate normal distribution with expectation $(\mu_1,\mu_2)^\top = (0,0)^\top$ and covariance matrix $[\sigma_{ij}]_{2\times 2}$ where $\sigma_{11} = \sigma_{22} = 1$ and $\sigma_{12} = \sigma_{21} = \rho$. Then for $0 < \rho < 1$ and $t > 0$,
	\begin{align} \label{Eq: Monhor1}
	\mP(\xi_1 \leq t, \xi_2 \leq t) \leq \Phi^2(t) + \frac{1}{2\pi} \exp \left( - \frac{t^2}{1+\rho} \right) \text{\emph{arcsin}}(\rho)
	\end{align}
	and
	\begin{align} \label{Eq: Monhor2}
	\mP( \xi_1 \leq t, \xi_2 \leq t) \geq  \Phi^2(t) + \frac{1}{2\pi} \exp \left( -t^2 \right) \text{\emph{arcsin}}(\rho).
	\end{align}
\end{lemma}
It is not difficult to see that a similar result can be obtained for $-1 < \rho \leq 0$ as
\begin{align} \label{Eq: Monhor3}
\mP(\xi_1 \leq t, \xi_2 \leq t) \leq \Phi^2(t) - \frac{1}{2\pi} \exp \left( - \frac{t^2}{1+\rho} \right) \text{{arcsin}}(-\rho)
\end{align}
and
\begin{align} \label{Eq: Monhor4}
\mP( \xi_1 \leq t, \xi_2 \leq t) \geq  \Phi^2(t) - \frac{1}{2\pi} \exp \left( -t^2 \right) \text{{arcsin}}(-\rho).
\end{align}
In fact, (\ref{Eq: Monhor1}), (\ref{Eq: Monhor2}), (\ref{Eq: Monhor3}) and (\ref{Eq: Monhor4}) hold for any $t$. By taking $t \rightarrow 0$ in the previous inequalities, we have
\begin{align} \label{Eq: orthant probability}
\mP(\xi_1 \leq 0, \xi_2 \leq 0) = \frac{1}{4} + \frac{1}{2\pi} \text{arcsin}(\rho) = \frac{1}{2} - \frac{1}{2\pi} \text{arccos}(\rho),
\end{align}
for any $-1 \leq \rho \leq 1$. The above identity is classical and can be found in different places \citep[e.g.][]{slepian1962one,childs1967reduction,xu2013comparative}.

\vskip 1em

Turning now to Lemma~\ref{Lemma: Integration over Unit Sphere (2 terms)}, let $\mathcal{Z}$ have a multivariate normal distribution with zero mean vector and identity covariance matrix. It is well-known that $\mathcal{Z}/\|\mathcal{Z}\|$ is uniformly distributed over $\mathbb{S}^{d-1}$ \citep[e.g. page 15 of][]{anderson2003introduction}. This leads to the key observation that
\begin{align} \label{Eq: Identity of beta and Z}
\int_{\mathbb{S}^{d-1}} \ind(\beta^\top U_1 \leq 0) \ind(\beta^\top U_2 \leq 0) d \lambda(\beta)  = \mE_{\mathcal{Z}} \left[ \ind(\mathcal{Z}^\top U_1 \leq 0) \ind(\mathcal{Z}^\top U_2 \leq 0) \right],
\end{align}
where $\mE_\mathcal{Z}[\cdot]$ is the expectation with respect to $\mathcal{Z}$. Note that $(\mathcal{Z}^\top U_1,\mathcal{Z}^\top U_2)^\top$ follows a bivariate normal distribution with correlation matrix $[\varrho_{ij}]_{2\times 2}$ where  $\varrho_{ij} =  U_i^\top U_j /\{ \|U_i\| \|U_j\|\}$. Using this connection and the equality (\ref{Eq: orthant probability}), we can obtain the closed-form expression for the left-hand side of (\ref{Eq: Identity of beta and Z}) and thus complete the proof.

\vskip 1em

\subsection{Proof of Theorem~\ref{Theorem: Multivaraite CvM Expression}}
Since $\beta^\top X$ and $\beta^\top Y$ are assumed to have continuous distribution functions, $\beta^\top X_1, \beta^\top X_2$ and $\beta^\top X_3$ have distinct values with probability one. This is also true for $\beta^\top Y_1, \beta^\top Y_2$ and $\beta^\top Y_3$. Therefore, the following identities hold for $\lambda$-almost all $\beta \in \mathbb{S}^{d-1}$.
\begin{equation} \label{Eq: CDF Identity (1)}
\begin{aligned}
& \int \left( F_{\beta^\top X} (t) \right)^2 d F_{\beta^\top X}(t) = \mP\left( \max\{\beta^\top X_1, \beta^\top X_2 \} \leq \beta^\top X_3 \right) =\frac{1}{3}, \\[.5em]
& \int \left( F_{\beta^\top Y} (t) \right)^2 d F_{\beta^\top Y}(t) = \mP\left( \max\{\beta^\top Y_1, \beta^\top Y_2 \} \leq \beta^\top Y_3 \right) =\frac{1}{3}, \\[.5em]
& \int \left( F_{\beta^\top X} (t) \right)^2 d F_{\beta^\top Y}(t) = \mP\left( \max\{\beta^\top X_1, \beta^\top X_2 \} \leq \beta^\top Y_1 \right), \\[.5em]
& \int \left( F_{\beta^\top Y} (t) \right)^2 d F_{\beta^\top X}(t) = \mP\left( \max\{\beta^\top Y_1, \beta^\top Y_2 \} \leq \beta^\top X_1 \right).
\end{aligned}
\end{equation}
Also note that 
\begin{align*}
& \mP\left(\max\{\beta^\top X_1, \beta^\top X_2 \} \leq \beta^\top Y_1 \right) + \mP\left(\max\{ \beta^\top X_1,  \beta^\top Y_1  \} \leq \beta^\top X_2 \right) \\[.5em]
+ ~  & \mP\left(\max\{\beta^\top X_2,\beta^\top Y_1 \} \leq \beta^\top X_1 \right) = 1
\end{align*}
and
\begin{align*}
&\mP\left(\max\{ \beta^\top X_1,  \beta^\top Y_1  \} \leq \beta^\top X_2 \right) = \mP\left(\max\{\beta^\top X_2,\beta^\top Y_1 \} \leq \beta^\top X_1 \right).
\end{align*}
These two identities give
\begin{equation}
\begin{aligned} \label{Eq: CDF Identity (2)}
\int F_{\beta^\top X}(t) F_{\beta^\top Y}(t) dF_{\beta^\top X}(t) & =  \mP\left(\max\{ \beta^\top X_1,  \beta^\top Y_1  \} \leq \beta^\top X_2 \right) \\[.5em]
& = \frac{1}{2} - \frac{1}{2} \mP \left( \max\{\beta^\top X_1, \beta^\top X_2 \}  \leq \beta^\top Y_1 \right).
\end{aligned}
\end{equation}
Similarly,
\begin{equation}
\begin{aligned} \label{Eq: CDF Identity (3)}
\int F_{\beta^\top X}(t) F_{\beta^\top Y}(t) dF_{\beta^\top Y}(t) & = \mP\left(\max\{ \beta^\top Y_1,  \beta^\top X_1  \} \leq \beta^\top Y_2 \right) \\[.5em] 
&= \frac{1}{2} - \frac{1}{2} \mP \left( \max\{\beta^\top Y_1, \beta^\top Y_2 \}  \leq \beta^\top X_1 \right).
\end{aligned}
\end{equation}
Now, combine (\ref{Eq: CDF Identity (1)}), (\ref{Eq: CDF Identity (2)}) and (\ref{Eq: CDF Identity (3)}) to have
\begin{align*}
& \int_{\mathbb{S}^{d-1}} \int_{\mathbb{R}} \left( F_{\beta^\top X} (t) - F_{\beta^\top Y}(t) \right)^2 d \{ \vartheta_X F_{\beta^\top X}(t) + \vartheta_Y F_{\beta^\top Y}(t)\} d\lambda(\beta) \\[.5em]
= ~ & \int_{\mathbb{S}^{d-1}} \mP\left( \max\{\beta^\top X_1,\beta^\top X_2 \} \leq \beta^\top Y_1  \right)  d\lambda(\beta) \\[.5em]
& + \int_{\mathbb{S}^{d-1}} \mP\left( \max\{\beta^\top Y_1,\beta^\top Y_2 \} \leq \beta^\top X_1  \right)  d\lambda(\beta) - \frac{2}{3}.
\end{align*}
Hence, 
\begin{align*}
W_d^2 & =~   \mE \left[ \ind(\beta^\top X_1 \leq \beta^\top Y_1, \beta^\top X_2 \leq \beta^\top Y_1 ) \right]  \\[.5em]
& ~ +  \mE \left[ \ind(\beta^\top Y_1 \leq \beta^\top X_1, \beta^\top Y_2 \leq \beta^\top X_1 ) \right] - \frac{2}{3}.
\end{align*}
Then apply Lemma~\ref{Lemma: Integration over Unit Sphere (2 terms)} to obtain the result.

\vskip 2em

\subsection{Proof of Theorem~\ref{Theorem: Asymptotic Null Distribution}}
We first show that $h$ is degenerate under $H_0$. Then apply the limit theorem for two-sample degenerate $U$-statistics~\citep{bhat1995theory}. 

\vskip 1em
\noindent \textbf{1. Degeneracy}

\noindent Recall the definition of the kernel $h_{\text{CvM}}$, i.e.
\begin{align*}
h_{\text{CvM}}(x_1,x_2;y_1,y_2) =  \frac{1}{3} - \frac{1}{2 \pi} \mathsf{Ang}(x_1 - y_1, x_2 - y_1)  - \frac{1}{2\pi} \mathsf{Ang} (y_1 - x_1, y_2 - x_1).
\end{align*}
Let us denote the symmetrized version of $h_{\text{CvM}}$ by $\widetilde{h}_{\text{CvM}}$ in the sense of (\ref{Eq: Symmetrization}), i.e.
\begin{align*}
\widetilde{h}_{\text{CvM}}(x_1,x_2; y_1,y_2) = \frac{1}{2} h_{\text{CvM}}(x_1,x_2;y_1,y_2) + \frac{1}{2} h_{\text{CvM}}(x_2,x_1; y_2,y_1).
\end{align*}

We first focus on the univariate case where $x_1,x_2,y_1,y_2 \in \mathbb{R}$ and make a connection to Lehmann's two-sample statistic \citep{lehmann1951consistency}. Let $\widetilde{h}^{(1)}_{\text{CvM}}$ denote the symmetrized $h_{\text{CvM}}$ for the univariate case, that can be written as
\begin{align*}
\widetilde{h}^{(1)}_{\text{CvM}} (x_1,x_2;y_1,y_2) := & \frac{1}{2} \Big\{ \ind(\max\{x_1,x_2\} \leq y_1 ) +  \ind(\max\{x_1,x_2\} \leq y_2 ) \\
& +  \ind(\max\{y_1,y_2\} \leq x_1 ) +  \ind(\max\{y_1,y_2\} \leq x_2 ) \Big\} - \frac{2}{3}.
\end{align*}
From the following identity, 
\begin{align*}
& \ind(\max\{x_1,x_2\} \leq \min \{y_1,y_2\}) + \ind(\max\{y_1,y_2\} \leq \min\{x_1,x_2\}) \\[.5em]
= ~ & \ind(\max\{x_1,x_2\} \leq y_1 ) +  \ind(\max\{x_1,x_2\} \leq y_2 ) \\[.5em]
+ ~ & \ind(\max\{y_1,y_2\} \leq x_1 ) +  \ind(\max\{y_1,y_2\} \leq x_2) - 1, 
\end{align*}
the univariate kernel has another expression as
\begin{align*}
2 \widetilde{h}^{(1)}_{\text{CvM}} (x_1,x_2;y_1,y_2) ~=~ & \ind(\max\{x_1,x_2\} \leq \min\{y_1,y_2\}) \\ 
+ ~ & \ind(\max\{y_1,y_2\} \leq \min\{x_1,x_2\}) - \frac{1}{3}.
\end{align*}
Thus $\widetilde{h}^{(1)}_{\text{CvM}}$ is equivalent to the kernel for Lehmann's two-sample statistic \citep{lehmann1951consistency}. Using this connection and the known results for Lehmann's two-sample statistic, we have 
\begin{equation} \label{Eq: Degeneracy of Lehmann Stat}
\begin{aligned}
& \widetilde{h}_{\text{CvM},{1,0}}^{(1)}( x_1) := \mE \left[ \widetilde{h}_{\text{CvM}}^{(1)}( x_1, X_2; Y_1, Y_2) \right] = 0,\\[.5em]
& \widetilde{h}_{\text{CvM},{0,1}}^{(1)}( y_1) := \mE \left[ \widetilde{h}_{\text{CvM}}^{(1)}( X_1, X_2; y_1, Y_2) \right] = 0,
\end{aligned}
\end{equation}
for any $x_1,y_1 \in \mathbb{R}$ under $H_0$. See Chapter 4 of \cite{bhat1995theory} for details.

\vskip 1em 

Let us now turn to multivariate cases where $x_1,x_2,y_1,y_2 \in \mathbb{R}^d$. By the definition of $\widetilde{h}_\text{CvM}$, we have
\begin{align*}
\widetilde{h}_{\text{CvM}}(x_1,x_2,y_1,y_2) = 
\int_{\mathbb{S}^{d-1}} \widetilde{h}_{\text{CvM}}^{(1)}(\beta^\top x_1, \beta^\top x_2; \beta^\top y_1,\beta^\top x_2) d\lambda(\beta).
\end{align*}
Now the Fubini's theorem combined with (\ref{Eq: Degeneracy of Lehmann Stat}) gives
\begin{align*}
\mE \left[ \widetilde{h}_{\text{CvM}}^{(1)}( \beta^\top x_1, \beta^\top X_2; \beta^\top Y_1, \beta^\top Y_2) \right]  = 
\mE \left[ \widetilde{h}_{\text{CvM}}^{(1)}( \beta^\top X_1, \beta^\top X_2; \beta^\top y_1, \beta^\top Y_2) \right] = 0,
\end{align*} 
for $\lambda$-almost all $\beta \in \mathbb{S}^{d-1}$. As a consequence, it is seen that
\begin{align*}
\widetilde{h}_{\text{CvM},{1,0}}(x_1) & := \mE \left[\widetilde{h}_{\text{CvM}} (x_1,X_2;Y_1,Y_2) \right] \\[.5em]
& = \int_{\mathbb{S}^{d-1}}\mE 
\left[ \widetilde{h}_{\text{CvM}}^{(1)}( \beta^\top x_1, \beta^\top X_2; \beta^\top Y_1, \beta^\top Y_2) \right]  d\lambda(\beta) = 0, \\[1em]
\widetilde{h}_{\text{CvM},{0,1}} (y_1) & := \mE \left[\widetilde{h}_{\text{CvM}} (X_1,X_2;y_1,Y_2) \right] \\[.5em]
& = \int_{\mathbb{S}^{d-1}} \mE \left[ \widetilde{h}_{\text{CvM}}^{(1)}
( \beta^\top X_1, \beta^\top X_2; \beta^\top y_1, \beta^\top Y_2) \right]  d\lambda(\beta) = 0.
\end{align*}

On the other hand, 
\begin{align*}
\widetilde{h}_{\text{CvM},2,0}(x_1,x_2) & :=  \mE\left[\widetilde{h}_{\text{CvM}}(x_1,x_2;Y_1,Y_2)\right] \\[.5em]
& = \frac{1}{2}  \int_{\mathbb{S}^{d-1}} \left( 1 - F_{\beta^\top X}(\max\{\beta^\top x_1,\beta^\top x_2\}) \right)^2 d\lambda(\beta) \\[.5em]
& ~+ \frac{1}{2}  \int_{\mathbb{S}^{d-1}} F_{\beta^\top X}^2(\min\{\beta^\top x_1,\beta^\top x_2\}) d\lambda(\beta) - \frac{1}{6}, \\[.5em]
\widetilde{h}_{\text{CvM},{0,2}}(y_1,y_2) & :=  \mE\left[\widetilde{h}_{\text{CvM}}(X_1,X_2;y_1,y_2)\right], \\[.5em]
& = \frac{1}{2} \int_{\mathbb{S}^{d-1}}  \left( 1 - F_{\beta^\top Y}(\max\{\beta^\top y_1,\beta^\top y_2\}) \right)^2 d\lambda(\beta) \\[.5em]
& ~ +  \frac{1}{2} \int_{\mathbb{S}^{d-1}} F_{\beta^\top Y}^2(\min\{\beta^\top y_1,\beta^\top y_2\}) d\lambda(\beta) - \frac{1}{6}, \\[.5em]
\widetilde{h}_{\text{CvM},{1,1}}(x_1,y_1) & :=   \mE\left[\widetilde{h}_{\text{CvM}}(x_1,X_2;y_1,Y_2)\right] \\[.5em]
& = -\frac{1}{2} \widetilde{h}_{\text{CvM},2,0}(x_1,y_1).
\end{align*}
Note that $\widetilde{h}_{\text{CvM},{2,0}}(x_1,x_2) \neq 0$ for some $(x_1,x_2)$. For example, when $x_1 = x_2$, it is seen that 
\begin{align*}
\frac{1}{2}\big\{ 1 - F_{\beta^\top X}( \beta^\top x_1 ) \big\}^2 + \frac{1}{2}F_{\beta^\top X}^2( \beta^\top x_1)  - \frac{1}{6}  ~ \geq ~ \frac{1}{12} \quad \text{for all $\beta \in \mathbb{S}^{d-1}$},
\end{align*}
which implies $\widetilde{h}_{\text{CvM},{2,0}}(x_1,x_1) \geq 1/12$. By the continuity of $\widetilde{h}_{\text{CvM},{2,0}}$ at $(x_1,x_1)$, there exist a set with nonzero measure such that $\widetilde{h}_{\text{CvM},{2,0}}(x_1,x_2)>0$. Therefore, we conclude that $\widetilde{h}_{\text{CvM}}$ (and $h_{\text{CvM}}$) has degeneracy of order one under $H_0$.

\vskip 1em

\noindent \textbf{2. Limiting distribution of the $U$-statistic}

\noindent To obtain the limiting null distribution of $U_{\text{CvM}}$, we apply the result given in Chapter 3 of \cite{bhat1995theory} to have
\begin{align*}
N U_{\text{CvM}} \convD  \frac{1}{\vartheta_X} \sum_{k=1}^{\infty} \lambda_k (\xi_k^2 - 1) + \frac{1}{\vartheta_Y} \sum_{k=1}^{\infty} \lambda_k (\xi_k^{\prime 2} - 1) - \frac{2}{\sqrt{\vartheta_X \vartheta_Y}} \sum_{k=1}^\infty \lambda_k \xi_k \xi_k^\prime,
\end{align*}
where $\xi_k, \xi_k^\prime \overset{i.i.d.}{\sim} N(0,1)$. Based on the observation that
\begin{align*}
\sqrt{\vartheta_Y} \xi_k  - \sqrt{\vartheta_X} \xi_k^\prime \sim N(0,1),
\end{align*}
the result follows.

\vskip 2em

\subsection{Proof of Theorem~\ref{Theorem: Asymptotic distribution under fixed alternatives}}
Let us write $\widetilde{h}_{\text{CvM},1,0}(x) = \mE [\widetilde{h}_{\text{CvM}}(x,X_1;Y_1,Y_2)]$ and $\widetilde{h}_{\text{CvM},0,1}(y) = \mE [\widetilde{h}_{\text{CvM}}(X_1,X_2;y,Y_1)]$. By Hoeffding's decomposition of a two-sample $U$-statistic \citep[e.g.~page 40 of][]{lee1990u}, the CvM-statistic can be approximated by 
\begin{align*}
U_{\text{CvM}} - W_d^2 =  \frac{2}{m}\sum_{i=1}^m \widetilde{h}_{\text{CvM},1,0}(X_i) + \frac{2}{n} \sum_{j=1}^n \widetilde{h}_{\text{CvM},0,1}(Y_j) + O_{\mathbb{P}}(N^{-1}).
\end{align*}
Then the result follows by the central limit theorem.

%\vskip 2em
%
%\subsection{Proof of Theorem~\ref{Theorem: Asymptotic distribution under contiguous alternatives}}
%
%This result is a direct consequence of Lemma~\ref{Lemma: Contiguous alternative} with $r=2$. 
%

\vskip 2em

\subsection{Proof of Theorem~\ref{Theorem: Critical value of Permutation test}}
Under the null hypothesis, we need to verify the conditions given in Theorem~\ref{Theorem: Two-Sample Degenerate Kernel}. Indeed, these conditions are satisfied with $r=2$ as in the proof of Theorem~\ref{Theorem: Asymptotic Null Distribution}. Hence, the result follows under $H_0$.

Next, we focus on the alternative hypothesis. The proof consists of two steps. In the first step, we show that (\ref{Eq: Coupling Goal}) is satisfied for the CvM-statistic. In the second step, we show that the two CvM-statistics --- one based on $i.i.d.$ samples from $\frac{m}{N}P_X + \frac{n}{N} P_Y$ and the other based on $i.i.d.$ samples from $\vartheta_X P_X + \vartheta_Y P_Y$ --- have the same limiting distribution under the given conditions.

\vskip 1em

\noindent \textbf{$\bullet$ Step 1.} 
\vskip .3em

\noindent For the first step, we use the coupling argument (Algorithm~\ref{Alg: Coupling}) to show that the difference between the two CvM-statistics --- one is based on the randomly permuted original samples and the other is based on the corresponding coupled $i.i.d.$ samples --- is asymptotically negligible. Formally, we state the result in the following lemma.

\begin{lemma}[Coupling for the CvM-statistic] \label{Lemma: Coupling for the CvM-statistic}
	Consider the two sets of samples $\{Z_1,\ldots, Z_N \}$ and $\{ \overline{Z}_{\varpi_0(1)}, \ldots, \overline{Z}_{\varpi_0(N)} \}$ from Algorithm~\ref{Alg: Coupling} and their random permutations $\{Z_{\varpi(1)}, \ldots, Z_{\varpi(N)}\}$ and $\{\overline{Z}_{\varpi(\varpi_0(1))}, \ldots, \overline{Z}_{\varpi(\varpi_0(N))} \}$. Then we have
	\begin{align} \label{Eq: Difference between the original and coupled U-statistics}
	NU_{\text{\emph{CvM}}}(Z_{\varpi(1)}, \ldots, Z_{\varpi(N)}) - N U_{\text{\emph{CvM}}}( \overline{Z}_{\varpi(\varpi_0(1))}, \ldots, \overline{Z}_{\varpi(\varpi_0(N))} ) \convP 0.
	\end{align}
	%where $U_{\text{CvM}}$ is the CvM-statistic in (\ref{Eq: Definition of U-statistic}). 
	\begin{proof}
		Using the result in Lemma~\ref{Lemma: another expression for the U-statistic}, we work with the third-order kernel $h^\star_{\text{CvM}}$ in (\ref{Eq: third-order kernel}). First notice that the expectations of both $U_{\text{CvM}}(Z_{\varpi(1)}, \ldots, Z_{\varpi(N)})$ and $U_{\text{CvM}}( \overline{Z}_{\varpi(\varpi_0(1))}, \ldots, \overline{Z}_{\varpi(\varpi_0(N))})$ are zero. To see this, putting $\mathcal{E} = \{ \beta, Z_1,\ldots,Z_N, \varpi(2), \varpi(3), \varpi(m+2) \}$, write
		\begin{align*}
		f(\mathcal{E}) = \mE_{\varpi(1),\varpi(m+1)} \big[ \{ \ind(\beta^\top Z_{\varpi(1)} \leq \beta^\top Z_{\varpi(3)}) - \ind(\beta^\top Z_{\varpi(m+1)} \leq \beta^\top Z_{\varpi(3)}) \} \big|~ \mathcal{E} \big]
		\end{align*}
		and note that $f(\mathcal{E})$ is zero for any $\mathcal{E}$. As a result, the law of total expectation gives
		\begin{align*}
		& \mE \big[ \{ \ind(\beta^\top Z_{\varpi(1)} \leq \beta^\top Z_{\varpi(3)}) - \ind(\beta^\top Z_{\varpi(m+1)} \leq \beta^\top Z_{\varpi(3)}) \}  \\[.5em]
		& ~~~~~~ \times \{ \ind(\beta^\top Z_{\varpi(2)} \leq \beta^\top Z_{\varpi(3)}) - \ind(\beta^\top Z_{\varpi(m+2)} \leq \beta^\top Z_{\varpi(3)}) \} \big] \\[.5em]
		=  ~ &  \mE \big[ f(\mathcal{E}) \times \{ \ind(\beta^\top Z_{\varpi(2)} \leq \beta^\top Z_{\varpi(3)}) - \ind(\beta^\top Z_{\varpi(m+2)} \leq \beta^\top Z_{\varpi(3)}) \} \big] = 0.
		\end{align*}
		By applying the same logic to the other terms, it is clear that the expectations of both test statistics are zero.

		Based on the previous observation, it now suffices to show that 
		\begin{align} \label{Eq: the expected value of the squared difference}
		\mE \left[ \{NU_{\text{CvM}}(Z_{\varpi(1)}, \ldots, Z_{\varpi(N)}) - N U_{\text{CvM}}( \overline{Z}_{\varpi(\varpi_0(1))}, \ldots, \overline{Z}_{\varpi(\varpi_0(N))}) \}^2 \right] = o(1)
		\end{align}
		to establish (\ref{Eq: Difference between the original and coupled U-statistics}). For simplicity, denote 
		\begin{align*}
		& v_{\varpi} (i_1,i_2,i_3;j_1,j_2,j_3) \\[.5em] 
		= ~& h^\star_{\text{CvM}}(Z_{\varpi(i_1)},Z_{\varpi(i_2)},Z_{\varpi(i_3)};Z_{\varpi(j_1+m)},Z_{\varpi(j_2+m)},Z_{\varpi(j_3+m)}) \\[.5em] 
		&- h^\star_{\text{CvM}}(\overline{Z}_{\varpi(\varpi_0(i_1))},\overline{Z}_{\varpi(\varpi_0(i_2))},\overline{Z}_{\varpi(\varpi_0(i_3))};\overline{Z}_{\varpi(\varpi_0(j_1+m))},\overline{Z}_{\varpi(\varpi_0(j_2+m))},\overline{Z}_{\varpi(\varpi_0(j_3+m))}).
		\end{align*}
		Then the square of $NU_{\text{CvM}}(Z_{\varpi(1)}, \ldots, Z_{\varpi(N)}) - N U_{\text{CvM}}( \overline{Z}_{\varpi(\varpi_0(1))}, \ldots, \overline{Z}_{\varpi(\varpi_0(N))} )$ can be written as
		\begin{align*}
		\mathcal{D}_{m,n} & := ~  \frac{N^2}{(m)_3^2 (n)_3^2} ~ \times \\[.5em]
		& ~ \sum_{i_1,i_2,i_3=1}^{m,\neq} \sum_{j_1,j_2,j_3=1}^{n,\neq} \sum_{i_1^\prime,i_2^\prime,i_3^\prime=1}^{m,\neq} \sum_{j_1^\prime,j_2^\prime,j_3^\prime=1}^{n,\neq} v_{\varpi} (i_1,i_2,i_3;j_1,j_2,j_3) v_{\varpi} (i_1^\prime,i_2^\prime,i_3^\prime;j_1^\prime,j_2^\prime,j_3^\prime).
		\end{align*}
		Further write
		\begin{align} \label{Eq: Definition of I and J}
		\mathcal{I}_3 = \{i_1,i_2,i_3\} \cap \{i_1^\prime, i_2^\prime, i_3^\prime \} \quad \text{and} \quad 	\mathcal{J}_3 = \{j_1,j_2,j_3\} \cap \{j_1^\prime, j_2^\prime, j_3^\prime\}.
		\end{align}
		By the law of total expectation, it can be seen that 
		\begin{align*}
		\mE \left[ v_{\varpi} (i_1,i_2,i_3;j_1,j_2,j_3) v_{\varpi} (i_1^\prime,i_2^\prime,i_3^\prime;j_1^\prime,j_2^\prime,j_3^\prime) | ~ \beta,Z_1,\ldots,Z_N,\overline{Z}_1,\ldots,\overline{Z}_N \right] = 0
		\end{align*}
		whenever $\#|\mathcal{I}_3| + \#|\mathcal{J}_3| \leq 1$. Thus the unconditional expectation is also zero in these cases. Next consider the cases where $\#|\mathcal{I}_3| + \#|\mathcal{J}_3| = 2$. More specifically, we split the cases into
		\begin{itemize}%[\hspace{1.2cm}]
		\item $\mathcal{C}_{a} = \{ i_1,\ldots ,i_3^\prime, j_1,\ldots,j_3^\prime :  \#|\mathcal{I}_3| = 2$ and $\#| \mathcal{J}_3 | = 0\}$,
		\item $\mathcal{C}_{b} = \{ i_1,\ldots,i_3^\prime, j_1,\ldots,j_3^\prime :  \#|\mathcal{I}_3| = 0$ and $\#| \mathcal{J}_3 | = 2\}$,
		\item $\mathcal{C}_{c} = \{ i_1,\ldots ,i_3^\prime, j_1,\ldots,j_3^\prime :  \#|\mathcal{I}_3| = 1$ and $\#| \mathcal{J}_3 | = 1\}$.
		\end{itemize}	
		Suppose there are $B_1$ different observations between 
		\begin{align*}
		\{Z_{\varpi(1)},\ldots, Z_{\varpi(m)} \} \quad  \text{and} \quad \{\overline{Z}_{\varpi(\varpi_0(1))},\ldots, \overline{Z}_{\varpi(\varpi_0{(m)})} \}
		\end{align*}
		and $B_2$ different observations between 
		\begin{align*}
		\{Z_{\varpi(m+1)},\ldots, Z_{\varpi(m+n)} \} \quad \text{and} \quad \{\overline{Z}_{\varpi(\varpi_0(m+1))},\ldots, \overline{Z}_{\varpi(\varpi_0{(m+n)})}\}.
		\end{align*}
		Hence, we have $D=B_1+B_2$ different observations in total between the original samples and the coupled samples. In these cases, it can be seen that 
		\begin{align*}
		&\#| \mathcal{C}_a | ~ \lesssim ~ B_1 m^3 n^6 + B_2 m^4 n^5, \\[.5em]
		&\#| \mathcal{C}_b | ~ \lesssim ~  B_1 m^5 n^4 + B_2 m^6 n^3, \\[.5em]
		&\#| \mathcal{C}_c | ~ \lesssim ~ B_1 m^4 n^5 + B_2 m^5n^4.
		\end{align*}
%		Then for Case (a), there are at least $\comb{3}{2} \times 3! \times (m-B_1)(m-B_1-1)(m-B_1-2)(m-B_1-3)(n-B_2)(n-B_2-1)(n-B_2-2)(n-B_2-3)(n-B_2-4)(n-B_2-5)$ numbers of zero product terms in $(\ast\ast)$ out of $\comb{3}{2} \times 3! \times m(m-1)(m-2)(m-3)n(n-1)(n-2)(n-3)(n-4)(n-5)$ total number of product terms. 
%		
%		Similarly for Case (b), there are at least $\comb{3}{2} \times 3! \times (m-B_1)(m-B_1-1)(m-B_1-2)(m-B_1-3)(m-B_1-4)(m-B_1-5)(n-B_2)(n-B_2-1)(n-B_2-2)(n-B_2-3)$ numbers of zero product terms in $(\ast\ast)$ out of $\comb{3}{2} \times 3! \times m(m-1)(m-2)(m-3)(m-4)(m-5)n(n-1)(n-2)(n-3)$ total number of product terms. 
%		
%		For Case (c), there are at least $\comb{3}{1} \times 3! \times (m-B_1)(m-B_1-1)(m-B_1-2)(m-B_1-3)(m-B_1-4) \times \comb{3}{1} \times 3! \times (n-B_2)(n-B_2-1)(n-B_2-2)(n-B_2-3)(n-B_2-4)$ numbers of zero product terms in $(\ast\ast)$ out of $\comb{3}{1} \times 3! \times m(m-1)(m-2)(m-3)(m-4) \times \comb{3}{1} \times 3! \times n(n-1)(n-2)(n-3)(n-4)$ total number of product terms.
		
%		To summarize it, the number of non-zero product terms for each case is bounded by
%		\begin{align*}
%		&\text{Case 2(a)} ~ \lesssim ~ B_1 m^3 n^6 + B_2 m^4 n^5, \\[.5em]
%		&\text{Case 2(b)} ~ \lesssim ~  B_1 m^5 n^4 + B_2 m^6 n^3, \\[.5em]
%		&\text{Case 2(c)} ~ \lesssim ~ B_1 m^4 n^5 + B_2 m^5n^4.
%		\end{align*}
		\noindent Also note that the number of the other cases such that $\#|\mathcal{I}_3| + \#|\mathcal{J}_3| > 2$ are at most $O(N^9)$. Since $\mE[B_1] = O(\sqrt{N}),\mE[B_2] = O(\sqrt{N})$ and the kernel $v_{\varpi}$ is bounded, we can conclude that 
		\begin{align*}
		\mE [\mathcal{D}_{m,n}] = O\left( \frac{1}{\sqrt{N}}\right) = o(1).
		\end{align*}
		This shows (\ref{Eq: the expected value of the squared difference}) and thus completes the proof. 
	\end{proof}
\end{lemma}

\vskip 1em

\noindent \textbf{$\bullet$ Step 2.} 

\vskip .3em

\noindent From Lemma~\ref{Lemma: Coupling for the CvM-statistic}, we have established that $NU_{\text{CvM}}(Z_{\varpi(1)}, \ldots, Z_{\varpi(N)})$ and $N U_{\text{CvM}}( \overline{Z}_{\varpi(\varpi_0(1))}, \\ \ldots, \overline{Z}_{\varpi(\varpi_0(N))})$ have the same limiting distribution. Note that $\overline{Z}_{\varpi(\varpi_0(1))}, \ldots, \overline{Z}_{\varpi(\varpi_0(N))}$ are sampled from $\frac{m}{N} P_X + \frac{n}{N} P_Y$. Next, we will further show that the limiting distribution of $N U_{\text{CvM}}$ based on samples from $\frac{m}{N} P_X + \frac{n}{N} P_Y$ and that based on samples from $\vartheta_X P_X + \vartheta_Y P_Y$ are equivalent when $\frac{m}{N} \rightarrow \vartheta_X$ and $\frac{n}{N} \rightarrow \vartheta_Y$ as $N \rightarrow \infty$ where $0 < \vartheta_X, \vartheta_Y < 1$. Since the limiting distribution of $N U_{\text{CvM}}$ is the weighted sum of independent chi-square statistics, the limiting distribution is decided by the weights, which are eigenvalues of the integral equation associated with the kernel. Using the symmetrized kernel $\widetilde{h}_{\text{CvM}}(x_1,x_2;y_1,y_2)$, define
\begin{align*}
\widetilde{h}_{\text{CvM},2,0}^{(m,n)}(x_1,x_2) = \int \widetilde{h}_{\text{CvM}}(x_1,x_2;y_1,y_2)  dH_{m,n}(y_1) dH_{m,n}(y_2)
\end{align*}
where $H_{m,n} = \frac{m}{N} P_X + \frac{n}{N} P_Y$. Similarly, define 
\begin{align*}
\widetilde{h}_{\text{CvM},2,0}(x_1,x_2) = \int \widetilde{h}_{\text{CvM}}(x_1,x_2;y_1,y_2) dH(y_1) dH(y_2)
\end{align*}
where $H = \vartheta_X P_X + \vartheta_Y P_Y$. Then it can be seen that 
\begin{align} \label{Eq: Upper bound of the difference}
|\widetilde{h}_{\text{CvM},2,0}^{(m,n)}(x_1,x_2) - \widetilde{h}_{\text{CvM},2,0}(x_1,x_2)| \leq \sum_{\substack{i=0,j=0 \\ i+j=4}}^4 \Big| \left(\frac{m}{N}\right)^i \left(\frac{n}{N}\right)^j - \vartheta_X^i \vartheta_Y^j \Big|,
\end{align}
by the boundedness of $\widetilde{h}_{\text{CvM}}$, i.e. $|\widetilde{h}_{\text{CvM}}| \leq 1$. Let $\{\lambda_i^{(m,n)} \}_{i=1}^\infty$ and $\{ \phi_i^{(m,n)}(\cdot) \}_{i=1}^\infty$ be eigenvalues and square integrable eigenfunctions of the integral equation
\begin{align} \label{Eq: Integral Equation 2}
\int \widetilde{h}_{\text{CvM},2,0}^{(m,n)} (x_1,x_2) \phi_i^{(m,n)}(x_2) d H_{m,n}(x_2) = \lambda_i^{(m,n)}  \phi_i^{(m,n)}(x_1).
\end{align}
Let us denote their limits by $\lambda_i^\ast = \lim_{N \rightarrow \infty}\lambda_i^{(m,n)}$ and $\phi_i^\ast (z) = \lim_{N \rightarrow \infty}\phi_i^{(m,n)}(z)$. In the next lemma, we will show that $\lambda_i^\ast$ and $\phi_i^\ast(z)$ satisfy the integral equation
\begin{align} \label{Eq: Integral Equation}
\int \widetilde{h}_{\text{CvM},2,0} (x_1,x_2) \phi_i^\ast(x_2) d H(x_2) = \lambda_i^\ast  \phi_i^\ast(x_1)
\end{align}
for all $x_1$. Thus the limits are the eigenvalues and the eigenfunctions of (\ref{Eq: Integral Equation}). 

\begin{lemma} \label{Lemma: Eigenvalues and Eigenfunctions}
	Let us denote the eigenvalues and the eigenfunctions of the integral equation in (\ref{Eq: Integral Equation 2}) by $\{\lambda_i^{(m,n)}\}_{i=1}^\infty$ and $\{\phi_i^{(m,n)}(\cdot) \}_{i=1}^\infty$, respectively. Further denote their limits by $\lambda_i^\ast = \lim_{N \rightarrow \infty}\lambda_i^{(m,n)}$ and $\phi_i^\ast (z) = \lim_{N \rightarrow \infty}\phi_i^{(m,n)}(z)$. Then $\{\lambda_i^\ast\}_{i=1}^\infty$ and $\{\phi_i^\ast(\cdot)\}_{i=1}^{\infty}$ are the eigenvalues and the eigenfunctions of the integral equation in (\ref{Eq: Integral Equation}). In addition, we have 
	\begin{align*}
	\sum_{i=1}^\infty \left( \lambda_i^{(m,n)} \right)^2 \rightarrow \sum_{i=1}^\infty \lambda_i^2 \quad \text{as $N \rightarrow \infty$.}
	\end{align*}
	\begin{proof}
		Note that 
		\begin{align*}
		& \Big| \int \widetilde{h}_{\text{CvM},2,0}^{(m,n)} (x_1,x_2) \phi_i^{(m,n)}(x_2) d H_{m,n}(x_2) - \int \widetilde{h}_{\text{CvM},2,0} (x_1,x_2)  \phi_i^{(m,n)}(x_2) d H(x_2) \Big| \\[.5em]
		\leq ~ &  \Big| \int \widetilde{h}_{\text{CvM},2,0}^{(m,n)} (x_1,x_2) \phi_i^{(m,n)}(x_2) d H_{m,n}(x_2) - \int \widetilde{h}_{\text{CvM},2,0}^{(m,n)}  (x_1,x_2)  \phi_i^{(m,n)}(x_2) d H(x_2) \Big|  \\[.5em]
		& +  \Big| \int \widetilde{h}_{\text{CvM},2,0}^{(m,n)} (x_1,x_2) \phi_i^{(m,n)}(x_2) d H(x_2) - \int \widetilde{h}_{\text{CvM},2,0} (x_1,x_2)  \phi_i^{(m,n)}(x_2) d H(x_2) \Big| \\[.5em]
		= ~ &  (I) + (II) \quad \text{(say).}
		\end{align*}
		For $(I)$, we have 
		\begin{align*}
		(I) = ~ & \Bigg| \left( \frac{m}{N}  - \vartheta_X \right) \int \widetilde{h}_{\text{CvM},2,0}^{(m,n)} (x_1,x_2) \phi_i^{(m,n)}(x_2) d P_X(x_2) \\[.5em]
		&  + \left( \frac{n}{N}  - \vartheta_Y \right) \int \widetilde{h}_{\text{CvM},2,0}^{(m,n)} (x_1,x_2) \phi_i^{(m,n)}(x_2) d P_Y(x_2)   \Bigg| \\[.5em]
		\leq ~ & \Big|\frac{m}{N}  - \vartheta_X \Big|  \int | \widetilde{h}_{\text{CvM},2,0}^{(m,n)} (x_1,x_2) \phi_i^{(m,n)}(x_2)  | d P_X(x_2) \\[.5em]
		&  + \Big| \frac{n}{N}  - \vartheta_Y \Big|  \int  | \widetilde{h}_{\text{CvM},2,0}^{(m,n)} (x_1,x_2) \phi_i^{(m,n)}(x_2) | d P_Y(x_2)  \\[.5em]
		\leq ~ &  \Big|\frac{m}{N}  - \vartheta_X \Big| \sqrt{ \int  \left(  \phi_i^{(m,n)}(x_2) \right)^2 d P_X(x_2)} + \Big|\frac{n}{N}  - \vartheta_Y \Big| \sqrt{ \int  \left(  \phi_i^{(m,n)}(x_2) \right)^2 d P_Y(x_2)}
		\end{align*}
		where the last inequality is due to Cauchy-Schwarz inequality and the boundedness of the kernel. Since $\phi_i^{(m,n)}$ is a normalized function, i.e.
		\begin{align*}
		&  \int \left(  \phi_i^{(m,n)}(x_2) \right)^2 d H_{m,n}(x_2)  \\[.5em]
		= ~ & \frac{m}{N}\int  \left(  \phi_i^{(m,n)}(x_2) \right)^2 d P_X(x_2) +	\frac{n}{N} \int  \left(  \phi_i^{(m,n)}(x_2) \right)^2 d P_Y(x_2)  = 1,
		\end{align*}
		we obtain the upper bound
		\begin{align} \label{Eq: Bound of the Sum}
		\int  \left(  \phi_i^{(m,n)}(x_2) \right)^2 d P_X(x_2) + \int  \left(  \phi_i^{(m,n)}(x_2) \right)^2 d P_Y(x_2)  \leq \frac{N}{\min\{m,n\}}.
		\end{align}
		Using this, $(I)$ is further bounded by 
		\begin{align*}
		(I) \leq ~ \sqrt{\frac{N}{\min\{m,n\}}} \left( \Big| \frac{m}{N}  - \vartheta_X \Big|  + \Big| \frac{n}{N}  - \vartheta_Y \Big|   \right).
		\end{align*}
		Next, focusing on $(II)$, we have 
		\begin{align*}
		(II) \leq ~ & \int \Big| \widetilde{h}_{\text{CvM},2,0}^{(m,n)} (x_1,x_2) - \widetilde{h}_{\text{CvM},2,0}(x_1,x_2)  \Big|  \phi_i^{(m,n)}(x_2) d H(x_2) \\[.5em]
		\leq ~ & \sum_{\substack{i=0,j=0 \\ i+j=4}}^4 \Big| \left(\frac{m}{N}\right)^i \left(\frac{n}{N}\right)^j - \vartheta_X^i \vartheta_Y^j \Big| \sqrt{\max(\vartheta_X,\vartheta_Y) \times \frac{N}{\min\{m,n\}}}.
		\end{align*}
		Since the upper bounds are uniform over $x_1$ and $m/N \rightarrow \vartheta_X, n/N \rightarrow \vartheta_Y$ as $N \rightarrow \infty$ by the assumption, we have
		\begin{align*}
		\lim_{N \rightarrow \infty}\sup_{x_1 \in \mathbb{R}^d}  \Big|& \int \widetilde{h}_{\text{CvM},2,0}^{(m,n)} (x_1,x_2) \phi_i^{(m,n)}(x_2) d H_{m,n}(x_2) \\[.5em]
		& ~~~ - \int \widetilde{h}_{\text{CvM},2,0} (x_1,x_2)  \phi_i^{(m,n)}(x_2) d H(x_2) \Big| = 0.
		\end{align*}
		In addition, 
		\begin{align*}
		0 = ~ & \lim_{N \rightarrow \infty}\sup_{x_1 \in \mathbb{R}^d}  \Big| \int \widetilde{h}_{\text{CvM},2,0}^{(m,n)} (x_1,x_2) \phi_i^{(m,n)}(x_2) d H_{m,n}(x_2) \\[.5em] 
		& ~~~~~~~~~~~~~~~~~~~ - \int \widetilde{h}_{\text{CvM},2,0} (x_1,x_2)  \phi_i^{(m,n)}(x_2) d H(x_2) \Big|, \\[.5em]
		\geq ~ & \sup_{x_1 \in \mathbb{R}^d}   \lim_{N \rightarrow \infty} \Big| \int \widetilde{h}_{\text{CvM},2,0}^{(m,n)} (x_1,x_2) \phi_i^{(m,n)}(x_2) d H_{m,n}(x_2) \\[.5em]
		&  ~~~~~~~~~~~~~~~~~~~  - \int \widetilde{h}_{\text{CvM},2,0} (x_1,x_2)  \phi_i^{(m,n)}(x_2) d H(x_2) \Big|, \\[.5em]
		= ~ & \sup_{x_1 \in \mathbb{R}^d}   \lim_{N \rightarrow \infty} \Big| \lambda_i^{(m,n)} \phi_i^{(m,n)}(x_1)  - \int \widetilde{h}_{\text{CvM},2,0} (x_1,x_2)  \phi_i^{(m,n)}(x_2) d H(x_2) \Big|, \\[.5em]
		= ~ & \sup_{x_1 \in \mathbb{R}^d}  \Big| \lambda_i^{\ast} \phi_i^{\ast}(x_1)   - \int \widetilde{h}_{\text{CvM},2,0} (x_1,x_2)  \phi_i^{\ast}(x_2) d H(x_2) \Big|,
		\end{align*}
		where the last equality is by the uniform integrability of $\widetilde{h}_{\text{CvM},2,0} (x_1,x_2)  \phi_i^{(m,n)}(x_2)$; hence we can interchange the order of the limit and the expectation. Specifically, it is seen that 
		\begin{align} \nonumber 
		& \int \left( \widetilde{h}_{\text{CvM},2,0} (x_1,x_2)  \phi_i^{(m,n)}(x_2) \right)^2 d H(x_2) \\[.5em]
		\leq ~ & \int \left(\phi_i^{(m,n)}(x_2) \right)^2 d H(x_2) \leq \max\{\vartheta_X,\vartheta_Y\} \times  \frac{N}{\min\{m,n\}} \label{Eq: Upper Bound -- Uniformly Integrability}
		\end{align} 
		based on (\ref{Eq: Bound of the Sum}). Since $N/\min\{m,n\} \rightarrow \max \{ \vartheta_X^{-1}, \vartheta_Y^{-1}  \}$ as $N \rightarrow \infty$ by the assumption, choose $N_0$ such that for all $N > N_0$, $|N/\min\{m,n\}  - \max \{ \vartheta_X^{-1}, \vartheta_Y^{-1}  \} | < 1$ and let $B_0 = \max\{N/\min\{m,n\}: N \leq N_0 \}$. Hence, (\ref{Eq: Upper Bound -- Uniformly Integrability}) is uniformly bounded by 
		\begin{align*}
		\max\{\vartheta_X,\vartheta_Y\} \times  \max \Bigg\{ 1 + \frac{1}{\min \{ \vartheta_X, \vartheta_Y  \}}, ~ B_0 \Bigg\}
		\end{align*}
		for all $N$. This implies the uniform integrability of $\widetilde{h}_{\text{CvM},2,0} (x_1,x_2)  \phi_i^{(m,n)}(x_2)$. Therefore, we conclude that the eigenvalues of (\ref{Eq: Integral Equation 2}) converge to those of (\ref{Eq: Integral Equation}).

		\vskip .8em

		In order to verify the second argument, note that
		\begin{align*}
		\int \int \left( \widetilde{h}_{\text{CvM},2,0} (x_1,x_2) \right)^2 dH(x_1) dH(x_2) =  \sum_{i=1}^\infty \lambda_i^2,
		\end{align*}
		where $\lambda_i$ are eigenvalues of (\ref{Eq: Integral Equation}) and 
		\begin{align*}
		\int \int \left( \widetilde{h}_{\text{CvM},2,0}^{(m,n)} (x_1,x_2) \right)^2 dH_{m,n}(x_1) dH_{m,n}(x_2) =  \sum_{i=1}^\infty \left( \lambda_i^{(m,n)} \right)^2.
		\end{align*}
		Based on (\ref{Eq: Upper bound of the difference}) and the boundedness of the kernel, we see that 
		\begin{align*}
		\Bigg| \sum_{i=1}^\infty \lambda_i^2  - \sum_{i=1}^\infty \left( \lambda_i^{(m,n)} \right)^2   \Bigg| \leq \Big| \frac{m}{N} - \vartheta_X \Big| + \Big| \frac{n}{N} - \vartheta_Y \Big| + 2 \sum_{\substack{i=0,j=0 \\ i+j=4}}^4 \Big| \left(\frac{m}{N}\right)^i \left(\frac{n}{N}\right)^j - \vartheta_X^i \vartheta_Y^j \Big|
		\end{align*}
		and thus
		\begin{align*}
		\lim_{N \rightarrow \infty } \sum_{i=1}^\infty \left( \lambda_i^{(m,n)} \right)^2 = \sum_{i=1}^\infty \lambda_i^2 .
		\end{align*}	
	\end{proof}
\end{lemma}

\begin{lemma}
	Let $NU_{\text{\emph{CvM}}}^{(1)}$ be the CvM-statistic based on $i.i.d.$ samples from $\frac{m}{N} P_X + \frac{n}{N} P_Y$. Similarly, let $NU_{\text{\emph{CvM}}}^{(2)}$ be the CvM-statistic based on $i.i.d.$ samples from $\vartheta_X P_X + \vartheta_Y P_Y$ where $m/N \rightarrow \vartheta_X$ and $n/N \rightarrow \vartheta_Y$. Then $NU_{\text{\emph{CvM}}}^{(1)}$ and $NU_{\text{\emph{CvM}}}^{(2)}$ have the same limiting distribution.
	\begin{proof}
		The proof proceeds by following the similar steps in Section~\ref{Section: Proof of Permutation Consistency}. Let us denote by $\widehat{U}_{\text{CvM},K}^{(1)}$, the truncated projection of $U_{\text{{CvM}}}^{(1)}$, which is similarly defined as (\ref{Eq: Truncated Projection}). Based on $i.i.d.$ samples $\{Z_1,\ldots, Z_{m+n}\}$ from $\frac{m}{N} P_X + \frac{n}{N} P_Y$, we can arrive at 
		\begin{align*}
		& N \widehat{U}_{\text{CvM},K}^{(1)}  \\
		= ~ &  \sum_{k=1}^K \lambda_k^{(m,n)} \left( \frac{ \sqrt{N}}{m}  \sum_{i=1}^m \phi_k^{(m,n)}(Z_{i})   - \frac{\sqrt{N}}{n} \sum_{i=m+1}^{m+n}  \phi_k^{(m,n)}(Z_{i})  \right)^2 -  \frac{1}{\vartheta_X \vartheta_Y}\sum_{k=1}^K \lambda_k  + o_{\mathbb{P}}(1).
		\end{align*}
		By the multivariate central limit theorem and Slutsky's theorem with $\lambda_i^{(m,n)} \rightarrow \lambda_i$,  $i=1,\ldots,K$ and $m/N \rightarrow \vartheta_X, n/N \rightarrow \vartheta_Y$, it can be seen that
		\begin{align*}
		N \widehat{U}_{\text{CvM},K}^{(1)}   \convD  \frac{1}{\vartheta_X \vartheta_Y } \sum_{k=1}^K \lambda_k (\xi_k^2 - 1),
		\end{align*}
		where $\xi_k^2$ are independent chi-square random variables with one degree of freedom. The remainder term can be similarly controlled by noting that  
		\begin{align*}
		\lim_{N \rightarrow \infty } \sum_{k=K+1}^\infty \left( \lambda_k^{(m,n)} \right)^2 = \sum_{k=K+1}^\infty \lambda_k^2
		\end{align*}
		from Lemma~\ref{Lemma: Eigenvalues and Eigenfunctions}. This shows that $NU_{\text{CvM}}^{(1)}$ has the same limiting distribution as $NU_{\text{CvM}}^{(2)}$.
	\end{proof}
\end{lemma}

\vskip 2em

\subsection{Proof of Proposition~\ref{Theorem: Asymptotic Power of Permutation Tests}}
The type I error control of the oracle test and the permutation test are obvious and well-known \citep[Chapter 15 of][]{lehmann2006testing}. Hence we focus on the asymptotic power of the tests. When $P_X$ and $P_Y$ are fixed, it is not difficult to show that both tests have asymptotic power equal to one; hence the result follows. In fact, we can prove a more general result that even if the CvM-distance between $P_X$ and $P_Y$ shrinks to zero as the sample size increases, the given tests can be consistent (see Theorem~\ref{Theorem: Upper Bound}). 

Next moving onto the contiguous alternative, we know from Theorem~\ref{Theorem: Asymptotic Null Distribution} that for some $\{\lambda_k\}_{k=1}^\infty$, the null distribution of $NU_{\text{CvM}}$ converges weakly to
\begin{align*}
N  U_{\text{CvM}} \convD  \vartheta_X^{-1} \vartheta_Y^{-1} \sum_{k=1}^\infty \lambda_k (\xi_k^2 - 1).
\end{align*}
Let us write the $(1-\alpha)$ quantile of $\vartheta_X^{-1} \vartheta_Y^{-1} \sum_{k=1}^\infty \lambda_k (\xi_k^2 - 1)$ by $q_{\alpha}$. Then under the null, $c_{\alpha,\text{CvM},s}^\ast \convP q_{\alpha}$, which further implies that $c_{\alpha, \text{CvM},s} \convP q_{\alpha}$ by Theorem~\ref{Theorem: Critical value of Permutation test}. By contiguity as described in the proof of Lemma~\ref{Lemma: Contiguous alternative}, $c_{\alpha,\text{CvM},s}^\ast \convP q_{\alpha}$ and $c_{\alpha, \text{CvM},s} \convP q_{\alpha}$ under the contiguous alternative as well. Then the result follows by Theorem~\ref{Theorem: Asymptotic distribution under contiguous alternatives} and Slutsky's theorem.

\vskip 2em

\subsection{Proof of Theorem~\ref{Theorem: Robustness}}
%The boundedness property of the multivariate CvM-statistic guarantees that the variance of $U_{\text{CvM}}$ is bounded regardless of the presence of contamination. On the other hand, $U_{\text{Energy}}$ can be adversely affected by outliers, which can make the variance of the energy statistic extremely high. Using this observation, we will show that as the same size tends to infinity, the energy statistic can be completely dominated by these outliers and the resulting test becomes powerless in the end. On the other hand, the power of the CvM test remains robust even in the presence of outliers.

To start, we present two lemmas: in Lemma~\ref{Lemma: Variance of U-Stat}, we bound the variance of $U_{\text{CvM}}$ and in Lemma~\ref{Lemma: Two moments under permutations}, we consider the two moments of $U_{\text{CvM}}$ under permutations. %These lemmas will also be used to prove Theorem~\ref{Theorem: Upper Bound}.

\begin{lemma}[Variance of $U_{\text{CvM}}$] \label{Lemma: Variance of U-Stat}
	Consider the CvM-statistic in (\ref{Eq: Definition of U-statistic}). Then there exist universal constants $C_1,C_2,C_3,C_4>0$ such that
	\begin{align*}
	\mV \left[ U_{\text{\emph{CvM}}} \right] ~ \leq ~ C_1 \mE \left[ U_{\text{\emph{CvM}}} \right] \left( \frac{1}{m} + \frac{1}{n} \right) + \frac{C_2}{m^2} + \frac{C_3}{n^2} + \frac{C_4}{mn}. 
	\end{align*}	
	\begin{proof}
		For this proof, it is more convenient to work with the third-order kernel given in (\ref{Eq: third-order kernel}). Let $\widetilde{h}^{\star}_{\text{CvM}}$ be the symmetrized kernel of $h^\star_{\text{CvM}}$ in the sense of (\ref{Eq: Symmetrization}) and define $\widetilde{h}^{\star}_{\text{CvM},c,d}$ in the sense of (\ref{Eq: definition of g_{c,d} }) for $0 \leq c,d, \leq 3$. Further denote the variance of $\widetilde{h}^{\star}_{\text{CvM},c,d}$ by $\sigma_{c,d}^2$ as in (\ref{Eq: Definition of xi_{cd}}). Then the variance of $U_{\text{CvM}}$ can be written as (Lemma~\ref{Lemma: Variance of two-sample U-statistic})
		\begin{align} \label{Eq: U-Stat Variance Formula}
		\mV\left( U_{\text{CvM}} \right) = \sum_{c=0}^{3} \sum_{d=0}^{3} \frac{\binom{3}{c} \binom{3}{d} \binom{m-3}{3-c} \binom{n-3}{3-d}}{\binom{m}{3} \binom{n}{3}} \sigma_{c,d}^2.
		\end{align}
		First we bound $\sigma_{1,0}^2$. After applying the law of total expectation repeatedly, we obtain that 
		\begin{align*}
		& \widetilde{h}^{\star}_{\text{CvM},1,0}(x_1) - \mE [\widetilde{h}^{\star}_{\text{CvM},1,0}(x_1) ] \\[.5em]
		=~ &  \mE \Big[ \Big\{ \ind(\beta^\top x_1 \leq \beta^\top X) - F_{\beta^\top X} (\beta^\top X) \Big\} \cdot \Big\{ F_{\beta^\top Y} (\beta^\top X) - F_{\beta^\top X} (\beta^\top X) \Big\} \Big]  \\[.5em]
		+~ &  \mE \Big[ \Big\{ \ind(\beta^\top x_1 \leq \beta^\top Y) - F_{\beta^\top X} (\beta^\top Y) \Big\} \cdot \Big\{ F_{\beta^\top Y} (\beta^\top Y) - F_{\beta^\top X} (\beta^\top Y) \Big\} \Big]  \\[.5em]
		+ ~ &  \frac{1}{2}\mE \Big[ \Big\{ F_{\beta^\top X} (\beta^\top x_1)  - F_{\beta^\top Y} (\beta^\top x_1) \Big\}^2 \Big] - \frac{1}{2}  \mE \Big[ \Big\{ F_{\beta^\top X} (\beta^\top X) - F_{\beta^\top Y} (\beta^\top X) \Big\}^2 \Big] \\[.5em]
		=~ & f_1(x_1) + f_2(x_1) + f_3(x_1) \quad \text{(say).}
		\end{align*}
		Using the basic inequality $\{f_1(x_1)+f_2(x_1)+f_3(x_1)\}^2 \leq 3f_1^2(x_1) + 3f_2^2(x_1) + 3f_3^2(x_1)$, we have
		\begin{align*}
		\sigma_{1,0}^2 =  & ~ \mE \big[ \big\{ \widetilde{h}^{\star}_{\text{CvM},1,0}(X) - \mE [\widetilde{h}^{\star}_{\text{CvM},1,0}(X) ] \big\}^2 \big] \\[.5em] 
		\leq & ~ 3 \mE \left[ f_1^2(X)  \right] +  3 \mE \left[ f_2^2(X)  \right] +  3 \mE \left[ f_3^2(X)  \right].
		\end{align*}
		By applying Cauchy-Schwarz inequality, the first two terms are bounded by
		\begin{align*}
		& \mE \left[ f_1^2(X)  \right] \leq ~ \mE \big[ \big\{ F_{\beta^\top X} (\beta^\top X) - F_{\beta^\top Y} (\beta^\top X) \big\}^2 \big], \\[.5em]
		& \mE \left[ f_2^2(X)  \right] \leq ~ \mE \big[ \big\{ F_{\beta^\top X} (\beta^\top Y) - F_{\beta^\top Y} (\beta^\top Y) \big\}^2 \big]. 
		\end{align*}
		Since $0 \leq  \mE \big[ \big\{ F_{\beta^\top X} (\beta^\top x_1)  - F_{\beta^\top Y} (\beta^\top x_1) \big\}^2 \big] \leq 1$ for all $x_1 \in \mathbb{R}^d$, the third term is also bounded by
		\begin{align*}
		\mE \left[ f_3^2(X) \right] \leq ~ & \frac{1}{4} \mE \Big[ \Big\{  \mE \big[ \big\{ F_{\beta^\top X} (\beta^\top X)  - F_{\beta^\top Y} (\beta^\top X) \big\}^2 \big]  \Big\}^2  \Big] \\[.5em]
		\leq ~ & \frac{1}{4} \mE \big[ \big\{ F_{\beta^\top X} (\beta^\top X) - F_{\beta^\top Y} (\beta^\top X) \big\}^2 \big].
		\end{align*}
		Thus the following fact (see Theorem~\ref{Theorem: Multivaraite CvM Expression})
		\begin{align*}
		\mE [U_{\text{CvM}}]  =  \frac{1}{2}  \mE \big[ \big\{ F_{\beta^\top X} (\beta^\top X) - F_{\beta^\top Y} (\beta^\top X) \big\}^2 \big] + \frac{1}{2} \mE \big[ \big\{ F_{\beta^\top X} (\beta^\top Y) - F_{\beta^\top Y} (\beta^\top Y) \big\}^2 \big],
		\end{align*}
		leads to $\sigma_{1,0}^2 \lesssim \mE[U_{\text{CvM}}]$. Similarly we have $\sigma_{0,1}^2 \lesssim \mE[U_{\text{CvM}}]$. The rest of $\sigma_{c,d}^2$ can be uniformly bounded due to the boundedness of $\widetilde{h}^{\star}_{\text{CvM}}$. Hence the result follows.
	\end{proof}
\end{lemma}

\begin{lemma}[Two moments under permutations] \label{Lemma: Two moments under permutations}
	The first and second moments of $U_{\text{\emph{CvM}}}$ under permutations are 
	\begin{align*}
	\mE_\varpi \left[ U_{\text{\emph{CvM}}} \right] = 0 \quad \text{and} \quad \mE_\varpi \left[ U_{\text{\emph{CvM}}}^2 \right]  \leq C  \left( \frac{1}{m} + \frac{1}{n} \right)^2,
	\end{align*}
	where $C$ is a universal constant.
	\begin{proof}
		Working directly with the kernel $h_{\text{CvM}}$ is less intuitive to understand the moments of $U_{\text{CvM}}$ under permutations. So we consider the third-order kernel $h^\star_{\text{CvM}}$ in (\ref{Eq: third-order kernel}). Then from Lemma~\ref{Lemma: another expression for the U-statistic}, we have
		\begin{align*}
		U_{\text{CvM}} & = \frac{1}{(m)_3 (n)_3} \sum_{i_1,i_2,i_3=1}^{m,\neq}
		\sum_{j_1,j_2,j_3=1}^{n,\neq} h^\star_{\text{CvM}}(X_{i_1},X_{i_2},X_{i_3}; Y_{j_1},Y_{j_2},Y_{j_3}).
		\end{align*}
		
		\vskip 1em
		
		\noindent \textbf{1. First moment}
		
		\vskip .3em
		
		\noindent Let $\{ Z_1,\ldots,Z_{m+n} \} =\{ X_1,\ldots,X_m, Y_1,\ldots,Y_n \}$ be the pooled samples. Then the first moment of $U_{\text{CvM}}$ becomes
		\begin{align*}
		\mE_\varpi \left[ U_{\text{CvM}} \right] = \mE_{\varpi} \left[ h^\star_{\text{CvM}}(Z_{\varpi(1)}, Z_{\varpi(2)}, Z_{\varpi(3)}; Z_{\varpi(m+1)}, Z_{\varpi(m+2)}, Z_{\varpi(m+3)} )  \right].
		\end{align*}

		Notice that $h^\star_{\text{CvM}}(x_1,x_2,x_3;y_1,y_2,y_3) = - h^\star_{\text{CvM}}(y_1,x_2,x_3;x_1,y_2,y_3)$. This observation shows that the conditional expectation of $h^\star_{\text{CvM}}$ given a subset of permutations $\mathcal{P}_{\varpi,4}  = \{\varpi(2),\varpi(3),\varpi(m+2),\varpi(m+3)\}$ becomes zero, i.e.
		\begin{align*}
		\mE_{\varpi(1),\varpi(m+1)} \left[ h^\star_{\text{CvM}}(Z_{\varpi(1)}, Z_{\varpi(2)}, Z_{\varpi(3)}; Z_{\varpi(m+1)}, Z_{\varpi(m+2)}, Z_{\varpi(m+3)} ) \big| \mathcal{P}_{\varpi,4} \right] = 0,
		\end{align*}
		for all $\mathcal{P}_{\varpi,4}$. Hence, $\mE_\varpi \left[  U_{\text{CvM}} \right] = 0$ by the law of total expectation.

		\vskip 1em
		
		\noindent \textbf{2. Second moment}
		
		\vskip .3em
		
		\noindent Next we calculate the second moment of $U_{\text{CvM}}$ under permutations where
		\begin{align*}
		U^2_{\text{CvM}} =  & \frac{1}{(m)_3^2 (n)_3^2 } \sum_{i_1,i_2,i_3=1}^{m,\neq} \sum_{j_1,j_2,j_3=1}^{n,\neq} \sum_{i_1^\prime,i_2^\prime,i_3^\prime=1}^{m,\neq} \sum_{j_1^\prime,j_2^\prime,j_3^\prime=1}^{n,\neq} \Big\{ \\[.5em]
		&  h^\star_{\text{CvM}}(Z_{i_1},Z_{i_2},Z_{i_3}; Z_{j_1 + m},Z_{j_2 + m},Z_{j_3 + m}) h^\star_{\text{CvM}}(Z_{i_1^\prime},Z_{i_2^\prime},Z_{i_3^\prime}; Z_{j_1^\prime + m },Z_{j_2^\prime + m},Z_{j_3^\prime + m}) \Big\}.
		\end{align*}
		Recall the definition of $\mathcal{I}_3$ and $\mathcal{J}_3$ given in (\ref{Eq: Definition of I and J}). When $\#| \mathcal{I}_3| + \#|\mathcal{J}_3| \leq 1$, we apply the law of total expectation as in the proof of Lemma~(\ref{Lemma: Coupling for the CvM-statistic}) to show that 
		\begin{equation} \label{Eq: Expectation of product f}
		\begin{aligned}
		\mE_\varpi  \big[ & h^\star_{\text{CvM}}(Z_{\varpi(i_1)},Z_{\varpi(i_2)},Z_{\varpi(i_3)}; Z_{\varpi(j_1 + m)},Z_{\varpi(j_2 + m)},Z_{\varpi(j_3 + m)}) \\[.5em]
		&  \times h^\star_{\text{CvM}}(Z_{\varpi(i_1^\prime)},Z_{\varpi(i_2^\prime)},Z_{\varpi(i_3^\prime)}; Z_{\varpi(j_1^\prime + m)},Z_{\varpi(j_2^\prime + m)},Z_{\varpi(j_3^\prime + m)} ) \big]  = 0.
		\end{aligned}
		\end{equation}		
		If $\#| \mathcal{I}_3| + \#|\mathcal{J}_3| > 1$, we use the fact that the kernel $h^\star_{\text{CvM}}$ is bounded by one in absolute value to have
		\begin{align*}
		\big| \mE_\varpi  \big[ & h^\star_{\text{CvM}}(Z_{\varpi(i_1)},Z_{\varpi(i_2)},Z_{\varpi(i_3)}; Z_{\varpi(j_1 + m)},Z_{\varpi(j_2 + m)},Z_{\varpi(j_3 + m)}) \\[.5em]
		&  \times h^\star_{\text{CvM}}(Z_{\varpi(i_1^\prime)},Z_{\varpi(i_2^\prime)},Z_{\varpi(i_3^\prime)}; Z_{\varpi(j_1^\prime + m)},Z_{\varpi(j_2^\prime + m)},Z_{\varpi(j_3^\prime + m)} ) \big] \big| \leq 1.
		\end{align*}		
		Based on the previous observations and the fact that the size of the cases where $\#| \mathcal{I}_3| + \#|\mathcal{J}_3| > 1$ is at most $\prod_{i=0}^4 (m-i) \times \prod_{j=0}^6 (n-j) + \prod_{i=0}^5 (m-i) \times \prod_{j=0}^5 (n-j) + \prod_{i=0}^6 (m-i) \times \prod_{j=0}^4 (n-j)$ up to scaling factors, we conclude that
		\begin{align*}
		\mE_{\varpi} \left[ U_{\text{CvM}}^2 \right]  \leq C \left( \frac{1}{m} + \frac{1}{n} \right)^2
		\end{align*}
		as desired.
	\end{proof} 
\end{lemma}

\vskip 1em

\noindent \textbf{1. Multivariate CvM-statistic}

\vskip .3em 

\noindent We follow the similar steps used in the proof of Theorem~\ref{Theorem: Upper Bound} to show the robustness of the CvM test. 
Since we assume that $Q_X \neq Q_Y$, there exists a positive constant $\delta_1$ such that $W_d(P_{X,N},P_{Y,N}) \geq (1-\epsilon) W_d(Q_X,Q_Y) \geq \delta_1$. Thus $\mE [U_{\text{CvM}}] \geq \delta_1^2$. We first upper bound the type II error as
\begin{align*}
\mP_1 \left( U_{\text{CvM}} \leq c_{\alpha,\text{CvM}} \right) & = ~ \mP_1 \left( U_{\text{CvM}} \leq c_{\alpha,\text{CvM}}, c_{\alpha,\text{CvM}} > \delta_1^2 /2 \right) \\[.5em]
& ~~ +\mP_1 \left( U_{\text{CvM}} \leq c_{\alpha,\text{CvM}}, c_{\alpha,\text{CvM}} \leq \delta_1^2/2\right) \\[.5em]
&  \leq ~ \mP_1 \left( c_{\alpha,\text{CvM}} > \delta_1^2 /2 \right)  + \mP_1 \left( U_{\text{CvM}} \leq \delta_1^2 /2 \right) \\[.5em]
& = ~ (I) + (II) \quad \text{(say).}
\end{align*}
For $(I)$, Lemma~\ref{Lemma: Two moments under permutations} and Chebyshev's inequality yield
\begin{align*}
\mP_\varpi \left( U_{\text{CvM}} \geq t  \right) \leq \frac{\mV_{\varpi}(U_{\text{CvM}})}{t^2} \leq \frac{C_0}{t^2} \cdot \left( \frac{1}{m} + \frac{1}{n} \right)^2
\end{align*}
where $C_0$ is some universal constant. This shows that the critical value of the permutation test is uniformly bounded by
\begin{align*}
c_{\alpha,\text{CvM}}  \leq \sqrt{\frac{C_0}{\alpha}}\left( \frac{1}{m} + \frac{1}{n} \right).
\end{align*}
Hence, we can bound $(I)$ by 
\begin{align*}
(I) = \mP_1 \left( c_{\alpha,\text{CvM}} > \delta_1^2 /2 \right) \leq ~ \frac{4}{\delta_1^4}\mE_1 \left[ c^2_{\alpha,\text{CvM}} \right]  \leq ~ \frac{4C_0}{\alpha \delta_1^4} \left( \frac{1}{m} + \frac{1}{n} \right)^2. 
\end{align*}
Next, 
\begin{align*}
(II) = \mP_1 \left( U_{\text{CvM}} \leq \delta_1^2 /2 \right) = ~  & \mP_1 \left( \frac{U_{\text{CvM}} - \mE_1  [U_{\text{CvM}}]}{\sqrt{\mV_1(U_{\text{CvM}})}}  \leq \frac{\delta_1^2 /2 - \mE_1  [U_{\text{CvM}}]}{\sqrt{\mV_1(U_{\text{CvM}})}} \right) \\[.5em]
\overset{(i)}{\leq} ~ & \mP_1 \left( \frac{U_{\text{CvM}} - \mE_1  [U_{\text{CvM}}]}{\sqrt{\mV_1(U_{\text{CvM}})}}  \leq \frac{-\delta_1^2 /2}{\sqrt{\mV_1(U_{\text{CvM}})}} \right) \\[.5em]
= ~ & \mP_1 \left( \frac{-U_{\text{CvM}} + \mE_1  [U_{\text{CvM}}]}{\sqrt{\mV_1(U_{\text{CvM}})}}  \geq \frac{\delta_1^2 /2}{\sqrt{\mV_1(U_{\text{CvM}})}} \right)  \\[.5em]
\overset{(ii)}{\leq} ~ & \frac{4\mV_1(U_{\text{CvM}})}{\delta_1^4} \\[.5em]
\overset{(iii)}{\leq} ~ & \frac{C_1}{\delta_1^2} \left( \frac{1}{m} + \frac{1}{n} \right) + \frac{C_2}{\delta_1^4} \left( \frac{1}{m} + \frac{1}{n} \right)^2
\end{align*}
where $(i)$ uses $\mE[U_{\text{CvM}}] \geq \delta_1^2$, $(ii)$ is by Chebyshev's inequality and $(iii)$ uses Lemma~\ref{Lemma: Variance of U-Stat} with universal constants $C_1$ and $C_2$. In the end, we have
\begin{align*}
\lim_{m,n \rightarrow \infty} \inf_{G_N} \mE_1 [\phi_{\text{CvM}}] & \geq ~ 1 - \lim_{m,n \rightarrow \infty} \inf_{G_N} \Bigg\{ \frac{4C_0}{\alpha \delta_1^4} \left( \frac{1}{m} + \frac{1}{n} \right)^2   +   \frac{C_1}{\delta_1^2} \left( \frac{1}{m} + \frac{1}{n} \right)  \\[.5em] 
&~~~~~~~~~~~~~~~~~~~~~~~~~~~~ + \frac{C_2}{\delta_1^4} \left( \frac{1}{m} + \frac{1}{n} \right)^2 \Bigg\} = 1,
\end{align*}
which completes the proof of the first part. 

\vskip 1em 

\noindent \textbf{2. Energy statistic}

\vskip .3em

\noindent We continue our discussion from the main text (see the proof of Theorem~\ref{Theorem: Robustness} in the main text). Recall that we take $G_N$ to have a multivariate normal distribution with zero mean vector and covariance matrix $\sigma_N^2 I_d$. Also recall the truncated random vectors coupled with $X$ and $Y$ defined as
\begin{align*}
\widetilde{X} = 
\begin{cases}
(0,\ldots,0)^\top, \quad & \text{if} ~ X \sim Q_X, \\
X / \sigma_N, \quad & \text{if} ~ X \sim G_N,
\end{cases} 
\quad \text{and} \quad
\widetilde{Y} = 
\begin{cases}
(0,\ldots,0)^\top, \quad & \text{if} ~ Y \sim Q_Y, \\
Y / \sigma_N, \quad & \text{if} ~ Y \sim G_N.
\end{cases} 
\end{align*}
We shall first show that the energy statistic based on the original samples and the other energy statistic based on the truncated samples are asymptotically equivalent.

\begin{lemma} \label{Lemma: Asymptotic Equivalence bewteen energy statistics}
	Suppose $\sigma_N^2 \asymp N^q$ for some $q>2$. Let $\widetilde{U}_{\text{\emph{Energy}}}$ be the energy statistic based on $\{\widetilde{X}_1,\ldots,\widetilde{X}_m, \widetilde{Y}_1, \ldots, \widetilde{Y}_n \}$ coupled with the original samples $\{ X_1,\ldots, X_m, Y_1,\ldots, Y_n \}$ and $U_{\text{\emph{Energy}}}$ be the energy statistic based on the original samples. Then under the asymptotic regime in (\ref{Eq: Limit of Sample Proportion}),
	\begin{align*}
	N \sigma_N^{-1} U_{\text{\emph{Energy}}} - N  \widetilde{U}_{\text{\emph{Energy}}} \convP 0.
	\end{align*}
	\begin{proof}
		%Let us define the difference between the Euclidean distances associated with original samples $(X_1,X_2)$ and truncated samples $(\widetilde{X}_1, \widetilde{X}_2)$ by
		Let us denote
		\begin{align*}
		\Delta_{m,n}(X_1,X_2) & =  \sigma_N^{-1} \| X_1- X_2 \|  - \|\widetilde{X}_{1} - \widetilde{X}_{2}\|.
		\end{align*}
		Observe that there are four possible cases for $\Delta_{m,n}(X_1,X_2)$:
		\begin{align*}
		\Delta_{m,n}(X_1,X_2) = 
		\begin{cases}
		\text{{Case~(a):}} \quad \frac{1}{\sigma_N}\| X_1 - X_2\|, ~~ & \text{if} ~~ X_1, X_2 \sim Q_X, \\
		\text{{Case~(b):}} \quad \frac{1}{\sigma_N}\| X_1 - X_2  \| - \frac{1}{\sigma_N}\| X_2 \|, ~~ & \text{if} ~~ X_1 \sim Q_X, X_2 \sim G_N, \\
		\text{{Case~(c):}} \quad \frac{1}{\sigma_N}\| X_1 - X_2 \| - \frac{1}{\sigma_N}\| X_1 \|, ~~ & \text{if} ~~ X_1 \sim G_N, X_2 \sim Q_X, \\
		\text{{Case~(d):}} \quad 0, ~~ & \text{if} ~~ X_1, X_2 \sim H_m.
		\end{cases}
		\end{align*}
		In any case, one can verify under the finite second moment condition that
		\begin{align} \label{Eq: Expectaion of Delta}
		\mE \left[ \Delta_{m,n}^2(X_1,X_2)  \right] \lesssim \sigma_N^{-2}.
		\end{align}
		Similarly, it can be seen that
		$\mE \left[ \Delta_{m,n}^2(X_1,X_2)  \right] \lesssim \sigma_N^{-2}$, $\mE \left[ \Delta_{m,n}^2(Y_1,Y_2)  \right] \lesssim \sigma_N^{-2}$ and $\mE \left[ \Delta_{m,n}^2(X_1,Y_1)  \right] \lesssim \sigma_N^{-2}$.
		
		\vskip .8em
		
		Write the symmetrized kernel of the energy statistic as
		\begin{align*}
		& \widetilde{h}_{\text{Energy}}(x_1,x_2;y_1,y_2)  \\[.5em]
		= ~& \frac{1}{2}\|x_1 - y_1\| + \frac{1}{2}\|x_1 - y_1\| + \frac{1}{2} \|x_2 - y_1 \| + \frac{1}{2}\|x_2 - y_2\| - \|x_1 - x_2\| - \|y_1 - y_2\|.
		\end{align*}
		Then the energy statistic based on the truncated random samples can be written as
		\begin{align*}
		\widetilde{U}_{\text{Energy}} = \frac{1}{(m)_2 (n)_2} \sum_{i_1,i_2 =1}^{m,\neq} \sum_{j_1,j_2 =1}^{n,\neq} \widetilde{h}_{\text{Energy}}(\widetilde{X}_{i_1},\widetilde{X}_{i_2}; \widetilde{Y}_{j_1}, \widetilde{Y}_{j_2}).
		\end{align*}
		Now our goal is to show
		\begin{align} \nonumber 
		&N \big( \sigma_N^{-1}U_{\text{Energy}} - \widetilde{U}_{\text{Energy}} \big)  \\[.5em] \nonumber
		= ~ & \frac{N}{(m)_2 (n)_2} \sum_{i_1,i_2 =1}^{m,\neq} \sum_{j_1,j_2 =1}^{n,\neq} \Bigg\{ \frac{1}{\sigma_N}\widetilde{h}_{\text{Energy}}(X_{i_1}, X_{i_2}; Y_{j_1}, Y_{j_2}) - \widetilde{h}_{\text{Energy}}(\widetilde{X}_{i_1},\widetilde{X}_{i_2}; \widetilde{Y}_{j_1}, \widetilde{Y}_{j_2}) \Bigg\} \\[.5em] \label{Eq: Definition of D}
		:=~ & \frac{N}{(m)_2 (n)_2} \sum_{i_1,i_2 =1}^{m,\neq} \sum_{j_1,j_2 =1}^{n,\neq} h_D\{ (X_{i_1},\widetilde{X}_{i_1}), (X_{i_2},\widetilde{X}_{i_2}); (Y_{j_1},\widetilde{Y}_{j_1}), (Y_{j_2},\widetilde{Y}_{j_2})\}  ~ \convP ~ 0.
		\end{align}
		For simplicity we will write
		\begin{align*}
		h_D(i_1,i_2;j_1,j_2)  = h_D\{ (X_{i_1},\widetilde{X}_{i_1}), (X_{i_2},\widetilde{X}_{i_2}); (Y_{j_1},\widetilde{Y}_{j_1}), (Y_{j_2},\widetilde{Y}_{j_2})\}.
		\end{align*}
		To show (\ref{Eq: Definition of D}), we first apply Cauchy-Schwarz inequality to bound
		\begin{align*}
		\mE \big[ h_D(i_1,i_2;j_1,j_2) h_D(i_1^\prime,i_2^\prime;j_1^\prime,j_2^\prime) \big]  \leq ~ & \sqrt{\mE \big[ h_D^2(i_1,i_2;j_1,j_2) \big] } \sqrt{\mE \big[ h_D^2(i_1^\prime,i_2^\prime;j_1^\prime,j_2^\prime)\big]}, \\[.5em]
		\lesssim ~ & \sigma_N^{-2} ,
		\end{align*}
		which holds for any set of indices such that $i_1 \neq i_2, j_1 \neq j_2, i_1^\prime \neq i_2^\prime, j_1^\prime \neq j_2^\prime$. Note that for the second inequality, we used
		\begin{align*}
		\mE \big[ h_D^2(i_1,i_2;j_1,j_2)) \big]  \lesssim ~ & \mE [\Delta_{m,n}^2(X_{i_1},X_{i_2})] + \mE [\Delta_{m,n}^2(X_{i_1},Y_{j_1})] + \mE [\Delta_{m,n}^2(X_{i_1},Y_{j_2})] \\[.5em]
		+ ~ & \mE [\Delta_{m,n}^2(X_{i_2},Y_{i_1})] + \mE [\Delta_{m,n}^2(X_{i_2},Y_{j_2})] + \mE [\Delta_{m,n}^2(Y_{j_1},Y_{j_2})], \\[.5em]
		\lesssim ~ & \sigma_N^{-2},
		\end{align*}
		by (\ref{Eq: Expectaion of Delta}) and similarly for the other cases.  As a consequence, 
		\begin{align*}
		\mE \Big[ N^2 \left( \sigma_N^{-1} U_{\text{Energy}} - \widetilde{U}_{\text{Energy}} \right)^2 \Big] \lesssim ~ \sigma_N^{-2} N^2.
		\end{align*}
		Under the given assumptions that $\sigma_N^2 \asymp (m+n)^{q}$ with $q > 2$ and $m/N \rightarrow \vartheta_X \in (0,1)$, we obtain $N ( \sigma_N^{-1} U_{\text{Energy}} - \widetilde{U}_{\text{Energy}}) \convP 0$ as desired.	
	\end{proof}
\end{lemma}

\vskip 1.5em

Since $\widetilde{U}_{\text{Energy}}$ has degeneracy of order one, $N \widetilde{U}_{\text{Energy}}$ converges to an infinite weighted sum of chi-square random variables (Theorem~\ref{Theorem: Asymptotic Null Distribution}):
\begin{align*}
N \widetilde{U}_{\text{Energy}} \convD \sum_{k=1}^{\infty} \lambda_k (\xi_k^2 - 1),
\end{align*}
for some $\{\lambda_k\}_{k=1}^\infty$. Lemma~\ref{Lemma: Asymptotic Equivalence bewteen energy statistics} then implies that $N U_{\text{Energy}} / \sigma_N$ converges to the same distribution: 
\begin{align*}
\frac{N}{\sigma_N} U_{\text{Energy}}  \convD \sum_{k=1}^{\infty} \lambda_k (\xi_k^2 - 1).
\end{align*}
Furthermore, the permutation distribution of $N \sigma_N^{-1} U_{\text{Energy}}$ is asymptotically equivalent to the limiting distribution of $N \widetilde{U}_{\text{Energy}}$ as shown in the next lemma.

\begin{lemma} \label{Lemma: Permutation distribution of Energy Statistic}
	Consider the same assumptions and notation used in Lemma~\ref{Lemma: Asymptotic Equivalence bewteen energy statistics}. Let $R(t)$ be the cumulative distribution function of the limiting distribution of $N \widetilde{U}_{\text{\emph{Energy}}}$. Then the permutation distribution function of $N \sigma_N^{-1} U_{\text{\emph{Energy}}}$, denoted by $\widehat{R}_{m,n}(t)$, satisfies
	\begin{align} \label{Eq: Permutation Limit of Energy Statistic}
	\sup_{t \in \mathbb{R}} \Big| \widehat{R}_{m,n}(t) - R(t)  \Big| \convP 0.
	\end{align}
	\begin{proof}
		Let $\{Z_1,\ldots, Z_{m+n}\}$ be the pooled samples of $\{X_1,\ldots, X_m, Y_1,\ldots,Y_n \}$ and similarly $\{\widetilde{Z}_1,\ldots, \widetilde{Z}_{m+n}\}$ be the pooled samples of $\{\widetilde{X}_1,\ldots, \widetilde{X}_m, \widetilde{Y}_1,\ldots,\widetilde{Y}_n\}$. For any random permutation $\varpi = \{ \varpi(1), \ldots, \varpi(N)\}$ of $\{1,\ldots,N\}$, we will show that
		\begin{align} \label{Eq: Energy Statistic Difference}
		N \sigma_N^{-1} U_{\text{Energy}} (Z_{\varpi})- N \widetilde{U}_{\text{Energy}} (\widetilde{Z}_\varpi)\convP 0,
		\end{align}
		where $Z_{\varpi} = (Z_{\varpi(1)}, \ldots, Z_{\varpi(N)})$ and $\widetilde{Z}_\varpi= (\widetilde{Z}_{\varpi(1)}, \ldots, \widetilde{Z}_{\varpi(N)})$. If this is the case, then for two independent $\varpi$ and $\varpi^\prime$, the following result
		\begin{align} \label{Eq: Hoeffding's condition for the energy statistic}
		(N \widetilde{U}_{\text{Energy}} (\widetilde{Z}_\varpi), N \widetilde{U}_{\text{Energy}} (\widetilde{Z}_{\varpi^\prime}))  \convD (T, T^\prime)
		\end{align}
		implies
		\begin{align*}
		(N \sigma_N^{-1} U_{\text{Energy}} ({Z}_\varpi), N \sigma_N^{-1} U_{\text{Energy}}({Z}_{\varpi^\prime}))  \convD (T, T^\prime),
		\end{align*}
		by Slutsky's theorem. Here $T$ and $T^\prime$ are independent and identically distributed with the distribution function $R(t)$. Then Hoeffding's condition in Lemma (\ref{Lemma: Hoeffding's condition}) establishes (\ref{Eq: Permutation Limit of Energy Statistic}). Indeed, (\ref{Eq: Hoeffding's condition for the energy statistic}) holds from Theorem~\ref{Theorem: Two-Sample Degenerate Kernel}; hence it is enough to show (\ref{Eq: Energy Statistic Difference}) to complete the proof. 
		
		Note that 
		\begin{align*}
		& N \sigma_N^{-1} U_{\text{Energy}} (Z_{\varpi})- N \widetilde{U}_{\text{Energy}} (\widetilde{Z}_\varpi) = ~  \frac{N}{(m)_2 (n)_2} \sum_{i_1,i_2 =1}^{m,\neq} \sum_{j_1,j_2 =1}^{n,\neq}  \Big[ \\[.5em]
		& h_D\{ (Z_{\varpi(i_1)},\widetilde{Z}_{\varpi(i_1)}), (Z_{\varpi(i_2)},\widetilde{Z}_{\varpi(i_2)}); (Z_{\varpi(j_1+m)},\widetilde{Z}_{\varpi(j_1+m)}), (Z_{\varpi(j_2+m)},\widetilde{Z}_{\varpi(j_2+m)})\} \Big]
		\end{align*}
		where kernel $h_D$ is given in (\ref{Eq: Definition of D}). Note further by (\ref{Eq: Expectaion of Delta}) that 
		\begin{align*}
		& \mE \left[ h_D^2\{ (Z_{\varpi(i_1)},\widetilde{Z}_{\varpi(i_1)}), (Z_{\varpi(i_2)},\widetilde{Z}_{\varpi(i_2)}); (Z_{\varpi(j_1+m)},\widetilde{Z}_{\varpi(j_1+m)}), (Z_{\varpi(j_2+m)},\widetilde{Z}_{\varpi(j_2+m)})\} \right] \\[.5em]
		\lesssim ~ &   \mE\left[ \Delta_{m,n}^2 (Z_{\varpi(i_1)}, Z_{\varpi(i_2)}) \right] + \mE\left[ \Delta_{m,n}^2 (Z_{\varpi(i_1)}, Z_{\varpi(j_1 + m)}) \right] \\[.5em]
		+ ~ &  \mE\left[ \Delta_{m,n}^2 (Z_{\varpi(i_1)}, Z_{\varpi(j_2 + m)}) \right] + \mE\left[ \Delta_{m,n}^2 (Z_{\varpi(i_2)}, Z_{\varpi(j_1 + m)}) \right] \\[.5em]
		+ ~ & \mE\left[ \Delta_{m,n}^2 (Z_{\varpi(i_2)}, Z_{\varpi(j_2 + m)}) \right]   + \mE\left[ \Delta_{m,n}^2 (Z_{\varpi(j_1+m)}, Z_{\varpi(j_2 + m)}) \right]  \\[.5em]
		\lesssim ~ & \sigma_N^{-2}
		\end{align*}
		and similarly for the other cases. Then it is easy to see that
		\begin{align*}
		\mE \big[ \big( N \sigma_N^{-1} U_{\text{Energy}} (Z_{\varpi})- N \widetilde{U}_{\text{Energy}} (\widetilde{Z}_\varpi) \big)^2 \big] \lesssim  \sigma_{N}^{-2} N^2 = o(1)
		\end{align*}
		whenever $\sigma_N^2 \asymp N^{q}$ for some $q > 2$. This implies (\ref{Eq: Energy Statistic Difference}), which completes the proof.
	\end{proof}
\end{lemma}

\vskip 1em

Combining the previous results yields
\begin{align*}
\lim_{N \rightarrow \infty} \mP \left( U_{\text{Energy}} >  c_{\alpha,\text{Energy}} \right) = ~ & \lim_{N \rightarrow \infty} \mP \left( N \sigma_N^{-1} U_{\text{Energy}} >  N \sigma_N^{-1}  c_{\alpha,\text{Energy}} \right) \\[.5em] 
= ~ & \lim_{N \rightarrow \infty} \mP \left( N \widetilde{U}_{\text{Energy}} >  \widetilde{c}_{\alpha,\text{Energy}}  \right) \leq \alpha,
\end{align*}
where $\widetilde{c}_{\alpha,\text{Energy}}$ is the $(1-\alpha)$ quantile of the permutation distribution of $N \widetilde{U}_{\text{Energy}}$. Hence the result follows.

\vskip 2em

\subsection{Proof of Lemma~\ref{Lemma: Lower Bound of CvM-distance}}
Let $\beta^\top Z$ have the distribution function $F_{\beta^\top X}(t)/2 + F_{\beta^\top Y}(t)/2$. First notice from the definition of the multivariate CvM-distance that
\begin{align*}
W_d^2 ~ &  = ~ \mE \Big[ \Big\{ F_{\beta^\top X} (\beta^\top Z) - F_{\beta^\top Y} (\beta^\top Z )  \Big\}^2 \Big]  \geq ~ \Big\{ \mE \Big[\Big| F_{\beta^\top X} (\beta^\top Z) - F_{\beta^\top Y} (\beta^\top Z )  \Big|\Big]  \Big\}^{2},
\end{align*}
where we used Jensen's inequality. Let us denote the expectation with respect to $X_1,X_2,Y_1$ (and $X_1,Y_1,Y_2$) by $\mE_{X_1,X_2,Y_1}$ (and $\mE_{X_1,Y_1,Y_2}$). Then from the definition of $\beta^\top Z$, we have
\begin{align*}
& \mE \big[ \big| F_{\beta^\top X} (\beta^\top Z) - F_{\beta^\top Y} (\beta^\top Z )  \big| \big] \\[.5em]
= ~ & \frac{1}{2}\mE \big[\big| F_{\beta^\top X} (\beta^\top X_1) - F_{\beta^\top Y} (\beta^\top X_1 )\big|\big] +  \frac{1}{2}\mE \big[\big| F_{\beta^\top X} (\beta^\top Y_1) - F_{\beta^\top Y} (\beta^\top Y_1 ) \big| \big] \\[.5em]
\geq~  &  \frac{1}{2}\mE_{\beta} \Big[ \Big| \mE_{X_1,X_2,Y_1} \Big\{  \ind(\beta^\top X_1 \leq \beta^\top X_2)  -  \ind(\beta^\top Y_1 \leq \beta^\top X_2)\Big\} \Big| \Big] \\[.5em]
+ ~  & \frac{1}{2} \mE_{\beta} \Big[ \Big| \mE_{X_1,Y_1,Y_2} \Big\{  \ind(\beta^\top X_1 \leq \beta^\top Y_2)  -  \ind(\beta^\top Y_1 \leq \beta^\top Y_2)\Big\} \Big| \Big], 
\end{align*}
where we used Jensen's inequality once again to obtain the lower bound. The last expression can be simplified based on the observation that $\mP(\beta^\top X_1 \leq \beta^\top X_2) = \mP(\beta^\top Y_1 \leq \beta^\top Y_2) = 1/2$ as
\begin{align*}
\mE_{\beta} \Big[\Big| \frac{1}{2} - \mP \left( \beta^\top X \leq \beta^\top Y \right)  \Big| \Big].
\end{align*}
Therefore,
\begin{align*}
W_d^2 & ~ \geq ~  \bigg\{ \int_{\mathbb{S}^{d-1}} \Big| \frac{1}{2} - \mP \left( \beta^\top X \leq \beta^\top Y \right)  \Big| d\lambda(\beta)\bigg\}^2,
\end{align*}
which completes the proof.

\vskip 2em

\subsection{Proof of Theorem~\ref{Theorem: Lower Bound}}

The minimax lower bound is based on a standard application of Neyman-Pearson lemma \citep[see e.g.][]{baraud2002non}. Here we write the joint distributions of samples under the null and alternative hypotheses by $P_0^{m,n}$ and $P_1^{m,n}$, respectively. Then
\begin{align} \nonumber
\inf_{\phi \in \mathds{T}_{m,n}(\alpha)} \sup_{P_X,P_Y \in \mathcal{F}(\epsilon_{m,n}^\star)} \mP_1 \left( \phi = 0 \right) \geq & ~ 1 - \alpha - \sup_{A \in \mathcal{A}} \big| P_0^{m,n}(A) - P_1^{m,n}(A) \big| \\[.5em] 
\geq  & ~ 1 - \alpha - \sqrt{ \frac{1}{2} \mathsf{KL} \left( P_1^{m,n}, P_0^{m,n} \right) }, \label{Eq: Lower Bound for KL}
\end{align}
where the second inequality is by Pinsker's inequality \citep[e.g.~Lemma 2.5 of][]{tsybakov2009introduction}.

Recall the example considered in Lemma~\ref{Lemma: Least Favorable Distribution}:
\begin{align*}
X^\ast := (\xi_{1},0,\ldots,0)^\top \quad \text{and} \quad Y^\ast := (\xi_{2},0,\ldots, 0)^\top,
\end{align*}
where $\xi_{1} \sim N(\mu_{X^\ast},1)$ and $\xi_{2} \sim N(\mu_{Y^\ast},1)$. We let $\mu_{X^\ast} = \mu_{Y^\ast} = 0$ under the null and 
\begin{align*}
\mu_{X^\ast} = \frac{\sqrt{2}(1-\alpha-\zeta)}{\sqrt{m}} \quad \text{and} \quad \mu_{Y^\ast} = -\frac{\sqrt{2}(1-\alpha-\zeta)}{\sqrt{n}},
\end{align*}
under the alternative. Then from Lemma~\ref{Lemma: Least Favorable Distribution}, we have $P_{X^\ast}, P_{Y^\ast} \in \mathcal{F}(\epsilon^\ast_{m,n})$ for all $d$. In this case, the Kullback-Leibler divergence is calculated as
\begin{align*}
\mathsf{KL} \left(  P_1^{m,n}, P_0^{m,n} \right) = \frac{m}{2}\mu_{X^\ast}^2 + \frac{n}{2} \mu_{Y^\ast}^2 = 2(1-\alpha - \zeta)^2.
\end{align*}
By plugging this into (\ref{Eq: Lower Bound for KL}), we conclude that 
\begin{align*}
\inf_{\phi \in \mathds{T}_{m,n}(\alpha)} \sup_{P_X,P_Y \in \mathcal{F}(\epsilon_{m,n}^\star)} \mP_1 \left( \phi = 0 \right) \geq \zeta.
\end{align*}
Hence the result follows.

\vskip 2em

\subsection{Proof of Theorem~\ref{Theorem: Upper Bound}}

To finish the proof, we need to verify the condition in (\ref{Eq: Minimax Power Sufficient Condition}). Using Chebyshev's inequality and Lemma~\ref{Lemma: Two moments under permutations}, 
\begin{align*}
\mP_{\varpi} \left( U_{\text{CvM}} \geq t \right) \leq \frac{\mE_{\varpi} [U_{\text{CvM}}^2]  }{t^2} \leq \frac{C_0}{t^2} \left( \frac{1}{m} + \frac{1}{n} \right)^2.
\end{align*}
As a result, the permutation critical value $c_{\alpha,\text{CvM}}$ is upper bounded by $\sqrt{C_0/\alpha} (1/m + 1/n)$ with probability one. This implies that its $\zeta/2$ upper quantile $c_{\zeta/2}^\ast$ is also bounded by
\begin{align*}
c_{\zeta/2}^\ast ~ \leq ~ \sqrt{\frac{C_0}{\alpha}} \left( \frac{1}{\sqrt{m}} + \frac{1}{\sqrt{n}} \right)^2.
\end{align*}
From Lemma~\ref{Lemma: Variance of U-Stat}, we have 
\begin{align*}
\sqrt{ \frac{\zeta}{2} \text{{Var}}_1 \left[ U_{\text{CvM}} \right]} ~ \leq ~ & \sqrt{ \frac{\zeta}{2} \cdot \Bigg\{ C_1 \mE_1 \left[ U_{\text{CvM}} \right] \cdot \left( \frac{1}{m} + \frac{1}{n} \right) + \frac{C_2}{m^2} + \frac{C_3}{n^2} + \frac{C_4}{mn} \Bigg\}} \\[.5em]
~ \leq ~ & C_5 \left( \frac{1}{\sqrt{m}} + \frac{1}{\sqrt{n}} \right)^2.
\end{align*}
By choosing a sufficiently large $c>0$ in (\ref{Eq: Radius of interest}), we conclude that
\begin{align*}
\mE_1 [U_{\text{CvM}}] \geq c_{\zeta/2}^\ast  + \sqrt{ \frac{\zeta}{2} \text{{Var}}_1 \left[ U_{\text{CvM}} \right]} .
\end{align*}

\vskip 2em

\subsection{Proof of Proposition~\ref{Proposition: Power of Linear Statistic}}
Let $\sigma_0^2$ and $\sigma_1^2$ be the variance of 
\begin{align*}
\widetilde{h}_{\text{CvM}}(X_1,X_2;Y_1,Y_2) = \frac{1}{2} \{ h_{\text{CvM}}(X_{1},X_{2}; Y_{1},Y_{2}) + h_{\text{CvM}}(X_{2},X_{1}; Y_{2},Y_{1})\},
\end{align*}
under the null and alternative, respectively. From the boundedness of $h_{\text{CvM}}$, we have $0 < \sigma_0^2, \sigma_1^2 < \infty$. Then by the central limit theorem, the null distribution approximates
\begin{align*}
& \frac{\sqrt{M}  L_{\text{CvM}} }{\sigma_0} \convD N(0,1) \quad \text{under $H_0$},
\end{align*}
which implies that $\sqrt{M} \sigma_0^{-1} c_{\alpha, \text{linear}} \rightarrow  - z_{\alpha}$ where $z_{\alpha}$ is the $\alpha$ quantile of the standard normal distribution and $z_{\alpha} < 0 $ for $\alpha < 1/2$. Hence, the power function approximates
\begin{align*}
\lim_{N \rightarrow \infty} \mP_1 \left( L_{\text{CvM}} > c_{\alpha, \text{linear}} \right) & = \lim_{N \rightarrow \infty} \mP_1 \left( \frac{\sqrt{M} (L_{\text{CvM}} - W_d^2)}{\sigma_1} > \frac{\sqrt{M}c_{\alpha, \text{linear}}}{\sigma_1}  - \frac{\sqrt{M} W_d^2}{\sigma_1} \right) \\[.5em]
& = \lim_{N \rightarrow \infty} \mP_1 \left( \frac{\sqrt{M} (L_{\text{CvM}} - W_d^2)}{\sigma_1} > - \frac{\sigma_0}{\sigma_1}z_\alpha  - \frac{\sqrt{M} W_d^2}{\sigma_1} \right) \\[.5em]
& \leq \lim_{N \rightarrow \infty} \mP_1 \left( \frac{\sqrt{M} (L_{\text{CvM}} - W_d^2)}{\sigma_1} >  - \frac{\sqrt{M} W_d^2}{\sigma_1} \right) \\[.5em]
& = \frac{1}{2},
\end{align*}
where the last equality uses
\begin{align*}
\frac{\sqrt{M} (L_{\text{CvM}} - W_d^2)}{\sigma_1} \convD N(0,1) \quad \text{under $H_1$}
\end{align*}
and $\sqrt{M} W_d^2 \convP 0$ by the assumption. This completes the proof.

\vskip 2em

\subsection{Proof of Theorem~\ref{Theorem: HDLSS Consistency}}
The proof consists of two parts. In the first part, we will present some lemmas, which investigate the limiting behavior of $\widetilde{h}_{\text{CvM}}$ under the HDLSS setting, and in part two, we will prove the main result.

\vskip 1em 

\noindent \textbf{$\bullet$ Part 1.}

\vskip .3em

\noindent First define the five quantities
\begin{align*}
Q_1 := ~ & \frac{1}{3} - \frac{1}{2\pi} \text{arccos} \left( \frac{\overline{\delta}_{XY}^2 + \overline{\sigma}_{X}^2}{\overline{\delta}_{XY}^2 + \overline{\sigma}_{X}^2 + \overline{\sigma}_{Y}^2} \right) - \frac{1}{2\pi} \text{arccos} \left( \frac{\overline{\delta}_{XY}^2 + \overline{\sigma}_{Y}^2}{\overline{\delta}_{XY}^2 + \overline{\sigma}_{X}^2 + \overline{\sigma}_{Y}^2} \right), \\[1.5em]
Q_2 : = ~ & \frac{1}{3} - \frac{1}{2\pi} \text{{arccos}} \left(  \frac{\overline{\sigma}_{X}^2}{ ( 2\overline{\sigma}_{X}^2 )^{1/2} (\overline{\delta}_{XY}^2 + \overline{\sigma}_{X}^2 + \overline{\sigma}_{Y}^2 )^{1/2} } \right) \\[.5em]
& ~~- \frac{1}{2\pi} \text{{arccos}} \left(  \frac{\overline{\sigma}_{Y}^2}{( 2\overline{\sigma}_{Y}^2 )^{1/2} (\overline{\delta}_{XY}^2 + \overline{\sigma}_{X}^2 + \overline{\sigma}_{Y}^2 )^{1/2} } \right), \\[1.5em]
Q_3 : = ~ &  \frac{1}{3} - \frac{1}{4\pi}  \Bigg[ \text{arccos} \left( \frac{1}{2} \right)   + \text{arccos} \left(  \frac{\overline{\delta}_{XY}^2 + \overline{\sigma}_{Y}^2}{\overline{\delta}_{XY}^2 + \overline{\sigma}_{X}^2 + \overline{\sigma}_{Y}^2} \right) \\[.5em] 
& ~~~~~~~~~~~~~~~~~~~~~~~~~+ 2\text{arccos} \left(  \frac{\overline{\sigma}_{X}^2}{ (2 \overline{\sigma}_X^2)^{1/2} ( \overline{\delta}_{XY}^2 + \overline{\sigma}_{X}^2 + \overline{\sigma}_{Y}^2)^{1/2}} \right)  \Bigg], \\[1.5em]
Q_4 := ~ & \frac{1}{3} - \frac{1}{4\pi}  \Bigg[ \text{arccos} \left( \frac{1}{2} \right)   + \text{arccos} \left(  \frac{\overline{\delta}_{XY}^2 + \overline{\sigma}_{X}^2}{\overline{\delta}_{XY}^2 + \overline{\sigma}_{X}^2 + \overline{\sigma}_{Y}^2} \right) \\[.5em] 
& ~~~~~~~~~~~~~~~~~~~~~~~~~ + 2\text{arccos} \left(  \frac{\overline{\sigma}_{Y}^2}{ (2 \overline{\sigma}_Y^2)^{1/2} ( \overline{\delta}_{XY}^2 + \overline{\sigma}_{X}^2 + \overline{\sigma}_{Y}^2)^{1/2}} \right)  \Bigg], \\[1.5em]
Q_5 := ~ & 0.
\end{align*}
Then by the weak law of large number and the continuous mapping theorem under \textbf{(A1)} and \textbf{(A2)}, it is not difficult to see that for any distinct indices $1 \leq i_1,i_2,i_3,i_4 \leq m$ and $1 \leq j_1,j_2,j_3,i_4 \leq n$, 
\begin{align*}
&\widetilde{h}_{\text{CvM}}(X_{i_1},X_{i_2};Y_{j_1},Y_{j_2}) = \widetilde{h}_{\text{CvM}}(Y_{j_1},Y_{j_2}; X_{i_1},X_{i_2}) \convP Q_1, \\[.5em]
& \widetilde{h}_{\text{CvM}}(X_{i_1},Y_{j_1};X_{i_2},Y_{j_2}) = \widetilde{h}_{\text{CvM}}(Y_{j_1},X_{i_1};Y_{j_2},X_{i_2})\convP Q_2.
\end{align*}
Similarly, 
\begin{align*}
&\widetilde{h}_{\text{CvM}}(X_{i_1},X_{i_2};X_{i_3},Y_{j_1}) = \widetilde{h}_{\text{CvM}}(X_{i_1},X_{i_2};Y_{j_1}, X_{i_3}) \\[.5em]
= ~&  \widetilde{h}_{\text{CvM}}(X_{i_3},Y_{j_1}; X_{i_1},X_{i_2}) =  \widetilde{h}_{\text{CvM}}(Y_{j_1},X_{i_3}; X_{i_1},X_{i_2}) ~ \convP ~ Q_3,
\end{align*}
and
\begin{align*}
&\widetilde{h}_{\text{CvM}}(Y_{j_1},Y_{j_2};Y_{j_3},X_{i_1}) = \widetilde{h}_{\text{CvM}}(Y_{j_1},Y_{j_2};X_{i_1}, Y_{j_3}) \\[.5em]
= ~&  \widetilde{h}_{\text{CvM}}(Y_{j_3},X_{i_1}; Y_{j_1},Y_{j_2}) =  \widetilde{h}_{\text{CvM}}(X_{i_1},Y_{j_3}; Y_{j_1},Y_{j_2}) ~ \convP ~ Q_4.
\end{align*}
When all components are from the same distribution, then $\widetilde{h}_{\text{CvM}}(X_{i_1},X_{i_2};X_{i_3},X_{i_4}) \convP Q_5 = 0$  and $\widetilde{h}_{\text{CvM}}(Y_{j_1},Y_{j_2};Y_{j_3},Y_{j_4}) \convP Q_5 = 0$.

\vskip 1em

In the next lemmas, we show that $Q_1$ is strictly greater than any of $Q_2$, $Q_3$, $Q_4$ and $Q_5$ whenever $\overline{\delta}_{XY}^2 > 0$ or $\overline{\sigma}_X^2 \neq \overline{\sigma}_Y^2$. In addition they all become equivalent to each other only when $\overline{\delta}_{XY}^2 = 0$ and $\overline{\sigma}_X^2 = \overline{\sigma}_Y^2$. We start by proving that the inverse cosine function is concave on $x \in [0,1]$. 

%To prove these results, we will apply reverse Jensen's inequality based on the fact that the inverse cosine function is concave on $x \in [0,1]$.

\begin{lemma} \label{Lemma: Concave function}
	The inverse cosine function is concave on $x \in [0,1]$. 
	\begin{proof}
		The result follows by observing that
		\begin{align*} 
		\frac{d}{dx}\text{arccos}(x) = - \frac{1}{\sqrt{1 - x^2}} \quad \text{and} \quad \frac{d^2}{dx^2}\text{arccos}(x)  = - \frac{x}{(1-x^2)^{3/2}}. 
		\end{align*} 
	\end{proof}
\end{lemma}

\vskip 1em

\begin{lemma}
	Assume \textbf{\emph{(A1)}} and \textbf{\emph{(A2)}} hold. Then we have $Q_1 \geq Q_2$ and the equality holds if and only if $\overline{\delta}_{XY}^2 = 0$ or $\overline{\sigma}_X^2 = \overline{\sigma}_Y^2$.
	\begin{proof}
		From Lemma~\ref{Lemma: Concave function}, the inverse cosine function is concave on $x \in [0,1]$.  So we apply reverse Jensen's inequality to have 
		\begin{align*}
		\text{{arccos}} \left( \frac{\overline{\delta}_{XY}^2 + \overline{\sigma}_{X}^2}{\overline{\delta}_{XY}^2 + \overline{\sigma}_{X}^2 + \overline{\sigma}_{Y}^2} \right) + \text{{arccos}} \left( \frac{\overline{\delta}_{XY}^2 + \overline{\sigma}_{Y}^2}{\overline{\delta}_{XY}^2 + \overline{\sigma}_{X}^2 + \overline{\sigma}_{Y}^2} \right) \leq 2 \text{arccos} \left(  \frac{2\overline{\delta}_{XY}^2 + \overline{\sigma}_{X}^2 + \overline{\sigma}_{Y}^2}{2(\overline{\delta}_{XY}^2 + \overline{\sigma}_{X}^2 + \overline{\sigma}_{Y}^2)}\right).
		\end{align*}
		Then it is enough to show that 
		\begin{align} \nonumber
		& \text{{arccos}} \left(  \frac{\overline{\sigma}_{X}^2}{(2\overline{\sigma}_{X}^2)^{1/2} (\overline{\delta}_{XY}^2 + \overline{\sigma}_{X}^2 + \overline{\sigma}_{Y}^2 )^{1/2} } \right) + \text{{arccos}} \left(  \frac{\overline{\sigma}_{Y}^2}{( 2\overline{\sigma}_{Y}^2)^{1/2} (\overline{\delta}_{XY}^2 + \overline{\sigma}_{X}^2 + \overline{\sigma}_{Y}^2 )^{1/2} } \right) \\[.5em]
		\geq ~ & 2 \text{arccos} \left(  \frac{2\overline{\delta}_{XY}^2 + \overline{\sigma}_{X}^2 + \overline{\sigma}_{Y}^2}{2 (\overline{\delta}_{XY}^2 + \overline{\sigma}_{X}^2 + \overline{\sigma}_{Y}^2)}\right). \label{Eq: Middle Step 3}
		\end{align}
		
		\vskip 1em
		
		Before we proceed, we introduce the following quantities to simplify the expressions.
		\begin{align*}
		T_{XY} ~ = ~ & \frac{2\overline{\delta}_{XY}^2 + \overline{\sigma}_{X}^2 + \overline{\sigma}_{Y}^2}{2 (\overline{\delta}_{XY}^2 + \overline{\sigma}_{X}^2 + \overline{\sigma}_{Y}^2)}, \\[.5em]
		T_{X} ~ = ~ & \frac{\overline{\sigma}_{X}^2}{ (2\overline{\sigma}_{X}^2)^{1/2} (\overline{\delta}_{XY}^2 + \overline{\sigma}_{X}^2 + \overline{\sigma}_{Y}^2 )^{1/2}}, \\[.5em]
		T_{Y} ~ = ~ & \frac{\overline{\sigma}_{Y}^2}{( 2\overline{\sigma}_{Y}^2)^{1/2} (\overline{\delta}_{XY}^2 + \overline{\sigma}_{X}^2 + \overline{\sigma}_{Y}^2)^{1/2}}
		\end{align*}
		and
		\begin{align*}
		T_1 ~ = ~& \overline{\delta}_{XY}^2 (\overline{\sigma}_{X}^2+ 2\overline{\sigma}_{Y}^2 +2\overline{\delta}_{XY}^2)^{1/2} \{2\overline{\sigma}_{X}^2 + \overline{\sigma}_{Y}^2 +2\overline{\delta}_{XY}^2\}^{1/2}, \\[.5em]
		T_2 ~=~&  \overline{\delta}_{XY}^2 ( 2\overline{\delta}_{XY}^2  - \overline{\sigma}_{X}\overline{\sigma}_{Y}),  \\[.5em]
		T_3 ~=~ &  (\overline{\sigma}_{X}^2+ \overline{\sigma}_{Y}^2) (\overline{\sigma}_{X}^2+ 2\overline{\sigma}_{Y}^2 +2\overline{\delta}_{XY}^2)^{1/2} (2\overline{\sigma}_{X}^2 + \overline{\sigma}_{Y}^2 +2\overline{\delta}_{XY}^2)^{1/2}, \\[.5em]
		T_4 ~= ~& -(\overline{\sigma}_{X}^2+ \overline{\sigma}_{Y}^2 )(\overline{\sigma}_{X}^2+ \overline{\sigma}_{Y}^2 + \overline{\sigma}_{X}\overline{\sigma}_{Y}).
		\end{align*}
		Based on the monotonicity of the inverse cosine function and the basic identity 
		\begin{align*}
		\text{arccos}(x) + \text{arccos}(y) = \text{arccos}\big(xy - \sqrt{1-x^2} \sqrt{1-y^2}\big) \quad \text{ for $0 \leq x,y \leq 1$,}
		\end{align*}
		it can be seen that proving the inequality (\ref{Eq: Middle Step 3}) is equivalent to proving
		\begin{align} \label{Eq: Middle Step 4}
		2 T_{XY}^2 - 1 \geq T_X T_Y - (1- T_X^2)^{1/2} (1-T_Y^2)^{1/2}.
		\end{align}
		After rearrangement, it can be further seen that the inequality (\ref{Eq: Middle Step 4}) is equivalent to
		\begin{align} \label{Eq: Middle Step 5}
		T_1 + T_2 + T_3 + T_4  \geq 0.
		\end{align}
		The inequality (\ref{Eq: Middle Step 5}) is indeed true and the equality holds only when $\overline{\delta}_{XY} = 0$ and $\overline{\sigma}_{X}^2 = \overline{\sigma}_{Y}^2$ since
		\begin{align*}
		T_1 + T_2 \geq 0 \quad 	\text{\emph{if and only if}} \quad   \overline{\delta}_{XY}^4\{(6\overline{\sigma}_{X}^2+4\overline{\sigma}_{X}\overline{\sigma}_{Y} +6 \overline{\sigma}_{Y}^2)\overline{\delta}_{XY}^2 + 2( \overline{\sigma}_{X}^2 +   \overline{\sigma}_{Y}^2)^2\} \geq 0
		\end{align*}
		and
		\begin{align*}
		& T_3 + T_4 \geq 0 \quad  \text{\emph{if and only if}} \\[.5em]
		& (\overline{\sigma}_{X}^2 + \overline{\sigma}_{Y}^2)(\overline{\sigma}_{X} - \overline{\sigma}_{Y})^2 + 2\overline{\delta}_{XY}^2(2\overline{\sigma}_{X}^2 + \overline{\sigma}_{Y}^2 ) + 2\overline{\delta}_{XY}^2(\overline{\sigma}_{X}^2 +2 \overline{\sigma}_{Y}^2 ) \geq 0.
		\end{align*}
		This completes the proof. 
	\end{proof}
\end{lemma}

\vskip 1em

\begin{lemma} \label{Lemma: HDLSS inequality}
	Assume \textbf{\emph{(A1)}} and \textbf{\emph{(A2)}} hold. Then we have $Q_1 \geq Q_3$ and the equality holds if and only if $\overline{\delta}_{XY}^2 = 0$ or $\overline{\sigma}_X^2 = \overline{\sigma}_Y^2$.
	\begin{proof}
		Using reverse Jensen's inequality, we have
		\begin{align*}
		\text{{arccos}} \left( \frac{1}{2} \right) \geq \frac{1}{2}\text{{arccos}} \left( \frac{\overline{\delta}_{XY}^2 + \overline{\sigma}_{X}^2}{\overline{\delta}_{XY}^2 + \overline{\sigma}_{X}^2 + \overline{\sigma}_{Y}^2} \right) + \frac{1}{2} \text{{arccos}} \left( \frac{\overline{\delta}_{XY}^2 + \overline{\sigma}_{Y}^2}{\overline{\delta}_{XY}^2 + \overline{\sigma}_{X}^2 + \overline{\sigma}_{Y}^2} \right)
		\end{align*}
		where the equality holds only when $\overline{\delta}_{XY} = 0$ and $\overline{\sigma}_{X}^2 = \overline{\sigma}_{Y}^2$. Then it is enough to verify that
		\begin{equation}
		\begin{aligned}   \label{Eq: Middle Step 1}
		& \text{arccos} \left(  \frac{\overline{\sigma}_{X}^2}{ (2 \overline{\sigma}_X^2)^{1/2} ( \overline{\delta}_{XY}^2 + \overline{\sigma}_{X}^2 + \overline{\sigma}_{Y}^2)^{1/2}} \right)  \\[.5em]
		\geq ~ & \frac{3}{4}\text{{arccos}} \left( \frac{\overline{\delta}_{XY}^2 + \overline{\sigma}_{X}^2}{\overline{\delta}_{XY}^2 + \overline{\sigma}_{X}^2 + \overline{\sigma}_{Y}^2} \right) + \frac{1}{4} \text{{arccos}} \left( \frac{\overline{\delta}_{XY}^2 + \overline{\sigma}_{Y}^2}{\overline{\delta}_{XY}^2 + \overline{\sigma}_{X}^2 + \overline{\sigma}_{Y}^2} \right).
		\end{aligned}
		\end{equation}
		By applying reverse Jensen's inequality and by the monotonicity of the inverse cosine function, it is seen that the following statement
		\begin{align} 
		 \frac{4\overline{\delta}_{XY}^2 + 3\overline{\sigma}_{X}^2 +  \overline{\sigma}_{Y}^2}{4(\overline{\delta}_{XY}^2 + \overline{\sigma}_{X}^2 + \overline{\sigma}_{Y}^2)}  ~\geq ~ \frac{\overline{\sigma}_{X}^2}{(2 \overline{\sigma}_X^2)^{1/2} ( \overline{\delta}_{XY}^2 + \overline{\sigma}_{X}^2 + \overline{\sigma}_{Y}^2)^{1/2}} \label{Eq: Middle Step 2-a}
		\end{align}
		implies (\ref{Eq: Middle Step 1}). Since (\ref{Eq: Middle Step 2-a}) is true if and only if 
		\begin{align}
		16 \overline{\delta}_{XY}^4 + 16 \overline{\delta}_{XY}^2 \overline{\sigma}_{X}^2  + 8 \overline{\delta}_{XY}^2 \overline{\sigma}_{Y}^2  + (\overline{\sigma}_{X}^2 - \overline{\sigma}_{Y}^2 )^2 \geq 0 \label{Eq: Middle Step 2}
		\end{align}		
		and the equality of (\ref{Eq: Middle Step 2}) holds only if $\overline{\delta}_{XY} = 0$ and $\overline{\sigma}_{X}^2 = \overline{\sigma}_{Y}^2$, the result follows. 
	\end{proof}
\end{lemma}

\vskip 1em

\begin{lemma}
	Assume \textbf{\emph{(A1)}} and \textbf{\emph{(A2)}} hold. Then we have $Q_1 \geq Q_4$ and the equality holds if and only if $\overline{\delta}_{XY}^2 = 0$ or $\overline{\sigma}_X^2 = \overline{\sigma}_Y^2$.
	\begin{proof}
		The proof is similar to that of Lemma~\ref{Lemma: HDLSS inequality}. Hence we omit the proof. 
	\end{proof}
\end{lemma}

\vskip 1em

\begin{lemma}
	Assume \textbf{\emph{(A1)}} and \textbf{\emph{(A2)}} hold. Then we have $Q_1 \geq Q_5$ and the equality holds if and only if $\overline{\delta}_{XY}^2 = 0$ or $\overline{\sigma}_X^2 = \overline{\sigma}_Y^2$. 
	\begin{proof}
		Using reverse Jensen's inequality, we see that
		\begin{align*}
		& \frac{1}{\pi} \text{arccos} \left(  \frac{2\overline{\delta}_{XY}^2 + \overline{\sigma}_{X}^2 + \overline{\sigma}_{Y}^2}{2(\overline{\delta}_{XY}^2 + \overline{\sigma}_{X}^2 + \overline{\sigma}_{Y}^2)} \right) \\[.5em]
		\geq ~ &\frac{1}{2\pi} \text{arccos} \left( \frac{\overline{\delta}_{XY}^2 + \overline{\sigma}_{X}^2}{\overline{\delta}_{XY}^2 + \overline{\sigma}_{X}^2 + \overline{\sigma}_{Y}^2} \right) + \frac{1}{2\pi} \text{arccos} \left( \frac{\overline{\delta}_{XY}^2 + \overline{\sigma}_{Y}^2}{\overline{\delta}_{XY}^2 + \overline{\sigma}_{X}^2 + \overline{\sigma}_{Y}^2} \right). 
		\end{align*}
		In addition, the inverse cosine function is monotone decreasing. So
		\begin{align*}
		\frac{1}{\pi} \text{arccos} \left(  \frac{2\overline{\delta}_{XY}^2 + \overline{\sigma}_{X}^2 + \overline{\sigma}_{Y}^2}{2(\overline{\delta}_{XY}^2 + \overline{\sigma}_{X}^2 + \overline{\sigma}_{Y}^2)} \right) \leq \frac{1}{\pi} \text{arccos} \left(  \frac{\overline{\delta}_{XY}^2 + \overline{\sigma}_{X}^2 + \overline{\sigma}_{Y}^2}{2(\overline{\delta}_{XY}^2 + \overline{\sigma}_{X}^2 + \overline{\sigma}_{Y}^2)} \right)  = \frac{1}{3},
		\end{align*}
		where the last step uses 
		\begin{align*}
		\frac{1}{\pi}\text{arccos} \left( \frac{1}{2} \right) = \frac{1}{3}.
		\end{align*}
		Notice that the first inequality becomes the equality only when $\overline{\sigma}_{X}^2 = \overline{\sigma}_{Y}^2$. The second inequality becomes the equality only when $\overline{\delta}_{XY}^2 = 0$. This proves the result.
	\end{proof}
\end{lemma}

\vskip 1em

Combining the previous lemmas, we give a summary:
\begin{lemma} \label{Lemma: HDLSS kernel limit}
	Assume \textbf{\emph{(A1)}} and \textbf{\emph{(A2)}} hold. Then we have
	\begin{align*}
	Q_1 \geq \max\{Q_2,Q_3,Q_4,Q_5 \}
	\end{align*}
	and the equality holds as $Q_1=Q_2=Q_3=Q_4=Q_5$ if and only if $\overline{\delta}_{XY}^2 = 0$ or $\overline{\sigma}_X^2 = \overline{\sigma}_Y^2$. 
\end{lemma}

\vskip 1em 

\noindent \textbf{$\bullet$ Part 2.}

\vskip .3em

\noindent In this part, we prove Theorem~\ref{Theorem: HDLSS Consistency}. Notice that $U_{\text{CvM}}$ is a linear combination of kernel $\widetilde{h}_{\text{CvM}}$ evaluated on different samples. Hence from the previous observation made in Part 1, it is seen that
\begin{align*}
U_{\text{CvM}} \convP Q_1 \quad \text{under $H_1$.}
\end{align*}
For a given permutation $\varpi$ of $\{1,\ldots,N\}$, let us denote by $U_{\text{CvM}}^\varpi$, the $U$-statistic computed based on $\{Z_{\varpi(1)},\ldots,Z_{\varpi(N)}\}$, i.e.~$U_{\text{CvM}}(Z_{\varpi(1)},\ldots,Z_{\varpi(N)})$. Let $\varpi_0 = \{1,\ldots,N\}$ be the original permutation. Then $U_{\text{CvM}}^{\varpi_0}$ becomes $U_{\text{CvM}}(Z_1,\ldots,Z_N)$ computed based on the original samples. Let us define that the permutation $\varpi$ is a neighbor of $\varpi_0$ if $\#|\{ \varpi(1),\ldots,\varpi(m)\} \cap \{1,\ldots,m\} | = m$.

We first consider the unbalanced case where $m \neq n$.  Observe that $U_{\text{CvM}}^\varpi$ converges to $Q_{\varpi}$, which is a weighted average of $Q_1,\ldots,Q_5$. According to Lemma~\ref{Lemma: HDLSS kernel limit}, $Q_1 \geq Q_{\varpi}$ and it is not difficult to see that $Q_1 = Q_\varpi$ only if $\varpi$ is a neighbor of $\varpi_0$. This means that $U_{\text{CvM}}^{\varpi_0} > U_{\text{CvM}}^\varpi$ in the limit for all $\varpi$ but neighbors of $\varpi_0$ under $H_1$. Since there are $m!n!$ neighbors of $\varpi_0$ out of $N!$ permutations, if we choose $\alpha > 1/\{N!/(m!n!)\}$, then we have $\lim_{d \rightarrow \infty} \mE [\phi_{\text{CvM}}] = 1$.

For the balanced case where $m=n$, the result follows by a similar argument but now we also need to consider $\varpi$ that satisfies $\#|\{ \varpi(1),\ldots,\varpi(m)\} \cap \{m+1,\ldots,m+n\} | = n$ to be a neighbor of $\varpi_0$. This is because $U_{\text{CvM}}(Z_1,\ldots,Z_N) = U_{\text{CvM}}(Z_N,\ldots,Z_1)$ if $m=n$. Hence now we have $2m!n!$ neighbors of $\varpi_0$ out of $N!$ permutations and if we choose $\alpha > 2/\{N!/(m!n!)\}$, then we have $\lim_{d \rightarrow \infty} \mE [\phi_{\text{CvM}}] = 1$.

\vskip 2em

\subsection{Proof of Theorem~\ref{Theorem: HDLSS Equivalence}} \label{Section: Connection of different statistics to CQ}
Our strategy to prove the given result is to connect different statistics to the CQ statistic, which is relatively easy to handle. Each connection can be found in
\begin{itemize}
	\item Section~\ref{Section: Connection of CvM to CQ}: Connection of $U_{\text{CvM}}^\varpi$ to $U_{\text{CQ}}^\varpi$,
	\item Section~\ref{Section: Connection of WMW to CQ}: Connection of $U_{\text{WMW}}^\varpi$ to $U_{\text{CQ}}^\varpi$,
	\item Section~\ref{Section: Connection of Energy to CQ}: Connection of $U_{\text{Energy}}^\varpi$ to $U_{\text{CQ}}^\varpi$,
	\item Section~\ref{Section: Connection of MMD to CQ}: Connection of $U_{\text{MMD}}^\varpi$ to $U_{\text{CQ}}^\varpi$. 
\end{itemize}
For notational simplicity, we will denote $Z_i^\ast,Z_2^\ast,Z_3^\ast,Z_4^\ast$ by $Z_1,Z_2,Z_3,Z_4$ throughout this section.

%First we collect some lemmas and their proofs that will be used throughout this section. Since we assume that $\Sigma = \Sigma_X = \Sigma_Y$, we  denote $d^{-1} \tr(\Sigma) := \overline{\sigma}_d^2 \rightarrow \overline{\sigma}^2$.

\subsubsection{Connection of $U_{\text{CvM}}^\varpi$ to $U_{\text{CQ}}^\varpi$} \label{Section: Connection of CvM to CQ}
In this subsection, we connect $U_{\text{CvM}}^\varpi$ to $U_{\text{CQ}}^\varpi$ under the HDLSS setting. We first list some lemmas and their proofs. The final connection between $U_{\text{CvM}}^\varpi$ and $U_{\text{CQ}}^\varpi$ can be found in Proposition~\ref{Proposition: CvM and CQ}. 

\vskip 1em 

\begin{lemma} \label{Lemma: Basic approximation}
	Under \textbf{\emph{(A1)}}, \textbf{\emph{(A2)}} and \textbf{\emph{(A4)}}, we have
	\begin{align*}
	& \frac{1}{d} \| Z_1 - Z_2\|^2 - 2\overline{\sigma}_d^2 = O_{\mathbb{P}}(d^{-1/2}) \quad \text{and} \\[.5em]
	& \frac{1}{d}(Z_1 - Z_3)^\top (Z_2 - Z_3) =  \overline{\sigma}_d^2  + O_{\mathbb{P}}(d^{-1/2}). 
	\end{align*}
	\begin{proof}
 		Under the assumption that $\mV[\|Z_1 - Z_2\|^2] = O(d)$, we apply Chebyshev's inequality to obtain
 		\begin{align*}
 		\frac{1}{d} \| Z_1 - Z_2\|^2 - \frac{1}{d}\mE[ \| Z_1 - Z_2\|^2 ] =  O_{\mathbb{P}}(d^{-1/2}).
 		\end{align*}
 		Note that regardless of the distributions of $Z_1$ and $Z_2$, the expected value of $\| Z_1 - Z_2\|^2$ is bounded by
 		\begin{align*}
 		\mE [\| Z_1 - Z_2\|^2 ] \leq \|\mu_X - \mu_Y\|^2 + 2 \tr(\Sigma^2).
 		\end{align*}
 		Thus under \textbf{{(A4)}},
 		\begin{align*}
 		\frac{1}{d}\mE[ \| Z_1 - Z_2\|^2 ] -  2\overline{\sigma}_d^2 = O(d^{-1/2}).
 		\end{align*}
		By combining the results, we prove the first part. The second part follows similarly. 
	\end{proof}
\end{lemma}

\vskip 1em

\begin{lemma} \label{Lemma: 1/sqrt(x) expansion}
	Under \textbf{\emph{(A1)}}, \textbf{\emph{(A2)}} and \textbf{\emph{(A4)}}, we have
	\begin{align*}
	\frac{\sqrt{d}}{\|Z_1 - Z_2\|} = \frac{1}{(2\overline{\sigma}_d^2)^{1/2}} - \frac{1}{2(2\overline{\sigma}_d^2)^{3/2}} \left( d^{-1} \|Z_1 - Z_2\|^2 - 2\overline{\sigma}_d^2  \right) + O_{\mathbb{P}}(d^{-1}).
	\end{align*}
	\begin{proof}
		Consider $f(x) = 1/\sqrt{x}$ and represent 
		\begin{align*}
		f \big(d^{-1}\|Z_1 - Z_2\|^2 \big) = \frac{\sqrt{d}}{\|Z_1 - Z_2\|}.
		\end{align*}
		By using the second order Taylor expansion of $f(x)$ around $f(2\overline{\sigma}_d^2)$ with Lemma~\ref{Lemma: Basic approximation}, we obtain the result.
	\end{proof}
\end{lemma}

\vskip 1em

\begin{lemma}  \label{Lemma: inverse function approximation with two terms}
	Under \textbf{\emph{(A1)}}, \textbf{\emph{(A2)}} and \textbf{\emph{(A4)}}, we have
	\begin{align*}
	\frac{d}{\|Z_1 - Z_3\|\|Z_2-Z_3\|} = \frac{1}{2\overline{\sigma}_d^2} & - \frac{1}{8\overline{\sigma}_d^4} \left( d^{-1} \|Z_1 - Z_3\|^2 - 2\overline{\sigma}_d^2 \right) \\[.5em]
	& - \frac{1}{8\overline{\sigma}_d^4} \left( d^{-1} \|Z_2 - Z_3\|^2 - 2\overline{\sigma}_d^2 \right) + O_{\mathbb{P}}(d^{-1}).
	\end{align*}
	\begin{proof}
		Based on Lemma~\ref{Lemma: 1/sqrt(x) expansion}, we have
		\begin{align*}
		\frac{d}{\|Z_1 - Z_3\|\|Z_2-Z_3\|} = ~& \bigg\{ \frac{1}{(2\overline{\sigma}_d^2)^{1/2}} - \frac{1}{2(2\overline{\sigma}_d^2)^{3/2}} \left( d^{-1} \|Z_1 - Z_3\|^2 - 2\overline{\sigma}_d^2  \right) + O_{\mathbb{P}}(d^{-1}) \bigg\} \\[.5em]
		\times ~ & \bigg\{ \frac{1}{(2\overline{\sigma}_d^2)^{1/2}} - \frac{1}{2(2\overline{\sigma}_d^2)^{3/2}} \left( d^{-1} \|Z_2 - Z_3\|^2 - 2\overline{\sigma}_d^2  \right) + O_{\mathbb{P}}(d^{-1}) \bigg\}. 
		\end{align*}
		By expanding the right-hand side and the following observations made from Lemma~\ref{Lemma: Basic approximation},
		\begin{align*}
		& \frac{1}{2(2\overline{\sigma}_d^2)^{3/2}} \left( d^{-1} \|Z_1 - Z_3\|^2 - 2\overline{\sigma}_d^2  \right) = O_{\mathbb{P}}(d^{-1/2}), \\[.5em]
		& \frac{1}{2(2\overline{\sigma}_d^2)^{3/2}} \left( d^{-1} \|Z_2 - Z_3\|^2 - 2\overline{\sigma}_d^2  \right) = O_{\mathbb{P}}(d^{-1/2}),
		\end{align*}
		the result follows. 
	\end{proof}
\end{lemma}

\vskip 1em

\begin{lemma} \label{Lemma: CvM kernel 1}
	Under \textbf{\emph{(A1)}}, \textbf{\emph{(A2)}} and \textbf{\emph{(A4)}}, we have
	\begin{align*}
	& \text{\emph{arccos}} \Bigg\{ \frac{(Z_1 - Z_3)^\top (Z_2 - Z_3)}{\|Z_1 - Z_3\| \|Z_2 - Z_3\| } \Bigg\}  \\[.5em]
	= ~ & \text{\emph{arccos}}\left(\frac{1}{2}\right) - \frac{2}{\sqrt{3}} \Bigg\{ \frac{(Z_1 - Z_3)^\top (Z_2 - Z_3)}{\|Z_1 - Z_3\| \|Z_2 - Z_3\| }  - \frac{1}{2} \Bigg\} + O_{\mathbb{P}}(d^{-1}).
	\end{align*}
	\begin{proof}
		First note that 
		\begin{align*}
		\frac{(Z_1 - Z_3)^\top (Z_2 - Z_3)}{\|Z_1 - Z_3\|\|Z_2-Z_3\|} - \frac{1}{2} = O_{\mathbb{P}}(d^{-1/2}),
		\end{align*}
		which follows from Lemma~\ref{Lemma: Basic approximation} and Lemma~\ref{Lemma: inverse function approximation with two terms}. We then use the second order Taylor expansion of the inverse cosine function around $\text{arccos}(1/2)$ to obtain the result.
	\end{proof}
\end{lemma}

\vskip 1em

\begin{lemma} \label{Lemma: CvM kernel 2}
	Under \textbf{\emph{(A1)}}, \textbf{\emph{(A2)}} and \textbf{\emph{(A4)}}, we have
	\begin{align*}
	\frac{(Z_1 - Z_3)^\top (Z_2 - Z_3)}{\|Z_1 - Z_3\|\|Z_2-Z_3\|} - \frac{1}{2}  =  & \frac{(Z_1 - Z_3)^\top (Z_2 - Z_3) - d \overline{\sigma}_d^2 }{2d \overline{\sigma}_d^2} \\[.5em] 
	& - \frac{1}{8d \overline{\sigma}_d^2} \left( \|Z_1 - Z_3\|^2 + \|Z_2-Z_3\|^2 - 4d \overline{\sigma}_d^2 \right) + O_{\mathbb{P}}(d^{-1}).
	\end{align*}
	\begin{proof}
		We split the left-hand side into two terms:
		\begin{align*}
		\frac{(Z_1 - Z_3)^\top (Z_2 - Z_3)}{\|Z_1 - Z_3\|\|Z_2-Z_3\|} - \frac{1}{2}  = ~ &  \frac{(Z_1 - Z_3)^\top (Z_2 - Z_3)}{\|Z_1 - Z_3\|\|Z_2-Z_3\|} - \frac{(Z_1 - Z_3)^\top (Z_2 - Z_3)}{2d \overline{\sigma}_d^2} \\[.5em]
		+ ~ & \frac{(Z_1 - Z_3)^\top (Z_2 - Z_3)}{2d \overline{\sigma}_d^2} - \frac{1}{2}.
		\end{align*}
		Now it is enough to show that
		\begin{align*}
		& \frac{(Z_1 - Z_3)^\top (Z_2 - Z_3)}{\|Z_1 - Z_3\|\|Z_2-Z_3\|} - \frac{(Z_1 - Z_3)^\top (Z_2 - Z_3)}{2d \overline{\sigma}_d^2} \\[.5em]
		= & - \frac{1}{8d \overline{\sigma}_d^2} \left( \|Z_1 - Z_3\|^2 + \|Z_2-Z_3\|^2 - 4d \overline{\sigma}_d^2 \right) + O_{\mathbb{P}}(d^{-1}).
		\end{align*}
		Note that 
		\begin{align*}
		& \frac{(Z_1 - Z_3)^\top (Z_2 - Z_3)}{\|Z_1 - Z_3\|\|Z_2-Z_3\|} - \frac{(Z_1 - Z_3)^\top (Z_2 - Z_3)}{2d \overline{\sigma}_d^2}  \\[.5em]
		= ~ & (Z_1 - Z_3)^\top (Z_2 - Z_3) \times \bigg(  \frac{1}{\|Z_1 - Z_3\|\|Z_2-Z_3\|} - \frac{1}{2d \overline{\sigma}_d^2} \bigg) \\[.5em]
		= ~ & (I) \times (II) \quad \text{(say).}
		\end{align*}
		From Lemma~\ref{Lemma: Basic approximation} and Lemma~\ref{Lemma: inverse function approximation with two terms}, it is seen that
		\begin{align*}
		& (I) = d \overline{\sigma}_d^2 + O_{\mathbb{P}}(d^{1/2}), \\[.5em]
		&  (II) =  - \frac{1}{8d\overline{\sigma}_d^4} \bigg[ d^{-1} \|Z_1 -Z_3\|^2 + d^{-1}\|Z_1 -Z_3\|^2 - 4 \overline{\sigma}_d^2 + O_{\mathbb{P}}(d^{-2}) \bigg].
		\end{align*}
		Expanding the terms in $(I) \times (II)$, we obtain the result.  
	\end{proof}
\end{lemma}

\vskip 1em

Based on the previous lemmas, we prove the main result of this subsection.

\begin{proposition} \label{Proposition: CvM and CQ}
	Under \textbf{\emph{(A1)}}, \textbf{\emph{(A2)}} and \textbf{\emph{(A4)}}, we have
	\begin{equation}
	\begin{aligned} \label{Eq: kernel approximation between h_CvM and h_CQ}
	& \widetilde{h}_{\text{\emph{CvM}}}(Z_1,Z_2;Z_3,Z_4) \\[.5em]
	= ~ & \frac{1}{4\pi \sqrt{3} d \overline{\sigma}_d^2} \{ (Z_1-Z_3)^\top (Z_2-Z_4)  + (Z_1 - Z_4)^\top (Z_2 - Z_3)\} + O_{\mathbb{P}}(d^{-1})
	\end{aligned}
	\end{equation}
	and thus
	\begin{align*}
	U_{\text{\emph{CvM}}}^{\varpi} = \frac{1}{2\pi \sqrt{3} d \overline{\sigma}_d^2} U_{\text{\emph{CQ}}}^{\varpi} + O_{\mathbb{P}}(d^{-1}).
	\end{align*}
	\begin{proof}
		By Lemma~\ref{Lemma: CvM kernel 1} and Lemma~\ref{Lemma: CvM kernel 2},
		\begin{align*}
		& \text{{arccos}} \Bigg\{ \frac{(Z_1 - Z_3)^\top (Z_2 - Z_3)}{\|Z_1 - Z_3\| \|Z_2 - Z_3\| } \Bigg\} \\[.5em]
		= ~ & \text{arccos} \bigg(\frac{1}{2}\bigg)  -\frac{2}{\sqrt{3}}\Bigg\{ \frac{(Z_1-Z_3)^\top (Z_2-Z_3)}{2d\overline{\sigma}_d^2} - \frac{1}{2}  \\[.5em]
		& ~~~~~~~~~~~~~~~~~~~~~~~~ - \frac{1}{8d\overline{\sigma}_d^2} \bigg( \|Z_1 - Z_3\|^2 + \|Z_2 - Z_3\|^2 - 4d \overline{\sigma}_d^2 \bigg)  \Bigg\} + O_{\mathbb{P}}(d^{-1}).
		\end{align*}
		We can obtain (\ref{Eq: kernel approximation between h_CvM and h_CQ}) by first plugging the above approximation into $\widetilde{h}_{\text{CvM}}$ for each inverse cosine function and then simplifying the expression. The second result is trivial by noting that 
		\begin{align*}
		\widetilde{h}_{\text{CQ}}(x_1,x_2;y_1,y_2)= \frac{1}{2}(x_1-y_1)^\top (x_2-y_2)  + \frac{1}{2} (x_1 - y_2)^\top (x_2 - y_1)
		\end{align*}
		is the symmetrized kernel of the CQ statistic. 
	\end{proof}
\end{proposition}

\vskip 2em

\subsubsection{Connection of $U_{\text{WMW}}^\varpi$ to $U_{\text{CQ}}^\varpi$} \label{Section: Connection of WMW to CQ}
Note that the symmetrized kernel of the WMW statistic can be written as
\begin{align*}
\widetilde{h}_{\text{WMW}}(x_1,x_2;y_1,y_2) = \frac{1}{2}\frac{(x_1-y_1)^\top (x_2-y_2)}{\|x_1 - y_1\| \|x_2 - y_2\|} +  \frac{1}{2}\frac{(x_1-y_2)^\top (x_2-y_1)}{\|x_1 - y_2\| \|x_2 - y_1\|}.
\end{align*}
We first provide a couple of lemmas and their proofs. We then present the main result in Proposition~\ref{Proposition: WMW and CQ}.

\vskip 1em

\begin{lemma} \label{Lemma: WMW and CQ Middle step}
	Under \textbf{\emph{(A1)}}, \textbf{\emph{(A2)}}, \textbf{\emph{(A3)}} and \textbf{\emph{(A4)}}, we have
	\begin{align*}
	\frac{d}{\|Z_1 - Z_2\|\|Z_3-Z_4\|} = \frac{1}{2\overline{\sigma}_d^2} & - \frac{1}{8\overline{\sigma}_d^4} \left( d^{-1} \|Z_1 - Z_2\|^2 - 2\overline{\sigma}_d^2 \right) \\[.5em]
	&- \frac{1}{8\overline{\sigma}_d^4} \left( d^{-1} \|Z_3 - Z_4\|^2 - 2\overline{\sigma}_d^2 \right) + O_{\mathbb{P}}(d^{-1}).
	\end{align*}
	\begin{proof}
	The proof is similar to Lemma~\ref{Lemma: inverse function approximation with two terms}; hence omitted.
	\end{proof}
\end{lemma}

\vskip 1em

\begin{lemma} \label{Lemma: WMW and CQ Middle step 2}
	Under \textbf{\emph{(A1)}}, \textbf{\emph{(A2)}}, \textbf{\emph{(A3)}} and \textbf{\emph{(A4)}}, we have
	\begin{align*}
	\frac{(Z_1-Z_3)^\top (Z_2-Z_4)}{\|Z_1 - Z_3\| \|Z_2 - Z_4\|} = \frac{(Z_1-Z_3)^\top (Z_2-Z_4)}{2d\overline{\sigma}_d^2 } + O_{\mathbb{P}}(d^{-1}).
	\end{align*}
	\begin{proof}
		Under \textbf{(A3)}, it can be seen as similar to Lemma~\ref{Lemma: Basic approximation} that 
		\begin{align*}
		d^{-1}(Z_1 - Z_3)^\top (Z_2 - Z_4) =  O_{\mathbb{P}}(d^{-1/2}).
		\end{align*}
		Then combining the above with Lemma~\ref{Lemma: Basic approximation} and Lemma~\ref{Lemma: WMW and CQ Middle step},
		\begin{align*}
		& \frac{(Z_1-Z_3)^\top (Z_2-Z_4)}{\|Z_1 - Z_3\| \|Z_2 - Z_4\|}  - \frac{(Z_1-Z_3)^\top (Z_2-Z_4)}{2d\overline{\sigma}_d^2}  \\[.5em]
		= ~ & d^{-1}(Z_1 - Z_3)^\top (Z_2 - Z_4) \times \Bigg\{  \frac{d}{\|Z_1 - Z_3\| \|Z_2 - Z_4\|} - \frac{1}{2\overline{\sigma}_d^2}  \Bigg\} \\[.5em] 
		= ~ & O_{\mathbb{P}}(d^{-1/2}) \times O_{\mathbb{P}}(d^{-1/2}). 
		\end{align*}
		Hence the result follows. 
	\end{proof}
\end{lemma}

\vskip 1em

Based on the previous lemmas, we prove the main result of this subsection.

\begin{proposition} \label{Proposition: WMW and CQ}
	Under \textbf{\emph{(A1)}}, \textbf{\emph{(A2)}}, \textbf{\emph{(A3)}} and \textbf{\emph{(A4)}}, we have
	\begin{align*}
	& \widetilde{h}_{\text{\emph{WMW}}}(Z_1,Z_2;Z_3,Z_4) \\[.5em]
	= ~& \frac{1}{2d\overline{\sigma}_d^2}  \{ (Z_1-Z_3)^\top (Z_2-Z_4)  + (Z_1 - Z_4)^\top (Z_2 - Z_3)\} + O_{\mathbb{P}}(d^{-1})
	\end{align*}
	and thus 
	\begin{align*}
	U_{\text{\emph{WMW}}} =  \frac{1}{2d\overline{\sigma}_d^2} U_{\text{\emph{CQ}}} + O_{\mathbb{P}}(d^{-1}).
	\end{align*}
	\begin{proof}
		The result is a direct consequence of Lemma~\ref{Lemma: WMW and CQ Middle step 2}.
	\end{proof}
\end{proposition}

\vskip 2em

\subsubsection{Connection of $U_{\text{Energy}}^\varpi$ to $U_{\text{CQ}}^\varpi$} \label{Section: Connection of Energy to CQ}
Next we find a connection between $U_{\text{Energy}}^\varpi$ and $U_{\text{CQ}}^\varpi$. Note that the symmetrized kernel of the energy statistic can be written as
\begin{align*}
\widetilde{h}_{\text{Energy}}(x_1,x_2;y_1,y_2) = & \frac{1}{2} \|x_1-y_1\| +  \frac{1}{2}\|x_1-y_2\| +  \frac{1}{2}\|x_2 - y_1\| +  \frac{1}{2}\|x_2 -y_2\| \\[.5em]
&  - \|x_1 - x_2\| - \|y_1 - y_2\|.
\end{align*}
Using this kernel expression, we connect $U_{\text{Energy}}$ to $U_{\text{CQ}}$ in Proposition~\ref{Proposition: Energy and CQ}. 

\vskip 1em

We start with one lemma. 

\begin{lemma} \label{Lemma: Energy and CQ Middle Step}
	Under \textbf{\emph{(A1)}} and \textbf{\emph{(A2)}}, we have
	\begin{align*}
	\frac{1}{\sqrt{d}}\|Z_1 - Z_2\| = (2\overline{\sigma}_d^2)^{1/2} + \frac{1}{2(2\overline{\sigma}_d^2)^{1/2}} \left( d^{-1}\|Z_1 - Z_2\|^2 - 2\overline{\sigma}_d^2 \right) + O_{\mathbb{P}}(d^{-1}).
	\end{align*}
	\begin{proof}
		We use the second order Taylor expansion of $f(x) = \sqrt{x}$ around $f(2\overline{\sigma}_d^2)$ with Lemma~\ref{Lemma: Basic approximation} to prove this result.
	\end{proof}
\end{lemma}

\vskip 1em

The main result of this subsection is stated as follows. 
\begin{proposition} \label{Proposition: Energy and CQ}
	Under \textbf{\emph{(A1)}} and \textbf{\emph{(A2)}}, we have
	\begin{align*}
	& \widetilde{h}_{\text{\emph{Energy}}}(Z_1,Z_2;Z_3,Z_4) \\[.5em]
	= ~ & \frac{1}{2 (2d \overline{\sigma}_d^2)^{1/2}} \{ (Z_1-Z_3)^\top (Z_2-Z_4)  + (Z_1 - Z_4)^\top (Z_2 - Z_3)\} + O_{\mathbb{P}}(d^{-1/2})
	\end{align*}
	and thus 
	\begin{align*}
	U_{\text{\emph{Energy}}} = \frac{1}{2 (d \overline{\sigma}_d^2)^{1/2}} U_{\text{\emph{CQ}}} + O_{\mathbb{P}}(d^{-1/2}).
	\end{align*}
	\begin{proof}
		We use Lemma~\ref{Lemma: Energy and CQ Middle Step} to approximate $\widetilde{h}_{\text{Energy}}$ to $\widetilde{h}_{\text{CQ}}$ and simplify the expression to obtain the first result. The second result is trivial. 
	\end{proof}
\end{proposition}

\vskip 2em

\subsubsection{Connection of $U_{\text{MMD}}^\varpi$ to $U_{\text{CQ}}^\varpi$} \label{Section: Connection of MMD to CQ}
In this subsection, we find a connection between $U_{\text{MMD}}^\varpi$ and $U_{\text{CQ}}^\varpi$. The symmetrized kernel of the MMD statistic can be written as
\begin{align*}
\widetilde{h}_{\text{MMD}}(x_1,x_2;y_1,y_2) = & -\frac{1}{2} \exp\bigg(-\frac{1}{2\varsigma_d^2}\|x_1-y_1\|^2\bigg) - \frac{1}{2} \exp\bigg(-\frac{1}{2\varsigma_d^2}\|x_1-y_2\|^2\bigg) \\[.5em] 
& -\frac{1}{2} \exp\bigg(-\frac{1}{2\varsigma_d^2}\|x_2-y_1\|^2\bigg)-\frac{1}{2} \exp\bigg(-\frac{1}{2\varsigma_d^2}\|x_2-y_2\|^2\bigg) \\[.5em] 
& +  \exp\bigg(-\frac{1}{2\varsigma_d^2}\|x_1-x_2\|^2\bigg)+ \exp\bigg(-\frac{1}{2\varsigma_d^2}\|y_1-y_2\|^2\bigg)
\end{align*}

and we assume that $\varsigma_d^2 \asymp d$. We first provide an approximation of the Gaussian kernel. 

\begin{lemma} \label{Lemma: MMD and CQ Middle Step}
	Under \textbf{\emph{(A1)}}, \textbf{\emph{(A2)}} and $\varsigma_d^2 \asymp d$, we have
	\begin{align*}
	&\exp \left(  - \frac{1}{2\varsigma_d^2} \|Z_1 - Z_2 \|^2 \right) \\[.5em] 
	= ~& \exp \left( - \frac{d\overline{\sigma}_d^2}{\varsigma_d^2} \right) - \exp \left( - \frac{d\overline{\sigma}_d^2}{\varsigma_d^2} \right) \left[ \frac{1}{2\varsigma_d^2} \|Z_1 - Z_2\|^2 - \frac{d\overline{\sigma}_d^2}{\varsigma_d^2} \right] + O_{\mathbb{P}}(d^{-1}).
	\end{align*}
	\begin{proof}
		We consider the second order Taylor expansion of $f(x) = e^{-x}$ around $f(d\overline{\sigma}_d^2 /\varsigma_d^2)$. Notice that under $\varsigma_d^2 \asymp d$, we have $d\overline{\sigma}_d^2 /\varsigma_d^2= O(1)$ and
		\begin{align*}
		 \frac{1}{2\varsigma_d^2} \|Z_1 - Z_2\|^2 - \frac{d\overline{\sigma}_d^2}{\varsigma_d^2} = \frac{d}{2\varsigma_d^2} \big(d^{-1}\|Z_1 - Z_2\|^2 - 2\overline{\sigma}_d^2 \big) = O_{\mathbb{P}}(d^{-1/2})
		\end{align*}
		from Lemma~\ref{Lemma: Basic approximation}. Thus the result follows. 
	\end{proof}
\end{lemma}

\vskip 1em

The main result of this subsection is stated as follows. 
\begin{proposition}  \label{Proposition: MMD and CQ}
	Under \textbf{\emph{(A1)}}, \textbf{\emph{(A2)}} and $\varsigma_d^2 \asymp d$, we have
	\begin{align*}
	\widetilde{h}_{\text{{\emph{MMD}}}}(Z_1,Z_2;Z_3,Z_4) =  \frac{e^{-d\overline{\sigma}_d^2/\varsigma_d^2}}{2\varsigma_d^2} \{ (Z_1-Z_3)^\top (Z_2-Z_4)  + (Z_1 - Z_4)^\top (Z_2 - Z_3)\} + O_{\mathbb{P}}(d^{-1})
	\end{align*}
	and thus 
	\begin{align*}
	U_{\text{\emph{MMD}}} =  \varsigma_d^{-2} e^{-d\overline{\sigma}_d^2/\varsigma_d^2} U_{\text{\emph{CQ}}} + O_{\mathbb{P}}(d^{-1/2}).
	\end{align*}
	\begin{proof}
		We use Lemma~\ref{Lemma: MMD and CQ Middle Step} to approximate $\widetilde{h}_{\text{MMD}}$ to $\widetilde{h}_{\text{CQ}}$ and simplify the expression to obtain the first result. The second result is trivial. 
	\end{proof}
\end{proposition}

\vskip 1em 

\noindent \textbf{$\bullet$ Main proof of Theorem~\ref{Theorem: HDLSS Equivalence}.} 

\vskip .3em

\noindent By collecting the results in Proposition~\ref{Proposition: CvM and CQ}, Proposition~\ref{Proposition: WMW and CQ}, Proposition~\ref{Proposition: Energy and CQ} and Proposition~\ref{Proposition: MMD and CQ}, it is easily checked that Theorem~\ref{Theorem: HDLSS Equivalence} holds and thus we complete the proof.

%
%\vskip 2em 
%
%
%
%\subsection{Proof of Corollary~\ref{Corollary: HDLSS Equivalence Critical Values}}
%We will only show that $c_{\alpha,\text{CvM}} = c_{\alpha, \text{CQ}} + O_{\mathbb{P}}(d^{-1/2})$. The remaining results follow similarly. From Theorem~\ref{Theorem: HDLSS Equivalence}, we know that 
%\begin{align*}
%2\pi\sqrt{3d} \overline{\sigma}_d^2(U_{\text{CvM}}^{\varpi_1},\ldots,U_{\text{CvM}}^{\varpi_{N!}}) = d^{-1/2}(U_{\text{CQ}}^{\varpi_1},\ldots,U_{\text{CQ}}^{\varpi_{N!}}) + O_{\mathbb{P}}(d^{-1/2})
%\end{align*}
%where $\varpi_i$ is an element of $\mathcal{S}_N$ for $i=1,\ldots,N!$. For simplicity, let us write $2\pi\sqrt{3d} \overline{\sigma}_d^2 U_{\text{CvM}}^{\varpi_i} = U_{\text{CvM},s}^{\varpi_i}$ and $d^{-1/2}U_{\text{CQ}}^{\varpi_i} = U_{\text{CQ},s}^{\varpi_i}$. Then $c_{\alpha,\text{CvM}}$ and $c_{\alpha,\text{CQ}}$ are the $\ceil{N!(1-\alpha)}$th order statistic of $\{ U_{\text{CvM},s}^{\varpi_1},\ldots, U_{\text{CvM},s}^{\varpi_{N!}}\}$ and $\{ U_{\text{CQ},s}^{\varpi_1},\ldots, U_{\text{CQ},s}^{\varpi_{N!}}\}$, respectively. It is well-known that the order statistic is a Lipschitz function. More specifically, using Pigeonhole principle, it can be seen that 
%\begin{align*}
%|c_{\alpha,\text{CvM}} - c_{\alpha,\text{CQ}}| \leq \Bigg\{ \sum_{i=1}^{N!} (U_{\text{CvM},s}^{\varpi_i} - U_{\text{CQ},s}^{\varpi_i})^2 \Bigg\}^{1/2} = O_{\mathbb{P}}(d^{-1/2}).
%\end{align*}
%Hence the result follows. 

\vskip 2em

\subsection{Proof of Theorem~\ref{Corollary: HDLSS Power equivalence}}
Under the stated assumptions, Theorem 2.1 of \cite{chakraborty2017tests} is satisfied. Hence the results for the CQ and WMW tests follow. For the rest of the tests, we apply Slutsky's theorem combined with Theorem~\ref{Theorem: HDLSS Equivalence} to obtain the results. This completes the proof.

\vskip 2em

\subsection{Proof of Lemma~\ref{Lemma: Angle distance is of negative type}}

For given $w \in \mathbb{R}^{d}$, it is seen that
\begin{align} \label{Eq: Triangle Inequality}
& \int_{\mathbb{S}^{d-1}} \Big| \ind(\beta^\top z \leq \beta^\top w) - \ind(\beta^\top z^\prime \leq \beta^\top w) \Big| d \lambda(\beta) \\[.5em] \nonumber
=~ &  \int_{\mathbb{S}^{d-1}} \ind(\beta^\top z \leq \beta^\top w < \beta^\top z^\prime) + \ind(\beta^\top z^\prime \leq \beta^\top w < \beta^\top z) d \lambda(\beta) \\[.5em]  \nonumber
=~ &  \frac{1}{2} - \frac{1}{2\pi} \text{arccos} \Bigg\{ \frac{(z-w)^\top (w -z^\prime)}{\|z-w\| \| w- z^\prime\|} \Bigg\} + \frac{1}{2} - \frac{1}{2\pi} \text{arccos} \Bigg\{ \frac{(z^\prime-w)^\top (w - z)}{\|z^\prime - w\| \|w - z\|} \Bigg\} \\[.5em] \nonumber
=~ &  1 - \frac{1}{\pi}\text{arccos} \Bigg\{ \frac{(z-w)^\top (w - z^\prime)}{\|z - w\| \|w - z^\prime\|} \Bigg\} \\[.5em]  \nonumber
=~ &  \frac{1}{\pi} \left(  \pi - \text{arccos} \Bigg\{ \frac{(z-w)^\top (w - z^\prime)}{\|z - w\| \|w - z^\prime\|} \Bigg\} \right) \\[.5em]  \nonumber
\overset{(i)}{=}~ &  \frac{1}{\pi} \text{arccos} \Bigg\{ \frac{(z-w)^\top (z^\prime - w)}{\|z - w\| \|z^\prime - w\|} \Bigg\} :=~  \rho_{Angle}(z,z^\prime; w) , \nonumber 
\end{align}
where $(i)$ is due to $\text{arccos}(x) + \text{arccos}(-x) = \pi$. Then $\rho_{Angle}(z,z^\prime)$ is the expected value of $\rho_{Angle}(z,z^\prime;Z^\ast)$ over $Z^\ast \sim (1/2) P_X + (1/2) P_Y$, i.e.
\begin{align*}
\rho_{Angle}(z,z^\prime) & = \mE \left[ \rho_{Angle}(z,z^\prime;Z^\ast) \right] \\[.5em]
& = \frac{1}{\pi}\mE \left[ \text{arccos}\Bigg\{ \frac{(z-Z^\ast)^\top (z^\prime - Z^\ast)}{\|z-Z^\ast\| \|z^\prime - Z^\ast\|} \Bigg\}\right].
\end{align*}
Now, if $z = z^\prime$, it is trivial to see $\rho_{Angle}(z,z^\prime) = 0$. In addition, if $\rho_{Angle}(z,z^\prime) = 0$, then we have $z = z^\prime$. In order to show the second direction, note that $\text{arccos}(x)$ is positive and monotone decreasing over $x \in [-1,1]$ and so $\rho_{Angle}(z,z^\prime) = 0$ implies that 
\begin{align*}
\frac{(z-Z^\ast)^\top (z^\prime - Z^\ast)}{\|z-Z^\ast\| \|z^\prime - Z^\ast\|} = 1,
\end{align*}
almost surely with respect to $(1/2) P_X + (1/2) P_Y$. By Cauchy-Schwarz inequality, the inner product becomes one if and only if $(z-Z^\ast)$ or $(z^\prime - Z^\ast)$ is a multiple of the other. This is only possible when $z - Z^\ast = z^\prime - Z^\ast$ almost surely, which implies $z = z^\prime$. The symmetry property follows easily by the definition of $\rho_{{Angle}}$. In addition, from triangle inequality, we have
\begin{align*}
&\int_{\mathbb{S}^{d-1}} \Big| \ind(\beta^\top z \leq \beta^\top w) - \ind(\beta^\top z^\prime \leq \beta^\top w) \Big| d \lambda(\beta) \\[.5em] \leq ~ &\int_{\mathbb{S}^{d-1}} \Big| \ind(\beta^\top z \leq \beta^\top w) - \ind(\beta^\top z^{\prime \prime} \leq \beta^\top w) \Big| d \lambda(\beta) \\[.5em]
& + \int_{\mathbb{S}^{d-1}} \Big| \ind(\beta^\top z^{\prime\prime} \leq \beta^\top w) - \ind(\beta^\top z^\prime \leq \beta^\top w) \Big| d \lambda(\beta),
\end{align*}
and therefore by the equality in (\ref{Eq: Triangle Inequality}), we can establish
\begin{align*}
\rho_{Angle}(z,z^\prime; w) ~  \leq ~ \rho_{Angle}(z,z^{\prime\prime}; w)  + \rho_{Angle}(z^{\prime},z^{\prime\prime}; w).
\end{align*}
Now, by taking the expectation over $Z^\ast$, we conclude that 
\begin{align*}
\rho_{Angle}(z,z^\prime) ~ \leq ~ \rho_{Angle}(z,z^{\prime \prime}) + \rho_{Angle}(z^\prime,z^{\prime \prime}).
\end{align*}

Next, we will show that for $\forall n \geq 2$, $z_1,\ldots, z_n \in S$, and $\alpha_1,\ldots,\alpha_n \in \mathbb{R}$, with $\sum_{i=1}^n \alpha_i = 0$,
\begin{align*}
\sum_{i=1}^n \sum_{j=1}^n \alpha_i \alpha_j \rho_{Angle}(z_i,z_j) \leq 0.
\end{align*} 
The result follows from Section 6 of \cite{bogomolny2007distance} who showed that for each fixed $z^\ast$, 
\begin{align} \label{Eq: Result by Bogomolny}
\sum_{i=1}^n \sum_{j=1}^n \alpha_i \alpha_j \rho_{Angle}(z_i,z_j;z^\ast) \leq 0,
\end{align} 
for any $\alpha_1,\ldots,\alpha_n \in \mathbb{R}$, with $\sum_{i=1}^n \alpha_i = 0$. Therefore, by taking the expected value over $z^\ast$ in (\ref{Eq: Result by Bogomolny}), we conclude that $\rho_{Angle}$ is of negative-type. 

\vskip 1em

Regarding Remark~\ref{Remark: Generalization of Angular distance}, note that
\begin{align*}
& \int_{\mathbb{R}^d} \rho_{Angle}(z,z^\prime; t) dt \\[.5em] 
=~ &  \int_{\mathbb{S}^{d-1}} \int_{\mathbb{R}} I (\beta^\top z \leq \beta^\top t < \beta^\top z^\prime) + \ind(\beta^\top z^\prime \leq \beta^\top t < \beta^\top z) d \beta^\top t  d \lambda(\beta) \\[.5em] 
\overset{(i)}{=} ~ & \int_{\mathbb{S}^{d-1}} | \beta^\top \left( z -  z^\prime\right) |  d \lambda(\beta) \\[.5em] 
\overset{(ii)}{=} ~ &  \gamma_d \| z - z^\prime\|,
\end{align*}
where $(i)$ and $(ii)$ are due to Lemma 2.1 and Lemma 2.3 of \cite{baringhaus2004new} and 
\begin{align*}
\gamma_d =  \frac{\sqrt{\pi} (d-1) \Gamma \left( (d-2)/2  \right) }{2 \Gamma(d/2)}.
\end{align*}
Therefore, the generalized angular distance with Lebesgue measure corresponds to the Euclidean distance.

\vskip 2em

\subsection{Proof of Proposition~\ref{Proposition: Identity between the generalized energy distance and CvM-distance}} 
From the definition of $\rho_{{Angle}}$, it is seen that
\begin{align*}
& 2 \mE \left[ \rho_{Angle} (X_1,Y_1) \right] - \mE \left[ \rho_{Angle} (X_1,X_2) \right] - \mE \left[\rho_{Angle}(Y_1,Y_2)\right] \\[.5em] 
= ~ & \frac{1}{\pi} \mE\left[ \mathsf{Ang} (X_1 - X_2, Y_1 - X_2) \right] + \frac{1}{\pi} \mE\left[ \mathsf{Ang} (X_1 - Y_2, Y_1 - Y_2) \right] \\[.5em]
-&  \frac{1}{2\pi}  \mE \left[ \mathsf{Ang} (X_1-X_3, X_2 - X_3) \right] - \frac{1}{2\pi} \mE \left[ \mathsf{Ang} (X_1-Y_1, X_2 - Y_1) \right]  \\[.5em]
 -& \frac{1}{2\pi}  \mE \left[ \mathsf{Ang} (Y_1-X_2, Y_2 - X_2) \right] - \frac{1}{2\pi} \mE \left[ \mathsf{Ang} (Y_1-Y_3, Y_2 - Y_3) \right].
\end{align*}
Then the result follows by Lemma~\ref{Lemma: Multivaraite CvM Expression}.

\vskip 2em

\subsection{Proof of Theorem~\ref{Theorem: Kendall's tau}}
Given $\alpha \in \mathbb{S}^{p-1}, \beta \in \mathbb{S}^{q-1}$, expand the square term to have
\begin{align*}
& \Big\{ 4 \mP\left(\alpha^\top (X_1 - X_2) < 0 , \beta^\top (Y_1 - Y_2) < 0 \right) - 1 \Big\}^2  \\[.5em]
= ~ & 16\mE \Big[ \ind(\alpha^\top (X_1 - X_2) < 0, \alpha^\top (X_3 - X_4) < 0) \\[.5em]
& ~~~~~~~~\times  \ind(\beta^\top (Y_1 - Y_2) < 0, \beta^\top (Y_3 - Y_4) < 0)  \Big] \\[.5em]
~~ -  &8  \mE \left[  \ind(\alpha^\top (X_1 - X_2) < 0 ) \times \ind(\beta^\top (Y_1 - Y_2) < 0) \right] + 1.
\end{align*}
By applying Lemma~\ref{Lemma: Integration over Unit Sphere (2 terms)}, the first term becomes
\begin{align*}
\mE \left[ \left(2  -  \frac{2}{\pi}  \mathsf{Ang} \left( X_1 - X_2, X_3 - X_4 \right)\right) \cdot  \left( 2 -  \frac{2}{\pi}\mathsf{Ang} \left( Y_1 - Y_2, Y_3 - Y_4 \right)\right) \right]
\end{align*}
and the remainder terms become $-1$, which yields the expression.

\vskip 1em

\subsection{Proof of Theorem~\ref{Theorem: BKR coefficient}}
Given $\alpha \in \mathbb{S}^{p-1}$ and $\beta \in \mathbb{S}^{q-1}$, 
\begin{align*}
& \int_{\mathbb{R}^2} \Big[ F_{\alpha^\top X, \beta^\top Y}(u,v) - F_{\alpha^\top X}(u) F_{\beta^\top Y}(v) \Big]^2 dF_{\alpha^\top X}(u) dF_{\beta^\top Y}(v) \\[.5em]
=~  & \mE \Big[ \ind(\alpha^\top (X_1 - X_3) \leq 0, \alpha^\top (X_2 - X_3) \leq 0) \\
& ~~~~~~~~~~~~~ \times  \ind(\beta^\top (Y_1 - Y_4) \leq 0, \beta^\top (Y_2 - Y_4) \leq 0) \Big] \\[.5em]
+ ~ & \mE \Big[ \ind(\alpha^\top (X_1 - X_5) \leq 0, \alpha^\top (X_2 - X_5) \leq 0) \\
&  ~~~~~~~~~~~~~  \times   \ind(\beta^\top (Y_3 - Y_6) \leq 0, \beta^\top (Y_4 - Y_6) \leq 0) \Big] \\[.5em]
-2 & \mE \Big[ \ind(\alpha^\top (X_1 - X_4) \leq 0, \alpha^\top (X_2 - X_4) \leq 0) \\
&  ~~~~~~~~~~~~~  \times  \ind(\beta^\top (Y_1 - Y_5) \leq 0, \beta^\top (Y_3 - Y_5) \leq 0) \Big]. 
\end{align*}
Then apply Lemma~\ref{Lemma: Integration over Unit Sphere (2 terms)} to obtain the expression.

\vskip 1em

\subsection{Proof of Lemma~\ref{Lemma: Extension of Escanciano (2006) with three arguments}}
To prove the results, we apply the same argument used in Section~\ref{Section: Proof of Lemma: Integration over unit sphere}. Let $\mathcal{Z}$ have a multivariate normal distribution with zero mean vector and identity covariance matrix. Then as in Section~\ref{Section: Proof of Lemma: Integration over unit sphere},
\begin{equation}
\begin{aligned} \label{Eq: Identity of beta and Z in general}
& \int_{\mathbb{S}^{d-1}} \prod_{i=1}^3 \ind(\beta^\top U_i \leq 0) d \lambda(\beta)  = \mE_{\mathcal{Z}} \bigg[ \prod_{i=1}^3 \ind(\mathcal{Z}^\top U_i \leq 0) \bigg].
\end{aligned}
\end{equation}
Since $(\mathcal{Z}^\top U_1,\mathcal{Z}^\top U_2,\mathcal{Z}^\top U_3)^\top$ has a multivariate normal distribution with zero mean vector and correlation matrix $[\varrho_{ij}]_{3 \times 3}$ with $\varrho_{ij} = U_i^\top U_j /\{ \|U_i\| \|U_j\|\}$, the right-hand side of (\ref{Eq: Identity of beta and Z in general}) can be computed based on orthant probabilities for normal distributions \citep[e.g.][]{childs1967reduction,xu2013comparative}. This completes the proof. 

%\cite{childs1967reduction} also provided expressions for higher order integrations. We will provide this result in Section~\ref{Section: Generalization of Escanciano Lemma}. 

\vskip 2em

\subsection{Proof of Theorem~\ref{Theorem: Multivariate tau star}}

From \cite{bergsma2014consistent}, the univariate $\tau^\ast$ can be written as
\begin{align*}
\tau^\ast = ~ & 4 \mP \left( X_1 \vee X_2 < X_3 \wedge X_4, ~ Y_1 \vee Y_2 < Y_3 \wedge Y_4 \right) \\[.5em]
+ ~ & 4 \mP \left( X_1\vee X_2 < X_3 \wedge X_4, ~Y_1\wedge Y_2 > Y_3 \vee Y_4 \right) \\[.5em]
- ~ & 8 \mP \left( X_1 \vee X_2 < X_3 \wedge X_4,~ Y_1 \vee Y_3 < Y_2 \wedge Y_4 \right).
\end{align*}
Notice that 
\begin{align*}
& \ind(X_1 \vee X_2 < X_3 \wedge X_4) \\[.5em]
=~ & \ind(X_1 < X_2 < X_3 < X_4) + \ind(X_2 < X_1 < X_3 < X_4) \\[.5em]
+~ &  \ind(X_1 < X_2 < X_4 < X_3) + \ind(X_2 < X_1 < X_4 < X_3) \\[.5em]
=~ &  \ind(X_1 < X_2) \ind(X_2 < X_3) \ind(X_3<X_4) +  \ind(X_2 < X_1) \ind(X_1 < X_3) \ind(X_3<X_4) \\[.5em] 
+~ &   \ind(X_1 < X_2) \ind(X_2 < X_4) \ind(X_4<X_3) +  \ind(X_2 < X_1) \ind(X_1 < X_4) \ind(X_4<X_3). 
\end{align*}
Similarly, we have 
\begin{align*}
& \ind(Y_1 \vee Y_2 < Y_3 \wedge Y_4)  \\[.5em]
=~ & \ind(Y_1 < Y_2) \ind(Y_2 < Y_3) \ind(Y_3 < Y_4) +  \ind(Y_2 < Y_1) \ind(Y_1 < Y_3) \ind(Y_3<Y_4) \\[.5em] 
+~ &   \ind(Y_1 < Y_2) \ind(Y_2 < Y_4) \ind(Y_4 < Y_3) +  \ind(Y_2 < Y_1) \ind(Y_1 < Y_4) \ind(Y_4 < Y_3). 
\end{align*}
Therefore, the product $I(X_1 \vee X_2 < X_3 \wedge X_4) \ind(Y_1 \vee Y_2 < Y_3 \wedge Y_4)$ can be expressed as the linear combination of 
\begin{align*}
\ind(X_{i_1} < X_{i_2}) \ind(X_{i_2} < X_{i_3}) \ind(X_{i_3} < X_{i_4})  \ind(Y_{j_1} < Y_{j_2}) \ind(Y_{j_2} < Y_{j_3}) \ind(Y_{j_3} < Y_{j_4}). 
\end{align*}
Using Lemma~\ref{Lemma: Extension of Escanciano (2006) with three arguments}, 
\begin{align*}
& \int_{\mathbb{S}^{p-1}} \ind(\alpha^\top X_{i_1} < \alpha^\top X_{i_2}) \ind( \alpha^\top X_{i_2} < \alpha^\top  X_{i_3}) \ind(\alpha^\top X_{i_3} < \alpha^\top X_{i_4}) d\lambda(\alpha) \\[.5em]
%& = \int_{\mathbb{S}^{p-1}} \ind(\alpha^\top (X_{i_1}  - X_{i_2}) < 0) \ind( \alpha^\top (X_{i_2} - X_{i_3}) < 0 ) \ind(\alpha^\top (X_{i_3} - X_{i_4}) < 0) d\lambda(\alpha)  \\[.5em]
& = \frac{1}{2} -\frac{1}{4\pi} \left[ \mathsf{Ang} \left( U_1, U_2 \right) +  \mathsf{Ang} \left( U_1, U_3 \right) + \mathsf{Ang} \left( U_2, U_3 \right)  \right],
\end{align*}
where $U_1 = X_{i_1}  - X_{i_2}$, $U_2 = X_{i_2} - X_{i_3}$ and $U_3 = X_{i_3} - X_{i_4}$. 

\vskip 1em

Similarly,
\begin{align*}
&\int_{\mathbb{S}^{q-1}} \ind(\beta^\top Y_{j_1} < \beta^\top Y_{j_2}) \ind( \beta^\top Y_{j_2} < \beta^\top  Y_{j_3}) \ind(\beta^\top Y_{j_3} < \beta^\top Y_{j_4}) d\lambda(\beta) \\[.5em]
& =  \frac{1}{2} -\frac{1}{4 \pi} \left[ \mathsf{Ang}\left( V_1, V_2 \right) + \mathsf{Ang}\left( V_1, V_3 \right) + \mathsf{Ang}\left( V_2, V_3 \right) \right] ,
\end{align*}
where $V_1 = Y_{j_1}  - Y_{j_2}$, $V_2 = Y_{j_2} - Y_{j_3}$ and $V_3 = Y_{j_3} - Y_{j_4}$. 

\vskip .5em 

As a result, we have
\begin{align*}
& \int_{\mathbb{S}^{p-1}} \int_{\mathbb{S}^{q-1}} \mP(\alpha^\top X_1 \vee \alpha^\top X_2 < \alpha^\top X_3 \wedge \alpha^\top X_4, ~ \\
& ~~~~~~~~~~~~~~~~~~~ ~ \beta^\top Y_1 \vee \beta^\top Y_2 < \beta^\top Y_3 \wedge \beta^\top  Y_4) d\lambda(\alpha) d\lambda(\beta)  \\[.5em]
= ~ &  \mE \left[ h_p(X_1,X_2,X_3,X_4) h_q(Y_1,Y_2,Y_3,Y_4)\right].
\end{align*}

\noindent Applying the same argument to the rest, we can obtain the explicit expression for $\tau^\ast_{p,q}$ as in Theorem~\ref{Theorem: Multivariate tau star}.

\vskip 2em

\subsection{Proof of Theorem~\ref{Theorem: Two-Sample Degenerate Kernel}} \label{Section: Proof of Permutation Consistency}

Let us write
\begin{align*}
U^\ast_{m,n}(Z_{m,n}) & := U^\ast_{m,n}(Z_1,\ldots,Z_{N})  \\[.5em]
& = N \{  U_{m,n}(Z_1,\ldots,Z_{N}) - \mE\left[ U_{m,n}(Z_1,\ldots,Z_{N}) \right] \}
\end{align*}
and denote $U^\ast_{m,n}(Z_{\varpi(1)}, \ldots, Z_{\varpi(N)})$ by ${U}^{\ast}_{m,n}(Z_{\varpi})$. Our goal is to show that for two independent random permutations $\varpi, \varpi^\prime$,
\begin{align} \label{Eq: Goal of Permutation result}
\left( U^\ast_{m,n}(Z_{\varpi}), ~ U^\ast_{m,n}(Z_{\varpi^\prime})  \right) \convD  (T, T^\prime),
\end{align}
where $T,T^\prime$ are independent and identically distributed with the distribution function $R(t)$. Then the desired result follows by Lemma~\ref{Lemma: Hoeffding's condition}. The proof consists of several pieces and closely follows the proof of the limiting distribution of a two-sample degenerate $U$-statistic in Chapter 3 of \cite{bhat1995theory}. 

\vskip 1em

We start with the projection of the two-sample $U$-statistic via Hoffding's decomposition. Consider the projection of the two-sample degenerate $U$-statistic based on $Z_{m,n}$:
\begin{align*}
\widehat{U}_{m,n}(Z_{m,n}) = ~ & \frac{r(r-1)}{m(m-1)}\sum_{1 \leq i_1 < i_2 \leq m} g^\ast_{2,0}(Z_{i_1},Z_{i_2}) + \frac{r(r-1)}{n(n-1)} \sum_{1 \leq j_1 < j_2 \leq n} g^\ast_{0,2}(Z_{j_1+m},Z_{j_2+m})  \\
& + \frac{r^2}{mn} \sum_{i=1}^m \sum_{j=1}^n g^\ast_{1,1}(Z_i,Z_{j+m}).
\end{align*}

\noindent Then it can be seen that
\begin{align*}
\mE [ (U_{m,n}(Z_{m,n}) - \widehat{U}_{m,n}(Z_{m,n} )] = 0 ~ \text{and} ~ \mV[ U_{m,n}(Z_{m,n}) - \widehat{U}_{m,n}(Z_{m,n} )] = O(N^{-3}),
\end{align*}
which implies
\begin{align} \label{Eq: Asymptotic Equivalence between U- and Projected U-statistic}
N(U_{m,n}(Z_{m,n})- \theta )= N( \widehat{U}_{m,n}(Z_{m,n}) - \theta )+ o_{\mathbb{P}}(1).
\end{align}

Under the finite second moment of the kernel $g$, we may have the decompositions
\begin{align*}
& g^\ast_{2,0}(x,y) = \sum_{i=1}^{\infty} \lambda_i \phi_i(x) \phi_i(y), \\[.5em]
& g^\ast_{0,2}(x,y) = \sum_{i=1}^{\infty} \gamma_i \psi_i(x) \psi_i(y), \\[.5em]
& g^\ast_{1,1}(x,y) = \sum_{i=1}^{\infty} \alpha_i \phi^\ast_i(x) \psi^\ast_i(y),
\end{align*}
where $\{\phi_i (\cdot )\}$, $\{ \psi_i(\cdot) \}$, $\{\phi^\ast(\cdot), \psi^\ast(\cdot)\}$ are orthonormal eigenfunctions and the corresponding eigenvalues $\{\lambda_i\}, \{\gamma_i\}, \{\alpha_i\}$, associated with $g^\ast_{2,0}, g^\ast_{0,2}$ and $g^\ast_{1,1}$, respectively \citep[see e.g.][for details]{bhat1995theory}. From the given conditions of the theorem, the eigenvalues and the eigenfunctions are related as follows:
\begin{align*}
& \phi_i(z) = \psi_i(z) = \phi_i^\ast(z) = \psi_i^\ast (z), \\[.5em]
& \gamma_i = \lambda_i \quad \text{and} \quad \alpha_i = \frac{1-r}{r} \lambda_i.  
\end{align*}
Therefore, 
\begin{align*}
N \widehat{U}_{m,n}(Z_{m,n}) = ~ & \widehat{a}_1 \left[ \frac{1}{m} \sum_{1 \leq i_1 \neq i_2 \leq m} \sum_{i=1}^\infty \lambda_i \phi_i (Z_{i_1}) \phi_i (Z_{i_2}) \right] \\
+ ~ &  \widehat{a}_2 \left[ \frac{1}{n} \sum_{1 \leq j_1 \neq j_2 \leq n} \sum_{j=1}^\infty \lambda_j \phi_j (Z_{j_1+m}) \phi_j (Z_{j_2+m}) \right] \\[.5em]
+ ~ & \widehat{a}_3 \left[ \frac{1}{\sqrt{mn}} \sum_{i_1 = 1}^{m} \sum_{j_1=1}^{n} \sum_{k=1}^{\infty} \lambda_k \phi_k(Z_{i_1}) \phi_k(Z_{j_1+m}) \right] \\[.8em]
= ~ &  \widehat{a}_1 T_m + \widehat{a}_2 T_n^\prime + \widehat{a}_3 T_{mn}^{\prime \prime},
\end{align*}
where 
\begin{align*}
& \widehat{a}_1 = \frac{r(r-1)}{2} \frac{N}{m-1}, \quad \widehat{a}_2 = \frac{r(r-1)}{2} \frac{N}{n-1} \quad \text{and} \quad \widehat{a}_3 = -r(r-1) \frac{N}{\sqrt{mn}}.
\end{align*}

Denote the centered and scaled projection of the $U$-statistic by
\begin{align*}
\widetilde{U}_{m,n} := N (\widehat{U}_{m,n}({Z}_{\varpi}) - \theta) \quad \text{and} \quad \widetilde{U}^\prime_{m,n} := N( \widehat{U}_{m,n}( {Z}_{\varpi^\prime})-\theta).
\end{align*}
Then due to (\ref{Eq: Asymptotic Equivalence between U- and Projected U-statistic}), 
\begin{align*}
\left( U^\ast_{m,n}(Z_{\varpi}), ~ U^\ast_{m,n}(Z_{\varpi^\prime})  \right) = \left( \widetilde{U}_{m,n}(Z_\varpi), ~\widetilde{U}_{m,n}^\prime(Z_{\varpi^\prime}) \right)  + o_{\mathbb{P}}(1).
\end{align*}
Therefore it suffices to show
\begin{align*}
\left( \widetilde{U}_{m,n}, ~\widetilde{U}_{m,n}^\prime \right) \convD (T,T^\prime)
\end{align*}
to complete the main proof. Having this goal in mind, we start with a truncation of the degenerate $U$-statistic.

\vskip 1em

\noindent \textbf{$\bullet$ Truncation of the $U$-statistics.}

\vskip 1em

\noindent Now, define a truncated version of $N(\widehat{U}_{m,n}(Z_{m,n}) - \theta)$ by
\begin{equation}
\begin{aligned} \label{Eq: Truncated Projection}
N(\widehat{U}_{m,n,K}(Z_{m,n}) - \theta)  = ~& \widehat{a}_1 \left[ \frac{1}{m} \sum_{1 \leq i_1 \neq i_2 \leq m} \sum_{i=1}^K \lambda_i \phi_i (Z_{i_1}) \phi_i (Z_{i_2}) \right]  \\
+~ &  \widehat{a}_2 \left[ \frac{1}{n} \sum_{1 \leq j_1 \neq j_2 \leq n} \sum_{j=1}^K \lambda_j \phi_j (Z_{j_1+m}) \phi_j (Z_{j_2+m}) \right] \\[.5em]
+~  & \widehat{a}_3 \left[ \frac{1}{\sqrt{mn}} \sum_{i_1 = 1}^{m} \sum_{j_1=1}^{n} \sum_{k=1}^{K} \lambda_k \phi_k(Z_{i_1}) \phi_k(Z_{j_1+m}) \right] \\[.8em]
=~  &  \widehat{a}_1 T_{mK} + \widehat{a}_2 T_{nK}^\prime + \widehat{a}_3 T_{mnK}^{\prime \prime}.
\end{aligned}
\end{equation}

\noindent Write
\begin{align*}
& \widehat{a}_1 T_{mK} + \widehat{a}_2 T_{nK}^\prime + \widehat{a}_3 T_{mnK}^{\prime \prime}  \\[.5em]
= ~ & \widehat{a}_1 \left[ \sum_{k=1}^K \lambda_k \left( W_{km}^2 -  V_{km} \right) \right] +\widehat{a}_2 \left[ \sum_{k=1}^K \lambda_k \left(  W_{kn}^{\prime 2} -  V_{kn}^\prime  \right) \right] +  \widehat{a}_3\left[ \sum_{k=1}^K \lambda_k W_{km} W_{kn}^\prime \right] \\[.5em]
= ~ & \frac{r(r-1)}{2} \Bigg\{ \sum_{k=1}^K \lambda_k \left( \sqrt{\frac{N}{m}} W_{km} - \sqrt{\frac{N}{n}}W_{kn}^\prime \right)^2 - \sum_{k=1}^K \lambda_k \left( \frac{N}{m} V_{km} + \frac{N}{n} V_{kn}^\prime \right) \Bigg\},
\end{align*} 
where
\begin{align*}
&W_{km} = \frac{1}{\sqrt{m}} \sum_{i_1=1}^m \phi_k (Z_{i_1}), \quad W_{kn}^\prime = \frac{1}{\sqrt{n}} \sum_{j_1=1}^n \phi_k (Z_{j_1+m}), \\[.5em]
&V_{km} = \frac{1}{m} \sum_{i_1=1}^m \phi_k^2 (Z_{i_1}), \quad V_{kn}^\prime = \frac{1}{n} \sum_{j_1=1}^n \phi_k^2 (Z_{j_1+m}),
\end{align*}
for $k=1,\ldots, K$. 

\vskip 1em

By strong law of large numbers, 
\begin{align*}
V_{mn}^{\ast \top} := (V_{1m}, \ldots, V_{Km}, V_{1n}^\prime, \ldots, V_{Kn}^\prime)^\top \overset{a.s.}{\longrightarrow} V^{\ast \top} = (V_{1}, \ldots, V_{K}, V_{1}^\prime, \ldots, V_{K}^\prime)^\top
\end{align*}		
and by the assumption that $m/N \rightarrow \vartheta_X$, $n/N \rightarrow \vartheta_Y$,
\begin{align*}
& N(\widehat{U}_{m,n, K} - \theta)  \\
=~  & \frac{r(r-1)}{2} \Bigg\{ \sum_{k=1}^K \lambda_k \left( \sqrt{\frac{N}{m}} W_{km} - \frac{r(r-1)}{2}  \sqrt{\frac{N}{n}}W_{kn}^\prime \right)^2 - \frac{1}{\vartheta_X \vartheta_Y}\sum_{k=1}^K \lambda_k \Bigg\}+ o_{\mathbb{P}}(1) \\[.5em]
=~ & \frac{r(r-1)}{2} \Bigg\{ N\sum_{k=1}^K \lambda_k \left( \frac{1}{m} \sum_{i=1}^m \phi_k (Z_i) - \frac{1}{n} \sum_{j=1}^n \phi_k (Z_{j+m}) \right)^2 - \frac{1}{\vartheta_X \vartheta_Y}\sum_{k=1}^K \lambda_k \Bigg\} + o_{\mathbb{P}}(1) \\[.5em]
=~ &  \frac{r(r-1)}{2} \Bigg\{  N\sum_{k=1}^K \lambda_k \left( \sum_{i=1}^N \epsilon_i \phi_k(Z_i) \right)^2 - \frac{1}{\vartheta_X \vartheta_Y}\sum_{k=1}^K \lambda_k \Bigg\} + o_{\mathbb{P}}(1)
\end{align*}
where 
\begin{align*}
(\epsilon_1, \ldots, \epsilon_m, \epsilon_{m+1}, \ldots, \epsilon_{m+n}) = (\underbrace{m^{-1}, \ldots, m^{-1}}_{\text{$m$ terms}}, \underbrace{-n^{-1}, \ldots, -n^{-1}}_{\text{$n$ terms}}).
\end{align*}

\vskip 1em
\noindent \textbf{$\bullet$ Proving independence of the truncated $U$-statistics.}
\vskip 1em

\noindent Consider the truncated permutation statistics
\begin{align*}
\widetilde{U}_{m,n,K} & := N(\widehat{U}_{m,n,K}(Z_\varpi) - \theta)  \\[.5em]
& = \frac{r(r-1)}{2} \Bigg\{ N \sum_{k=1}^K \lambda_k \left( \sum_{i=1}^N \epsilon_{\varpi(i)} \phi_k(Z_{i}) \right)^2 - \frac{1}{\vartheta_X \vartheta_Y}\sum_{k=1}^K \lambda_k \Bigg\} + o_{\mathbb{P}}(1) \\[.5em]
\widetilde{U}_{m,n,K}^\prime &  := N(\widehat{U}_{m,n,K}(Z_{\varpi^\prime}) - \theta)  \\[.5em]
& = \frac{r(r-1)}{2} \Bigg\{ N \sum_{k=1}^K \lambda_k \left( \sum_{i=1}^N \epsilon_{\varpi^\prime(i)} \phi_k(Z_{i}) \right)^2 - \frac{1}{\vartheta_X \vartheta_Y}\sum_{k=1}^K \lambda_k  \Bigg\} + o_{\mathbb{P}}(1). 
\end{align*}
Note that $\epsilon_{\varpi(i)}$ and $\epsilon_{\varpi^\prime(i)}$ are independent random variables by the assumption having either $1/m$ or $-1/n$ with $m/N$ and $n/N$ probabilities; hence
\begin{align*}
\text{Cov}\left(  \epsilon_{ \varpi (i)} \phi_k(Z_{i}),  \epsilon_{\varpi^\prime (i)}  \phi_k(Z_{i}) \right) = \mE \left[ \epsilon_{\varpi(i)}  \right] \mE \left[ \epsilon_{\varpi^\prime(i)}\right] \mE \left[ \phi_k^2(Z_{i}) \right] =0.
\end{align*}
By the Cram{\'e}r-Wold device and the Lindeberg condition, we see that
\begin{align*}
& \sqrt{N} \left( \sum_{i=1}^N \epsilon_{\varpi(i)} \phi_1(Z_{i}), \ldots, \sum_{i=1}^N \epsilon_{\varpi(i)} \phi_K(Z_{i}), \sum_{i=1}^N \epsilon_{\varpi^\prime(i)} \phi_1(Z_{i}) , \ldots, \sum_{i=1}^N \epsilon_{\varpi^\prime(i)} \phi_K(Z_{i}) \right)^\top \\[.5em]
&~~~~~~~~~~~~~  \convD N(0, {\vartheta_X}^{-1} {\vartheta_Y}^{-1} I_{2K}).
\end{align*}
Thus the components of the vector are asymptotically independent to each other. Then apply the continuous mapping theorem together with Slutsky's theorem to have
\begin{align}  \label{Eq: Multivariate Normal Approximation}
(\widetilde{U}_{m,n,K}, \widetilde{U}_{m,n,K}^\prime) \convD (T_K, T_K^\prime)
\end{align}
where $T_K$ and $T_K^\prime$ are independent and have the same distribution as
\begin{align*}
\frac{r(r-1)}{2\vartheta_X \vartheta_Y} \sum_{k=1}^{K} \lambda_k (\xi_k^2 - 1),
\end{align*}
where $\xi_k \overset{i.i.d.}{\sim} N(0,1)$. 

\vskip 1em

\noindent \textbf{$\bullet$ Bounding the difference between characteristic functions.}

\vskip 1em

\noindent We will use the characteristic functions to show 
\begin{align*}
\left( \widetilde{U}_{m,n}, ~\widetilde{U}_{m,n}^\prime  \right) \convD  (T, T^\prime).
\end{align*}
More specifically, we will show that for any $x,y \in \mathbb{R}$ and any $\epsilon > 0$ and sufficiently large $N$,
\begin{align*}
& \Big| \mE \left[ e^{i(x \widetilde{U}_{m,n} + y \widetilde{U}_{m,n}^\prime)} \right] -  \mE \left[ e^{i(x T + y T^\prime)} \right] \Big|  \leq ~ (I) + (II)+ (III) ~ < \epsilon	
\end{align*}
where 
\begin{align*}
& (I) = \Big| \mE \left[ e^{i(x \widetilde{U}_{m,n} + y \widetilde{U}_{m,n}^\prime)} \right] -  \mE \left[ e^{i(x \widetilde{U}_{m,n,K} + y \widetilde{U}_{m,n,K}^\prime)} \right] \Big|, \\[.5em]
& (II) = \Big| \mE \left[ e^{i(x \widetilde{U}_{m,n,K} + y \widetilde{U}_{m,n,K}^\prime)} \right] -  \mE \left[ e^{i(x T_{K} + y T_{K}^\prime)} \right] \Big|, \\[.5em]
& (III) = \Big| \mE \left[ e^{i(x T_{K} + y T_{K}^\prime)} \right] -  \mE \left[ e^{i(x T + y T^\prime)} \right] \Big|.
\end{align*}
We bound these terms in sequence.

\vskip 1em

\noindent \textbf{1. Bounding $(I)$.}
\vskip .5em

\noindent Based on $|e^{iz}| = 1$ and $|e^{iz} - 1| \leq |z|$, we bound $(I)$ by
\begin{align*}
(I) ~=~ & \Big| \mE \left[ e^{i(x \widetilde{U}_{m,n} + y \widetilde{U}^\prime_{m,n} )} \right] -  \mE \left[ e^{i(x \widetilde{U}_{m,n,K} + y \widetilde{U}^\prime_{m,n,K} ))} \right] \Big|  \\[.5em]
~ \leq ~ & ~ |x| \left[ \mE \left( \widetilde{U}_{m,n,K} - \widetilde{U}_{m,n} \right)^2 \right]^{1/2} +  |y| \left[ \mE \left( \widetilde{U}_{m,n,K}^\prime - \widetilde{U}_{m,n}^\prime  \right)^2 \right]^{1/2} \\[.5em]
~ \leq ~ & ~ \left( |x| + |y| \right) \Bigg\{  \frac{r(r-1)}{2\widehat{\vartheta}_1} \left(2\sum_{k=K+1}^\infty \lambda_k^2\right)^{1/2}  +  \frac{r(r-1)}{2\widehat{\vartheta}_2} \left( 2\sum_{k=K+1}^\infty \lambda_k^2 \right)^{1/2} \\
& ~~~~~~~~~~~~~~ - \frac{r(r-1)}{\sqrt{\widehat{\vartheta}_1 \widehat{\vartheta}_2}}  \left( \sum_{k=K+1}^\infty \lambda_k^2 \right)^{1/2}  \Bigg\} \\[.5em]
~ = ~ & ~  \left( |x| + |y| \right) \frac{r(r-1)}{\sqrt{2}} \left( \frac{1}{\sqrt{\widehat{\vartheta}_1}} - \frac{1}{\sqrt{\widehat{\vartheta}_2}} \right)^2 \left(\sum_{k=K+1}^\infty \lambda_k^2 \right)^{1/2}  \\[.5em]
~ \leq ~ & ~ \left( |x| + |y| \right) \frac{r(r-1)}{\sqrt{2}\widehat{\vartheta}_1\widehat{\vartheta}_2} \left(\sum_{k=K+1}^\infty \lambda_k^2 \right)^{1/2} 
\end{align*}
where $\widehat{\vartheta}_1 = m/N$ and $\widehat{\vartheta}_2 = n/N$. 

\vskip 1em

Now, for fixed $x$ and $y$ and any given $\epsilon>0$, we choose $K$ large enough to bound 
\begin{align} \label{Eq: Choice of K, epsilon}
\left( |x| + |y| \right) \frac{r(r-1)}{\sqrt{2}\vartheta_X \vartheta_Y} \left(\sum_{k=K+1}^\infty \lambda_k^2 \right)^{1/2}  < \frac{\epsilon}{3}.
\end{align}
Since $\widehat{\vartheta}_1 \rightarrow \vartheta_X$ and $\widehat{\vartheta}_2 \rightarrow \vartheta_Y$ as $N \rightarrow \infty$, we have 
\begin{align*}
(I) ~ \leq & ~ \left( |x| + |y| \right) \frac{r(r-1)}{\sqrt{2}\widehat{\vartheta}_1\widehat{\vartheta}_2} \left(\sum_{k=K+1}^\infty \lambda_k^2 \right)^{1/2} < \frac{\epsilon}{3},
\end{align*}
for all sufficiently large $N$.

\vskip 1em

\noindent \textbf{2. Bounding $(II)$.} 
\vskip .5em

\noindent From the result established in (\ref{Eq: Multivariate Normal Approximation}), we have 
\begin{align*}
(II) ~=~ \Big| \mE \left[ e^{i(x \widetilde{U}_{m,n,K} + y \widetilde{U}^\prime_{m,n,K} ))} \right]  -  \mE \left[ e^{i(x T_K + y T_K^\prime)} \right] \Big| < \frac{\epsilon}{3}  \quad \text{for all sufficiently large $N$}.
\end{align*}

\vskip 1em

\noindent \textbf{3. Bounding $(III)$.} 
\vskip .5em

\noindent From Chapter 3 of \cite{bhat1995theory} with the conditions given on the kernel, the asymptotic distribution of a degenerate $U$-statistic converges to
\begin{equation}
\begin{aligned} \label{Eq: the limiting distribution of a degenerate U-statistic}
N \left( U_{m,n}  - \theta \right) ~ \convD ~ & \frac{r(r-1)}{2\vartheta_X} \sum_{k=1}^{\infty} \lambda_k(\xi_k^2 - 1) + \frac{r(r-1)}{2\vartheta_Y} \sum_{k=1}^{\infty} \lambda_k (\xi_k^{\prime 2}-1) \\[.5em]
& - \frac{r(r-1)}{\sqrt{\vartheta_X \vartheta_Y}} \sum_{k=1}^{\infty} \lambda_k \xi_k \xi_k^\prime
\end{aligned}
\end{equation}
where $\{\xi_k\}$ and $\{ \xi_k^\prime \}$ are independent standard normal random variables and $\{\lambda_k\}$ are eigenvalues associated with the kernel. Note that the right-side of (\ref{Eq: the limiting distribution of a degenerate U-statistic}) can be re-written as
\begin{align*}
\frac{r(r-1)}{2\vartheta_X \vartheta_Y} \sum_{k=1}^{\infty} \lambda_k \left[(\sqrt{\vartheta_Y}\xi_k - \sqrt{\vartheta_X} \xi_k^\prime )^2 - 1\right],
\end{align*}
where $\sqrt{\vartheta_Y}\xi_k - \sqrt{\vartheta_X} \xi_k^\prime \sim N(0,1)$. Therefore, $T, T^\prime$ are identically distributed as
\begin{align*}
\frac{r(r-1)}{2\vartheta_X \vartheta_Y} \sum_{k=1}^{\infty} \lambda_k (\xi_k^2 - 1). 
\end{align*}
Recall that $T_K, T_K^\prime$ have the same distribution as 
\begin{align*}
\frac{r(r-1)}{2\vartheta_X \vartheta_Y} \sum_{k=1}^{K} \lambda_k (\xi_k^2 - 1). 
\end{align*}
Consequently, 
\begin{align*}
\Big| \mE \left[ e^{i(x T_K + y T_K^\prime)} \right] -  \mE \left[ e^{i(x T + y T^\prime)} \right] \Big| & \leq |x| \left[ \mE \left( T_K - T \right)^2 \right]^{1/2} +  |y| \left[ \mE \left( T_K^\prime - T^\prime \right)^2 \right]^{1/2} \\[.5em]
& \leq \left( |x| + |y| \right) \frac{r(r-1)}{\sqrt{2}\vartheta_X \vartheta_Y} \left(\sum_{k=K+1}^{\infty} \lambda_k^2 \right)^{1/2} < \frac{\epsilon}{3},
\end{align*}		
with the same choice of $x,y,\epsilon,K$ in (\ref{Eq: Choice of K, epsilon}).

\vskip 1em 

\noindent \textbf{$\bullet$ Combining the bounds.}
\vskip 1em
\noindent From the previous results, we conclude that for any $x, y \in \mathbb{R}$ and any $\epsilon >0$ with sufficiently large $N$,
\begin{align*}
\Big| \mE \left[ e^{i(x \widetilde{U}_{m,n} + y \widetilde{U}_{m,n}^\prime)} \right] -  \mE \left[ e^{i(x T + y T^\prime)} \right] \Big|  < \epsilon,
\end{align*}
and therefore
\begin{align*}
\left( \widetilde{U}_{m,n}, ~\widetilde{U}_{m,n}^\prime  \right) \convD  (T, T^\prime).
\end{align*}
This completes the proof.

\vskip 2em

\section{Additional Results} \label{Section: Additional Results}

In this section, we provide details on Equation~(\ref{Eq: CQ assumption}), Remark~\ref{Remark: Extension} and Remark~\ref{Remark: Equivalence of PA to SP} in the main text.

\subsection{Verification of (\ref{Eq: CQ assumption}) in the main text}
First we state the distributional assumptions made in \cite{bai1996effect} and \cite{chen2010two}:
\begin{align} \label{Eq: BSCQ Model}
X = \Gamma_X V_X + \mu_X \quad \text{and} \quad Y = \Gamma_Y V_Y + \mu_Y, 
\end{align}
where $V_X$ and $V_Y$ are independent random vectors in $\mathbb{R}^u$ for some $u \geq d$ such that $\mE(V_X)=\mE(V_Y)= 0$ and $\mV(V_X)= \mV(V_Y) = I_u$, the $u \times u$ identity matrix. $\Gamma_X$ and $\Gamma_Y$ are non-random $d \times u$ matrices such that $\Sigma_X = \Gamma_X \Gamma_X^\top$ and $\Sigma_Y = \Gamma_Y \Gamma_Y^\top$ are positive definite and $\mu_X$ and $\mu_Y$ are non-random $d$-dimensional vectors. Write $V_X=(V_{X,1},\ldots,V_{X,m})$ and $V_Y=(V_{Y,1},\ldots,V_{Y,m})$. Assume that $\mE (V_{X,i}^4) = \mE (V_{Y,i}) = 3 + \varDelta < \infty$ for $i=1,\ldots,m$ where $\varDelta$ is the difference between the fourth moment of $V_{X,i}$ and $N(0,1)$. In addition assume that
\begin{align*}
\mE (V_{X,l_1}^{\alpha_1} V_{X,l_2}^{\alpha_2} \cdots V_{X,l_q}^{\alpha_q}) = \prod_{i=1}^q  \mE (V_{X,l_i}^{\alpha_i}) \quad \text{and} \quad
\mE (V_{Y,l_1}^{\alpha_1} V_{Y,l_2}^{\alpha_2} \cdots V_{Y,l_q}^{\alpha_q})  = \prod_{i=1}^q  \mE (V_{Y,l_i}^{\alpha_i}) 
\end{align*}
for a positive integer $q$ such that $\sum_{l=1}^q \alpha_l \leq 8$ and $l_1 \neq l_2 \neq \cdots \neq l_q$. 

\vskip 1em 

Our goal here is to show that $\mV(\|Z_1 - Z_2\|^2) = O(d)$ and $\mV\{(Z_1-Z_3)^\top (Z_2 - Z_3)\} = O(d)$ are implied by 
\begin{align*}
(\mu_X - \mu_Y)^\top (\Sigma_X + \Sigma_Y) (\mu_X - \mu_Y)  = O(d) \quad \text{and} \quad  \tr \{ (\Sigma_X + \Sigma_Y )^2\} = O(d).
\end{align*}
where $Z_1,Z_2,Z_3$ are independent and each $Z_i$ is identically distributed as either $X$ or $Y$ in (\ref{Eq: BSCQ Model}). First let us focus on $\mV(\|Z_1 - Z_2\|^2)$. Denote $\overline{Z}_1 = Z_1  - \mE(Z_1)$, $\overline{Z}_2 = Z_2  - \mE(Z_2)$ and $\delta_{12} = \mE(Z_1) - \mE(Z_2)$. Based on the basic inequality,
\begin{align*}
\mV \Big(\sum_{i=1}^k X_i \Big) \leq k \sum_{i=1}^k \mV(X_i) \quad \text{for any $k \geq 1$,}
\end{align*}
we have
\begin{align*}
\mV(\|Z_1 - Z_2\|^2) & = ~ \mV\{(\overline{Z}_1 - \overline{Z}_2)^\top(\overline{Z}_1 - \overline{Z}_2) + 2 \delta_{12}^\top (\overline{Z}_1 - \overline{Z}_2)\} \\[.5em]
& \leq ~2 \mV\{(\overline{Z}_1 - \overline{Z}_2)^\top(\overline{Z}_1 - \overline{Z}_2)\} + 8 \mV\{\delta_{12}^\top (\overline{Z}_1 - \overline{Z}_2)\} \\[.5em]
& \leq ~ 8 \mV (\overline{Z}_1^\top \overline{Z}_1) + 8 \mV (\overline{Z}_2^\top \overline{Z}_2) + 16 \mV (\overline{Z}_1^\top \overline{Z}_2) + 8 \delta_{12}^\top \mV (\overline{Z}_1 - \overline{Z}_2) \delta_{12}.
\end{align*}
Now using Proposition A.1 of \cite{chen2010tests}, we have that $\mV (\overline{Z}_1^\top \overline{Z}_1) \leq (2+\varDelta) \tr(\Sigma_{Z_1}^2)$ and  $\mV (\overline{Z}_2^\top \overline{Z}_2) \leq (2+\varDelta) \tr(\Sigma_{Z_2}^2)$ where $\Sigma_{Z_i} = \mV(Z_i)$ for $i=1,2$. Additionally we know that $\mV (\overline{Z}_1^\top \overline{Z}_2 )  \leq \mE \{ (\overline{Z}_1^\top \overline{Z}_2)^2 \} = \tr(\Sigma_{Z_1} \Sigma_{Z_2})$. Combining the results,
\begin{align*}
\mV (\|Z_1 - Z_2\|^2)  \lesssim \tr \{ (\Sigma_X + \Sigma_Y )^2\}  + (\mu_X - \mu_Y)^\top (\Sigma_X + \Sigma_Y) (\mu_X - \mu_Y).
\end{align*}
Hence $\mV (\|Z_1 - Z_2\|^2) = O(d)$ under (\ref{Eq: CQ assumption}).

\vskip 1em

Next moving onto $\mV \{(Z_1-Z_3)^\top (Z_2 - Z_3)\}$, write $\overline{Z}_3 = Z_3 - \mE(Z_3)$, $\delta_{13} = \mE(Z_1) - \mE (Z_3)$ and $\delta_{23} = \mE (Z_2) - \mE (Z_3)$. Then 
\begin{align*}
& \mV \{(Z_1-Z_3)^\top (Z_2 - Z_3) \} \\[.5em]
= ~ & \mV\{ (\overline{Z}_1-\overline{Z}_3)^\top (\overline{Z}_2 - \overline{Z}_3) + \delta_{13}^\top (\overline{Z}_2 - \overline{Z}_3) + (\overline{Z}_1 - \overline{Z}_3)^\top \delta_{23} \} \\[.5em]
\leq ~ & 3\mV\{ (\overline{Z}_1-\overline{Z}_3)^\top (\overline{Z}_2 - \overline{Z}_3)\} + 3 \mV \{ \delta_{13}^\top (\overline{Z}_2 - \overline{Z}_3)\}  + 3 \mV \{(\overline{Z}_1 - \overline{Z}_3)^\top \delta_{23} \} \\[.5em]
\leq ~ & 12 \mV ( \overline{Z}_1^\top \overline{Z}_2) + 12\mV(\overline{Z}_1^\top \overline{Z}_3) + 12 \mV(\overline{Z}_3^\top \overline{Z}_2) + 12 \mV(\overline{Z}_3^\top \overline{Z}_3)  \\[.5em]
& + 3 \delta_{13}^\top \mV (\overline{Z}_2 - \overline{Z}_3) \delta_{13} + 3 \delta_{23}^\top \mV ( \overline{Z}_1 - \overline{Z}_3 ) \delta_{23}.
\end{align*}
Now similarly as before,
\begin{align*}
\mV \{(Z_1-Z_3)^\top (Z_2 - Z_3) \}   \lesssim  \tr \{ (\Sigma_X + \Sigma_Y )^2\}  + (\mu_X - \mu_Y)^\top (\Sigma_X + \Sigma_Y) (\mu_X - \mu_Y).
\end{align*}
Hence $\mV \{(Z_1-Z_3)^\top (Z_2 - Z_3) \} = O(d)$ under (\ref{Eq: CQ assumption}).

\vskip 2em

\subsection{Generalization of Lemma~\ref{Lemma: Integration over Unit Sphere (2 terms)}} \label{Section: Generalization of Escanciano Lemma}

In Lemma~\ref{Lemma: Extension of Escanciano (2006) with three arguments}, we provided the explicit formula for the integration involving three indicator functions. Here we extend the result to the integration involving four indicator functions.

\begin{lemma} \label{Lemma: Extension of Escanciano (2006) with four arguments}
	For arbitrary vectors $U_1, U_2, U_3, U_4 \in \mathbb{R}^d$, let us denote $\varrho_{ij} = U_i U_j / \{\|U_i\| \|U_j\|\}$ for $i,j \in \{1,2,3,4\}$. Then
	\begin{align} \label{Eq: integration over four indicators}
	\int_{\mathbb{S}^{d-1}}\prod_{i=1}^4 \ind(\beta^\top U_i \leq 0)  d\lambda (\beta) = \frac{7}{16} + \frac{1}{8\pi} \sum_{i=1}^3 \sum_{j=i+1}^4 \mathsf{Ang}\left(U_i, U_j \right) + Q
	\end{align}
	where 
	\begin{align*}
	Q = \frac{1}{4\pi^2} \sum_{\ell=1}^4 \int_{0}^1 \frac{\varrho_{1\ell}}{(1 - \varrho_{1 \ell}^2u^2)^{1/2}} \text{\emph{arcsin}} \Bigg\{ \frac{\gamma_{1,\ell}(u) }{\gamma_{2,\ell}(u) \gamma_{3,\ell}(u) } \Bigg\} du
	\end{align*}
	with
	\begin{align*}
	& \gamma_{1,2} = \varrho_{34} -\varrho_{23} \varrho_{24}  -[ \varrho_{13}\varrho_{14}+ \varrho_{12}(\varrho_{12}\varrho_{34} - \varrho_{14} \varrho_{23}- \varrho_{13} \varrho_{24} ) ]u^2  \\[.5em]
	& \gamma_{1,3} = \varrho_{24} -\varrho_{23} \varrho_{34}  -[ \varrho_{12}\varrho_{14}+ \varrho_{13}(\varrho_{13}\varrho_{24} - \varrho_{14} \varrho_{23}- \varrho_{12} \varrho_{34} ) ]u^2   \\[.5em]
	& \gamma_{1,4} = \varrho_{23} -\varrho_{24} \varrho_{34}  -[ \varrho_{12}\varrho_{13}+ \varrho_{14}(\varrho_{14}\varrho_{23} - \varrho_{13} \varrho_{24}- \varrho_{12} \varrho_{34} ) ]u^2  \\[.5em]
	& \gamma_{2,2} = \gamma_{2,3} = [1 - \varrho_{23}^2 -(\varrho_{12}^2 + \varrho_{13}^2 - 2\varrho_{12}\varrho_{13} \varrho_{23})u^2 ]^{1/2} \\[.5em]
	& \gamma_{3,2} = \gamma_{2,4} = [1 - \varrho_{24}^2 -(\varrho_{12}^2 + \varrho_{14}^2 - 2\varrho_{12}\varrho_{14} \varrho_{24})u^2 ]^{1/2} \\[.5em]
	& \gamma_{3,3} = \gamma_{3,4} = [1 - \varrho_{34}^2 -(\varrho_{13}^2 + \varrho_{14}^2 - 2\varrho_{13}\varrho_{14} \varrho_{34})u^2 ]^{1/2}.
	\end{align*}
	\begin{proof}
		To prove the results, we apply the same argument used in Section~\ref{Section: Proof of Lemma: Integration over unit sphere}. Let $\mathcal{Z}$ have a multivariate normal distribution with zero mean vector and identity covariance matrix. Then as in Section~\ref{Section: Proof of Lemma: Integration over unit sphere}, we have
		\begin{equation}
		\begin{aligned} \label{Eq: Identity of beta and Z in general2}
		\int_{\mathbb{S}^{d-1}} \prod_{i=1}^4 \ind(\beta^\top U_i \leq 0) d \lambda(\beta)  = \mE_{\mathcal{Z}} \bigg[ \prod_{i=1}^4 \ind(\mathcal{Z}^\top U_i \leq 0) \bigg].
		\end{aligned}
		\end{equation}
		Since $(\mathcal{Z}^\top U_1,\mathcal{Z}^\top U_2, \mathcal{Z}^\top U_3, \mathcal{Z}^\top U_4)^\top$ has a multivariate normal distribution with zero mean vector and correlation matrix $[\varrho_{ij}]_{4 \times 4}$ with $\varrho_{ij} = U_i^\top U_j /\{ \|U_i\| \|U_j\|\}$, the right-hand side of (\ref{Eq: Identity of beta and Z in general2}) can be computed based on orthant probabilities for normal distributions \citep[e.g.][]{childs1967reduction,xu2013comparative}. This completes the proof. 
	\end{proof}
\end{lemma}

\begin{remark}
	Although the explicit formula given in Lemma~\ref{Lemma: Extension of Escanciano (2006) with four arguments} looks complicated, it reduces the integral over $\mathbb{S}^{d-1}$ to a more tractable single integral over the unit interval. Hence it would help significantly improve computational time and efficiency in practical applications. 
\end{remark}

\begin{remark}
	\cite{childs1967reduction} also provided expressions for higher order integrations. Using the same argument as before, it is possible to further generalize Lemma~\ref{Lemma: Extension of Escanciano (2006) with four arguments}.
\end{remark}

\vskip 2em 

\subsection{Asymptotic Equivalences between Projection-Averaging and Spatial-Sign Statistics} \label{Section: Asymptotic Equivalences between Projection-Averaging and Spatial-Sign Statistics}
In this section, we provide details on Remark~\ref{Remark: Equivalence of PA to SP}. Based on $U$-statistics, the multivariate one-sample sign test statistic and the two-sample WMW test statistic via projection-averaging can be defined as
\begin{align*}
& U_{\text{Sign-Proj}} = \frac{1}{(m)_2} \sum_{i,j=1}^{m,\neq} h_{\text{Sign-Proj}}(X_i,X_j), \\[.5em]
& U_{\text{WMW-Proj}} = \frac{1}{(m)_2 (n)_2} \sum_{i_1,i_2=1}^{m,\neq} \sum_{j_1,j_2=1}^{n,\neq} h_{\text{WMW-Proj}}(X_{i_1},X_{i_2};Y_{j_1},Y_{j_2}),
\end{align*}
where
\begin{align*}
& h_{\text{Sign-Proj}}(x,y) = \frac{1}{4} - \frac{1}{2\pi} \mathsf{Ang}(x,y) \quad \text{and} \\[.5em] 
& h_{\text{WMW-Proj}}(x_1,x_2;y_1,y_2) = \frac{1}{4} - \frac{1}{2\pi} \mathsf{Ang}(x_1-y_1,x_2-y_2).
\end{align*}
On the other hand, the multivariate one-sample sign test statistic and two-sample WMW test statistic based on the spatial sign are 
\begin{align*}
& U_{\text{Sign-SS}} = \frac{1}{(m)_2} \sum_{i,j=1}^{m,\neq} \frac{X_i^\top X_j}{\|X_i\| \|X_j\|}, \\[.5em]
& U_{\text{WMW-SS}} = \frac{1}{(m)_2 (n)_2} \sum_{i_1,i_2=1}^{m,\neq} \sum_{j_1,j_2=1}^{n,\neq} \frac{(X_{i_1}-Y_{j_1})^\top(X_{i_2}-Y_{j_2})}{\|X_{i_1}-Y_{j_1}\| \|X_{i_2}-Y_{j_2}\|}.
\end{align*}

\vskip 1em

We provide the following proposition for the one-sample case where we prove the asymptotic equivalence between $U_{\text{Sign-Proj}}$ and $U_{\text{Sign-SS}}$.

\begin{proposition} \label{Proposition: One-Sample Equivalence}
	Suppose that ${\mV}[X_1^\top X_2] = O(d)$ and ${\mV}[\|X_1\|^2] = O(d)$. Let us write and assume that
	\begin{align*}
	& \eta_{X,d} = \frac{\|\mu_X\|^2}{\|\mu_X\|^2 + \emph{\tr}(\Sigma_X)} \rightarrow \eta_X \in [0,1), \\[.5em]
	& \delta_{X,d} =  \frac{1}{4} - \frac{1}{2\pi} \text{\emph{arccos}}(\eta_{X,d}) - \frac{\eta_{X,d}}{2\pi(1-\eta_{X,d}^2)^{1/2}}.
	\end{align*}
	Then under the HDLSS setting,
	\begin{align*}
	U_{\text{\emph{Sign-Proj}}} = \delta_{X,d} + \frac{1}{2\pi(1-\eta_{X,d}^2)^{1/2}} U_{\text{\emph{Sign-SS}}} + O_{\mathbb{P}}(d^{-1}).
	\end{align*}
	When $\mu_X = 0$, the expression can be simplified as
	\begin{align*}
	U_{\text{\emph{Sign-Proj}}}  = \frac{1}{\sqrt{2\pi}} U_{\text{\emph{Sign-SS}}} + O_{\mathbb{P}}(d^{-1}).
	\end{align*}
	\begin{proof}
		Similarly as in Section~\ref{Section: Connection of different statistics to CQ}, we use the Taylor expansion and the weak law of large numbers to obtain
		\begin{align*}
		\frac{X_1^\top X_2}{\|X_1\|\|X_2\|} = \eta_{X,d} + O_{\mathbb{P}}(d^{-1/2}).
		\end{align*}
		Next applying the second order Taylor expansion of $f(x) = \text{arccos}(x)$ around $f(\eta_{X,d})$ yields
		\begin{align*}
		\text{arccos} \bigg\{ \frac{X_1^\top X_2}{\|X_1\|\|X_2\|} \bigg\}   = \text{arccos} (\eta_{X,d})- \frac{1}{(1 - \eta_{X,d}^2)^{1/2}} \bigg( \frac{X_1^\top X_2}{\|X_1\|\|X_2\|}  - \eta_{X,d} \bigg) + O_{\mathbb{P}}(d^{-1}).
		\end{align*}
		We finish the proof by plugging this approximation into $U_{\text{{Sign-Proj}}}$.
	\end{proof}
\end{proposition}

\vskip 1em

For the two-sample case, we present the following result. 
\begin{proposition}
	Suppose that $ {\mV}[(X_1-Y_1)^\top (X_2-Y_2)] = O(d)$, ${\mV}[\|X_1 - Y_1\|^2] = O(d)$. Let us write and assume that 
	\begin{align*}
	& \eta_{XY,d} = \frac{\|\mu_X - \mu_Y\|^2}{\|\mu_X -\mu_Y\|^2 + \emph{\tr}(\Sigma_X) + \emph{\tr}(\Sigma_Y)} \rightarrow \eta_{XY} \in [0,1). \\[.5em]
	& \delta_{XY,d} = \frac{1}{4} - \frac{1}{2\pi} \text{\emph{arccos}}(\eta_{XY,d}) - \frac{\eta_{XY,d}}{2\pi(1-\eta_{XY,d}^2)^{1/2}}.
	\end{align*}
	Then under the HDLSS setting,
	\begin{align*}
	U_{\text{\emph{WMW-Proj}}} = \delta_{XY,d}   + \frac{1}{2\pi(1-\eta_{XY,d}^2)^{1/2}} U_{\text{\emph{WMW-SS}}} + O_{\mathbb{P}}(d^{-1}).
	\end{align*}
	When $\mu_X = \mu_Y$, the expression can be simplified as
	\begin{align*}
	U_{\text{\emph{WMW-Proj}}}  = \frac{1}{\sqrt{2\pi}} U_{\text{\emph{WMW-SS}}} + O_{\mathbb{P}}(d^{-1}).
	\end{align*}
	\begin{proof}
		The proof is similar to that of Proposition~\ref{Proposition: One-Sample Equivalence}; hence omitted. 
	\end{proof}
\end{proposition}

\vskip 2em

%\clearpage 

\section{Additional Simulations} \label{Section: Additional Simulations}

This section provides additional simulation results under the setting where the component variables are strongly dependent. Specifically, we assume that $X$ has a multivariate $t$-distribution with the location parameter $\mu_X=(0,\ldots,0)^\top$, the degrees of freedom $\upsilon$ and the $d \times d$ shape matrix $S$ where $[S]_{ij} = 1$ if $i=j$ and $[S]_{ij} = 0.9$ otherwise. Note that when $\upsilon > 2$, the covariance matrix of $X$ is given by $\frac{\upsilon}{\upsilon-2} S$. Similarly, we assume that $Y$ has a multivariate $t$-distribution with the location parameter $\mu_X= (0.2,\ldots,0.2)^\top$, the degrees of freedom $\upsilon$ and the shape matrix $S$. Under the given setting, we generated $m=n=20$ random samples from each distribution with $d=200$ and carried out the permutation tests as in Section~\ref{Section: Simulations}. We increased the degrees of freedom from $\upsilon=1$ to $\upsilon=\infty$ to vary the moment conditions. As shown in Table~\ref{Table: Additional result}, the WMW test performs the best when $\upsilon \leq 7$ closely followed by the CvM test. When $\upsilon$ is large (e.g. $\upsilon \geq 20$) meaning that  $X$ and $Y$ have relatively light-tailed distributions, the power of the five tests (CvM, Energy, MMD, CQ, WMW) are very similar as observed in Section~\ref{Section: Simulations}. These empirical results provide evidence that the findings in Section~\ref{Section: High Dimension, Low Sample Size Analysis} may hold under even more general settings where the component variables are strongly dependent. %We hope that further studies will be carried out to confirm this observation. 

\begin{table}[h!]
	\begin{center}
		\renewcommand{\tabcolsep}{7pt}
		\renewcommand\arraystretch{1}
		\small
		\caption{\small Empirical power of the considered tests at $\alpha=0.05$ against the location models when the component variables are strongly dependent.}
		\label{Table: Additional result}	
	\begin{tabular}{ccccccccc}
					\toprule
		$m=20,n=20$ & $\upsilon=1$ & $\upsilon=3$ & $\upsilon=5$ & $\upsilon=7$ & $\upsilon=9$ & $\upsilon=11$ & $\upsilon=20$ & $\upsilon=\infty$ \\ \midrule
		CvM    & 0.118	& \textbf{0.653}	&\textbf{0.823}	&\textbf{0.880}&	\textbf{0.907}	&\textbf{0.918}	& \textbf{0.943}	& \textbf{0.943}    \\
		Energy & 0.053	&0.332	&\textbf{0.642}	&\textbf{0.808}	& \textbf{0.865}	& \textbf{0.887}	& \textbf{0.937} &	\textbf{0.945}    \\
		MMD    & 0.075&	0.162	&0.363	&0.595	&0.755	&0.810	& \textbf{0.923}	& \textbf{0.945}    \\
		CQ     & 0.063 &	0.470	&\textbf{0.692}	&\textbf{0.815}	&\textbf{0.842}	& \textbf{0.892}	&\textbf{0.920}	& \textbf{0.943}    \\
		WMW    & \textbf{0.340}	& \textbf{0.767}	&\textbf{0.865}	&\textbf{0.892}	&\textbf{0.892}	& \textbf{0.930}	&\textbf{0.942}	& \textbf{0.943}    \\ \midrule
		NN     & \textbf{0.293}	& \textbf{0.490}	&0.528	&0.532	&0.528	&0.533	&0.577	&0.583    \\
		FR     & \textbf{0.225}	& 0.322	&0.305	&0.313	&0.307	&0.293	&0.283	&0.378    \\
		MBG    & 0.047	& 0.062	&0.053	&0.043	&0.048	&0.052	&0.050	&0.100    \\
		Ball   & 0.063	& 0.050	&0.057	&0.053	&0.070	&0.070	&0.075	&0.620    \\
		CM     & 0.052 &	0.067	&0.057	&0.057	&0.065	&0.075	&0.093	&0.125    \\
		BG     & 0.040	& 0.045 &	0.047	&0.040	&0.065	&0.048	&0.058	&0.185    \\
		Run    & 0.112	&0.112	&0.155	&0.152	&0.167	&0.187	&0.198	&0.325  \\
					\bottomrule    
	\end{tabular}
	\end{center}
\end{table}

\clearpage

\end{document}